%% file: Ratner.tex
\newif\ifarXiv \arXivfalse        \newif\iftwocolumn \twocolumntrue
\newif\ifportable \portablefalse   \newif\ifEPS \EPStrue
\newif\ifhyperref \hyperreftrue    \newif\ifLucida \Lucidatrue
\newif\ifcancopy \cancopytrue   
 \newif\ifupdatenotation \updatenotationtrue

 \arXivtrue
 \twocolumnfalse
 \hyperreffalse
\cancopyfalse
\ifarXiv
 \twocolumnfalse
 \portabletrue
 \EPSfalse
 \hyperreffalse
 \Lucidafalse
 \cancopytrue
 \updatenotationfalse %% would like true, but doesn't seem to work
 \fi

\def\wlog#1{}
 \makeatletter \renewcommand{\@latex@info}[1]{} \makeatother
 \renewcommand{\PackageInfo}[2]{}

 \makeatletter \renewcommand{\@font@info}[1]{} \makeatother

 \documentclass[reqno]{amsbook}

\newcommand{\includeall}{\renewcommand{\includeonly}[1]{}}
 \includeall

 \sfcode`\.=1250
 \sfcode`\?=1250
 \sfcode`!=1250
 \sfcode`:=1250
 \sfcode`;=1250
 \sfcode`,=1250

\ifLucida
 \usepackage[LM1]{fontenc}
 \usepackage[expert,
 LM1,
 nofontinfo,
 nolucidascale]{lucidabr}
 \renewcommand{\notin}{\not\in}
 \fi

\usepackage{amssymb, graphicx, ifthen}

\ifhyperref
 \usepackage[textures]{hyperref}
 \else
 \def\href#1{}
 
 \fi

 \ifhyperref
 \makeatletter
 \renewcommand{\Hy@Warning}[1]{}
 \makeatother
 \fi

\newcommand{\twoonone}
 {\twocolumn
\setlength{\textwidth}{9 in}
\setlength{\textheight}{7.25 in}
\setlength{\columnsep}{0.5in}
\hoffset=-0.75in
\voffset=-1in  %% -0.875in
 }
\iftwocolumn \twoonone \fi

 \newcommand{\eps}{\ifEPS eps\else epsf\fi}

\ifportable \def\MR#1{MR~#1} %% had trouble with arXiv on this???
  \else
 \renewcommand{\MR}[1]%
 {\href{http://www.ams.org/mathscinet-getitem?mr=#1}%
 {MR #1}}
 \fi

\makeatletter
\def\plainref#1{\expandafter\@setref\csname r@#1\endcsname\@firstoftwo{#1}}
\makeatother %???

\newcommand{\pref}[1]{\mbox{\upshape(\ref{#1})}}
\newcommand{\refp}[1]{\mbox{\upshape\ref{#1}$'$}}
\newcommand{\prefp}[1]{\mbox{\upshape(\refp{#1})}}
\newcommand{\fullref}[2]{\ref{#1}\pref{#1-#2}}
\renewcommand{\see}[1]{{\upshape(}see~\ref{#1}{\upshape)}}
\newcommand{\fullsee}[2]{(see~\ref{#1}\pref{#1-#2})}
\newcommand{\fullseerem}[2]{(see Rem.~\ref{#1}\pref{#1-#2})}
\newcommand{\fullcf}[2]{(cf.~\ref{#1}\pref{#1-#2})}
\newcommand{\seeand}[2]{(see~\ref{#1} and~\ref{#2})}
\newcommand{\cf}[1]{(cf.~\ref{#1})} \newcommand{\cfdefn}[1]{(cf.\
Defn.~\ref{#1})} \newcommand{\cfcor}[1]{(cf.\ Cor.~\ref{#1})}
\newcommand{\cfprop}[1]{(cf.\ Prop.~\ref{#1})}
\newcommand{\cfand}[2]{(cf.~\ref{#1} and~\ref{#2})}
\newcommand{\seeex}[1]{{\upshape(}see Exer.~\ref{#1}{\upshape)}}
\newcommand{\seeprop}[1]{{\upshape(}see Prop.~\ref{#1}{\upshape)}}
\newcommand{\seethm}[1]{{\upshape(}see Thm.~\ref{#1}{\upshape)}}
\newcommand{\seecor}[1]{{\upshape(}see Cor.~\ref{#1}{\upshape)}}
\newcommand{\seeeg}[1]{{\upshape(}see Eg.~\ref{#1}{\upshape)}}
\newcommand{\seecors}[2]{(see Cors.~\ref{#1} and~\ref{#2})}
\newcommand{\seeprops}[2]{(see Props.~\ref{#1} and~\ref{#2})}
\newcommand{\seelem}[1]{{\upshape(}see Lem.~\ref{#1}{\upshape)}}
\newcommand{\seerem}[1]{{\upshape(}see Rem.~\ref{#1}{\upshape)}}
\newcommand{\seeassump}[1]{{\upshape(}see Assump.~\ref{#1}{\upshape)}}
\newcommand{\cfex}[1]{{\upshape(}cf.\ Exer.~\ref{#1}{\upshape)}} 
\newcommand{\seeexs}[2]{{\upshape(}see Exers.~\ref{#1}
and~\ref{#2}{\upshape)}}
\newcommand{\seeNot}[1]{{\upshape(}see Notn.~\ref{#1}{\upshape)}} 
\newcommand{\seeDefn}[1]{{\upshape(}see Defn.~\ref{#1}{\upshape)}}
\newcommand{\seefig}[1]{(see Fig.~\ref{#1})}
\newcommand{\seeSect}[1]{(see \S\ref{#1})}
\renewcommand{\eqref}[1]{Eq.~\pref{#1}}

\ifLucida %% Lucida goofs up some symbols
 \renewcommand{\S}{\text{\char'244}} 
 \fi

\newcommand{\defit}[2][\temp]{\def\temp{#2}{\it\bfseries#2\/}\index{#1}}
\newcommand{\term}[2][\temp]{\def\temp{#2}{#2}\index{#1}}

\newcommand{\bigset}[2]{\left\{\, #1
 \mathrel{\left| \vphantom {\left\{ #1 \mid #2 \right\} } \right.}
 #2 \,\right\} }

\newcommand{\emphit}{} %% make sure the name isn't already used
 \let\emphit=\emph
\renewcommand{\emph}[1]{{\it\bfseries#1\/}}

\itemsep\smallskipamount
\newcommand{\olditemize}{}
\let\olditemize=\itemize
 \makeatletter
 \renewcommand{\itemize}
 {\ifnum \@itemdepth >0 \smallskip \fi
 \olditemize}
 \makeatother

\newcommand{\oldenumerate}{}
\let\oldenumerate=\enumerate
 \makeatletter
 \renewcommand{\enumerate}
 {\ifnum \@enumdepth >0  \vskip-\lastskip \smallskip\fi
 \oldenumerate}
 \makeatother

\renewcommand{\labelenumi}{\theenumi)}

\makeatletter
\def\@listI{\leftmargin\leftmargini
  \parsep\z@skip
  \topsep\listisep %\itemsep\z@skip
  \listparindent\normalparindent}
\let\@listi\@listI
\def\@listii{\leftmargin\leftmarginii
  \labelwidth\leftmarginii \advance\labelwidth-\labelsep
  \topsep\z@skip \parsep\z@skip \partopsep\z@skip}
 \makeatother

\makeatletter
 \AtBeginDocument{%
  \labelsep=3pt\relax
  \setcounter{enumi}{13}\setcounter{enumii}{13}%
  \setcounter{enumiii}{13}\setcounter{enumiv}{13}%
  \settowidth\leftmargini{\labelenumi\hskip\labelsep}%
  \advance\leftmargini by 0.5\normalparindent
  \settowidth\leftmarginii{\labelenumii\hskip\labelsep}%
  \settowidth\leftmarginiii{\labelenumiii\hskip\labelsep}%
  \settowidth\leftmarginiv{\labelenumiv\hskip\labelsep}%
  \setcounter{enumi}{0}\setcounter{enumii}{0}%
  \setcounter{enumiii}{0}\setcounter{enumiv}{0}%
  \leftmarginv=10pt  \leftmarginvi=\leftmarginv
  \leftmargin=\leftmargini
  \labelwidth=\leftmargini \advance\labelwidth-\labelsep
  \@listi}
 \makeatother

\DeclareMathOperator{\SL}{SL}
\DeclareMathOperator{\GL}{GL}

\DeclareMathOperator{\Mat}{Mat}
\DeclareMathOperator{\SO}{SO}
\DeclareMathOperator{\Ad}{Ad}
\DeclareMathOperator{\Rad}{Rad}
\DeclareMathOperator{\ad}{ad}
\DeclareMathOperator{\Stab}{Stab}
\DeclareMathOperator{\graph}{graph}
\DeclareMathOperator{\vol}{vol}
\DeclareMathOperator{\diam}{diam}
\DeclareMathOperator{\supp}{supp}
\DeclareMathOperator{\Prob}{Prob}
\DeclareMathOperator{\Var}{Var} %% the variety of an ideal
\DeclareMathOperator{\Lyap}{Lyap}
\DeclareMathOperator{\trace}{trace}
\renewcommand{\Re}{\mathop{\rm Re}}
\newcommand{\hatpartA}{\skew6\widehat\partA}

\DeclareFontFamily{U}{euF}{}
\DeclareFontShape{U}{euF}{m}{n}{%
 <-5>*[1.2]eufm5%
 <5-10>gen*[1.2]eufm%
 <10->*[1.2]eufm10%
 }{}
\DeclareFontShape{U}{euF}{b}{n}{%
 <-5>*[1.2]eufb5%
 <5-10>gen*[1.2]eufb%
 <10->*[1.2]eufb10%
 }{}
 \DeclareMathAlphabet{\AMSfrak}{U}{euF}{m}{n}
 \SetMathAlphabet{\AMSfrak}{bold}{U}{euF}{b}{n}

\newcommand{\real}{{\mathord{\mathbb R}}}
\newcommand{\rational}{{\mathord{\mathbb Q}}}
\newcommand{\complex}{{\mathord{\mathbb C}}}
\newcommand{\torus}{{\mathord{\mathbb T}}}
\newcommand{\integer}{{\mathord{\mathbb Z}}}
\renewcommand{\natural}{\mathbb{N}}
\newcommand{\Lie}[1]{\lowercase{\AMSfrak{#1}}}  %% could use MathPiTwo ???
\newcommand{\LieSL}{\mathop{\Lie{SL}}}
\newcommand{\LieGL}{\mathop{\Lie{GL}}}
\newcommand{\so}{\mathop{\Lie{SO}}}
\newcommand{\lie}[1]{\mathord{\underline{#1}}}
\newcommand{\Zar}[1]{\overline{\overline{#1}}}
\newcommand{\fund}{\mathord{\mathcal{F}}}
\newcommand{\iso}{\cong}
\newcommand{\Id}{I}
\newcommand{\past}[1]{\mathord{{#1}^+\mkern-4mu}}
\newcommand{\symmdiff}{\mathbin{\triangle}}
\newcommand{\transpose}{{\kern-1.5pt\mathsf{T}}}
\newcommand{\RP}[1]{\mathord{\real \mathrm{P}}^{#1}}
\newcommand{\measures}{\mathcal{M}}
\newcommand{\Func}{\mathcal{F}}
 \newcommand{\hilbert}{\mathcal{H}}
 \newcommand{\coho}{\mathcal{H}}

\newcommand{\Bern}{T_{\text{\upshape Bern}}}
\newcommand{\Baker}{T_{\text{\upshape Bake}}}
\newcommand{\COIN}{\mathcal{C}}
\newcommand{\T}{\mathcal{T}}
\newcommand{\E}{\mathcal{E}}
\newcommand{\cover}{\widetilde}
\newcommand{\closure}[1]{\overline{#1}}

\newcommand{\head}{\mathsf{H}}
\newcommand{\tail}{\mathsf{T}}

 \ifportable \else
 \input MathPiTwo
 \font\MPT=MathPi2 at 10 pt %% ??? shouldn't hardcode the sizes???
 \font\MPTs=MathPi2 at 7 pt %% ???
 \font\MPTss=MathPi2 at 5 pt %% ???
 \fi
 \newcommand{\MathPiTwo}[3][\MathPiTemp]{\def\MathPiTemp{#3}
 \ifportable \newcommand{#2}{\mathcal{#1}}
 \else
 \newcommand{#2}{\mathord{\mathchoice
 {\mbox{\MPT \csname MathPiTwo#3\endcsname}}
 {\mbox{\MPT \csname MathPiTwo#3\endcsname}}
 {\mbox{\MPTs \csname MathPiTwo#3\endcsname}}
 {\mbox{\MPTss \csname MathPiTwo#3\endcsname}}}}\fi}
 \MathPiTwo[H]{\hyperbolic}{FrakH} %% Fraktur H
 \MathPiTwo{\poly}{P}
 \MathPiTwo{\partA}{A}
 \MathPiTwo{\partB}{B}
 \MathPiTwo{\partR}{R}
 \MathPiTwo{\headset}{H}
    \ifportable \renewcommand{\headset}{\mathbf{H}} \fi
 \MathPiTwo{\tailset}{T}
    \ifportable \renewcommand{\tailset}{\mathbf{T}} \fi
 \MathPiTwo{\module}{W}
 \MathPiTwo{\neigh}{O}
 \MathPiTwo{\ideal}{I}
 \MathPiTwo{\singular}{S}

\DeclareMathSymbol{\easybackslash}{\mathord}{symbols}{"6E}

\makeatletter
\newcommand\@dotsep{4.5}
\def\@tocline#1#2#3#4#5#6#7{\relax
  \ifnum #1>\c@tocdepth % then omit
  \else
    \par \addpenalty\@secpenalty\addvspace{#2}%
    \begingroup \hyphenpenalty\@M
    \@ifempty{#4}{%
      \@tempdima\csname r@tocindent\number#1\endcsname\relax
    }{%
      \@tempdima#4\relax
    }%
    \parindent\z@ \leftskip#3\relax \advance\leftskip\@tempdima\relax
    \rightskip\@pnumwidth plus1em \parfillskip-\@pnumwidth
  #5\leavevmode\hskip-\@tempdima #6\relax
 \ifnum #1>0 % we want dots
    \leaders\hbox{$\m@th
      \mkern \@dotsep mu\raise2pt\hbox{.}\mkern \@dotsep mu$}\fi\hfill
  \hbox to\@pnumwidth{\@tocpagenum{%\ifnum #1 = 0  \bfseries\fi
#7}}\par
    \nobreak
    \endgroup
  \fi}

\renewcommand{\tocchapter}[3]{\@ifnotempty{#2}{\bfseries}%
  \indentlabel{\@ifnotempty{#2}{\ignorespaces#1 #2.\quad}}#3}

\makeatother

 \makeatletter
 \def\swappedhead#1#2#3{%
 \thmnumber{\@upn{\normalfont\@ifnotempty{#2}{(}#2\@ifnotempty{#2}{)~}}}%
  \thmname{#1}%
  \thmnote{ {\the\thm@notefont#3}}}

\def\@thm#1#2#3{%
  \ifhmode\unskip\unskip\par\fi
  \normalfont
  \trivlist
  \let\thmheadnl\relax
  \let\thm@swap\@gobble
  \let\thm@indent\noindent % no indent
  \thm@headfont{\bfseries}% heading font bold
  \thm@notefont{\fontseries\mddefault\upshape}%
  \thm@headpunct{.}% add period after heading
  \thm@headsep 5\p@ plus\p@ minus\p@\relax
  \thm@space@setup
  #1% style overrides
  \@topsep \thm@preskip               % used by thm head
  \@topsepadd \thm@postskip           % used by \@endparenv
  \def\@tempa{#2}\ifx\@empty\@tempa
    \def\@tempa{\@oparg{\@begintheorem{#3}{}}[]}%
  \else
    \refstepcounter{#2}%
    \def\@tempa{\@oparg{\@begintheorem{#3}{\csname the#2\endcsname}}[]}%
  \fi
  \@tempa
}
 \makeatother

\swapnumbers  %% Prefer to have the theorem numbers at the margin.

\theoremstyle{plain}
\newtheorem{thm}[equation]{Theorem}
\newtheorem{lem}[equation]{Lemma}
\newtheorem{cor}[equation]{Corollary}
\newtheorem{prop}[equation]{Proposition}
\newtheorem{conj}[equation]{Conjecture}

\numberwithin{figure}{section}

\numberwithin{section}{chapter}

\numberwithin{equation}{section}

\theoremstyle{definition}
\newtheorem{defn}[equation]{Definition}
\newtheorem{assump}[equation]{Assumption}
\newtheorem{eg}[equation]{Example}
\newtheorem{ceg}[equation]{Counterexample}
\newtheorem{notation}[equation]{Notation}

\newtheorem{techassump}[equation]{Technical assumption}

\theoremstyle{remark}

\newtheorem{rem}[equation]{Remark}
\newtheorem{warn}[equation]{Warning}

\newenvironment{thmref}{\thmrefer}{}
\newcommand{\thmrefer}[1]{\renewcommand\theequation
  {\protect\ref{#1}$'$}\addtocounter{equation}{-1}}

\newcounter{saveenumi}

\makeatletter
 \def\recordsection#1{\@bsphack
  \protected@write\@auxout{}%
         {\string\newsectionref{#1}{\thesection}}%
  \@esphack}
 \def\newsectionref#1#2{\global\@namedef{section@#1}{#2}}
 \def\sectionref#1{\csname section@#1\endcsname}

\newcommand{\shortref}[1]{\ref{short@#1}}

 \makeatother

\makeatletter
 
 \newcommand{\l@exersection}{\@tocline{3}{0pt}{1pc}{7pc}{}}
 \let\exersectionname\@empty
 \let\subsubsectionmark\@gobble
 
 \def\toclevel@exersection{3}
 \newcommand{\exertitle}{Exercises for \S\thesection}

\newenvironment{exercises}
 {\medbreak\@startsection{exersection}{3}
  \z@\z@{-.5em}%
  {\bf}*{\exertitle}
 \begin{enumerate}
 \renewcommand{\theenumi}
    {{\thesection\#\arabic{enumi}}}
 \renewcommand{\labelenumi}{\#\arabic{enumi}.}
 \let\stupidlabel=\label
 \def\label##1{\stupidlabel{##1}
 \recordsection{##1}
 \begingroup
 \ifnum\@enumdepth= 1
 \renewcommand{\@currentlabel}{\arabic{enumi}}
 \else \ifnum\@enumdepth= 2
 \renewcommand{\@currentlabel}{\arabic{enumi}\alph{enumii}}
 \else \message{Error: No exercise number here!!!}
 \fi \fi
 \stupidlabel{short@##1}
 \endgroup}}
 {\end{enumerate}
 }
 \makeatother

 \newcommand{\onlyoneexercise}
 {\begingroup
 \renewcommand{\exertitle}{Exercise for \S\thesection}
 }

\newcommand{\hint}[1]{{\ifvmode\else\newline\fi}[{\smaller\emphit{Hint:}
#1}]}
\newcommand{\nolinehint}[1]{[{\smaller\emphit{Hint:}
#1}]}

\makeatletter
 \newbox\chapwidth
 \setbox\chapwidth\hbox{\bf Chapter 1.}
\renewcommand{\tocsection}[3]{%
  \indentlabel{{\hbox to
\wd\chapwidth{\hfil\@ifnotempty{#2}{\S#2.}}\quad}}#3}

 \def\l@section{\@tocline{1}{2pt}{0pt}{}{}}

 \def\l@chapter{\filbreak\@tocline{0}{8pt plus1pt}{0pt}{}{}}

 \def\l@schapter{\filbreak\@tocline{1}{8pt plus1pt}{0pt}{}{}}
 
 \def\toclevel@schapter{0}

\renewcommand{\tocchapter}[3]{\@ifnotempty{#2}{\bf}%
  \indentlabel{{\hbox to
\wd\chapwidth{\hfil\@ifnotempty{#2}{Chapter #2.}}\quad}}#3}

\def\@schapter#1{\typeout{#1}%
  \let\@secnumber\@empty
  \def\@toclevel{1}%  %% changed from 0
  \ifx\chaptername\appendixname \@tocwriteb\tocappendix{schapter}{#1}%
  \else \@tocwriteb\tocchapter{schapter}{#1}\fi % changed from chapter
  \chaptermark{#1}%
  \addtocontents{lof}{\protect\addvspace{10\p@}}%
  \addtocontents{lot}{\protect\addvspace{10\p@}}%
  \@makeschapterhead{#1}\@afterheading}

 \makeatother

\newenvironment{notes}{\section*{Notes}\addtocounter{section}{1}}{}

\newcommand{\notesect}[1]{\subsection*{\S\ref{#1}}}

\makeatletter
 \renewenvironment{proof}[1][\proofname]{\par\medbreak
  \pushQED{\qed}%
  \normalfont \topsep6\p@\@plus6\p@\relax
  \noindent{\bf
    #1\@addpunct{.} }\ignorespaces
}{%
  \popQED\endtrivlist\@endpefalse\medbreak
}
 \makeatother

\makeatletter
 \renewcommand{\@captionheadfont}{\normalfont}
 
 \makeatother

 \newenvironment{claim}[1][\unskip]{\em
 \medskip \noindent Claim.\ }{\unskip\upshape}

 \newcounter{step}
 
 \renewcommand{\thestep}{\arabic{step}}

 \newcounter{case}
 \newenvironment{case}[1][\unskip]{\refstepcounter{case}\em
 \medskip \noindent Case \thecase\ #1.\ }{\unskip\upshape}
 \renewcommand{\thecase}{\arabic{case}}

 \newcounter{subcase}
 
 \numberwithin{subcase}{case}

 \newcounter{subsubcase}
 
 \numberwithin{subsubcase}{subcase}

 \makeatletter
 \newenvironment{numproof}[1][\proofname]{\par \normalfont
 \refstepcounter{equation}
  \topsep6\p@\@plus6\p@ \trivlist \itemindent\z@
  \item[\hskip\labelsep(\theequation) \bfseries
    #1\@addpunct{.}]\ignorespaces
}{%
  \qed\endtrivlist
}% 
 \makeatother

\renewcommand{\bibname}{References}

\makeatletter
 \renewenvironment{thebibliography}[1]{%
  \section*{\bibname}%
  \normalfont\labelsep .5em\relax
  \renewcommand\theenumiv{\arabic{enumiv}}\let\p@enumiv\@empty
  \list{\@biblabel{\theenumiv}}{\settowidth\labelwidth{\@biblabel{#1}}%
    \leftmargin\labelwidth \advance\leftmargin\labelsep
    \usecounter{enumiv}}%
  \sloppy \clubpenalty\@M \widowpenalty\clubpenalty
  \sfcode`\.=\@m
 \itemsep=\medskipamount
}{%
  \def\@noitemerr{\@latex@warning{Empty `thebibliography' environment}}%
  \endlist
}
 \makeatother

 \makeatletter
 \def\@cite#1#2{{%
 \m@th\upshape\mdseries[{#1}{\if@tempswa, #2\fi}]}}
 \makeatother

\makeatletter
 \def\@starttoc#1#2{
  \begingroup
  \setTrue{#1}
  \let\secdef\@gobbletwo \chapter
  \let\@secnumber\@empty  % for \@tocwrite and \chaptermark
  \ifx\contentsname#2
  \else \@tocwrite{chapter}{#2}\fi
  \typeout{#2}\@xp\chaptermark\@xp{#2}
  \@makeschapterhead{#2}\@afterheading
  \parskip\z@skip
  \makeatletter
  \@input{\jobname.#1}
  \if@filesw
    \@xp\newwrite\csname tf@#1\endcsname
    \immediate\@xp\openout\csname tf@#1\endcsname \jobname.#1\relax
  \fi
  \global\@nobreakfalse \endgroup
  \newpage
}
 \makeatother

\makeatletter
\def\makenotationnindex{%
  \newwrite\@notationfile
  \immediate\openout\@notationfile=\jobname.ntn
  \def\nindex{\@bsphack\begingroup
                \@sanitize
                \@wrnotation}
   \def\nchap{\@bsphack\begingroup
                \@sanitize
                \@wrnchap}\typeout
    {Writing notation index file \jobname.ntn }%
  \let\makenotationnindex\@empty
}
\@onlypreamble\makenotationnindex
\def\@wrnotation#1{%
   \protected@write\@notationfile{}%
      {\string\notationentry{#1}{\thepage}}%
 \endgroup
 \@esphack}
\def\@wrnchap#1{%
   \protected@write\@notationfile{}%
      {\string\notationchapter{#1}{\thechapter}}%
 \endgroup
 \@esphack}
\def\nindex{\@bsphack\begingroup\@sanitize\@index}
\def\nchap{\@bsphack\begingroup\@sanitize\@index}

 \newcommand{\notationentry}[2]{\item #1, #2 \smallskip}
 \newcommand{\notationindexname}{List of Notation}
 \newcommand{\notationchapter}[2]{\bigbreak
 \item \hbox{\hskip-\itemindent \bf \chaptername\ #2. #1} \medskip}

\newcommand{\NotationIndex}
 {\refstepcounter{section} % to correct the PDF bookmarks???
 \chapter*{\notationindexname}
 \begingroup
  \let\item\@idxitem
  \parindent\z@  \parskip\z@\@plus.3\p@\relax
 \vskip-2\bigskipamount
 \input \jobname.ntn
 \clearpage
 \endgroup}
 \makeatother

 \makeatletter
 \newcommand{\mychapter}[1]{\chapter{#1}\nchap{#1}{\def\@thefnmark{}
  \ifcancopy
\@footnotetext{
 Copyright \copyright\ 2003--2005 Dave Witte Morris. All rights
reserved. \\
 Permission to make copies of these lecture
notes for educational or scientific use, including
multiple copies for classroom or seminar teaching, is
granted (without fee), provided that any fees charged
for the copies are only sufficient to recover the
reasonable copying costs, and that all copies include
the title page and this copyright notice. Specific
written permission of the author is required to
reproduce or distribute this book (in whole or in
part) for profit or commercial advantage.}
 \fi
 }}
 \makeatother

\makeatletter
\def\partrunhead#1#2#3{%
  \@ifnotempty{#2}{\emphit{\ignorespaces#1
#2\unskip}\@ifnotempty{#3}{. }}%
  \def\@tempa{#3}%
  \ifx\@empty\@tempa\else
    \begingroup \def\\{ \ignorespaces}% defend against questionable usage
     \emphit{\@tempa}%\@tempa
    \endgroup
  \fi
}
\let\chapterrunhead\partrunhead
\let\sectionrunhead\partrunhead

\def\ps@headings{\ps@empty
  \def\@evenhead{%
    \setTrue{runhead}%
    \normalfont\small %% instead of \scriptsize ???
    \rlap{\thepage}\hfil \leftmark{}{}\hfil}%
  \def\@oddhead{%
    \setTrue{runhead}%
    \normalfont\small %% instead of \scriptsize ???
  \hfil
    \rightmark{}{}\hfil \llap{\thepage}}%
  \let\@mkboth\markboth
  \def\partmark{\@secmark\markboth\partrunhead\partname}%
  \def\chaptermark{%
    \@secmark\markboth\chapterrunhead{}}%
  \def\sectionmark{%
    \@secmark\markright\sectionrunhead\sectionname}%
}

\makeatother

\makeatletter

\def\@citex[#1]#2{%
  \let\@citea\@empty
  \@cite{\@for\@citeb:=#2\do
    {\@citea\def\@citea{,\penalty\@m\ }%
    \edef\@citeb{<Chap\thechapter>\expandafter\@firstofone\@citeb\@empty}%
     \if@filesw\immediate\write\@auxout{\string\citation{\@citeb}}\fi
     \@ifundefined{b@\@citeb}{\mbox{\reset@font\bfseries ?}%
       \G@refundefinedtrue
       \@latex@warning
         {Citation `\@citeb' on page \thepage \space undefined}}%
       {\hbox{\csname b@\@citeb\endcsname}}}}{#1}}

\AtBeginDocument{

\newcommand{\old@lbibitem}{} \let\old@lbibitem=\@lbibitem
\renewcommand{\@lbibitem}{}
 \def\@lbibitem[#1]#2{\old@lbibitem[#1]{<Chap\thechapter>#2}}

\newcommand{\old@bibitem}{} \let\old@bibitem=\@bibitem
\renewcommand{\@bibitem}{}
 \def\@bibitem#1{\old@bibitem{<Chap\thechapter>#1}}

 }

\makeatother

\makeindex

 \ifupdatenotation
 \makenotationnindex
 \fi

\begin{document}

%%% set up shortrefs for Exercises 
%%% need to do after \begin{document} for compatibility with hyperref

 \makeatletter
 \global\let\presmartref=\ref
 \gdef\ref#1{\@ifundefined{section@#1}{\presmartref{#1}}%
 {\ifthenelse{\equal{\thesection}{\sectionref{#1}}}%
 {\shortref{#1}}{\presmartref{#1}}}}
 \makeatother

%% make references upright
%%% need to do after \begin{document} for compatibility with hyperref
 \newcommand{\notupref}{}
 \let\notupref=\ref
 \renewcommand{\ref}[1]{{\upshape\notupref{#1}}}

\title{{\fontsize{0.6in}{0.7in}\selectfont
 Ratner's Theorems \\[-0.15in]
 on
 \\ Unipotent Flows}}

\author{Dave Witte Morris}

 \address{\hskip-\parindent\vbox{
 \centerline{\href{http://www.cs.uleth.ca/}%
 {Department of Mathematics and Computer Science}}
 \centerline{\href{http://www.uleth.ca/}{University of Lethbridge}}
 \centerline{Lethbridge, Alberta, T1K 3M4, Canada}
 \medskip
 \centerline{\href{mailto:Dave.Morris@uleth.ca}%
 {\tt Dave.Morris@uleth.ca}}
 \centerline{\href{http://people.uleth.ca/~dave.morris/}%
 {\tt http://people.uleth.ca/$\sim$dave.morris/}}
 \vskip 0.05in
 }}

% \address{\hskip-\parindent\vbox{
% \centerline{\href{http://www.uleth.ca/fas/mcs/}%
% {Department of Mathematics}}
% \centerline{\href{http://osu.okstate.edu/}%
% {Oklahoma State University}}
% \centerline{Stillwater, Oklahoma 74078, USA}
% \medskip
% \vskip 0.25in
% }}

\dedicatory{\upshape \normalsize \vfill Copyright
\copyright\ 2003--2005 Dave Witte Morris.
  All rights reserved. \\
 \mbox{ } \\
 \ifcancopy
 \small Permission to make copies of this book for educational or
scientific use, including multiple copies for classroom or seminar
teaching, is granted (without fee), provided that any fees charged
for the copies are only sufficient to recover the
reasonable copying costs, and that all copies include
this title page and its copyright notice. Specific
written permission of the author is required to
reproduce or distribute this book (in whole or in
part) for profit or commercial advantage.
 \\ \mbox{ } \\
 \fi
 { \ifarXiv
 {\bfseries ArXiv Final Version 1.1}
 \else
 {\bfseries Final version 1.2 for publication}
 %{\bfseries Online Final Version 1.1}
 \fi
  (February 11, 2005)
 % {\bfseries Comments welcome!}
 } %%???
 \\ \mbox{ } \\
 to appear in: \\  
 \href{http://www.press.uchicago.edu/Complete/Series/CLM.html}%
 {\emphit{Chicago Lectures in Mathematics Series}} \\
 \href{http://www.press.uchicago.edu/}{University of Chicago Press}
 } %% ???

{\sffamily
\maketitle
 }

\newpage

 \pagenumbering{roman}

\setcounter{page}{5}

%% Dedication
 \thispagestyle{empty}
 \text{ } \vfil
 \begin{center}
      \it To Joy, \\ my wife and friend
 \end{center}
 \vfil \vfil
 \newpage

 \tableofcontents

\makeatletter
  %% don't leave blank page between abs, sched, and ack
 \@openrightfalse
 \makeatother

\include{RatnerAbsAck}

\makeatletter
 \@openrighttrue %% RatnerAbsAck turned this off
 \makeatother

 %% to get page 1 on the right-hand side
\newpage \text{ } \thispagestyle{empty} \newpage

\setcounter{page}{1}
 \pagenumbering{arabic}

\include{RatnerIntro}

\include{RatnerEntropy}

\include{RatnerErgodic}

\include{RatnerAlgGrps}

\include{RatnerProof}

 \index{finite fibers|indseealso{self-joining,~finite~fibers}}

 \NotationIndex

\refstepcounter{section} %% to correct the PDF bookmarks???

\makeatletter  %% to eliminate Latex error
 \raggedright
 \input{Ratner.ind}
 \makeatother

\end{document}

%% file: RatnerAbsAck.tex
\thispagestyle{plain}

 \chapter*{Abstract}

 Unipotent flows are well-behaved dynamical systems.  In particular, Marina
Ratner has shown that the closure of every orbit for such a flow is of a
nice algebraic (or geometric) form. This is known as the Ratner
Orbit Closure Theorem; the Ratner Measure-Classification Theorem and the
Ratner Equidistribution Theorem are closely related results. After
presenting these important theorems and some of their consequences, the
lectures explain the main ideas of the proof. Some algebraic technicalities
will be pushed to the background.

Chapter~\ref{IntroChap} is the main part of the book.  It is intended for a
fairly general audience, and provides an elementary introduction to the
subject, by presenting examples that illustrate the theorems, some of their
applications, and the main ideas involved in the proof.

 Chapter~\ref{EntropyChap} gives an elementary introduction to the theory of
entropy, and proves an estimate used in the proof of Ratner's Theorems. It is
of independent interest.

Chapters~\ref{ErgodicChap} and~\ref{AlgGrpsChap} are utilitarian. They
present some basic facts of ergodic theory and the theory of algebraic
groups that are needed in the proof. The reader (or lecturer) may wish to
skip over them, and refer back as necessary.

Chapter~\ref{ProofChap} presents a fairly complete (but not entirely
rigorous) proof of Ratner's Measure-Classification Theorem. Unlike the other
chapters, it is rather technical. The entropy argument that finishes our
presentation of the proof is due to G.~A.~Margulis and G.~Tomanov. Earlier
parts of our argument combine ideas from Ratner's original proof with the
approach of G.~A.~Margulis and G.~Tomanov.

The first four chapters can be read independently, and are intended to be
largely accessible to second-year graduate students. All four are needed for
Chapter~\ref{ProofChap}. A reader who is familiar with ergodic theory and
algebraic groups, but not unipotent flows, may skip Chaps.~\ref{EntropyChap},
\ref{ErgodicChap}, and~\ref{AlgGrpsChap} entirely, and read only
\S\ref{Poly+ShearSect}--\S\ref{JoiningSect} of Chap.~\ref{IntroChap} before
beginning Chap.~\ref{ProofChap}.

\newpage \thispagestyle{plain}

\chapter*{Possible lecture schedules}

It is quite reasonable to stop anywhere after
\S\ref{Poly+ShearSect}. In particular, a single lecture (1--2 hours) can
cover the main points of \S\ref{WhatisRatnerSect}--\S\ref{Poly+ShearSect}.

\medskip

A good selection for a moderate series of lectures would be
 \S\ref{WhatisRatnerSect}--\S\ref{JoiningSect} and \S\ref{OutlinePfSect},
adding \S\ref{TwoEgsSect}--\S\ref{h(g)Sect} if the audience is not familiar
with entropy. For a more logical presentation, one should briefly
discuss \S\ref{PtwiseErgSect} (the Pointwise Ergodic Theorem) before starting
\S\ref{Poly+ShearSect}--\S\ref{JoiningSect}.

\medskip

Here are suggested guidelines for a longer course:

 \begin{enumerate} 

 \item [\S\ref{WhatisRatnerSect}--\S\ref{MeasVersSect}:]
 Introduction to Ratner's Theorems
 (0.5--1.5 hours)

 \item [\S\ref{ApplSect}:]
 Applications of Ratner's Theorems
 (optional, 0--1 hour)

 \item [\S\ref{Poly+ShearSect}--\S\ref{GenShearSect}:]
 Shearing and polynomial divergence
 (1--2 hours)

 \item [\S\ref{EntropyIntroSect}--\S\ref{JoiningSect}:]
 Other basic ingredients of the proof
 (1--2 hours)

 \item [\S\ref{MeasuresToOrbitsSect}:]
 From measures to orbit closures
 (optional, 0--1 hour)

\medskip

 \item [\S\ref{TwoEgsSect}--\S\ref{DefnEntopySect}:]
 What is entropy?
 (1--1.5 hours)

 \item [\S\ref{CalcEntropySect}--\S\ref{h(g)Sect}:]
 How to calculate entropy
 (1--2 hours)

\item [\S\ref{EntropyEstimateSect}:]
 Proof of the entropy estimate
 (optional, 1--2 hours)

\medskip

\item [\S\ref{PtwiseErgSect}:]
 Pointwise Ergodic Theorem
 (0.5--1.5 hours)

\item [\S\ref{MautnerSect}:]
 Mautner Phenomenon
 (optional, 0.5--1.5 hours)

\item [\S\ref{ErgDecompSect}:]
 Ergodic decomposition
 (optional, 0.5--1.5 hours)

\item [\S\ref{AvgSetSect}:]
 Averaging sets
 (0.5--1.5 hours)

\medskip

\item[\S\ref{AlgicGrpsSect}--\S\ref{LieGrpSect}:] Algebraic groups
(optional, 0.5--3 hours)

\medskip

\item [\S\ref{OutlinePfSect}:]
 Outline of the proof
 (0.5--1.5 hours)

\item [\S\ref{ShearPolyPfSect}--\S\ref{EndOfPfSect}:]
 A fairly complete proof
 (3--5 hours)

\item [\S\ref{PreciseSect}--\S\ref{RatnerPfAssumeASect}:]
 Making the proof more rigorous
 (optional, 1--3 hours)

 \end{enumerate}

\newpage \thispagestyle{plain}

\chapter*{Acknowledgments}

I owe a great debt to many people, including the audiences of my
lectures, for their many comments and helpful discussions that added to my
understanding of this material and improved its presentation. A few of the
main contributors are S.~G.~Dani, Alex Eskin, Bassam Fayad, David Fisher,
Alex Furman, Elon Lindenstrauss, G.~A.~Margulis, Howard Masur, Marina Ratner,
Nimish Shah, Robert J.~Zimmer, and three anonymous referees.
 I am also grateful to Michael Koplow, for pointing out many, many minor
errors in the manuscript.

Major parts of the book were written while I was visiting the Tata Institute
of Fundamental Research (Mumbai, India), the Federal Technical Institute
(ETH) of Zurich, and the University of Chicago. I am grateful to my
colleagues at all three of these institutions for their hospitality and for
the aid they gave me in my work on this project. Some financial support was
provided by a research grant from the National Science Foundation
(DMS--0100438). 

I gave a series of lectures on this material at the ETH of Zurich and at the
University of Chicago.
 Chapter~\ref{IntroChap} is an expanded version of a lecture that was first
given at Williams College in 1990, and has been repeated at several other
universities.
 Chapter~\ref{EntropyChap} is based on talks for the University of Chicago's
Analysis Proseminar in 1984 and Oklahoma State University's Lie Groups
Seminar in 2002.

I thank my wife, Joy Morris, for her emotional support and unfailing patience
during the writing of this book.

All author royalties from sales of this book will go to charity.

%% file: RatnerIntro.tex
\mychapter{Introduction to Ratner's Theorems} \label{IntroChap}

\section{What is Ratner's Orbit Closure Theorem?} \label{WhatisRatnerSect}

We begin by looking at an elementary example.

\begin{eg} \label{RatnerOrbForTn}
 For convenience, let us use $[x]$ to denote the image of a point $x \in
\real^n$ in the \term[torus!Tn*$\torus^n$]{$n$-torus}
 \nindex{$\torus^n = \real^n/\integer^n$ = $n$-torus}
 $\torus^n = \real^n/\integer^n$;
that is,
 $$ [x] = x + \integer^n .$$
 Any vector $v \in \real^n$ determines a $C^\infty$~flow~$\varphi_t$
on~$\torus^n$, by
 \begin{equation} \label{DefnTflow}
  \mbox{$\varphi_t \bigl( [x] \bigr) = [x + t v] $ 
 \qquad for $x \in \real^n$ and $t \in \real$}
 \end{equation}
 \seeex{isaflowonT}.
 It is well known that the closure of the orbit of each point
of~$\torus^n$ is a subtorus of~$\torus^n$ (see Exer.~\ref{ThomogOrbits},
or see Exers.~\ref{T2homogOrbits} and~\ref{T3homogOrbEg} for
examples). More precisely, for each $x \in \real^n$, there is a vector
subspace~$S$ of~$\real^n$, such that
 \begin{enumerate} \renewcommand{\theenumi}{S\arabic{enumi}}
 \item \label{RatT-vinS}
 $v \in S$ (so the entire $\varphi_t$-orbit of~$[x]$ is contained in
$[x+S]$),
 \item \label{RatT-Scpct}
 the image $[x+S]$ of $x + S$ in $\torus^n$ is compact (hence,
the image is diffeomorphic to~$\torus^k$, for some $k \in 
\{0,1,2,\ldots,n\}$),
 and
 \item \label{RatT-dense}
 the $\varphi_t$-orbit of~$[x]$  is dense in $[x+S]$ (so $[x+S]$ is
the closure of the orbit of~$[x]$).
 \end{enumerate}
 In short, the closure of every orbit is a nice, geometric subset
of~$\torus^n$.
 \end{eg}

\index{Ratner!Theorem|indsee{Ratner's~Theorems}}
 \index{Theorem!Ratner|indsee{Ratner's~Theorems}}
\term[Ratner's Theorems!Orbit Closure|(]{Ratner's Orbit Closure}
Theorem is a far-reaching generalization of Eg.~\ref{RatnerOrbForTn}. Let
us examine the building blocks of that example.
 \begin{itemize}
 \item Note that $\real^n$ is a \defit{Lie group}. That is, it is a group
(under vector addition) and a \index{manifold}manifold, and the group
operations are~$C^\infty$ functions.
 \item The subgroup $\integer^n$ is
 \index{discrete subgroup|indsee{subgroup,~discrete}}
 \defit[subgroup!discrete]{discrete}.
(That is, it has no \term[point!accumulation]{accumulation points}.)
Therefore, the quotient space $\real^n/\integer^n = \torus^n$ is a
\index{manifold}manifold.
 \item The quotient space $\real^n/\integer^n$ is compact.
 \item The map $t \mapsto tv$ (which appears in the formula
\pref{DefnTflow}) is a \defit[subgroup!one-parameter]{one-parameter
subgroup} of~$\real^n$; that is, it is a $C^\infty$ group
\index{homomorphism!of Lie groups}homomorphism from~$\real$ to~$\real^n$.
 \end{itemize}
 Ratner's Theorem allows:
 \begin{itemize}
 \item the Euclidean space~$\real^n$ to be replaced by any
 \nindex{$G$ = Lie group}
 \defit{Lie group}~$G$;
 \item the subgroup~$\integer^n$ to be replaced by any
 \nindex{$\Gamma$ = lattice in $G$}
 \term[subgroup!discrete]{\emph{discrete subgroup}}~$\Gamma$ of~$G$, such
that the quotient space $\Gamma \backslash G$ is compact;
 and
 \item the map $t \mapsto tv$ to be replaced by any
 \index{subgroup!one-parameter!unipotent|indsee{unipotent~subgroup,
\text{\quad one-parameter}}}
\term[unipotent!subgroup!one-parameter]{\emph{unipotent} one-parameter
subgroup}~$u^t$ of~$G$.
 \nindex{$u^t$ = unipotent one-parameter subgroup}
  (The definition of ``unipotent" will be explained later.)
 \end{itemize}
 Given $G$, $\Gamma$, and~$u^t$, we may define a $C^\infty$~flow
$\varphi_t$ on~$\Gamma \backslash G$ by
 \nindex{$\varphi_t$ = $u^t$-flow on $\Gamma \backslash G$}
  \begin{equation} \label{DefnUnipflow}
  \mbox{$\varphi_t ( \Gamma x ) = \Gamma x u^t $ 
 \qquad for $x \in G$ and $t \in \real$}
 \end{equation}
 (cf.~\ref{DefnTflow} and see Exer.~\ref{isaflow}).
 We may also refer to~$\varphi_t$ as the \defit[unipotent!flow]{\pmb{$u^t$}-flow}
on $\Gamma \backslash G$.
 \index{ut-flow*$u^t$-flow|indsee{unipotent~flow}}
 Ratner proved that the closure of every
$\varphi_t$-orbit is a nice, geometric subset of $\Gamma \backslash G$.
More precisely (note the direct analogy with the conclusions of
Eg.~\ref{RatnerOrbForTn}), if we write
 \nindex{$[x]$ = image of~$x$ in~$\Gamma \backslash G$}
 $[x]$ for the image of~$x$ in~$\Gamma \backslash G$, then, for each
$x \in G$, there is a closed, \term[subgroup!connected]{connected
subgroup}~$S$ of~$G$, such that 
 \begin{enumerate} \renewcommand{\theenumi}{S\arabic{enumi}$'$}
 \item \label{utInS}
 $\{u^t\}_{t \in \real} \subset S$ (so the entire $\varphi_t$-orbit
of~$[x]$ is contained in $[x S]$),
 \item \label{xScpct}
 the image $[x S]$ of $x S$ in $\Gamma \backslash G$ is compact
(hence, diffeomorphic to the \term[homogeneous!space]{homogeneous space}
$\Lambda \backslash S$, for some \term[subgroup!discrete]{discrete
subgroup}~$\Lambda$ of~$S$),
 and
 \item \label{utDenseInS}
 the $\varphi_t$-orbit of~$[x]$  is dense in $[x S]$ (so $[x S]$ is
the closure of the orbit).
 \end{enumerate}

\begin{rem} \ 
 \begin{enumerate}
 \item Recall that 
 \nindex{$\Gamma \backslash G$ = $\{\, \Gamma x \mid x \in G\,\}$}
 $\Gamma \backslash G = \{\, \Gamma x \mid x \in G\,\}$ is the set of
\defit[coset!right]{right} cosets of~$\Gamma$ in~$G$. We will consistently
use right cosets~$\Gamma x$, but all of the results can easily be
translated into the language of \defit[coset!left]{left} cosets $x
\Gamma$. For example, a $C^\infty$ flow~$\varphi_t'$ can be defined on
$G/\Gamma$ by $\varphi_t'(x \Gamma) = u^t x \Gamma$. 
 \item It makes no difference whether we write $\real^n /\integer^n$ or
$\integer^n \backslash \real^n$ for~$\torus^n$, because
\term[coset!right]{right} cosets and \term[coset!left]{left} cosets are
the same in an \term[group!abelian]{abelian group}.
 \end{enumerate}
 \end{rem}

\begin{notation} \label{SL2etcNotation}
 For a very interesting special case, which will be the main topic of most
of this chapter,
 \index{SL(2,R)*$\SL(2,\real)$|(}
 \begin{itemize}
 \item let 
 $$ G = \SL(2,\real) $$
 \nindex{$\SL(2,\real)$ = group of $2 \times 2$ real matrices of
determinant one}%
 be the group of $2 \times 2$ real matrices of
determinant\index{determinant} one; that is
  $$ \SL(2,\real) = \bigset{
 \begin{bmatrix}
 \mathsf{a} & \mathsf{b} \\
 \mathsf{c} & \mathsf{d}
 \end{bmatrix}
 }{
 \begin{matrix}
 \mathsf{a}, \mathsf{b}, \mathsf{c}, \mathsf{d} \in \real, \\
 \mathsf{a} \mathsf{d} - \mathsf{b} \mathsf{c} = 1
 \end{matrix}
 } ,$$
 and
 \item define $u,a  \colon \real \to \SL(2,\real)$ by
 \nindex{$u^t = 
 \begin{bmatrix}
 1&0\\
 t&1
 \end{bmatrix}$}
 $$ \text{
 $u^t = 
 \begin{bmatrix}
 1&0\\
 t&1
 \end{bmatrix}$ 
 \qquad and \qquad 
 $a^t =
 \begin{bmatrix}
 e^{t}&0 \\ 0&e^{-t}
 \end{bmatrix}$.
 } $$
 \end{itemize}
  \nindex{$a^t =
 \begin{bmatrix}
 e^{t}&0 \\ 0&e^{-t}
 \end{bmatrix}$}
 Easy calculations show that 
 \index{commutation relations}
 $$ \mbox{$u^{s+t} = u^s \, u^t$ \quad and \quad $a^{s+t} = a^s \, a^t$} $$
 \seeex{u+a=homo}, so $u^t$ and $a^t$ are
\term[subgroup!one-parameter]{one-parameter subgroups} of~$G$.
 For any subgroup~$\Gamma$ of~$G$, define flows $\eta_t$ and $\gamma_t$
on $\Gamma \backslash G$, by
 \nindex{$\eta_t$ = horocycle flow = $u^t$-flow on $\Gamma \backslash
\SL(2,\real)$}
 \nindex{$\gamma_t$ = geodesic flow = $a^t$-flow on $\Gamma \backslash
\SL(2,\real)$}
% \begin{align*}
% \eta_t( \Gamma x) &= \Gamma x u^t , \\
% \gamma_t( \Gamma x) &= \Gamma x a^t 
% .
% \end{align*}
  $$
 \text{$\eta_t( \Gamma x) = \Gamma x u^t$
 \quad and \quad
 $\gamma_t( \Gamma x) = \Gamma x a^t$.}
 $$

 \end{notation}

\begin{rem}
 Assume (as usual) that $\Gamma$ is \defit[subgroup!discrete]{discrete}
and that $\Gamma \backslash G$ is compact. If $G =\SL(2,\real)$, then, in
geometric terms, 
 \begin{enumerate}
 \item $\Gamma \backslash G$ is (essentially) the \term{unit tangent
bundle} of a compact \term{surface of constant negative curvature}
\seeex{H2exer},
 \item $\gamma_t$ is called the \defit[geodesic!flow]{geodesic flow} on
$\Gamma \backslash G$ \seeex{geodhoroonX},
 and
 \item $\eta_t$ is called the \defit[horocycle!flow]{horocycle flow} on
$\Gamma \backslash G$ \seeex{geodhoroonX}.
 \end{enumerate}
 \index{flow!horocycle|indsee{horocycle~flow}}
 \index{flow!geodesic|indsee{geodesic~flow}}
 \end{rem}
 
\begin{defn} \label{UnipMatDefn}
 A square matrix~$T$ is \defit[unipotent!matrix]{unipotent} if $1$ is the
only (complex) \term{eigenvalue} of~$T$; in other words, $(T-1)^n = 0$,
where $n$~is the number of rows (or columns) of~$T$.
 \end{defn}

\begin{eg}
 Because $u^t$ is a \term[unipotent!matrix]{unipotent matrix} for
every~$t$, we say that $u^t$ is a
\defit[unipotent!subgroup!one-parameter]{unipotent one-parameter subgroup}
of~$G$. Thus, \term[Ratner's Theorems!Orbit Closure]{Ratner's Theorem}
applies to the \term[horocycle!flow]{horocycle flow} $\eta_t$: the closure
of every $\eta_t$-orbit is a nice, geometric subset of $\Gamma \backslash
G$.

More precisely, algebraic calculations, using properties
(\ref{utInS}, \ref{xScpct}, \ref{utDenseInS}) show that $S = G$
\seeex{Cocpct->S=SL2}. Thus, the closure of every orbit is $[G] = \Gamma
\backslash G$. In other words, every $\eta_t$-orbit is \term[dense
orbit]{dense} in the entire space $\Gamma \backslash G$. 
 \end{eg}

\begin{ceg} \label{FractalOrbit}
\index{counterexample}
 In contrast, $a^t$ is \term[unipotent!not]{\emph{not} a unipotent matrix}
(unless $t = 0$), so $\{a^t\}$ is \term[unipotent!not]{\emph{not}} a
unipotent one-parameter subgroup. Therefore, \term[Ratner's Theorems!Orbit
Closure]{Ratner's Theorem} does \emph{not} apply to the
\term[geodesic!flow]{geodesic flow} $\gamma_t$.

Indeed, although we omit the proof, it can be shown that the closures of
some orbits of~$\gamma_t$ are very far from being nice, geometric subsets
of $\Gamma \backslash G$. For example, the closures of some orbits are
\term[fractal]{fractals} (nowhere close to being a \term{submanifold} of $\Gamma
\backslash G$). Specifically, for some orbits, if $C$ is the closure of
the orbit, then some neighborhood (in~$C$) of a point in~$C$ is
homeomorphic to $C' \times \real$, where $C'$~is a \term{Cantor set}.

When we discuss some ideas of Ratner's proof (in \S\ref{Poly+ShearSect}),
we will see, more clearly, why the flow generated by this
\term[subgroup!one-parameter!diagonal]{diagonal one-parameter subgroup}
behaves so differently from a \term[unipotent!flow]{unipotent flow}.
 \end{ceg}

\begin{rem}
 It can be shown fairly easily that \emph{almost} every orbit of the
\term[horocycle!flow]{horocycle flow} $\eta_t$ is \term[dense orbit]{dense}
in~$[G]$, and the same is true for the \term[geodesic!flow]{geodesic
flow}~$\gamma_t$ \cfand{HoroFlowErg}{GeodFlowErg}. Thus, for both of these
flows, it is easy to see that the closure of \emph{almost every} orbit is
$[G]$, which is certainly a nice \index{manifold}manifold. (This means
that the \term{fractal} orbits of \pref{FractalOrbit} are exceptional;
they form a set of measure zero.) The point of \term[Ratner's
Theorems]{Ratner's Theorem} is that it replaces ``almost every" by
``every."
 \end{rem}

Our assumption that $\Gamma \backslash G$ is compact can be relaxed. 

\begin{defn}
 Let $\Gamma$ be a subgroup of a \term{Lie group}~$G$.
 \begin{itemize}
 \item A measure $\mu$ on $G$ is \defit[measure!left invariant]{left
invariant} if $\mu(g A) = \mu(A)$ for all $g \in G$ and all measurable $A
\subset G$. Similarly, $\mu$ is \defit[measure!right invariant]{right
invariant} if $\mu(A g) = \mu(A)$ for all $g$ and~$A$.
 \item Recall that any \term{Lie group}~$G$ has a (left)
\defit[measure!Haar]{Haar measure}; that is, there exists a
\term[measure!left invariant]{left-invariant}
(\term[measure!regular]{regular}) Borel \term{measure}~$\mu$ on~$G$.
Furthermore, $\mu$ is unique up to a scalar multiple. (There is also a
measure that is \term[measure!right invariant]{right invariant}, but the
right-invariant measure may not be the same as the left-invariant measure.)
 \nindex{$\mu$ = measure on $G$ or on $\Gamma \backslash G$}
 \item A \defit{fundamental domain} for a subgroup~$\Gamma$ of a group~$G$
is a measurable subset~$\fund$ of~$G$, such that 
 \begin{itemize}
 \nindex{$\fund$ = fundamental domain for $\Gamma$ in $G$}
 \item $\Gamma \fund = G$,
 and
 \item $\gamma \fund \cap \fund$ has measure~$0$, for all $\gamma \in
\Gamma \smallsetminus \{e\}$.
 \end{itemize}
 \item A subgroup~$\Gamma$ of a \term{Lie group}~$G$ is a \defit{lattice}
if
 \begin{itemize}
 \item $\Gamma$ is \term[subgroup!discrete]{discrete},
 and
 \item some (hence, every) \term{fundamental domain} for $\Gamma$ has
finite measure \seeex{FundDomsSame}.
 \end{itemize}
 \end{itemize}
 \end{defn}

\begin{defn}
 If $\Gamma$ is a \term{lattice} in~$G$, then there is a unique
\term[measure!invariant]{$G$-invariant probability measure}~\nindex{$\mu_G$
= $G$-invariant (``Haar") probability measure on $\Gamma \backslash
G$}$\mu_G$ on $\Gamma \backslash G$ (see Exers.~\ref{LattUnimod->FinVol},
\ref{Latt->FinVol}, and~\ref{LattMeasUnique}). It turns out that $\mu_G$
can be represented by a smooth volume form on the \index{manifold}manifold
$\Gamma \backslash G$. Thus, we may say that $\Gamma \backslash G$ has
\defit{finite volume}. We often refer to $\mu_G$ as the
\defit[measure!Haar]{Haar measure} on $\Gamma \backslash G$.
 \end{defn}

\begin{eg} \label{SL2Zeg}
 Let
 \begin{itemize}
 \item $G = \SL(2,\real)$
 and
 \item \index{SL(2,Z)*$\SL(2,\integer)$}$\Gamma = \SL(2,\integer )$.
 \end{itemize}
 It is well known that $\Gamma$ is a \term{lattice} in~$G$. For example, a
\term{fundamental domain}~$\fund$ is illustrated in Fig.~\ref{sl2(R)}
\seeex{FundForSL2Z}, and an easy calculation shows that the
\term[hyperbolic!measure]{(hyperbolic) measure} of this set is finite
\seeex{SL2R/SL2ZFinVol}.
 \end{eg}

 \begin{figure}
 \begin{center}
 \includegraphics[scale=0.5]{EPS/sl2R.\eps}
 \caption{When $\SL(2,\real)$ is identified with (a double cover of the
unit tangent bundle of) the upper half plane~$\hyperbolic$, the shaded
region is a \term{fundamental domain} for $\SL(2,\integer)$.}
 \label{sl2(R)}
 \end{center}
 \end{figure}

Because compact sets have finite measure, one sees that if $\Gamma
\backslash G$ is compact (and $\Gamma$ is
\defit[subgroup!discrete]{discrete}!), then $\Gamma$ is a \term{lattice}
in~$G$ \seeex{Cocpct->Latt}. Thus, the following result generalizes our
earlier description of Ratner's Theorem. Note, however, that the subspace
$[x S]$ may no longer be compact; it, too, may be a noncompact space of
finite volume.

\begin{thm}[{(\term[Ratner's Theorems!Orbit Closure]{Ratner Orbit Closure
Theorem})}] \label{RatnerOrb}
 If
 \begin{itemize}
 \item $G$ is any \term{Lie group},
 \item $\Gamma$ is any \term{lattice} in~$G$,
 and
 \item $\varphi_t$ is any \term[unipotent!flow]{unipotent flow} on~$\Gamma
\backslash G$,
 \end{itemize}
 then the closure of every $\varphi_t$-orbit is
\term[homogeneous!subset]{homogeneous}.
 \end{thm}

\begin{rem} \label{RatnerOrbPrecise}
 Here is  a more precise statement of the conclusion of \term[Ratner's
Theorems!Orbit Closure]{Ratner's Theorem} \pref{RatnerOrb}. 
 \begin{itemize}
 \item Use $[x]$ to denote the image in $\Gamma \backslash G$ of an
element~$x$ of~$G$.
 \item Let $u^t$ be the \term[unipotent!subgroup!one-parameter]{unipotent
one-parameter subgroup} corresponding to~$\varphi_t$, so $\varphi_t \bigl(
[x] \bigr) = [\Gamma x u^t ]$.
 \end{itemize} 
 Then, for each $x \in G$, there is a \index{subgroup!connected}connected,
closed subgroup~$S$ of~$G$, such that
 \begin{enumerate}
 \item $\{u^t\}_{t \in \real} \subset S$,
 \item the image $[x S]$ of $x S$ in $\Gamma \backslash G$ is closed, and
has \term[finite volume]{finite $S$-invariant volume} (in other words,
$(x^{-1}\Gamma x) \cap S$ is a \term{lattice} in~$S$
\seeex{FinVolOrb->Lattice}),
 and
 \item the $\varphi_t$-orbit of~$[x]$  is dense in $[x S]$.
 \end{enumerate}
 \end{rem}

\begin{eg}
 \ 
 \begin{itemize}
 \item Let $G = \SL(2,\real)$ and \index{SL(2,Z)*$\SL(2,\integer)$}$\Gamma
= \SL(2,\integer)$ as in Eg.~\ref{SL2Zeg}.
 \item Let $u^t$ be the usual
\term[unipotent!subgroup!one-parameter]{unipotent one-parameter subgroup}
of~$G$ (as in Notn.~\ref{SL2etcNotation}).
 \end{itemize}
 Algebraists have classified all of the \term[subgroup!connected]{connected
subgroups} of~$G$ that contain~$u^t$. They are:
 \begin{enumerate}
 \item $\{u^t\}$,
 \item the \term[lower-triangular matrices]{lower-triangular group}
$\left\{ 
 \begin{bmatrix}
 * & 0 \\ * & *
 \end{bmatrix} \right\}$,
 and
 \item $G$.
 \end{enumerate}
 It turns out that the \term[lower-triangular matrices]{lower-triangular
group} does not have a \term{lattice} \cfex{Cocpct->S=SL2}, so we conclude
that the subgroup~$S$ must be either $\{u^t\}$ or~$G$.

In other words, we have the following \term{dichotomy}: 
 \index{SL(2,Z)*$\SL(2,\integer)$}
 $$ \begin{matrix}
 \mbox{each orbit of the \term[unipotent!flow]{$u^t$-flow} on
$\SL(2,\integer) \backslash \SL(2,\real)$} \\
 \mbox{is either closed or \term[dense orbit]{dense}.}
 \end{matrix}$$
 \end{eg}

\begin{eg}
 Let 
 \begin{itemize}
 \item \index{SL(3,R)*$\SL(3,\real)$}$G = \SL(3,\real)$,
 \item \index{SL(3,Z)*$\SL(3,\integer)$}$\Gamma = \SL(3,\integer)$,
 and
 \item $u^t = \begin{bmatrix}
 1 & 0 & 0 \\
 t & 1 & 0 \\
 0 & 0 & 1
 \end{bmatrix}$.
 \end{itemize}
 Some orbits of the \term[unipotent!flow]{$u^t$-flow} are closed, and some are
dense, but there are also intermediate possibilities.  
 For example, $\SL(2,\real)$ can be embedded in the top left corner of
$\SL(3,\real)$:
 $$ \SL(2,\real) \iso 
 \left\{ 
 \begin{bmatrix}
 * & * & 0 \\
 * & * & 0 \\
 0 & 0 & 1
 \end{bmatrix}
 \right\}
 \subset \SL(3,\real) .$$
 \index{SL(3,R)*$\SL(3,\real)$}
 This induces an embedding \index{SL(2,Z)*$\SL(2,\integer)$}
 \begin{equation} \label{EmbedSL2R/SL2Z}
 \SL(2,\integer)\backslash \SL(2,\real) \hookrightarrow 
 \SL(3,\integer)\backslash \SL(3,\real) .
 \end{equation}
 \index{SL(3,Z)*$\SL(3,\integer)$}%
 The image of this embedding is a \term{submanifold}, and it is the
closure of certain orbits of the $u^t$-flow \seeex{topcorner=closure}.
 \end{eg}

\begin{rem} \label{RatOrbH}
 \term[Ratner's Theorems!Orbit Closure]{Ratner's Theorem} \pref{RatnerOrb}
also applies, more generally, to the orbits of any subgroup~$H$ that is
\term{generated by unipotent elements}, not just a
one-dimensional\index{dimension!of a Lie group} subgroup. (However, if
the subgroup is disconnected, then the subgroup~$S$ of
Rem.~\ref{RatnerOrbPrecise} may also be disconnected. It is true, though,
that every connected \index{component!connected}component of~$S$ contains
an element of~$H$.)
 \end{rem}

\begin{exercises}

\item Show that, in general, the closure of a \term{submanifold} may be a bad
set, such as a \term{fractal}. (\term[Ratner's Theorems!Orbit Closure]{Ratner's Theorem} shows that this pathology cannot
not appear if the \term{submanifold} is an orbit of a
\term[unipotent!flow]{``unipotent" flow}.) More precisely, for any closed
subset~$C$ of~$\torus^2$, show there is an injective $C^\infty$ function
$f \colon \real \to \torus^3$, such that 
 $$ \closure{f(\real)} \cap \bigl( \torus^2 \times \{0\} \bigl) = C \times
\{0\} ,$$
 where $\closure{f(\real)}$ denotes the closure of the image of~$f$.
 \hint{Choose a countable, dense subset~$\{c_n\}_{n = -\infty}^\infty$
of~$C$, and choose~$f$ (carefully!) with $f(n) = c_n$.}

\item \label{isaflowonT}
 Show that \pref{DefnTflow} defines a $C^\infty$~\defit{flow}
on~$\torus^n$; that is,
 \begin{enumerate}
 \item $\varphi_0$ is the identity map,
 \item $\varphi_{s+t}$ is equal to the composition $\varphi_s \circ
\varphi_t$, for all $s,t \in \real$;
 and
 \item the map $\varphi \colon \torus^n \times \real \to \torus^n$,
defined by $\varphi(x,t) = \varphi_t(x)$ is~$C^\infty$.
 \end{enumerate}

\item \label{T2homogOrbits}
 Let $v = (\alpha, \beta) \in \real^2$. Show, for each $x \in \real^2$,
that
 the closure of $[x + \real v]$ is
 $$\begin{cases}
 \hfil \bigl\{[x]\bigr\} & \mbox{if $\alpha = \beta = 0$}, \\
 \hfil [x + \real v] & \mbox{if $\alpha/\beta \in \rational$ (or $\beta =
0$)}, \\
 \hfil \torus^2 & \mbox{if $\alpha/\beta \notin \rational$ (and $\beta \neq
0$)}.
 \end{cases}
 $$

\item \label{T3homogOrbEg}
 Let $v = (\alpha,1,0) \in \real^3$, with
$\alpha$~\term[irrational!number]{irrational}, and let $\varphi_t$ be the
corresponding flow on~$\torus^3$ \see{DefnTflow}. Show that the subtorus
$\torus^2  \times \{0\}$ of~$\torus^3$ is the closure of the
$\varphi_t$-orbit of $(0,0,0)$.

\item \label{ThomogOrbits}
 Given $x$ and~$v$ in~$\real^n$, show that there is a vector subspace~$S$
of~$\real^n$, that satisfies \pref{RatT-vinS}, \pref{RatT-dense},
and~\pref{RatT-dense} of Eg.~\ref{RatnerOrbForTn}.

\item \label{SIndepofx}
 Show that the subspace~$S$ of Exer.~\ref{ThomogOrbits} depends only
on~$v$, not on~$x$. (This is a special property of
\term[group!abelian]{abelian groups}; the analogous statement is \emph{not}
true in the general setting of \term[Ratner's Theorems!Orbit Closure|)]{Ratner's Theorem}.)

\item \label{isaflow}
 Given
 \begin{itemize}
 \item a \term{Lie group}~$G$,
 \item a closed subgroup~$\Gamma$ of~$G$,
 and
 \item a \term[subgroup!one-parameter]{one-parameter subgroup}~$g^t$
of~$G$,
 \end{itemize}
 show that $\varphi_t(\Gamma x) = \Gamma x g^t$ defines a \term{flow} on
$\Gamma \backslash G$.

\item \label{u+a=homo}
 For $u^t$ and~$a^t$ as in Notn.~\ref{SL2etcNotation}, and all $s,t \in
\real$, show that
 \index{commutation relations}
 \begin{enumerate}
 \item $u^{s+t} = u^s u^t$,
 and
 \item $a^{s+t} = a^s a^t$.
 \end{enumerate}

\item \label{ANormsU}
 Show that the subgroup $\{a^s\}$ of $\SL(2,\real)$
\term[normalizer]{normalizes} the subgroup $\{u^t\}$. That is, $a^{-s}
\{u^t\} a^s = \{u^t\}$ for all~$s$.

\item \label{H2exer}
 Let 
 \nindex{$\hyperbolic$ = hyperbolic plane}
 $ \hyperbolic = \{\, x + iy \in \complex \mid y > 0 \,\} $
 be the \defit{upper half plane} (or \defit[hyperbolic!plane]{hyperbolic
plane}), with Riemannian metric $\langle \cdot \mid \cdot \rangle$ defined
by
 $$ \langle v \mid w \rangle_{x+iy} = \frac{1}{y^2} (v \cdot w) ,$$
 for tangent vectors $v,w \in T_{x+iy} \hyperbolic$, where $v \cdot w$ is
the usual Euclidean inner product on $\real^2 \iso \complex$.
 \begin{enumerate}
 \item Show that the formula
 $$\mbox{$\displaystyle g z = \frac{ \mathsf{a} z + \mathsf{c} }{\mathsf{b} z + \mathsf{d}
}$
 \qquad for
 $z \in \hyperbolic$ and $g = \begin{bmatrix}
 \mathsf{a} & \mathsf{b} \\
 \mathsf{c} & \mathsf{d} 
 \end{bmatrix}
 \in \SL(2,\real)
 $}$$
 defines an action of $\SL(2,\real)$ by \term[isometry]{isometries}
on~$\hyperbolic$.
 \item Show that this action is \term[transitive action]{transitive}
on~$\hyperbolic$.
 \item Show that the \index{stabilizer!of a point}stabilizer
$\Stab_{\SL(2,\real)}(i)$ of the point~$i$ is 
 \index{SO(Q)*$\SO(Q)$}
 $$\SO(2) = 
 \bigset{
 \begin{bmatrix}
 \cos \theta, & \sin \theta \\
 -\sin \theta & \cos \theta
 \end{bmatrix}
 }{
 \theta \in \real} .$$
 \item The \defit{unit tangent bundle} $T^1 \hyperbolic$ consists of the
tangent vectors of length~$1$. By differentiation, we obtain an action of
$\SL(2,\real)$ on $T^1 \hyperbolic$. Show that this action is
\term[transitive action]{transitive}.
 \item \label{H2exer-T1=G}
 For any unit tangent vector $v \in T^1 \hyperbolic$, show
\index{stabilizer!of a tangent vector}
 $$\Stab_{\SL(2,\real)}(v) = \pm \Id .$$
 Thus, we may identify $T^1 \hyperbolic$ with $\SL(2,\real)/\{\pm \Id\}$.
 \item \label{H2exer-geod}
 It is well known that the \term[geodesic!in the hyperbolic
plane]{geodesics} in~$\hyperbolic$ are semicircles (or lines) that are
orthogonal to the real axis. Any $v \in T^1 \hyperbolic$ is tangent to a
unique geodesic. The \defit[geodesic!flow]{geodesic flow}~$\hat\gamma_t$
on~$T^1 \hyperbolic$ moves the unit tangent vector~$v$ a distance~$t$
along the geodesic it determines. Show, for some vector~$v$ (tangent to the
imaginary axis), that, under the identification of
Exer.~\ref{H2exer-T1=G}, the \term[geodesic!flow]{geodesic flow}
$\hat\gamma_t$ corresponds to the flow $x \mapsto x a^t$ on
$\SL(2,\real)/\{\pm \Id\}$, for some $c \in \real$.
 \item The \defit[horocycle!in the hyperbolic plane]{horocycles}
in~$\hyperbolic$ are the circles that are tangent to the real axis (and
the lines that are parallel to the real axis). Each $v \in T^1\hyperbolic$
is an inward unit \index{normal!vector}normal vector to a unique
horocycle~$H_v$.  The \defit[horocycle!flow]{horocycle flow}~$\hat\eta_t$
on~$T^1 \hyperbolic$ moves the unit tangent vector~$v$ a distance~$t$
(counterclockwise, if $t$~is positive) along the corresponding
horocycle~$H_v$. Show, for the identification in Exer.~\ref{H2exer-geod},
that the horocycle flow corresponds to the flow $x \mapsto x u^t$ on
$\SL(2,\real)/\{\pm \Id\}$.
 \end{enumerate}

\item \label{geodhoroonX}
 Let $X$ be any compact, connected \term[surface of constant
negative curvature]{surface of (constant) negative curvature~$-1$}. We use
the notation and terminology of Exer.~\ref{H2exer}. It is known that there
is a \term[map!covering]{covering map} $\rho \colon \hyperbolic \to X$ that
is a local isometry. Let 
 $$\Gamma = \{\, \gamma \in \SL(2,\real) \mid \mbox{$\rho(\gamma z) =
\rho(z)$ for all $z \in \hyperbolic$} \,\} .$$
 \vskip-\smallskipamount
 \begin{enumerate}
 \item Show that 
 \begin{enumerate}
 \item $\Gamma$ is \term[subgroup!discrete]{discrete},
 and
 \item $\Gamma \backslash G$ is compact.
 \end{enumerate}
 \item Show that the \term{unit tangent bundle} $T^1 X$ can be identified
with  $\Gamma \backslash G$, in such a way that
 \begin{enumerate}
 \item the \term[geodesic!flow]{geodesic flow} on $T^1 X$ corresponds to
the flow $\gamma_t$ on $\Gamma \backslash \SL(2,\real)$,
 and
  \item the \term[horocycle!flow]{horocycle flow} on $T^1 X$ corresponds
to the flow $\eta_t$ on $\Gamma \backslash \SL(2,\real)$.
 \end{enumerate}
 \end{enumerate}

 \begin{figure}
 \begin{center}
 \includegraphics[scale=0.44035]{EPS/geodhyper.\eps}
 \caption{The geodesic flow on $\hyperbolic$.}
 \label{geodhyper}
 \end{center}
 \end{figure}

 \begin{figure}
 \begin{center}
 \includegraphics[scale=0.44035]{EPS/horohyper.\eps}
 \caption{The horocycle flow on $\hyperbolic$.}
 \label{horocycles}
 \end{center}
 \end{figure}

\item Suppose $\Gamma$ and~$H$ are subgroups of a group~$G$. For $x \in
G$, let \nindex{$\Stab_H(\Gamma x)$ = stabilizer of $\Gamma x$
in~$H$}
 $$\Stab_H(\Gamma x) = \{\, h \in H \mid \Gamma x h = \Gamma x \,\} $$
 be the \defit[stabilizer!of a point]{stabilizer} of~$\Gamma x$
in~$H$.
 Show $\Stab_H(\Gamma x) = x^{-1} \Gamma x \cap H$.

\item \label{Cocpct->S=SL2}
 Let 
 \begin{itemize}
 \item $G = \SL(2,\real)$
 \item $S$ be a \term[subgroup!connected]{connected subgroup} of
$G$ containing $\{u^t\}$,
 and 
 \item $\Gamma$ be a \term[subgroup!discrete]{discrete subgroup} of
$G$, such that $\Gamma \backslash G$ is compact.
 \end{itemize}
 It is known (and you may assume) that
 \begin{enumerate} \renewcommand{\theenumi}{\alph{enumi}}
 \item if $\dim S = 2$, then $S$ is conjugate to the
\term[lower-triangular matrices]{lower-triangular group}~$B$,
 \item if there is a \term[subgroup!discrete]{discrete subgroup}~$\Lambda$
of~$S$, such that $\Lambda \backslash S$ is compact, then $S$ is
\defit[Lie group!unimodular]{unimodular}, that is, the
\index{determinant}determinant of the linear transformation $\Ad_S g$
is~$1$, for each $g \in S$,
 and
 \item $\Id$ is the only \term[unipotent!matrix]{unipotent matrix}
in~$\Gamma$.
 \end{enumerate}
 Show that if there is a \term[subgroup!discrete]{discrete
subgroup}~$\Lambda$ of~$S$, such that
 \begin{itemize}
 \item $\Lambda \backslash S$ is compact, 
 and
 \item $\Lambda$ is conjugate to a subgroup of~$\Gamma$,
 \end{itemize}
 then $S = G$.
 \hint{If $\dim S \in \{1,2\}$, obtain a contradiction.}

\item \label{FundDomsSame}
 Show that if $\Gamma$ is a \term[subgroup!discrete]{discrete subgroup}
of~$G$, then all \term[fundamental domain]{fundamental domains}
for~$\Gamma$ have the same measure.
 In particular, if one fundamental domain has finite measure, then all do.
 \hint{$\mu(\gamma A) = \mu(A)$, for all $\gamma \in \Gamma$, and every
subset~$A$ of~$\fund$.}

\item \label{LattUnimod->FinVol}
 Show that if $G$~is \defit[Lie group!unimodular]{unimodular} (that is, if
the left \term[measure!Haar]{Haar measure} is also
\term[measure!invariant]{invariant} under right translations) and $\Gamma$
is a \term{lattice} in~$G$, then there is a
\term[measure!invariant]{$G$-invariant probability measure} on $\Gamma
\backslash G$. \hint{For $A \subset \Gamma \backslash G$, define 
 $ \mu_G(A) = \mu \bigl( \{\, g \in \fund \mid \Gamma g \in A \,\}
\bigr) $.}

\item \label{Latt->FinVol}
 Show that if $\Gamma$ is a \term{lattice} in~$G$, then there is a
\term[measure!invariant]{$G$-invariant probability measure} on $\Gamma
\backslash G$. \hint{Use the uniqueness of \term[measure!Haar]{Haar
measure} to show, for $\mu_G$ as in Exer.~\ref{LattUnimod->FinVol} and $g
\in G$, that there exists $\Delta(g) \in \real^+$, such that $\mu_G(A g) =
\Delta(g) \, \mu_G(A)$ for all $A \subset \Gamma \backslash G$. Then show
$\Delta(g) = 1$.}

\item \label{LattMeasUnique}
 Show that if $\Gamma$ is a \term{lattice} in~$G$, then the
\term[measure!invariant]{$G$-invariant probability measure}~$\mu_G$ on
$\Gamma \backslash G$ is unique.
 \hint{Use $\mu_G$ to define a \term[measure!invariant]{$G$-invariant
measure} on~$G$, and use the uniqueness of Haar measure.}

\item \label{FundForSL2Z}
 Let 
 \begin{itemize}
 \item $G = \SL(2,\real)$,
 \item \index{SL(2,Z)*$\SL(2,\integer)$}$\Gamma = \SL(2,\integer)$,
 \item \index{fundamental domain}$\fund = \{\, z \in \hyperbolic \mid
\mbox{$|z| \ge 1$ and $-1/2 \le \Re z \le 1/2$} \, \}$,
 and
 \item $e_1 = (1,0)$ and $e_2 = (0,1)$,
 \end{itemize}
 and define
 \begin{itemize}
 \item $B \colon G \to \real^2$ by $B(g) = (g^\transpose e_1,
g^\transpose e_2)$, where $g^\transpose$ denotes the
transpose of~$g$,
 \item $C \colon \real^2 \to \complex$ by $C(x,y) = x + i y$,
 and
 \item $\zeta \colon G \to \complex$ by 
 $$ \zeta(g) = \frac{C( g^\transpose e_2)}{C( g^\transpose e_1)} .$$
 \end{itemize}
 Show: 
 \begin{enumerate}
 \item \term[hyperbolic!plane]{$\zeta(G) = \hyperbolic$},
 \item $\zeta$ induces a \index{homeomorphism}homeomorphism
$\overline{\zeta} \colon \hyperbolic \to \hyperbolic$, defined by
$\overline{\zeta}(g i) = \zeta(g)$,
 \item $\zeta(\gamma g) = \gamma \, \zeta(g)$, for all $g \in G$ and $\gamma \in
\Gamma$,
 \item for $g,h \in G$, there exists $\gamma \in \Gamma$, such
that $\gamma g = h$ if and only if $\langle g^\transpose e_1, g^\transpose
e_2 \rangle_{\integer} = \langle h^\transpose e_1, h^\transpose e_2
\rangle_{\integer}$, where $\langle v_1, v_2 \rangle_{\integer}$ denotes
the \index{group!abelian}abelian group consisting of all integral linear
combinations of~$v_1$ and~$v_2$,
 \item for $g \in G$, there exist $v_1,v_2 \in \langle g^\transpose e_1, g^\transpose
e_2 \rangle_{\integer}$, such that
 \begin{enumerate}
 \item $\langle v_1,v_2 \rangle_{\integer} =  \langle g^\transpose e_1, g^\transpose
e_2 \rangle_{\integer}$,
 and
 \item ${C(v_2)}{C(v_1)} \in \fund$,
 \end{enumerate}
 \item $\Gamma \fund = \hyperbolic$,
 \item if $\gamma \in \Gamma \smallsetminus \{\pm\Id\}$, then $\gamma \fund
\cap \fund$ has measure~$0$,
 and
 \item $\{\, g \in G \mid g i \in \fund\,\}$ is is a fundamental domain
for~$\Gamma$ in~$G$.
 \end{enumerate}
 \hint{Choose $v_1$ and~$v_2$ to be a nonzero vectors of minimal length in
$\langle g^\transpose e_1, g^\transpose e_2 \rangle_{\integer}$ and
$\langle g^\transpose e_1, g^\transpose e_2 \rangle_{\integer}
\smallsetminus \integer v_1$, respectively.}

\item \label{SL2R/SL2ZFinVol}
 Show:
 \begin{enumerate}
 \item the area element on the \term[hyperbolic!plane]{hyperbolic
plane}~$\hyperbolic$ is $dA = y^{-2} \, dx \, dy$,
 and
 \item the \term{fundamental domain}~$\fund$ in Fig.~\ref{sl2(R)} has
finite hyperbolic area.
 \end{enumerate}
 \hint{We have $\int_a^\infty \int_b^c y^{-2} \, dx \, dy <
\infty$.}

\item \label{FinMeasToFund}
 Show that if
 \begin{itemize}
 \item $\Gamma$ is a \term[subgroup!discrete]{discrete subgroup} of a
\term{Lie group}~$G$,
 \item $F$ is a measurable subset of~$G$,
 \item $\Gamma F = G$,
 and
 \item $\mu(F) < \infty$,
 \end{itemize}
 then $\Gamma$ is a \term{lattice} in~$G$.

\item \label{Cocpct->Latt}
 Show that if
 \begin{itemize}
 \item $\Gamma$ is a \term[subgroup!discrete]{discrete subgroup} of a
\term{Lie group}~$G$,
 and
 \item $\Gamma \backslash G$ is compact,
 \end{itemize}
 then $\Gamma$ is a \term{lattice} in~$G$.
 \hint{Show there is a compact subset~$C$ of~$G$, such that $\Gamma
C = G$, and use Exer.~\ref{FinMeasToFund}.}

\item \label{FinVolOrb->Lattice}
 Suppose
 \begin{itemize}
 \item $\Gamma$ is a \term[subgroup!discrete]{discrete subgroup} of a
\term{Lie group}~$G$,
 and
 \item $S$ is a closed subgroup of~$G$.
 \end{itemize}
 Show that if the image $[x S]$ of $x S$ in $\Gamma \backslash G$ is
closed, and has \term[finite volume]{finite $S$-invariant volume}, then
$(x^{-1}\Gamma x) \cap S$ is a \term{lattice} in~$S$.

\item \label{DivInG/Gamma}
 Let
 \begin{itemize}
 \item $\Gamma$ be a \term{lattice} in a \term{Lie group}~$G$,
 \item $\{x_n\}$ be a sequence of elements of~$G$.
 \end{itemize}
 Show that $[x_n]$ has no subsequence that converges in $\Gamma \backslash
G$ if and only if there is a sequence $\{\gamma_n\}$ of nonidentity
elements of~$\Gamma$, such that $x_n^{-1} \gamma_n x_n \to e$ as $n \to
\infty$.
 \hint{($\Leftarrow$)~Contrapositive. If $\{x_{n_k}\} \subset \Gamma C$,
where $C$ is compact, then $x_n^{-1} \gamma_n x_n$ is bounded away
from~$e$.
 ($\Rightarrow$)~Let $\neigh$ be a small open subset of~$G$. By passing to
a subsequence, we may assume $[x_m \neigh] \cap [x_n \neigh] = \emptyset$,
for $m \neq n$. Since $\mu(\Gamma \backslash G) < \infty$, then $\mu([x_n
\neigh]) \neq \mu(\neigh)$, for some~$n$. So the natural map $x_n \neigh
\to [x_n \neigh]$ is not injective. Hence, $x_n^{-1} \gamma x_n \in \neigh
\neigh^{-1}$ for some $\gamma \in \Gamma$.}

\item \label{Lattice->ClosedOrbit}
 Prove the converse of Exer.~\ref{FinVolOrb->Lattice}. That is, if
$(x^{-1}\Gamma x) \cap S$ is a \term{lattice} in~$S$, then the image $[x
S]$ of $x S$ in $\Gamma \backslash G$ is closed (and has \term[finite
volume]{finite $S$-invariant volume}).
 \hint{Exer.~\ref{DivInG/Gamma} shows that the inclusion of $\bigl(
(x^{-1}\Gamma x) \cap S \bigr) \backslash S$ into $\Gamma \backslash G$ is
a \index{map!proper}proper map.}

\item \label{topcorner=closure}
 Let $C$ be the image of the embedding \pref{EmbedSL2R/SL2Z}. Assuming
that $C$ is closed, show that there is an orbit of the
\term[unipotent!flow]{$u^t$-flow} on
\index{SL(3,Z)*$\SL(3,\integer)$}$\SL(3,\integer) \backslash \SL(3,\real)$
whose closure is~$C$.

\item {[\emphit{Requires some familiarity with
\term[hyperbolic!geometry]{hyperbolic geometry}}]}
 Let $M$ be a compact,
\term[hyperbolic!n-manifold*$n$-manifold]{hyperbolic $n$-manifold}, so $M
= \Gamma \backslash \hyperbolic^n$, for some
\defit[subgroup!discrete]{discrete} group~$\Gamma$ of isometries of
\term[hyperbolic!n-space*$n$-space]{hyperbolic $n$-space}~$\hyperbolic^n$.
For any $k \le n$, there is a natural embedding $\hyperbolic^k
\hookrightarrow \hyperbolic^n$. Composing this with the
\index{map!covering}covering map to~$M$ yields a $C^\infty$ immersion $f
\colon \hyperbolic^k \to M$. Show that if $k \neq 1$, then there is a
compact \index{manifold}manifold~$N$ and a $C^\infty$ function $\psi
\colon N \to M$, such that the closure $\closure{f(\hyperbolic^k)}$ is
equal to $\psi(N)$.

\item 
 Let \index{SL(2,Z)*$\SL(2,\integer)$}$\Gamma = \SL(2,\integer)$ and $G =
\SL(2,\real)$. Use \term[Ratner's Theorems!Orbit Closure]{Ratner's Orbit
Closure Theorem} (and Rem.~\ref{RatOrbH}) to show, for each $g \in G$,
that $\Gamma g \Gamma$ is either dense in~$G$ or
\term[subgroup!discrete]{discrete}.
 \hint{You may assume, without proof, the fact that if $N$ is any
\term[subgroup!connected]{connected subgroup} of~$G$ that is
\term[normalizer]{normalized} by~$\Gamma$, then either $N$~is trivial, or
$N = G$. (This follows from the Borel Density Theorem \pref{BDT}.}

 \end{exercises}

\section{Margulis, Oppenheim, and quadratic forms} \label{MargOppSect}

\index{quadratic form|(}

\term{Ratner's Theorems} have important applications in
number theory. In particular, the following result was a major motivating
factor. It is often called the ``\term{Oppenheim Conjecture}," but that
terminology is no longer appropriate, because it was proved more than 15
years ago, by G.~A.~Margulis. See \S\ref{ApplSect} for other (more recent)
applications.

\begin{defn} \label{QuadFormDefn} \ 
 \begin{itemize}
 \item A (real) \defit{quadratic form} is a homogeneous polynomial of
degree~$2$ (with real coefficients), in any number of variables.  For
example, 
 $$ Q(x,y,z,w) = x^2 - 2 xy + \sqrt{3} yz - 4 w^2 $$
 is a quadratic form (in 4 variables).
 \item A quadratic form~$Q$ is \defit[quadratic
form!indefinite]{indefinite} if $Q$ takes both positive and negative
values. For example, $x^2 - 3xy + y^2$ is indefinite, but $x^2 - 2xy +
y^2$ is definite \seeex{IndefEgExer}.
 \item A quadratic form~$Q$ in $n$~variables is \defit[quadratic
form!nondegenerate]{nondegenerate} if there does \emph{not} exist a
\emphit{nonzero} vector $x \in \real^n$, such that 
 $Q(v + x) = Q(v - x)$, for all $v \in \real^n$ \cfex{DegenExer}.
 \goodbreak
 \end{itemize}
 \end{defn}

 \begin{thm}[(Margulis)]
\label{Oppenheim}
 \index{Theorem!Margulis (on quadratic forms)}
 Let $Q$ be a real, indefinite, non-degenerate \term{quadratic form} in
$n \ge 3$~variables.

If $Q$ is not a scalar multiple of a form with integer coefficients, then
$Q(\integer^n)$ is dense in~$\real$.
 \end{thm}

\begin{eg} \label{OppAppl123}
 If $Q(x,y,z) = x^2 - \sqrt{2}xy + \sqrt{3}z^2$, then $Q$ is not a
scalar multiple of a form with integer coefficients \seeex{QnotZcoeffs}, so
Margulis' Theorem tells us that $Q(\integer^3)$ is dense in~$\real$.
 That is, for each $r \in \real$ and $\epsilon > 0$, there exist $a,b,c
\in \integer$, such that $|Q(a,b,c) - r| < \epsilon$.
 \end{eg}

\begin{rem} \label{OppRems}
 \ 
 \begin{enumerate}
 \item The hypothesis that $Q$ is \term[quadratic
form!indefinite]{indefinite} is necessary. If, say, $Q$ is positive
definite, then $Q(\integer^n) \subset \real^{\ge0}$ is \emph{not} dense in
all of~$\real$. In fact, if $Q$ is definite, then $Q(\integer^n)$ is
\emph{discrete} \seeex{Def->Disc}.
 \item \label{OppRems-2Var}
 There are counterexamples when $Q$~has only two variables
\seeex{OppEg2VarNotDense}, so the assumption that there are at least
3~variables cannot be omitted in general.
 \item A quadratic form is \term[quadratic form!degenerate]{degenerate}
if (and only if) a change of basis turns it into a form with less
variables. Thus, the \term[counterexample]{counterexamples} of
\pref{OppRems-2Var} show the assumption that $Q$~is \term[quadratic
form!nondegenerate]{nondegenerate} cannot be omitted in general
\seeex{OppEg3VarNotDense}.
 \item The converse of Thm.~\ref{Oppenheim} is true: if $Q(\integer^n)$ is
dense in~$\real$, then $Q$ cannot be a scalar multiple of a form with
integer coefficients \seeex{OppenConverse}.
 \end{enumerate}
 \end{rem}

Margulis' Theorem~\pref{Oppenheim} can be related to \term[Ratner's Theorems]{Ratner's Theorem} by considering the orthogonal group of the quadratic form~$Q$.

\begin{defn} \label{SOQDefn} \ 
  \begin{enumerate}
 \item If $Q$ is a quadratic form in $n$~variables, then
 \nindex{$\SO(Q)$ = orthogonal group of quadratic form~$Q$}$\SO(Q)$ is the
\defit[group!orthogonal (special)|indsee{$\SO(Q)$}]{orthogonal group}
(or \defit[group!isometry]{isometry group}) of~$Q$. That is,
 \index{SO(Q)*$\SO(Q)$|(}
 $$ \SO(Q) = \bigset{ h \in \SL(n,\real) }{ \text{$Q(v h) = Q(v)$ for all
$v \in \real^n$} } .$$
 \index{SL(l,R)*$\SL(\ell,\real)$}
 (Actually, this is the \emph{special} orthogonal group, because we are
including only the matrices of \index{determinant}determinant one.)
 \item As a special case,
 \nindex{$\SO(m,n)$ = orthogonal group of quadratic
form~$Q_{m,n}$}$\SO(m,n)$ is a shorthand for the orthogonal group
$\SO(Q_{m,n})$, where%
 \nindex{$Q_{m,n}$ = quadratic form $x_1^2 + \cdots + x_m^2 - x_{m+1}^2 -
\cdots - x_{m+n}^2$}
 $$Q_{m,n}(x_1,\ldots,x_{m+n}) = x_1^2 + \cdots + x_m^2 - x_{m+1}^2 -
\cdots - x_{m+n}^2 .$$
 \item Furthermore, we use \index{SO(Q)*$\SO(Q)$}$\SO(m)$ to denote
$\SO(m,0)$ (which is equal to $\SO(0,m)$).
 \end{enumerate}
 \end{defn}

\begin{defn}
 We use \nindex{$H^\circ$ = identity component of~$H$} $H^\circ$ to denote
the \defit[component!identity]{identity component} of a subgroup~$H$ of
$\SL(\ell,\real)$; that is, $H^\circ$ is the connected
\index{component!connected}component of~$H$ that contains the identity
element~$e$. It is a closed subgroup of~$H$.

Because \index{SO(Q)*$\SO(Q)$}$\SO(Q)$ is a real \index{algebraic
group!real}algebraic group \fullsee{EgsOfAlgGrps}{SO(Q)},
\term[Theorem!Whitney]{Whitney's Theorem} \pref{Zar->AlmConn} implies that
it has only finitely many \index{component!connected}components. (In fact,
it has only one or two components \seeexs{SOQconn}{SOQcomps}.) Therefore,
the difference between $\SO(Q)$ and $\SO(Q)^\circ$ is very minor, so it
may be ignored on a first reading.
 \end{defn}

\begin{proof}[{\bf Proof of \term[Theorem!Margulis (on quadratic
forms)]{Margulis' Theorem on values of quadratic forms}}]
 Let
 \begin{itemize}
 \item $G = \SL(3,\real)$, 
 \item \index{SL(3,Z)*$\SL(3,\integer)$}$\Gamma = \SL(3,\integer)$,
 \item $Q_0(x_1,x_2,x_3) = x_1^2 + x_2^2 - x_3^2$,
 and
 \item \index{SO(Q)*$\SO(Q)$}$H = \SO(Q_0)^\circ = \SO(2,1)^\circ$.
 \end{itemize}

  Let us assume $Q$ has exactly three variables (this causes no loss of
generality --- see Exer.~\ref{Opp3Var}). Then, because $Q$ is
\term[quadratic form!indefinite]{indefinite}, the
\term[signature (of a quadratic form)]{signature} of~$Q$ is either $(2,1)$
or $(1,2)$ \cfex{QFSignature}; hence, after a change of coordinates,
$Q$~must be a scalar multiple of~$Q_0$; thus, there exist $g \in
\SL(3,\real)$\index{SL(3,R)*$\SL(3,\real)$} and $\lambda \in
\real^{\times}$, such that
 $$Q = \lambda \, Q_0 \circ g .$$

Note that \index{SO(Q)*$\SO(Q)$}$\SO(Q)^\circ = g H g^{-1}$
\seeex{SOQ=gHg}. Because $H \approx \SL(2,\real)$ is \term{generated by
unipotent elements} \seeex{SO21=SL2R} and
$\SL(3,\integer)$\index{SL(3,Z)*$\SL(3,\integer)$} is a \term{lattice} in
$\SL(3,\real)$\index{SL(3,R)*$\SL(3,\real)$} \see{SLZLatt}, we can apply
\term[Ratner's Theorems!Orbit Closure]{Ratner's Orbit Closure Theorem}
\see{RatOrbH}. The conclusion is that there is a
\term[subgroup!connected]{connected subgroup}~$S$ of~$G$, such that
 \begin{itemize}
 \item $H \subset S$,
 \item the closure of $[g H]$ is equal to $[g S]$,
 and
 \item there is an \term[measure!invariant]{$S$-invariant probability
measure} on $[g S]$.
 \end{itemize}
 Algebraic calculations show that the only closed,
\term[subgroup!connected]{connected subgroups} of~$G$ that contain~$H$ are
the two obvious subgroups: $G$ and~$H$ \seeex{ContainSO21}. Therefore, $S$
must be either $G$ or~$H$. We consider each of these possibilities
separately.

\setcounter{case}{0}

\begin{case}
 Assume $S = G$.
  \end{case}
 This implies that 
 \begin{equation} \label{OppenPf-GamgHDenseG}
 \mbox{$\Gamma g H$ is dense in~$G$.}
 \end{equation}
 We have
 \begin{align*}
  Q(\integer^3) &= Q_0( \integer^3 g) && \mbox{(definition of~$g$)} \\
  & = Q_0(\integer^3 \Gamma g) && \mbox{($\integer^3 \Gamma =
\integer^3$)} \\
 &= Q_0(\integer^3 \Gamma g H) && \mbox{(definition of~$H$)} \\
 &\simeq Q_0(\integer^3 G) && \mbox{(\pref{OppenPf-GamgHDenseG} and $Q_0$
is continuous)} \\
 & = Q_0(\real^3 \smallsetminus \{0\}) && \mbox{($vG = \real^3
\smallsetminus \{0\}$ for $v \neq 0$)} \\
 &= \real
 , \end{align*}
 where ``$\simeq$" means ``is dense in."

\begin{case}
 Assume $S = H$.
 \end{case}
 This is a degenerate case; we will show that $Q$ is a scalar multiple of
a form with integer coefficients. To keep the proof short, we will apply
some of the theory of \term[algebraic group!theory of]{algebraic groups}.
The interested reader may consult Chapter~\ref{AlgGrpsChap} to fill in the
gaps.

Let $\Gamma_g = \Gamma \cap (g H g^{-1})$. Because the orbit $[g H] = [g
S]$ has finite $H$-\term[measure!invariant]{invariant measure}, we know
that $\Gamma_g$ is a \term{lattice} in~\index{SO(Q)*$\SO(Q)$}$g H g^{-1} =
\SO(Q)^\circ$. So the \term[Theorem!Borel Density]{Borel Density Theorem}
\pref{BDT} implies $\SO(Q)^\circ$ is contained in the
\term[Zariski!closure]{Zariski closure} of~$\Gamma_g$. Because $\Gamma_g
\subset \Gamma = \SL(3,\integer)$\index{SL(3,Z)*$\SL(3,\integer)$}, this
implies that the (almost) \index{algebraic group!real}algebraic group
$\SO(Q)^\circ$ is \index{algebraic group!defined over Q*defined over
$\rational$}defined over~$\rational$ \seeex{Zpts->def/Q}. Therefore, up to
a scalar multiple, $Q$ has integer coefficients \seeex{EgsOverQ-SO(Q)Ex}.
 \end{proof}

\begin{exercises}

\item 
 Suppose $\alpha$ and~$\beta$ are nonzero real numbers, such that
$\alpha/\beta$ is \term[irrational!number]{irrational}, and define $L(x,y)
= \alpha x + \beta y$. Show $L(\integer^2)$ is dense in~$\real$.
(\term[Theorem!Margulis (on quadratic forms)]{Margulis' Theorem}
\pref{Oppenheim} is a generalization to \term[quadratic form]{quadratic
forms} of this rather trivial observation about \term[linear!form]{linear
forms}.)

\item \label{IndefEgExer}
 Let $Q_1(x,y) = x^2 - 3xy + y^2$ and $Q_2(x,y) = x^2 - 2xy + y^2$. Show
\begin{enumerate}
 \item $Q_1(\real^2)$ contains both positive and negative numbers, 
 but
 \item $Q_2(\real^2)$ does not contain any negative numbers. 
 \end{enumerate}

\item \label{DegenExer}
 Suppose $Q(x_1,\ldots,x_n)$ is a \term{quadratic form}, and let $e_n =
(0,\ldots,0,1)$ be the $n^{\text{th}}$~standard basis vector. Show 
 $$ \text{$Q(v + e_n) = Q(v - e_n)$ for all $v \in \real^n$} $$
 if and only  if there is a quadratic form $Q'(x_1,\ldots,x_{n-1})$ in
$n-1$~variables, such that 
 $Q(x_1,\ldots,x_n) = Q'(x_1,\ldots,x_{n-1})$ for all $x_1,\ldots,x_n \in
\real$.

\item \label{QnotZcoeffs}
 Show that the \term[quadratic form]{form}~$Q$ of Eg.~\ref{OppAppl123} is
not a scalar multiple of a form with integer coefficients; that is, there
does not exist $k \in \real^\times$, such that all the coefficients of $k
Q$ are integers.

\item Suppose $Q$ is a quadratic form in $n$~variables. Define 
 $$
 \text{$B \colon \real^n \times \real^n \to \real$ by }
 B(v,w) = \frac{1}{4} \bigl( Q(v + w) - Q(v - w) \bigr) .$$
 \vskip-\smallskipamount
 \begin{enumerate}
 \item Show that $B$ is a \term[bilinear form!symmetric]{symmetric bilinear
form} on~$\real^n$. That is, for $v, v_1,v_2, w \in \real^n$ and $\alpha
\in \real$, we have:
 \begin{enumerate}
 \item $B(v,w) = B(w,v)$
 \item $B(v_1 + v_2, w) = B(v_1,w) + B(v_2,w)$,
 and
 \item $B(\alpha v, w) = \alpha \, B(v,w)$.
 \end{enumerate}
 \item For $h \in \SL(n,\real)$\index{SL(l,R)*$\SL(\ell,\real)$}, show $h
\in \SO(Q)$ if and only if $B( vh, wh ) = B(v,w)$ for all $v,w \in
\real^n$.
 \item  We say that the bilinear form $B$~is \defit[bilinear
form!nondegenerate]{nondegenerate} if for every nonzero $v \in \real^n$,
there is some nonzero $w \in \real^n$, such that $B(v,w) \neq 0$.
 Show that $Q$ is \term[quadratic form!nondegenerate]{nondegenerate}
if and only if $B$~is nondegenerate.
 \item For $v \in \real^n$, let 
 \nindex{$v^\perp$ = orthogonal complement of vector~$v$}
 $v^\perp = \{\, w \in \real^n \mid B(v,w) = 0 \,\} $.
 Show:
 \begin{enumerate}
 \item $v^\perp$ is a subspace of~$\real^n$, 
 and
 \item if $B$ is nondegenerate and $v \neq 0$, then $\real^n = \real v \oplus
v^\perp$.
 \end{enumerate}
 \end{enumerate}

\item \label{QFSignature}
  \begin{enumerate}
 \item Show that $Q_{k,n-k}$ is a \term[quadratic
form!nondegenerate]{nondegenerate quadratic form} (in $n$~variables).
 \item Show that $Q_{k,n-k}$ is \term[quadratic
form!indefinite]{indefinite} if and only if $i \notin \{0,n\}$.
 \item A subspace~$V$ of~$\real^n$ is \defit[totally isotropic
subspace]{totally isotropic} for a quadratic form~$Q$ if $Q(v) = 0$ for
all $v \in V$. Show that $\min(k,n-k)$ is the maximum
dimension\index{dimension!of a vector space} of a totally isotropic
subspace for~$Q_{k,n-k}$.
 \item \label{QFSignature-sign}
 Let $Q$ be a \term[quadratic form!nondegenerate]{nondegenerate quadratic
form} in $n$~variables. Show there exists a unique $k \in
\{0,1,\ldots,n\}$, such that there is an invertible linear
transformation~$T$ of~$\real^n$ with $Q
= Q_{k,n-k} \circ T$.
 We say that the \defit[signature (of a quadratic form)]{signature} of~$Q$
is $(k, n-k)$.
 \end{enumerate}
 \hint{\pref{QFSignature-sign}~Choose $v \in \real^n$ with $Q(v) \neq 0$.
By induction on~$n$, the restriction of~$Q$ to~$v^\perp$ can be
transformed to $Q_{k', n-1-k'}$.}

\item \label{Def->Disc}
 Let $Q$ be a real \term{quadratic form} in $n$~variables. 
 Show that if $Q$ is \defit[quadratic form!positive definite]{positive
definite} (that is, if $Q(\real^n) \ge 0$), then $Q(\integer^n)$ is a
\term[set!discrete]{\emph{discrete} subset of~$\real$}. 

\item \label{OppEg2VarNotDense}
 Show:
 \begin{enumerate}
 \item If $\alpha$ is an \term[irrational!number]{irrational} root of a
\index{function!polynomial!quadratic}quadratic polynomial (with integer
coefficients), then there exists $\epsilon > 0$, such that 
 $$|\alpha - (p/q)| > \epsilon/pq ,$$
 for all $p,q \in \integer$ (with $p,q \neq 0$).
 \hint{$k(x - \alpha)(x - \beta)$ has integer coefficients, for some $k
\in \integer^+$ and some $\beta \in \real \smallsetminus \{\alpha\}$.}
 \item The \term{quadratic form} $Q(x,y) = x^2 - (3 + 2 \sqrt{2}) y^2$ is
real, indefinite, and nondegenerate, and is \emph{not} a scalar multiple
of a form with integer coefficients.
 \item $Q(\integer, \integer)$ is not dense in~$\real$.
 \hint{$\sqrt{3 + 2 \sqrt{2}} = 1 + \sqrt{2}$ is a root of a quadratic
polynomial.}
 \end{enumerate}

\item \label{OppEg3VarNotDense}
 Suppose $Q(x_1,x_2)$ is a real, indefinite \term{quadratic form} in two
variables, and that $Q(x,y)$ is \emph{not} a scalar multiple of a
form with integer coefficients, and define $Q^*(y_1,y_2,y_3) = Q(y_1, y_2
- y_3)$.
 \begin{enumerate}
 \item Show that $Q^*$ is a real, indefinite quadratic form in two
variables, and that $Q^*$ is \emph{not} a scalar multiple of a
form with integer coefficients.
 \item Show that if $Q(\integer^2)$ is not dense in~$\real$, then
$Q^*(\integer^3)$ is not dense in~$\real$.
 \end{enumerate}

\item \label{OppenConverse}
 Show that if $Q(x_1,\ldots,x_n)$ is a \term{quadratic form}, and
$Q(\integer^n)$ is dense in~$\real$, then $Q$ is \emph{not} a scalar
multiple of a form with integer coefficients.

\item \label{SOQconn}
 Show that \index{SO(Q)*$\SO(Q)$}$\SO(Q)$ is
\index{subgroup!connected|(}connected if $Q$~is definite.
 \hint{Induction on~$n$. There is a natural embedding of $\SO(n-1)$ in
$\SO(n)$, such that the vector $e_n = (0,0,\ldots,0,1)$ is
\index{stabilizer!of a vector}fixed by $\SO(n-1)$. For $n \ge 2$, the map
$\SO(n-1) g \mapsto e_n g$ is a \index{homeomorphism}homeomorphism from
$\SO(n-1) \backslash \SO(n)$ onto the $(n-1)$-sphere~$S^{n-1}$.}

\item \label{WittThm}
 (\term[Theorem!Witt]{Witt's Theorem})
 Suppose $v,w \in \real^{m+n}$ with $Q_{m,n}(v) = Q_{m,n}(w) \neq 0$, and
assume $m+n \ge 2$. Show there exists \index{SO(Q)*$\SO(Q)$}$g \in
\SO(m,n)$ with $v g = w$.
 \hint{There is a linear map $T \colon v^\perp \to w^\perp$ with
$Q_{m,n}(x T) = Q_{m,n}(x)$ for all~$x$ \seeex{QFSignature}. (Use the
assumption $m + n \ge 2$ to arrange for $g$ to have
\index{determinant}determinant~$1$, rather than~$-1$.)}

\item \label{SOQcomps}
 Show that \index{SO(Q)*$\SO(Q)$}$\SO(m,n)$ has no more than two
\index{component!connected}components if $m,n \ge 1$. (In fact, although
you do not need to prove this, it has exactly two components.)
 \hint{Similar to Exer.~\ref{SOQconn}. (Use Exer.~\ref{WittThm}.) If $m >
1$, then $\{\, v \in \real^{m+n} \mid Q_{m,n} = 1 \,\}$ is connected. The
base case $m = n = 1$ should be done separately.}

\item \label{SOQ=gHg}
 In the notation of the proof of Thm.~\ref{Oppenheim}, show
\index{SO(Q)*$\SO(Q)$}$\SO(Q)^\circ = g H g^{-1}$.

\item \label{Opp3Var}
 Suppose $Q$ satisfies the hypotheses of Thm.~\ref{Oppenheim}. Show there
exist $v_1,v_2,v_3 \in \integer^n$, such that the quadratic form $Q'$
on~$\real^3$, defined by
 $ Q'(x_1,x_2,x_3) = Q(x_1 v_1 + x_2 v_2 + x_3 v_3 )$, 
 also satisfies the hypotheses of Thm.~\ref{Oppenheim}.
 \hint{Choose any $v_1,v_2$ such that $Q(v_1)/Q(v_2)$ is negative and
irrational. Then choose $v_3$ generically (so $Q'$ is nondegenerate).}

\index{quadratic form|)}

\item \label{SO21=SL2R}
 (\emphit{Requires some Lie theory}) Show:
 \begin{enumerate}
 \item The \index{determinant}determinant function $\det$ is a quadratic
form on $\LieSL(2,\real)$ of signature $(2,1)$.
 \item The adjoint \index{representation!adjoint}representation
$\Ad_{\SL(2,\real)}$ maps $\SL(2,\real)$ into
\index{SO(Q)*$\SO(Q)$}$\SO(\det)$.
 \item $\SL(2,\real)$ is locally isomorphic to
\index{SO(Q)*$\SO(Q)$}$\SO(2,1)^\circ$.
 \item $\SO(2,1)^\circ$ is generated by unipotent elements.
 \end{enumerate}

\item \label{ContainSO21}
 (\emphit{Requires some Lie theory})
 \begin{enumerate}
 \item Show that $\so(2,1)$ is a \term{maximal subalgebra} of the \term{Lie
algebra} $\LieSL(3,\real)$. That is, there does not exist a subalgebra
$\Lie H$ with $\so(2,1) \subsetneq \Lie H \subsetneq \LieSL(3,\real)$.
 \item Conclude that if $S$ is any closed,
\term[subgroup!connected]{connected subgroup} of
$\SL(3,\real)$\index{SL(3,R)*$\SL(3,\real)$} that contains
\index{SO(Q)*$\SO(Q)$|)}$\SO(2,1)$, then 
 $$ \text{either $S = \SO(2,1)$ or $S =
\SL(3,\real)$.} $$
 \end{enumerate}
 \hint{$ \lie u = \begin{bmatrix} 0 & 1 & 1 \\ -1 & 0 & 0 \\ 1 & 0 & 0
\end{bmatrix}$
 is a \index{matrix!nilpotent}nilpotent element of $\so(2,1)$, and the
kernel of $\ad_{\LieSL(3,\real)} \lie u$ is only
$2$-dimensional\index{dimension!of a vector space}. Since $\Lie H$ is a
submodule of $\LieSL(3,\real)$, the conclusion follows
\seeex{HighWtGen-VinW}.}

\end{exercises}

\section{Measure-theoretic versions of Ratner's Theorem}
\label{MeasVersSect}

For unipotent flows, \term[Ratner's Theorems!Orbit Closure]{Ratner's
Orbit Closure Theorem} \pref{RatnerOrb} states that the closure of each
orbit is a nice, geometric subset~$[x S]$ of the space $X = \Gamma
\backslash G$.
 This means that the orbit is \term[dense orbit]{dense} in~$[x S]$; in
fact, it turns out to be \defit{uniformly distributed} in~$[x S]$.
 Before making a precise statement, let us look at a simple example.

\begin{eg} \label{T2UnifDistEg}
  As in Eg.~\ref{RatnerOrbForTn}, let $\varphi_t$ be the flow
 $$ \varphi_t \bigl( [x] \bigr) = [x + t v] $$
 on~\term[torus!Tn*$\torus^n$]{$\torus^n$} defined by a vector $v \in
\real^n$.
 Let $\mu$ be the \term[measure!Lebesgue]{Lebesgue measure} on~$\torus^n$,
normalized to be a \defit[measure!probability]{probability measure} (so
$\mu(\torus^n) = 1$).
 \begin{enumerate}
 \item \label{T2UnifDistEg-measure}
 Assume $n = 2$, so we may write $v = (a,b)$. If $a/b$ is
\term[irrational!number]{irrational}, then every orbit of~$\varphi_t$ is
dense in~$\torus^2$ \seeex{T2homogOrbits}. In fact, every orbit is
\term{uniformly distributed} in~$\torus^2$: if $B$ is any nice open subset
of~$\torus^2$ (such as an open ball), then the amount of time that each
orbit spends in~$B$ is proportional to the area of~$B$. More precisely,
for each $x \in \torus^2$, and letting~$\lambda$ be the
\term[measure!Lebesgue]{Lebesgue measure} on~$\real$, we have
 \begin{equation} \label{T2UnifDistEg-Open}
 \frac{ \lambda \left(
 \bigset{t \in [0,T]}{\varphi_t(x) \in B}
 \right)}{T}
 \to \mu(B)
 \qquad \mbox{as $T \to \infty$}
 \end{equation}
 \seeex{T2UnifDistOpen}.
  \item \label{T2UnifDistEg-integral}
 Equivalently, if 
 \begin{itemize}
 \item $v = (a,b)$ with $a/b$ irrational, 
 \item $x \in \torus^2$,
 and 
 \item $f$~is any \term[function!continuous]{continuous function}
on~$\torus^2$,
 \end{itemize}
 then 
 \begin{equation} \label{T2UnifDistEg-Func}
 \lim_{T \to \infty} \frac{ \int_0^T f \bigl( \varphi_t(x) \bigr) \, dt }{
T } \to \int_{\torus^2} f\, d \mu
 \end{equation}
 \seeex{T2UnifDistFunc}.

 \item \label{T2UnifDistEg-NotUnif}
 Suppose now that $n = 3$, and assume $v = (a,b,0)$, with $a/b$
\term[irrational!number]{irrational}. Then the orbits of $\varphi_t$ are
\emph{not} dense in~$\torus^3$, so they are \emph{not} \term{uniformly
distributed} in~$\torus^3$ (with respect to the usual
\term[measure!Lebesgue]{Lebesgue measure} on~$\torus^3$). Instead, each
orbit is uniformly distributed in some subtorus of~$\torus^3$: given $x =
(x_1,x_2,x_3) \in \torus$, let $\mu_2$ be the \term[measure!Haar]{Haar
measure} on the horizontal $2$-torus $\torus^2 \times \{x_3\}$ that
contains~$x$. Then
 $$ \frac{1}{T} \, \int_0^T f \bigl( \varphi_t(x) \bigr) \, dt \to
\int_{\torus^2 \times \{x_3\}} f\, d \mu_2 \text{\qquad as $T \to \infty$}
$$ \seeex{T3NotUnifDistFunc}. 
 \item \label{T2UnifDistEg-Tn}
 In general, for any~$n$ and~$v$, and any $x \in \torus^n$, there is
a subtorus~$S$ of~$\torus^n$, with
\term[measure!Haar]{Haar measure}~$\mu_S$, such that
 $$ \int_0^T f \bigl( \varphi_t(x) \bigr) \, dt \to \int_{S} f\, d \mu_S$$
as $T \to \infty $
 \seeex{RatnerForTn}. 
 \end{enumerate}
 \end{eg}

The above example generalizes, in a natural way, to all
\term[unipotent!flow]{unipotent flows}:

\begin{thm}[{(\term[Ratner's Theorems!Equidistribution]{Ratner
Equidistribution Theorem})}] \label{RatnerUnifDist}
 If
 \begin{itemize}
 \item $G$ is any \term{Lie group},
 \item $\Gamma$ is any \term{lattice} in~$G$,
 and
 \item $\varphi_t$ is any \term[unipotent!flow]{unipotent flow} on~$\Gamma
\backslash G$,
 \end{itemize}
 then each $\varphi_t$-orbit is \term{uniformly distributed} in its
closure.
 \end{thm}

\begin{rem} \label{RatnerUnifDistPrecise}
 Here is a more precise statement of Thm.~\ref{RatnerUnifDist}.
 For any fixed $x \in G$, Ratner's Theorem \pref{RatnerOrb} provides a
\term[subgroup!connected|)]{connected, closed subgroup}~$S$ of~$G$
\see{RatnerOrbPrecise}, such that
 \begin{enumerate}
 \item $\{u^t\}_{t \in \real} \subset S$,
 \item the image $[x S]$ of $x S$ in $\Gamma \backslash G$ is closed, and
has \term[finite volume]{finite $S$-invariant volume},
 and
 \item the $\varphi_t$-orbit of~$[x]$  is dense in $[x S]$.
 \end{enumerate}
 Let \nindex{$\mu_S$ = $S$-invariant probability measure on an
$S$-orbit}$\mu_S$ be the (unique) $S$-\term[measure!invariant]{invariant
probability measure} on $[x S]$. Then \term[Ratner's Theorems!Equidistribution]{Thm.~\ref{RatnerUnifDist}} asserts, for every
\term[function!continuous]{continuous function}~$f$ on~$\Gamma \backslash
G$ with compact \term[support!of a function]{support}, that
 $$ \frac{1}{T} \int_0^T f \bigl( \varphi_t(x) \bigr) \, dt \to \int_{[x
S]} f\, d \mu_S
 \qquad \mbox{as $T \to \infty$.} $$
 \end{rem}

This theorem yields a classification of the
$\varphi_t$-\term[measure!invariant]{invariant probability measures}.

\begin{defn}
 Let
 \begin{itemize}
 \item $X$ be a metric space,
 \item $\varphi_t$ be a \term[flow!continuous]{continuous flow} on~$X$,
 and
 \item $\mu$ be a measure on~$X$.
 \end{itemize}
 We say:
 \begin{enumerate}
 \item $\mu$ is \defit[measure!invariant]{$\varphi_t$-invariant} if $\mu
\bigl( \varphi_t(A) \bigr) = \mu(A)$, for every Borel subset~$A$ of~$X$,
and every $t \in \real$.
 \item $\mu$ is \defit[ergodic!measure|(]{ergodic} if $\mu$ is
$\varphi_t$-\term[measure!invariant]{invariant}, and every
$\varphi_t$-\term[function!invariant]{invariant} Borel function on~$X$ is
essentially constant (w.r.t.~$\mu$). (A function~$f$ is
\defit[function!constant (essentially)]{essentially constant} on~$X$ if
there is a set~$E$ of measure~$0$, such that $f$~is constant on $X
\smallsetminus E$.)
 \end{enumerate}
 \end{defn}

Results of \term{Functional Analysis} (such as
\term[Theorem!Choquet]{Choquet's Theorem}) imply that every
\term[measure!invariant]{invariant probability measure} is a
\term[measure!convex combination]{convex combination} (or, more generally,
a \term[measure!direct integral]{direct integral}) of
\term[ergodic!measure]{ergodic probability measures} \seeex{ChoquetEx}.
(See \S\ref{ErgDecompSect} for more discussion of the relationship between
arbitrary measures and ergodic measures.) Thus, in order to understand all
of the \term[measure!invariant]{invariant} measures, it suffices to
classify the \term[ergodic!measure]{ergodic} ones. Combining
\term[Ratner's Theorems!Equidistribution]{Thm.~\ref{RatnerUnifDist}} with
the \term[Theorem!Pointwise Ergodic]{Pointwise Ergodic Theorem}
\pref{PtwiseErgThmFlow} implies that these \term[ergodic!measure]{ergodic}
measures are of a nice geometric form \seeex{RatnerMeasEx}:

 \begin{cor}[{(\term[Ratner's Theorems!Measure Classification]{Ratner
Measure Classification Theorem})}] \label{RatnerMeas}
 \index{Ratner's Theorems!Measure Classification|(}
 If
 \begin{itemize}
 \item $G$ is any \term{Lie group},
 \item $\Gamma$ is any \term{lattice} in~$G$,
 and
 \item $\varphi_t$ is any \term[unipotent!flow]{unipotent flow} on~$\Gamma
\backslash G$,
 \end{itemize}
 then every \term[ergodic!measure]{ergodic}
$\varphi_t$-\term[measure!invariant]{invariant probability measure} on
$\Gamma \backslash G$ is \term[measure!homogeneous]{homogeneous}.

 That is, every \term[ergodic!measure]{ergodic}
$\varphi_t$-\term[measure!invariant]{invariant} probability measure is of
the form~$\mu_S$, for some~$x$ and some subgroup~$S$ as in
Rem.~\ref{RatnerUnifDistPrecise}.
 \end{cor}

A logical development (and the historical development) of the material
proceeds in the opposite direction: instead of deriving
Cor.~\ref{RatnerMeas} from Thm.~\ref{RatnerUnifDist}, the main goal of
these lectures is to explain the main ideas in a direct proof of
\term[Ratner's Theorems!Measure Classification]{Cor.~\ref{RatnerMeas}}. Then
Thms.~\term[Ratner's Theorems!Orbit Closure]{\ref{RatnerOrb}}
and~\term[Ratner's Theorems!Equidistribution]{\ref{RatnerUnifDist}} can be
obtained as corollaries. As an illustrative example of this opposite
direction --- how knowledge of \term[measure!invariant]{invariant
measures} can yield information about closures of orbits --- let us prove
the following classical fact. (A more complete discussion appears in
Sect.~\ref{MeasuresToOrbitsSect}.)

\begin{defn} \label{UniqErgDefn}
 Let $\varphi_t$ be a \term[flow!continuous]{continuous flow} on a metric
space~$X$.
 \begin{itemize}
 \item $\varphi_t$ is \defit[flow!minimal]{minimal} if every orbit is
dense in~$X$.
 \item $\varphi_t$ is \defit[ergodic!uniquely]{uniquely ergodic} if there
is a \emph{unique} $\varphi_t$-\term[measure!invariant]{invariant}
probability measure on~$X$.
 \end{itemize}
 \end{defn}

\begin{prop} \label{UniqErg->Minl}
  Suppose
 \begin{itemize}
 \item $G$ is any \term{Lie group},
 \item $\Gamma$ is any \term{lattice} in~$G$, such that $\Gamma \backslash
G$ is \emph{compact},
 and
 \item $\varphi_t$ is any \term[unipotent!flow]{unipotent flow} on~$\Gamma
\backslash G$.
 \index{flow!unipotent|indsee{unipotent~flow}}
 \end{itemize}
 If $\varphi_t$ is \term[ergodic!uniquely]{uniquely ergodic}, then
$\varphi_t$ is \term[flow!minimal]{minimal}.
 \end{prop}

\begin{proof}
 We prove the contrapositive: assuming that some orbit $\varphi_{\real}(x)$
is not dense in $\Gamma \backslash G$, we will show that the
$G$-\term[measure!invariant]{invariant measure} $\mu_G$ is not the only
$\varphi_t$-\term[measure!invariant]{invariant} probability measure on
$\Gamma \backslash G$.

Let $\Omega$ be the closure of $\varphi_{\real}(x)$. Then $\Omega$ is a
compact $\varphi_t$-\term[set!invariant]{invariant subset} of $\Gamma
\backslash G$ \seeex{InvariantClosure}, so there is a
$\varphi_t$-\term[measure!invariant]{invariant} probability measure~$\mu$
on $\Gamma \backslash G$ that is \term[support!of a measure]{supported} on~$\Omega$
\seeex{MeasOnInvSet}. Because 
 $$ \supp \mu \subset \Omega \subsetneq \Gamma \backslash G = \supp \mu_G,
$$
 we know that $\mu \neq \mu_G$. Hence, there are (at least) two different
$\varphi_t$-\term[measure!invariant]{invariant} probability measures on
$\Gamma \backslash G$, so $\varphi_t$ is not
\term[ergodic!uniquely]{uniquely ergodic}.
 \end{proof}

\begin{rem} \label{RatMeasThmRems} \ 
 \begin{enumerate}
 \item \label{RatMeasThmRems-notlatt}
 There is no need to assume $\Gamma$ is a \term{lattice} in
Cor.~\ref{RatnerMeas} --- the conclusion remains true when $\Gamma$ is any
closed subgroup of~$G$. However, to avoid confusion, let us point out that
this is \emph{not} true of the \term[Ratner's Theorems!Orbit Closure]{Orbit Closure Theorem} --- there are
\term[counterexample]{counterexamples} to \pref{RatOrbH} in some cases
where $\Gamma \backslash G$ is not assumed to have finite volume.
 For example, a \term{fractal} orbit closure for $a^t$ on $\Gamma
\backslash G$ yields a \term{fractal} orbit closure for $\Gamma$ on
$\{a^t\} \backslash G$, even though the \term{lattice}~$\Gamma$ may be
generated by \term[unipotent!element]{unipotent elements}.
 \item An appeal to \term[Ratner's Theorems]{``Ratner's Theorem"} in the
literature could be referring to any of Ratner's three major
theorems: her Orbit Closure Theorem \pref{RatnerOrb}, her
Equidistribution Theorem \pref{RatnerUnifDist}, or her
Measure Classification Theorem \pref{RatnerMeas}.
 \item There is not universal agreement on the names of these three major
\term[Ratner's Theorems]{theorems of Ratner}. For example, the
Measure Classification Theorem is also known as ``Ratner's
Measure-Rigidity Theorem" or ``Ratner's Theorem on Invariant Measures," and
the Orbit Closure Theorem is also known as the ``topological version" of
her theorem.
 \item Many authors (including M.~Ratner) use the adjective
\defit[measure!algebraic]{algebraic}, rather than
\term[measure!homogeneous]{\emph{homogeneous}}, to describe
measures~$\mu_S$ as in \pref{RatnerUnifDistPrecise}. This is because
$\mu_S$ is defined via an algebraic (or, more precisely, group-theoretic)
construction.
 \end{enumerate}
 \end{rem}

\begin{exercises}

\item \label{T2UnifDistOpen}
 Verify Eg.~\fullref{T2UnifDistEg}{measure}.
 \hint{It may be easier to do Exer.~\ref{T2UnifDistFunc} first. The
characteristic function of~$B$ can be approximated by
\term[function!continuous]{continuous functions}.}

\item \label{T2UnifDistFunc}
 Verify Eg.~\fullref{T2UnifDistEg}{integral}; show that if $a/b$ is
\term[irrational!number]{irrational}, and $f$ is any
\term[function!continuous]{continuous function} on~$\torus^2$, then
\pref{T2UnifDistEg-Func} holds.
 \hint{ Linear combinations of functions of the form
 $$f(x,y) = \exp2\pi (m x + n y) i$$
 are dense in the space of
\term[function!continuous]{continuous functions}.
 \\ \emphit{Alternate solution:} If $T_0$ is sufficiently large, then, for
every $x \in \torus^2$, the segment $\{\varphi_t(x)\}_{t=0}^{T_0}$ comes
within~$\delta$ of every point in~$\torus^2$ (because $\torus^2$ is
compact and \index{group!abelian}abelian, and the orbits of~$\varphi_t$ are
dense). Therefore, the uniform continuity of~$f$ implies that if $T$ is
sufficiently large, then the value of $(1/T) \, \int_0^T f \bigl(
\varphi_t(x) \bigr) \, dt$ varies by less than~$\epsilon$ as $x$~varies
over~$\torus^2$.} 

\item \label{T3NotUnifDistFunc}
 Verify Eg.~\fullref{T2UnifDistEg}{NotUnif}.

\item \label{RatnerForTn}
 Verify Eg.~\fullref{T2UnifDistEg}{Tn}.

\item \label{erg<>extpt}
 Let
 \begin{itemize}
 \item $\varphi_t$ be a \term[flow!continuous]{continuous flow} on a
\index{manifold}manifold~$X$,
 \item $\Prob(X)_{\varphi_t}$ be
the set of $\varphi_t$-\term[measure!invariant]{invariant Borel probability
measures} on~$X$,
 and
 \item $\mu \in \Prob(X)_{\varphi_t}$.
 \end{itemize}
 Show that the following are equivalent:
 \begin{enumerate}
 \item \label{erg<>extpt-erg}
 $\mu$ is \term[ergodic!measure]{ergodic};
 \item \label{erg<>extpt-subset}
 every $\varphi_t$-\term[set!invariant]{invariant Borel subset} of~$X$ is
either null or conull;
 \item \label{erg<>extpt-extpt}
 $\mu$ is an \defit[point!extreme]{extreme point} of
$\Prob(X)_{\varphi_t}$, that is, $\mu$ is \emph{not} a
\term[measure!convex combination]{convex combination} of two other
measures in the space $\Prob(X)_{\varphi_t}$.
 \end{enumerate}
 \hint{(\ref{erg<>extpt-erg}$\Rightarrow$\ref{erg<>extpt-extpt})~If $\mu
= a_1 \mu_1 + a_2\mu_2$, consider the \term[Radon-Nikodym
derivative]{Radon-Nikodym derivatives} of $\mu_1$ and~$\mu_2$
(w.r.t.~$\mu$).
 (\ref{erg<>extpt-extpt}$\Rightarrow$\ref{erg<>extpt-subset})~If $A$~is
any subset of~$X$, then $\mu$ is the sum of two measures, one
\term[support!of a measure]{supported} on $A$, and the other
\term[support!of a measure]{supported} on the complement of~$A$.}

\item \label{ChoquetEx}
 \term[Theorem!Choquet]{Choquet's Theorem} states that if $C$ is any
compact subset of a Banach space, then each point in~$C$ is of the form
$\int_C c \, d\mu(c)$, where $\nu$ is a probability measure
\term[support!of a measure]{supported} on the \term[point!extreme]{extreme
points} of~$C$. Assuming this fact, show that every
$\varphi_t$-\term[measure!invariant]{invariant probability measure} is an
\term[measure!direct integral]{integral} of
\term[ergodic!measure|)]{ergodic}
$\varphi_t$-\term[measure!invariant]{invariant measures}.

\item \label{RatnerMeasEx}
 Prove \term[Ratner's Theorems!Measure
Classification|)]{Cor.~\ref{RatnerMeas}}.
 \hint{Use \pref{RatnerUnifDist} and \pref{PtwiseErgThmFlow}.}

\item \label{InvariantClosure}
 Let
 \begin{itemize}
 \item $\varphi_t$ be a \term[flow!continuous]{continuous flow} on a
metric space~$X$,
 \item $x \in X$,
 and
 \item $\varphi_{\real}(x) = \{\, \varphi_t(x) \mid t \in \real \,\}$ be
the \defit{orbit} of~$x$.
 \end{itemize}
 Show that the closure $\closure{\varphi_{\real}(x)}$ of
$\varphi_{\real}(x)$ is $\varphi_t$-\defit[set!invariant]{invariant}; that
is, $\varphi_t \bigl( \closure{\varphi_{\real}(x)} \bigr) =
\closure{\varphi_{\real}(x)}$, for all $t \in \real$.

\item \label{MeasOnInvSet}
  Let
 \begin{itemize}
 \item $\varphi_t$ be a \term[flow!continuous]{continuous flow} on a
metric space~$X$,
 and
 \item $\Omega$ be a nonempty, compact,
$\varphi_t$-\term[set!invariant]{invariant subset} of~$X$.
 \end{itemize}
 Show there is a $\varphi_t$-\term[measure!invariant]{invariant probability
measure}~$\mu$ on~$X$, such that $\supp(\mu) \subset \Omega$. (In other
words, the complement of~$\Omega$ is a null set, w.r.t.~$\mu$.)
 \hint{Fix $x \in \Omega$. For each $n \in \integer^+$, $(1/n)
\, \int_0^n f \bigl( \varphi_t(x) \bigr) \, dt$ defines a probability
measure~$\mu_n$ on~$X$. The limit of any convergent subsequence is
$\varphi_t$-\term[measure!invariant]{invariant}.} 

\item Let 
 \begin{itemize}
 \item $S^1 = \real \cup \{\infty\}$ be the one-point
compactification of~$\real$,
 and
 \item $\varphi_t(x) = x + t$ for $t \in \real$ and $x \in S^1$.
 \end{itemize}
 Show $\varphi_t$ is a flow on~$S^1$ that is
\term[ergodic!uniquely]{uniquely ergodic} (and
\term[flow!continuous]{continuous}) but not \term[flow!minimal]{minimal}. 

\item Suppose $\varphi_t$ is a uniquely \term[ergodic!flow]{ergodic},
continuous flow on a compact metric space~$X$.
 Show $\varphi_t$ is \term[flow!minimal]{minimal} if and only if there is a
$\varphi_t$-\term[measure!invariant]{invariant} probability measure~$\mu$
on~$X$, such that the \term[support!of a measure]{support} of~$\mu$ is all of~$X$.

\item \label{NoMeasOnInvSet}
 Show that the conclusion of
Exer.~\ref{MeasOnInvSet} can fail if we omit the hypothesis that $\Omega$
is compact.
 \hint{Let $\Omega = X = \real$, and define $\varphi_t(x) = x + t$.}

 \end{exercises}

\section{Some applications of Ratner's Theorems} \label{ApplSect}

This section briefly describes a few of the many results that rely
crucially on \term{Ratner's Theorems} (or the methods behind them).
 Their proofs require substantial new ideas, so, although we will
emphasize the role of Ratner's Theorems, we do not mean to imply that any
of these theorems are merely corollaries.

\subsection{Quantitative versions of Margulis' Theorem on values of
quadratic forms}

As discussed in \S\ref{MargOppSect}, \term[Theorem!Margulis
(on quadratic forms)]{G.~A.~Margulis} proved, under appropriate hypotheses
on the \term{quadratic form}~$Q$, that the values of~$Q$ on~$\integer^n$
are dense in~$\real$. By a more sophisticated argument, it can be shown
(except in some small cases) that the values are \term{uniformly
distributed} in~$\real$, not just dense:

\begin{thm} \label{QuantOppThm}
 Suppose
 \begin{itemize}
 \item $Q$ is a real, \term[quadratic form!nondegenerate]{nondegenerate
quadratic form},
 \item $Q$~is not a scalar multiple of a form with integer coefficients,
 and
 \item the \term[signature (of a quadratic form)]{signature} $(p,q)$ of~$Q$
satisfies $p \ge 3$ and $q \ge 1$.
 \end{itemize}
  Then, for any interval $(\mathsf{a},\mathsf{b})$ in~$\real$, we have
  $$ \frac{
 \# \bigset{v \in \integer^{p+q}}
 { \begin{matrix}
 \mathsf{a} < Q(v) < \mathsf{b} , \\
 \|v\| \le N 
 \end{matrix}
 }
 }{
  \vol \bigset{v \in \real^{p+q}}
 { \begin{matrix}
 \mathsf{a} < Q(v) < \mathsf{b} , \\
 \|v\| \le N 
 \end{matrix}
 }
 }
 \to 1
 \mbox{\qquad as $N \to \infty$}. $$
 \end{thm}

\begin{rem} \ 
 \begin{enumerate}
 \item By calculating the appropriate volume, one finds a constant~$C_Q$,
depending only on~$Q$, such that, as $N \to \infty$,
 $$ \# \bigset{v \in \integer^{p+q}}
 { \begin{matrix}
 \mathsf{a} < Q(v) < \mathsf{b} , \\
 \|v\| \le N 
 \end{matrix}
 }
 \sim 
 (\mathsf{b} - \mathsf{a}) C_Q N^{p + q - 2}
 . $$
 \item The restriction on the \term[signature (of a quadratic
form)]{signature} of~$Q$ cannot be eliminated; there are
\term[counterexample]{counterexamples} of signature $(2,2)$ and $(2,1)$.
 \end{enumerate}
 \end{rem}

\emphit{Why \term[Ratner's Theorems]{Ratner's Theorem} is relevant.} 
 We provide only an indication of the direction of attack, not an
actual proof of the theorem.
 \begin{enumerate}

 \item Let \index{SO(Q)*$\SO(Q)$}$K = \SO(p) \times \SO(q)$, so $K$ is a
\index{Lie group!compact}compact subgroup of $\SO(p,q)$.

 \item For $\mathsf{c}, \mathsf{r} \in \real$, it is not difficult to see
that $K$ is \term[transitive action]{transitive} on 
 $$ \{\, v \in \real^{p+q} \mid Q_{p,q}(v) = \mathsf{c}, \ \|v\| =
\mathsf{r} \,\}$$
 (unless $q = 1$, in which case $K$ has two orbits).

 \item
 Fix $g \in \SL(p+q,\real)$\index{SL(l,R)*$\SL(\ell,\real)$}, such that $Q
= Q_{p,q} \circ g$. (Actually, $Q$ may be a scalar multiple of $Q_{p,q}
\circ g$, but let us ignore this issue.)

 \item Fix a nontrivial one-parameter
\term[unipotent!subgroup!one-parameter]{unipotent subgroup} $u^t$ of
\index{SO(Q)*$\SO(Q)$}$\SO(p,q)$.

 \item Let $\neigh$ be a bounded open set that 
 \begin{itemize}
 \item intersects $Q_{p,q}^{-1}(c)$, for all $\mathsf{c} \in
(\mathsf{a},\mathsf{b})$,
 and
 \item does not contain any \index{fixed point}fixed points of~$u^t$ in
its closure.
 \end{itemize}
 By being a bit more careful in the choice of $\neigh$ and~$u^t$, we may
arrange that $\| w u^t \|$ is within a constant factor of~$t^2$ for all $w
\in \neigh$ and all large $t \in \real$.

 \item If $v$ is any large element of~$\real^{p+q}$, with $Q_{p,q}(v) \in
(\mathsf{a},\mathsf{b})$, then there is some $w \in \neigh$, such that
$Q_{p,q}(w) = Q_{p,q}(v)$. If we choose $t \in \real^+$ with $\|w u^{-t}\|
= \|v\|$ (note that $t < C \sqrt{\|v\|}$, for an appropriate
constant~$C$), then $w \in v K u^t$. Therefore
 \begin{equation} \label{OppQuant-Non0}
 \int_0^{C \sqrt{\|v\|}} \int_K \chi_{\neigh}( v k u^t) \, dk \, dt \neq 0
, \end{equation}
 where $\chi_{\neigh}$ is the characteristic function of~$\neigh$.

 \item  We have
 \begin{align*}
 & \bigset{v \in \integer^{p+q}}
 { \begin{matrix}
 \mathsf{a} < Q(v) < \mathsf{b} , \\
 \|v\| \le N 
 \end{matrix}
 }
 \\& \qquad = 
 \bigset{v \in \integer^{p+q}}
 { \begin{matrix}
 \mathsf{a} < Q_{p,q}(v g) < \mathsf{b} , \\
 \|v\| \le N 
 \end{matrix}
 }
 .
 \end{align*}
 From~\pref{OppQuant-Non0}, we see that the cardinality of the right-hand
side can be approximated by
 $$ \sum_{v \in \integer^{p+q}} \int_0^{C \sqrt{N}} \int_K \chi_{\neigh}(
v g k u^t) \, dk \, dt
 .$$

\item By
 \begin{itemize}
 \item bringing the sum inside the integrals,
 and
 \item defining $\widetilde{\chi_{\neigh}} \colon \Gamma \backslash G \to
\real$ by 
 $$ \widetilde{\chi_{\neigh}} (\Gamma x) = 
 \sum_{v \in \integer^{p+q}} \chi_{\neigh}( v x) , $$
 where $G = \SL(p+q, \real)$\index{SL(l,R)*$\SL(\ell,\real)$} and $\Gamma =
\SL(p+q,\integer)$\index{SL(l,Z)*$\SL(\ell,\integer)$},
 \end{itemize}
 we obtain
 \begin{equation} \label{OppEquiIntegToEstimate}
  \int_0^{C \sqrt{N}} \int_K \widetilde{\chi_{\neigh}}( \Gamma g k u^t)
\, dk \, dt
 .
 \end{equation}
 The outer integral is the type that can be calculated from
\term[Ratner's Theorems!Equidistribution]{Ratner's Equidistribution
Theorem} \pref{RatnerUnifDist} (except that the integrand is not
continuous and may not have compact \term[support!of a function]{support}).

 \end{enumerate}

\begin{rem} \ 
 \begin{enumerate}
 \item Because of technical issues, it is actually a more precise version
\pref{DaniMargUnifDist} of equidistribution that is used to estimate the
integral \pref{OppEquiIntegToEstimate}. In fact, the issues are so serious
that the above argument actually yields only a lower bound on the
integral. Obtaining the correct upper bound requires additional difficult
arguments. 
 \item Furthermore, the conclusion of Thm.~\ref{QuantOppThm} fails for
some forms of signature $(2,2)$ or $(2,1)$; the limit may be $+\infty$.
 \end{enumerate}
 \end{rem}

\subsection{Arithmetic Quantum Unique Ergodicity}
 \index{Quantum!Unique Ergodicity}

 Suppose $\Gamma$ is a \term{lattice} in $G = \SL(2,\real)$, such that
$\Gamma \backslash G$ is compact. Then $M = \Gamma \backslash \hyperbolic$
is a compact \index{manifold}manifold. (We should assume here that $\Gamma$
has no elements of finite order.) The \term[hyperbolic!metric]{hyperbolic
metric} on~$\hyperbolic$ yields a Riemannian metric on~$M$, and there is a
corresponding \term{Laplacian}~$\Delta$ and volume measure~$\vol$
(normalized to be a probability measure). Let 
 $$ 0 = \lambda_0 < \lambda_1 \le \lambda_2 \le\cdots $$
 be the \term[eigenvalue]{eigenvalues} of~$\Delta$ (with multiplicity). For
each $\lambda_n$, there is a corresponding \term{eigenfunction} $\phi_n$,
which we assume to be normalized (and real valued), so that $\int_M
\phi_n^2 \, d \vol = 1$.

In \index{Quantum!Mechanics}Quantum Mechanics, one may think of $\phi_n$ as
a possible state of a particle in a certain system; if the particle is in
this state, then the probability of finding it at any particular location
on~$M$ is represented by the probability distribution $\phi_n^2 \, d
\vol$. It is natural to investigate the limit as $\lambda_n \to \infty$,
for this describes the behavior that can be expected when there is enough
energy that quantum effects can be ignored, and the laws of classical
mechanics can be applied. 

It is conjectured that, in this \term{classical limit}, the particle
becomes \term{uniformly distributed}:

\begin{conj}[(Quantum Unique Ergodicity)] \label{QUEConj}
 $$\lim_{n \to \infty} \phi_n^2 \, d \vol = d \vol .$$
 \end{conj}
 \index{Quantum!Unique Ergodicity}

This conjecture remains open, but it has been proved in an important
special case.

\begin{defn} \ 
 \begin{enumerate}
 \item If $\Gamma$ belongs to a certain family of \term[lattice]{lattices}
(constructed by a certain method from an algebra of
\term[quaternion]{quaternions} over~$\rational$) then we say that $\Gamma$
is a \defit[lattice!congruence]{congruence} \term{lattice}. Although these
are very special \term[lattice]{lattices}, they arise very naturally in 
many applications in number theory and elsewhere.
 \item If the \term{eigenvalue} $\lambda_n$ is simple (i.e, if $\lambda_n$
is not a repeated eigenvalue), then the corresponding \term{eigenfunction}
$\phi_n$ is uniquely determined (up to a sign). If $\lambda_n$ is not
simple, then there is an entire space of possibilities for~$\phi_n$, and
this ambiguity results in a serious difficulty. 

Under the assumption that $\Gamma$ is a congruence \term{lattice}, it is
possible to define a particular orthonormal \index{basis (of a vector
space)}{basis} of each \index{eigenspace}eigenspace; the elements of this
basis are well defined (up to a sign) and are called \defit[Hecke
eigenfunction]{Hecke eigenfunctions}, (or \defit[Hecke-Maass cusp
form]{Hecke-Maass cusp forms}).
 \end{enumerate}
  We remark that if $\Gamma$ is a \term[lattice!congruence]{congruence}
\term{lattice}, and there are no repeated \term{eigenvalue}s, then each
$\phi_n$ is automatically a \term{Hecke eigenfunction}.
 \end{defn}

\begin{thm} \label{AQUE} \index{Quantum!Unique Ergodicity}
 If
 \begin{itemize}
 \item $\Gamma$ is a congruence \term{lattice},
 and
 \item each $\phi_n$ is a Hecke \term{eigenfunction},
 \end{itemize}
 then 
 $\lim_{n \to \infty} \phi_n^2 \, d \vol = d \vol$.
 \end{thm}

\emphit{Why \term[Ratner's Theorems]{Ratner's Theorem} is relevant.}
 Let $\mu$ be a limit of some subsequence of $\phi_n^2 \, d \vol$. Then
$\mu$ can be lifted to an $a^t$-\term[measure!invariant]{invariant}
probability measure~$\widehat\mu$ on $\Gamma \backslash G$. Unfortunately,
$a^t$ is not unipotent, so Ratner's Theorem does not immediately apply. 

Because each $\phi_n$ is assumed to be a Hecke \term{eigenfunction}, one is
able to further lift~$\mu$ to a measure~$\widetilde\mu$ on a certain
\term[homogeneous!space]{homogeneous space} $\widetilde\Gamma \backslash
\bigl( G \times \SL(2,\rational_p) \bigr)$, where $\rational_p$~denotes the
field of \term[p-adic*$p$-adic]{$p$-adic} numbers for an appropriate
prime~$p$. There is an additional action coming from the factor
$\SL(2,\rational_p)$. By combining this action with the
 ``\index{shearing property}Shearing Property"
 \index{property!shearing|indsee{shearing property}}
 of the \term[unipotent!flow]{$u^t$-flow}, much as in the proof of
\pref{RatMeas-S=U} below, one shows that $\widetilde\mu$ is
$u^t$-\term[measure!invariant]{invariant}. (This argument requires one to
know that the \index{entropy!of a flow}entropy $h_{\hat\mu}(a^t)$ is
nonzero.) Then a version of \term[Ratner's Theorems]{Ratner's Theorem}
generalized to apply to $p$-adic groups implies that $\widetilde\mu$ is
$\SL(2,\real)$-\term[measure!invariant]{invariant}.

\subsection{Subgroups generated by lattices in opposite horospherical
subgroups}
 \label{ApplOh}
 \index{generated by unipotent elements}
 \index{lattice|(}
 \index{subgroup!horospherical!opposite}

\begin{notation}
 For $1 \le k < \ell$, let
 \begin{itemize}
 \item
 $\mathbb{U}_{k,\ell} = 
 \bigset{ g \in \SL(\ell,\real) }{
 \text{$g_{i,j} = \delta_{i,j}$ if $i > k$ or $j \le k$}
 }$\index{SL(l,R)*$\SL(\ell,\real)$|(},
 and
 \item $\mathbb{V}_{k.\ell} = 
 \bigset{ g \in \SL(\ell,\real) }{
 \text{$g_{i,j} = \delta_{i,j}$ if $j > k$ or $i \le k$}
 }$.
 \end{itemize}
 (We remark that $\mathbb{V}_\ell$ is the \term{transpose}
of~$\mathbb{U}_\ell$.)
 \end{notation}

\begin{eg}
 $$ \mathbb{U}_{3,5} = \left\{ \begin{bmatrix}
 1 & 0 & 0 & * & * \\
 0 & 1 & 0 & * & * \\
 0 & 0 & 1 & * & * \\
 0 & 0 & 0 & 1 & 0 \\
 0 & 0 & 0 & 0 & 1 
 \end{bmatrix} \right\}
 \text{ and }
  \mathbb{V}_{3,5} = \left\{ \begin{bmatrix}
 1 & 0 & 0 & 0 & 0 \\
 0 & 1 & 0 & 0 & 0 \\
 0 & 0 & 1 & 0 & 0 \\
 * & * & * & 1 & 0 \\
 * & * & * & 0 & 1 
 \end{bmatrix} \right\}
 . $$
 \end{eg}

\begin{thm}
 Suppose
 \begin{itemize}
 \item $\Gamma_U$ is a \term{lattice} in $\mathbb{U}_{k,\ell}$, 
 and
 \item $\Gamma_V$ is a \term{lattice} in $\mathbb{V}_{k,\ell}$, 
 \item the subgroup $\Gamma = \langle \Gamma_U, \Gamma_V \rangle$ is
\defit[subgroup!discrete]{discrete},
 and
 \item $\ell \ge 4$.
 \end{itemize}
 Then $\Gamma$ is a \term{lattice} in
$\SL(\ell,\real)$\index{SL(l,R)*$\SL(\ell,\real)$}.
 \end{thm}

 \emphit{Why \term[Ratner's Theorems]{Ratner's Theorem} is relevant.} Let
$\mathcal{U}_{k,\ell}$ be the space of \term[lattice]{lattices} in
$\mathbb{U}_{k,\ell}$ and $\mathcal{V}_{k,\ell}$ be the space of
\term[lattice]{lattices} in $\mathbb{V}_{k,\ell}$. (Actually, we consider
only \term[lattice]{lattices} with the same ``covolume" as~$\Gamma_U$
or~$\Gamma_V$, respectively.) The block-diagonal subgroup $\SL(k,\real)
\times \SL(\ell-k,\real)$ \index{normalizer}normalizes
$\mathbb{U}_{k,\ell}$ and $\mathbb{V}_{k,\ell}$, so it acts by conjugation
on $\mathcal{U}_{k,\ell} \times \mathcal{V}_{k,\ell}$. There is a natural
identification of this with an action by translations on a
\term[homogeneous!space]{homogeneous space} of $\SL(k \ell, \real) \times
\SL(k \ell,\real)$, so  \term[Ratner's Theorems]{Ratner's Theorem} implies that the closure
of the orbit of $(\Gamma_U, \Gamma_V)$ is homogeneous \see{RatOrbH}. This
means that there are very few possibilities for the closure. By combining
this conclusion with the \index{subgroup!discrete}discreteness of~$\Gamma$
(and other ideas), one can establish that the orbit itself is closed. This
implies a certain compatibility between $\Gamma_U$ and~$\Gamma_V$, which
leads to the desired conclusion.

\begin{rem}
 For simplicity, we have stated only a very special case of the above
theorem. More generally, one can replace
$\SL(\ell,\real)$\index{SL(l,R)*$\SL(\ell,\real)$} with another
\term[group!simple (or almost)]{simple Lie group} of real rank at
least~$2$, and replace $\mathbb{U}_{k,\ell}$ and~$\mathbb{V}_{k,\ell}$
with a pair of opposite horospherical subgroups. The conclusion should be
that $\Gamma$ is a \term{lattice} in~$G$, but this has only been proved
under certain additional technical assumptions.
 \end{rem}
 \index{Lie group!simple|indsee{group,~simple}}
 \index{lattice|)}
 \index{SL(l,R)*$\SL(\ell,\real)$|)}

\subsection{Other results}

For the interested reader, we list some of the many additional
publications that put \term{Ratner's Theorems} to good
use in a variety of ways.

 \begin{itemize}

\item S. Adams:
 Containment does not imply Borel reducibility, in:
 S. Thomas, ed.,
 \emphit{Set theory (Piscataway, NJ, 1999)}, pages 1--23.
 % DIMACS Ser. Discrete Math. Theoret. Comput. Sci., 58, 
 Amer.\ Math.\ Soc., Providence, RI, 2002. 
 \MR{2003j:03059} 

\item A.~Borel and G.~Prasad:
 Values of isotropic quadratic forms at $S$-integral points,
 \emphit{Compositio Math.} 83 (1992), no.~3, 347--372.
 \MR{93j:11022} %MR1175945

 \item N.~Elkies and C.~T.~McMullen:
 Gaps in ${\sqrt n}\bmod 1$ and ergodic theory,
 \emphit{Duke Math.~J.}  123  (2004),  no.~1, 95--139.
 \MR{2060024}

\item A.~Eskin, H.~Masur, and M.~Schmoll:
 Billiards in rectangles with barriers,
 \emphit{Duke Math. J.}  118  (2003),  no.~3, 427--463.
 \MR{2004c:37059}

 \item A.~Eskin, S.~Mozes, and N.~Shah:
 Unipotent flows and counting lattice points on homogeneous varieties,
 \emphit{Ann. of Math.}  143  (1996),  no.~2, 253--299.
 \MR{97d:22012}

 \item   A.~Gorodnik:
 Uniform distribution of orbits of lattices on spaces of frames,
 \emphit{Duke Math.~J.}  122  (2004),  no.~3, 549--589.
 \MR{2057018}

 \item J.~Marklof:
 Pair correlation densities of inhomogeneous quadratic forms,
 \emphit{Ann. of Math.}  158  (2003),  no.~2, 419--471.
 \MR{2018926}

\item T.~L.~Payne:
 Closures of totally geodesic immersions into locally symmetric spaces of
noncompact type,
 \emphit{Proc.\ Amer.\ Math.\ Soc.} 127  (1999),  no.~3, 829--833.
 \MR{99f:53050}

 \item  V.~Vatsal:
 Uniform distribution of Heegner points, 
 \emphit{Invent.\ Math.}  148  (2002),  no.~1, 1--46.
 \MR{2003j:11070}

 \item R.~J.~Zimmer:
 Superrigidity, Ratner's theorem, and fundamental groups, 
 \emphit{Israel J.\ Math.} 74  (1991),  no.~2-3, 199--207.
 \MR{93b:22019} 

 \end{itemize}

\section{Polynomial divergence and shearing}
 \index{polynomial!divergence}
 \index{shearing property}
 \index{polynomial|indsee{function,~polynomial}}
 \label{Poly+ShearSect}

In this section, we illustrate some basic ideas that are used in Ratner's
proof that \term[ergodic!measure]{ergodic measures} are
\term[measure!homogeneous]{homogeneous}\index{Ratner's Theorems!Measure Classification}\pref{RatnerMeas}. This will be done
by giving direct proofs of some statements that follow easily from her
theorem. Our focus is on the group $\SL(2,\real)$.

\begin{notation}
 Throughout this section,
 \begin{itemize}
 \item $\Gamma$ and $\Gamma'$ are \term[lattice]{lattices} in $\SL(2,\real)$,
 \item $u^t$ is the one-parameter
\term[unipotent!subgroup!one-parameter]{unipotent subgroup} of
$\SL(2,\real)$ defined in \pref{SL2etcNotation},
 \item $\eta_t$ is the corresponding \term[unipotent!flow]{unipotent flow}
on $\Gamma \backslash \SL(2,\real)$,
 and
 \item $\eta'_t$ is the corresponding \term[unipotent!flow]{unipotent
flow} on $\Gamma' \backslash \SL(2,\real)$.
 \end{itemize}
 Furthermore, to provide an easy source of
\term[counterexample]{counterexamples},
 \begin{itemize}
 \item $a^t$ is the one-parameter
\term[subgroup!one-parameter!diagonal]{diagonal subgroup} of $\SL(2,\real)$
defined in \pref{SL2etcNotation},
 \item $\gamma_t$ is the corresponding \term[geodesic!flow]{geodesic flow}
on $\Gamma \backslash \SL(2,\real)$,
 and
 \item $\gamma'_t$ is the corresponding \term[geodesic!flow]{geodesic flow}
on $\Gamma' \backslash \SL(2,\real)$.
 \end{itemize}
 For convenience,
 \begin{itemize}
 \item we sometimes write $X$ for $\Gamma \backslash \SL(2,\real)$,
 and
 \item we sometimes write $X'$ for $\Gamma' \backslash \SL(2,\real)$.
 \end{itemize}
 \end{notation}

Let us begin by looking at one of Ratner's first major results in the
subject of \term[unipotent!flow]{unipotent flows}. 

\begin{eg}
 Suppose $\Gamma$ is \defit[subgroup!conjugate]{conjugate} to~$\Gamma'$.
That is, suppose there exists $g \in \SL(2,\real)$, such that $\Gamma =
g^{-1} \Gamma' g$.

Then $\eta_t$ is \defit[isomorphic, measurably]{measurably isomorphic}
to~$\eta'_t$. That is, there is a (measure-preserving) bijection $\psi
\colon \Gamma \backslash \SL(2,\real) \to \Gamma' \backslash
\SL(2,\real)$, such that $\psi \circ \eta_t = \eta'_t \circ \psi$ (a.e.).

Namely, $\psi(\Gamma x) = \Gamma' g x$ \seeex{ConjLat->IsoFlows}. One may
note that $\psi$~is continuous (in fact,~$C^\infty$), not just measurable.
 \end{eg}

The example shows that if $\Gamma$ is conjugate to~$\Gamma'$, then
$\eta_t$ is measurably isomorphic to~$\eta'_t$. (Furthermore, the
isomorphism is obvious, not some complicated measurable function.) Ratner
proved the converse. As we will see, this is now an easy consequence of
\term[Ratner's Theorems!Measure Classification]{Ratner's
Measure Classification Theorem}~\pref{RatnerMeas}, but it was once an
important theorem in its own right.

\begin{cor}[{(\term[Ratner!Rigidity Theorem]{Ratner Rigidity Theorem})}]
\label{RatnerRigThm}
  If $\eta_t$ is \term[isomorphic, measurably]{measurably isomorphic}
to~$\eta'_t$, then $\Gamma$ is \term[subgroup!conjugate]{conjugate}
to~$\Gamma '$.
 \end{cor}

This means that if $\eta_t$ is measurably isomorphic to~$\eta'_t$, then it
is obvious that the two flows are isomorphic, and an isomorphism can be
taken to be a nice, $C^\infty$~map. This is a very special property of
\term[unipotent!flow]{unipotent flows}; in general, it is difficult to
decide whether or not two flows are measurably isomorphic, and measurable
isomorphisms are usually not~$C^\infty$. For example, it can be shown that
$\gamma_t$ is always \term[isomorphic, measurably]{measurably isomorphic}
to~$\gamma'_t$ (even if $\Gamma$ is not conjugate to~$\Gamma'$), but there
is usually no $C^\infty$~isomorphism. (For the experts: this is because
\term[geodesic!flow]{geodesic flows} are \term[Bernoulli
shift]{Bernoulli}.) 

\begin{rem} \label{RatRigRem}
 \ 
 \begin{enumerate}
 \item \label{RatRigRem-Witte}
 A version of \term[Ratner!Rigidity Theorem]{Cor.~\ref{RatnerRigThm}}
remains true with any \term{Lie group}~$G$ in the place of $\SL(2,\real )$,
and any (\term[ergodic!flow]{ergodic}) \term[unipotent!flow]{unipotent
flows}.
 \item In contrast, the conclusion \term[counterexample]{fails miserably}
for some subgroups that are not unipotent.  For example, choose 
 \begin{itemize}
 \item \label{RatRigRem-Bernoulli}
 any $n,n' \ge 2$,
 and 
 \item any \term[lattice]{lattices} $\Gamma$ and $\Gamma'$ in $G =
\SL(n,\real)$ and $G' = \SL(n',\real)$\index{SL(l,R)*$\SL(\ell,\real)$},
respectively.
 \end{itemize}
 By embedding $a^t$ in the top left corner of $G$
and~$G'$, we obtain (\term[ergodic!flow]{ergodic}) flows $\varphi_t$ and
$\varphi'_t$ on $\Gamma \backslash G$ and $\Gamma' \backslash
G'$, respectively.

There is obviously no $C^\infty$~isomorphism between $\varphi_t$ and
$\varphi'_t$, because the \term[homogeneous!space]{homogeneous spaces}
$\Gamma \backslash G$ and $\Gamma' \backslash G'$ do not have the same
dimension\index{dimension!of a manifold} (unless $n = n'$). Even so,
it turns out that the two flows are
\term[isomorphic, measurably]{measurably isomorphic} (up to a change in
speed; that is, after replacing $\varphi_t$ with $\varphi_{ct}$ for some
$c \in \real^\times$). (For the experts: this is because the flows are
\term[Bernoulli shift]{Bernoulli}.) 
 \end{enumerate}
 \end{rem}

\begin{proof}[{\bf Proof of \term[Ratner!Rigidity
Theorem]{Cor.~\ref{RatnerRigThm}}}]
 Suppose $\psi  \colon (\eta_t,X) \to (\eta'_t,X')$ is a measurable
isomorphism. Consider the \index{graph}graph of~$\psi$:
 $$ \graph(\psi) = \bigset{ \bigl( x, \psi(x) \bigr)}{ x \in X} \subset X
\times X' .$$
 Because $\psi$ is measure preserving and \term{equivariant}, we see that
the measure~$\mu_G$ on~$X$ pushes to an \term[ergodic!measure]{ergodic}
$\eta_t \times \eta'_t$-\term[measure!invariant]{invariant
measure}~$\mu_\times$ on $X \times X'$  \seeex{DiagPushInvMeas}.
 \begin{itemize}
 \item Because $\eta_t \times \eta'_t$ is a \term[unipotent!flow]{unipotent
flow} \seeex{UnipInSL2xSL2}, \term[Ratner's Theorems!Measure Classification]{Ratner's Measure Classification Theorem}
\pref{RatnerMeas} applies, so we conclude that the \term[support!of a measure]{support}
of~$\mu_\times$ is a single orbit of a subgroup~$S$ of $\SL(2,\real)
\times \SL(2,\real)$. 
 \item On the other hand, \index{graph}$\graph(\psi)$ is the
\term[support!of a measure]{support} of~$\mu_\times$.
 \end{itemize}
 We conclude that the \index{graph}graph of~$\psi$ is a single~$S$-orbit
(a.e.). This implies that $\psi$ is equal to an \defit[map!affine]{affine
map} (a.e.); that is, $\psi$~the composition of a group
\index{homomorphism!of Lie groups}homomorphism and a translation
\seeex{GrfOrbit->Aff}. So $\psi$~is of a purely algebraic nature, not a
terrible measurable map, and this implies that $\Gamma$ is
\term[subgroup!conjugate]{conjugate} to~$\Gamma'$ \seeex{AffMap->ConjLatt}.
 \end{proof}

We have seen that Cor.~\ref{RatnerRigThm} is a consequence of Ratner's
Theorem \pref{RatnerMeas}. It can also be proved directly, but the proof
does not help to illustrate the ideas that are the main goal of this
section, so we omit it. Instead, let us consider another consequence of
\term[Ratner's Theorems]{Ratner's Theorem}.

\begin{defn}
 A flow $(\varphi_t,\Omega)$ is a \defit[quotient!of a flow]{quotient}
(or \defit[factor (of a flow)]{factor}) of $(\eta_t,X)$ if there is a
measure-preserving Borel function $\psi \colon X \to \Omega$, such that
 \begin{equation} \label{EquiEq}
 \mbox{$\psi \circ \eta_t = \varphi_t \circ \psi$
(a.e.).}
 \end{equation}
 For short, we may say $\psi$ is (essentially) \defit{equivariant} if
\pref{EquiEq} holds.

The function~$\psi$ is \emph{not} assumed to be injective. (Indeed,
quotients are most interesting when $\psi$ collapses substantial portions
of~$X$ to single points in~$\Omega$.) On the other hand, $\psi$ must be
essentially surjective \seeex{quot->surj}.
 \end{defn}

\begin{eg} \label{QuotEgs}
 \ 
 \begin{enumerate}
 \item \label{QuotEgs-latts}
 If $\Gamma \subset \Gamma'$, then the \term[horocycle!flow]{horocycle
flow} $(\eta'_t,X')$ is a \term[quotient!of a flow]{quotient} of
$(\eta_t,X)$ \seeex{QuotExLatts}.
 \item  \label{QuotEgs-tori}
 For $v \in \real^n$ and $v' \in \real^{n'}$, let $\varphi_t$ and
$\varphi'_t$ be the corresponding flows on
\term[torus!Tn*$\torus^n$]{$\torus^n$} and~$\torus^{n'}$. If
 \begin{itemize}
 \item $n' < n$, and
 \item $v'_i = v_i$ for $i = 1,\ldots, n'$,
 \end{itemize} then
$(\varphi'_t,X')$ is a \term[quotient!of a flow]{quotient} of
$(\varphi_t,X)$ \seeex{QuotExTori}.
 \item The one-point space $\{ *\}$ is a \term[quotient!of a flow]{quotient} of any flow. This is the \defit[trivial quotient]{trivial}
quotient.
 \end{enumerate}
 \end{eg}

\begin{rem}
 Suppose $(\varphi_t,\Omega)$ is a \term[quotient!of a flow]{quotient}
of $(\eta_t,X)$. Then there is a map $\psi \colon X \to \Omega$ that is
essentially \term{equivariant}. If we define
 $$ \mbox{$x \sim y$ when $\psi(x) = \psi(y)$} ,$$
 then $\sim$ is an \term{equivalence relation} on~$X$, and we may identify
$\Omega$ with the quotient space $X/{\sim}$.

For simplicity, let us assume $\psi$ is completely equivariant (not just
a.e.). Then the equivalence relation $\sim$ is
$\eta_t$-\term[equivalence relation!invariant]{invariant}; if $x \sim y$,
then $\eta_t(x) \sim \eta_t(y)$. Conversely, if $\equiv$ is an
$\eta_t$-\term[equivalence relation!invariant]{invariant} (measurable)
{equivalence relation} on~$X$, then $X/{\equiv}$ is a
\term[quotient!of a flow]{quotient} of $(\varphi_t,\Omega)$.
 \end{rem}

Ratner proved, for $G = \SL(2,\real)$, that \term[unipotent!flow]{unipotent
flows} are closed under taking \term[quotient!of a flow]{quotients}:

\begin{cor}[{(\term[Ratner!Quotients Theorem]{Ratner Quotients Theorem})}]
\label{RatnerQuotients}
 Each nontrivial \term[quotient!of a flow]{quotient} of $\bigl(
\eta_t,\Gamma \backslash \SL(2,\real ) \bigr)$ is isomorphic to a
\term[horocycle!flow]{unipotent flow} $\bigl( \eta'_t,\Gamma' \backslash
\SL(2,\real ) \bigr)$, for some \term{lattice}~$\Gamma'$.
 \end{cor}

One can derive this from \term[Ratner's Theorems!Measure Classification]{Ratner's Measure Classification Theorem}
\pref{RatnerMeas}, by putting an $(\eta_t \times
\eta_t)$-\term[measure!invariant]{invariant probability measure} on 
 $$ \bigset{ (x,y) \in X \times X }{ \psi(x) = \psi(y) } .$$
 We omit the argument (it is similar to the proof of \term[Ratner!Rigidity
Theorem]{Cor.~\ref{RatnerRigThm}} \seeex{RatnerQuotientsEx}), because it
is very instructive to see a direct proof that does not appeal to Ratner's
Theorem. However, we will prove only the following weaker statement. (The
proof of \pref{RatnerQuotients} can then be completed by applying
Cor.~\ref{FinFib->FinCov} below \seeex{QuotSmooth}.)

\begin{defn}
 A Borel function $\psi \colon X \to \Omega$ has
\defit{finite fibers} (a.e.) if there is a conull subset $X_0$ of~$X$,
such that $\psi^{-1}(\omega) \cap X_0$ is finite, for all $\omega \in
\Omega$.
 \end{defn}

\begin{eg}
 If $\Gamma \subset \Gamma'$, then  the natural quotient map $\psi \colon
X \to X'$ \fullcf{QuotEgs}{latts} has \term{finite fibers}
\seeex{AffQuotFinFib}.
 \end{eg}

\begin{cor}[(Ratner)] \label{QuotFinFib}
 If $\bigl( \eta_t,\Gamma \backslash \SL(2,\real ) \bigr) \to
(\varphi_t,\Omega)$ is any quotient map {\rm(}and $\Omega$ is
nontrivial{\rm)}, then $\psi$~has  \term{finite fibers} {\rm(}a.e.{\rm)}. 
 \end{cor}

In preparation for the direct proof of this result, let us develop some
basic properties of \term[unipotent!flow]{unipotent flows} that are also
used in the proofs of \term[Ratner's Theorems]{Ratner's general theorems}.

Recall that
 $$  \mbox{$u^t = \begin{bmatrix}
 1&0\\ t&1
 \end{bmatrix}$ 
  \qquad and \qquad
 $a^t =
 \begin{bmatrix} 
 e^{t}&0 \\
0&e^{-t}
 \end{bmatrix}$.}$$
 For convenience, let $G = \SL(2,\real)$.

\begin{defn}
 If $x$ and~$y$ are any two points of~$\Gamma \backslash G$, then there
exists $q \in G$, such that $y = xq$. If $x$ is close to~$y$ (which we
denote 
 $x \approx y$), then $q$~may be chosen close to the identity. Thus,
we may define a \term[metric on G*metric on $\Gamma \easybackslash
G$]{metric}~$d$ on $\Gamma \backslash G$ by
 \nindex{$d(x,y)$ = distance from $x$ to $y$}
 $$ d(x,y) = \min \bigset{ \|q - \Id \| }{
 \begin{matrix} q \in G, \\ xq = y \end{matrix}
  } ,$$
 where
 \begin{itemize}
 \item \nindex{$\Id$ = identity matrix}
 $\Id$ is the identity matrix,
 and
 \item  \nindex{$\| \cdot \|$ = matrix norm}
  $\| \cdot \|$ is any (fixed) \index{matrix!norm}matrix norm on
 \nindex{$\Mat_{k \times k}(\real)$ = $\{$ $k \times k$ real matrices $\}$}
 $\Mat_{2 \times 2}(\real)$. For example, one may take
 $$ \left\| \begin{bmatrix} \mathsf{a} & \mathsf{b} \\ \mathsf{c} &
\mathsf{d} \end{bmatrix} \right \|
 = \max \bigl\{ |\mathsf{a}|, |\mathsf{b}|, |\mathsf{c}|, |\mathsf{d}|
\bigr\} .$$
 \end{itemize}
 \end{defn}

\begin{figure}
 \begin{center}
 \includegraphics[scale=0.44035]{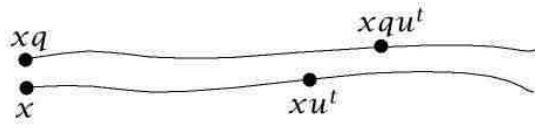}
 \caption{The $\eta_t$-orbits of two nearby orbits.}
 \label{movapart}
 \end{center}
 \end{figure}

A crucial part of \term[Ratner!method]{Ratner's method} involves looking at
what happens to two \term[point!nearby|(]{nearby points} as they move under
the flow~$\eta_t$. Thus, we consider two points $x$ and~$xq$, with $q
\approx \Id$, and we wish to calculate
 $ d \bigl( \eta_t(x), \eta_t(xq) \bigr) $,
 or, in other words, 
 $$ d( x u^t, xq u^t ) $$
 \seefig{movapart}. 
 \begin{itemize}
 \item To get from~$x$ to~$xq$, one multiplies by~$q$; therefore, $d(x,xq)
= \|q - \Id \|$.
 \item To get from $xu^t$ to~$xqu^t$, one multiplies by $u^{-t} q u^t$;
therefore 
 $$d( x u^t, xq u^t ) = \|u^{-t} q u^t - \Id \|$$
 (as long as this is small --- there are infinitely many elements~$g$
of~$G$ with $x u^t g = x q u^t$, and the distance is obtained by choosing
the smallest one, which may not be $u^{-t} q u^t$ if $t$ is large).
 \end{itemize}
 Letting
 $$ q - \Id = \begin{bmatrix}\mathsf{a} &\mathsf{b} \\ \mathsf{c} & \mathsf{d} \end{bmatrix}
,$$
 a simple matrix calculation \seeex{ConjByUEx} shows that
 \begin{equation} \label{ConjByU}
 u^{-t} q u^t - \Id = \begin{bmatrix}
 \mathsf{a} +\mathsf{b} t & \mathsf{b} \\
 \mathsf{c}  - (\mathsf{a} -\mathsf{d})t - \mathsf{b} t^2 & \mathsf{d} - \mathsf{b} t
 \end{bmatrix}
 .
  \end{equation}

All the entries of this matrix are \term[function!polynomial]{polynomials}
(in~$t$), so we have following obvious conclusion:

\begin{prop}[{(\term[polynomial!divergence]{Polynomial divergence})}]
\label{PolyDivSL2}
 Nearby points of $\Gamma \backslash G$ move apart at
\term[speed!polynomial]{polynomial speed}.
 \end{prop}

 In contrast, \term[point!nearby]{nearby points} of the
\term[geodesic!flow]{geodesic flow} move apart at
\term[speed!exponential]{exponential speed}:
 \begin{equation} \label{ConjByA}
 a^{-t} q a^t - \Id 
 =
 \begin{bmatrix}
 \mathsf{a} & \mathsf{b} e^{-2t} \\
 \mathsf{c} e^{2 t} & \mathsf{d} \end{bmatrix}
 \end{equation}
 \seeex{ConjByA}.
 Intuitively, one should think of \term[speed!polynomial]{polynomial
speed} as ``slow" and \term[speed!exponential]{exponential speed} as
``fast." Thus, 
 \begin{itemize}
 \item \term[point!nearby]{nearby points} of a
\term[unipotent!flow]{unipotent flow} drift slowly apart, but
 \item \term[point!nearby]{nearby points} of the
\term[geodesic!flow]{geodesic flow} jump apart rather suddenly.
 \end{itemize}
 More precisely, note that
 \begin{enumerate}
 \item if a \term[function!polynomial]{polynomial} (of bounded degree)
stays small for a certain length of time, then it must remain fairly
small for a proportional length of time \seeex{PolyProportion}:
 \begin{itemize}
 \item if the polynomial is small for a minute, then it must stay
fairly small for another second (say);
 \item if the polynomial is small for an hour, then it must stay
fairly small for another minute;
 \item if the polynomial is small for a year, then it must stay fairly
small for another week;
 \item if the polynomial is small for several thousand years, then it must
stay fairly small for at least a few more decades;
 \item if the polynomial has been small for an
\index{infinitely!long}infinitely long time, then it must stay small
forever (in fact, it is constant).
 \end{itemize}
 \item In contrast, the \term[exponential!function]{exponential function}
$e^t$ is fairly small ($< 1$) \index{infinitely!long}infinitely far into
the past (for $t < 0$), but it becomes arbitrarily large in finite time.
 \end{enumerate}
 Thus, 
 \begin{enumerate}
 \item If two points of a \term[unipotent!flow]{unipotent flow} stay close
together 90\% of the time, then they must stay fairly close together all
of the time.
 \item In contrast, two points of a \term[geodesic!flow]{geodesic flow}
may stay close together 90\% of the time, but spend the remaining 10\% of
their lives \term[wandering around the manifold]{wandering quite freely}
(and independently) around the \index{manifold}manifold.
 \end{enumerate}
 The upshot is that if we can get good bounds on a
\term[unipotent!flow]{unipotent flow} most of the time, then we have bounds
that are nearly as good all of the time:

\begin{notation}
 For convenience, let \nindex{$x_t$ = $x u^t$ (image of~$x$ under
unipotent flow)}$x_t = x u^t$ and $y_t = y u^t$.
 \end{notation}

\begin{cor} \label{PolyOrbsStayClose}
 For any $\epsilon > 0$, there is a $\delta > 0$, such that if $d(x_t,y_t)
< \delta$ for 90\% of the times~$t$ in an interval $[a, b]$, then 
 $d(x_t,y_t) < \epsilon$ for \emph{all} of the times~$t$ in the interval
$[a, b]$.
 \end{cor}

 \begin{figure}
 \begin{center}
 \includegraphics[scale=0.44035]{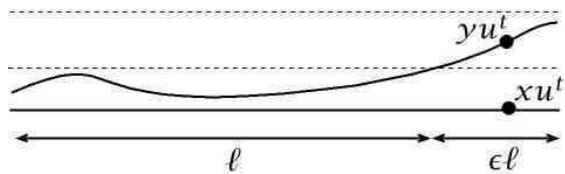}
 \caption{Polynomial divergence: Two points that stay close together for a
period of time of length~$\ell$ must stay fairly close for an additional
length of time~$\epsilon \ell$ that is proportional to~$\ell$.}
 \label{polydiv}
 \end{center}
 \end{figure}

\begin{rem}
 \term[babysitting]{Babysitting} provides an analogy that illustrates this
difference between unipotent flows and \term[geodesic!flow]{geodesic
flows}.
 \begin{enumerate}
 \item A \term[unipotent!child]{unipotent child} is easy to watch over. If
she sits quietly for an hour, then we may leave the room for a few minutes,
knowing that she will not get into trouble while we are away. Before she
leaves the room, she will start to make little motions, squirming in her
chair. Eventually, as the motions grow, she may get out of the chair, but
she will not go far for a while. It is only after giving many \term{warning
signs} that she will start to walk slowly toward the door.
 \item A \term[geodesic!child]{geodesic child}, on the other hand, must be
watched almost constantly. We can take our attention away for only a few
seconds at a time, because, even if she has been sitting quietly in her
chair all morning (or all week), the child might suddenly jump up and run
out of the room while we are not looking, getting into all sorts of
\term{mischief}. Then, before we notice she left, she might go back to her
chair, and sit quietly again. We may have no idea there was anything amiss
while we were not watching.
 \end{enumerate}
 \end{rem}

Consider the RHS of Eq.~\ref{ConjByU}, with $\mathsf{a}$, $\mathsf{b}$,
$\mathsf{c}$, and~$\mathsf{d}$ very small. Indeed, let us say they are
\index{infinitesimal}infinitesimal; too small to see. As $t$~grows, it is
the the bottom left corner that will be the first matrix entry to attain
macroscopic size \seeex{BottomLeftLarge}. Comparing with the definition
of~$u^t$ \see{SL2etcNotation}, we see that this is exactly the direction
of the $u^t$-orbit \seefig{shearing}. Thus:

 \begin{figure}
 \begin{center}
 \includegraphics[scale=0.44035]{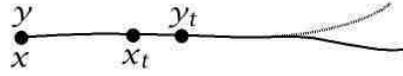}
 \caption{Shearing: If two points start out so close together that we
cannot tell them apart, then the first difference we see will be that one
gets ahead of the other, but (apparently) following the same path. It is
only much later that we will detect any difference between their paths.}
 \label{shearing}
 \end{center}
 \end{figure}

\begin{prop}[{(\term[shearing property]{Shearing Property})}]
\label{ShearingSL2R}
 The fastest \term{relative motion} between two
\term[point!nearby|)]{nearby points} is parallel to the orbits of the flow.
 \end{prop}

The only exception is that if $q \in \{u^t\}$, then $u^{-t} q u^t = q$ for
all~$t$; in this case, the points $x_t$ and~$y_t$ simply move along
together at exactly the same speed, with no \term{relative motion}.

\begin{cor} \label{ShearingSL2R(1step)}
 If $x$ and $y$ are \term[point!nearby]{nearby points}, then either
 \begin{enumerate}
 \item there exists $t > 0$, such that $y_t \approx x_{t \pm 1}$,
 or
 \item $y = x_{\epsilon}$, for some $\epsilon \approx 0$.
 \end{enumerate}
 \end{cor}

\begin{rem}[{(\term[infinitesimal]{Infinitesimals})}]
 Many theorems and proofs in these notes are presented in terms of
\term[infinitesimal]{infinitesimals}. (We write 
  \nindex{$x \approx y$ = $x$ is infinitesimally close to~$y$}
 $x \approx y$ if the distance from~$x$ to~$y$ is infinitesimal.)
There are two main reasons for this:
 \begin{enumerate}
 \item Most importantly, these lectures are intended more to communicate
ideas than to record rigorous proofs, and the terminology of
infinitesimals is very good at that. It is helpful to begin by
pretending that points are infinitely close together, and see what
will happen.  If desired, the reader may bring in epsilons and deltas
after attaining an intuitive understanding of the situation.
	\item \term{Nonstandard Analysis} is a theory that provides a rigorous
foundation to infinitesimals --- almost all of the infinitesimal proofs
that are sketched here can easily be made rigorous in these terms. For
those who are comfortable with it, the infinitesimal approach is often
simpler than the classical notation, but we will provide
\term[non-infinitesimal version]{non-infinitesimal versions} of the main
results in Chap.~\ref{ProofChap}.
 \end{enumerate}
 \end{rem}

\begin{figure}
 \begin{center}
 \includegraphics[scale=0.44035]{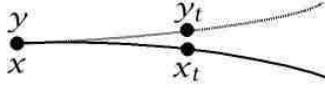}
 \caption{Exponential divergence: when two points start out so close
together that we cannot tell them apart, the first difference we see
may be in a direction transverse to the orbits.}
 \label{expdivfig}
 \end{center}
 \end{figure}

\begin{rem}
 In contrast to the above discussion of~$u^t$, 
 \begin{itemize}
 \item the matrix $a^t$ is diagonal, 
 but
 \item the \term[largest!term]{largest entry} in the RHS of
Eq.~\ref{ConjByA} is an off-diagonal entry,
 \end{itemize}
 so points in the \term[geodesic!flow]{geodesic flow} move apart (at
exponential speed)  in a direction transverse to the orbits
\seefig{expdivfig}. 
 \end{rem}

Let us now illustrate how to use the \term[shearing property]{Shearing
Property}.

\begin{proof}[{\bf Proof of Cor.~\ref{QuotFinFib}}]
 To bring the main ideas to the foreground, let us first consider a special
case with some (rather drastic) simplifying assumptions. We will then
explain that the assumptions are really not important to the argument.
  \begin{enumerate} \renewcommand{\theenumi}{A\arabic{enumi}}
 \item \label{AssXCpct}
 Let us assume that $X$ is compact (rather than merely having finite
volume).
 \item \label{AssNoFP}
 Because $(\varphi_t,\Omega)$ is \term[ergodic!dynamical system]{ergodic}
\seeex{QuotErgodic} and nontrivial, we know that the set of
\index{fixed point}fixed points has measure zero; let us assume that
$(\varphi_t,\Omega)$ has no \term[fixed point]{fixed points} at all.
Therefore, 
 \begin{equation} \label{d(pw,w)bdd}
 \hskip0.75in %% didn't center properly???
 \begin{matrix}
 \text{$d \bigl( \varphi_1(\omega), \omega \bigr)$ is bounded away
from~$0$,} \\
 \text{as $\omega$ ranges over~$\Omega$}
 \end{matrix}
 \end{equation}
 \seeex{NoFP->d>0}.
 \item \label{AssPsiUnifCont}
 Let us assume that the \term[quotient!map]{quotient map} $\psi$ is
\term[function!continuous!uniformly]{uniformly continuous} (rather than
merely being measurable). This may seem unreasonable, but
\term[Theorem!Lusin]{Lusin's Theorem} (Exer.~\ref{LusinThmEx}) tells us
that $\psi$ is \term[function!continuous!uniformly]{uniformly continuous}
on a set of measure $1 - \epsilon$, so, as we shall see, this is
actually not a major issue.
 \end{enumerate}
 Suppose some fiber $\psi^{-1} (\omega_0)$ is infinite. (This will lead
to a contradiction.) 

Because $X$ is compact, the infinite set $\psi^{-1} (\omega_0)$ must have
an \term[point!accumulation]{accumulation point}. Thus, there exist  $x
\approx y$  with  $\psi(x) = \psi(y)$. Because $\psi$ is
\term{equivariant}, we have
 \begin{equation} \label{psi(xt)=psi(yt)}
 \mbox{$\psi(x_t) = \psi(y_t)$ for all~$t$} .
 \end{equation}
 Flow along the orbits until the points $x_t$ and~$y_t$ have
diverged to a reasonable distance; say, $d(x_t,y_t) = 1$, and let 
 \begin{equation} \label{omega=yt}
 \omega = \psi(y_t) .
 \end{equation}
Then the \term[shearing property]{Shearing Property} implies
\see{ShearingSL2R(1step)} that 
 \begin{equation} \label{ytapproxP(xt)}
 y_t \approx \eta_1(x_t)
 . \end{equation}
 Therefore 
 \begin{align*}
 \omega
 &= \psi(y_t) 
 && \mbox{\pref{omega=yt}} \\
 &\approx \psi \bigl( \eta_1(x_t) \bigr)
 && \mbox{(\pref{ytapproxP(xt)} and $\psi$ is \term[function!continuous!uniformly]{uniformly continuous})} \\
 &= \varphi_1 \bigl( \psi (x_t) \bigr)
 && \mbox{($\psi$ is equivariant)} \\
 &= \varphi_1(\omega)
 && \mbox{(\pref{psi(xt)=psi(yt)} and \pref{omega=yt})} 
 .
 \end{align*}
 This contradicts \pref{d(pw,w)bdd}.

To complete the proof, we now indicate how to eliminate the assumptions
\pref{AssXCpct}, \pref{AssNoFP}, and~\pref{AssPsiUnifCont}. 

First, let us point out that \pref{AssXCpct} was not necessary. The proof
shows that $\psi^{-1}(\omega)$ has no
\term[point!accumulation]{accumulation point} (a.e.); thus,
$\psi^{-1}(\omega)$ must be countable. Measure theorists can show that a
\term{countable-to-one} \term{equivariant} map between
\term[ergodic!dynamical system]{ergodic spaces} with
\term[measure!invariant]{invariant probability measure} must actually be
\term[finite fibers]{finite-to-one} (a.e.) \seeex{CtblFibers}. Second,
note that it suffices to show, for each $\epsilon > 0$, that there is a
subset $\hat X$ of~$X$, such that
 \begin{itemize}
 \item $\mu(\hat X) > 1 - \epsilon$ and
 \item $\psi^{-1}(\omega) \cap \hat X$ is countable, for a.e.\ $\omega \in
\Omega$.
 \end{itemize}

Now, let $\hat\Omega$ be the complement of the set of \term[fixed point]{fixed points}
in~$\Omega$. This is conull, so $\psi^{-1}(\hat\Omega)$ is conull in~$X$.
Thus, by \term[Theorem!Lusin]{Lusin's Theorem}, $\psi^{-1}(\hat\Omega)$
contains a compact set~$K$, such that
 \begin{itemize}
 \item $\mu_G(K) > 0.99$, 
 and
 \item $\psi$ is \term[function!continuous!uniformly]{uniformly
continuous} on~$K$.
 \end{itemize}
 Instead of making assumptions \pref{AssNoFP} and~\pref{AssPsiUnifCont},
we work inside of~$K$.
 Note that:
 \begin{enumerate}
 \item[\prefp{AssNoFP}] $d \bigl( \varphi_1(\omega), \omega \bigr)$ is
bounded away from~$0$, for $\omega \in \psi(K)$;
 and
 \item[\prefp{AssPsiUnifCont}]
 $\psi$ is \term[function!continuous!uniformly]{uniformly continuous}
on~$K$.
 \end{enumerate}
 Let $\hat X$ be a \term[generic!set]{generic set} for~$K$; that is, points
in~$\hat X$ spend 99\% of their lives in~$K$. The \term[Theorem!Pointwise
Ergodic]{Pointwise Ergodic Theorem} \pref{PtwiseErgThmFlow} tells us that
the generic set is conull. (Technically, we need the points of~$\hat X$ to
be \defit[generic!uniformly]{uniformly generic}: there is a constant~$L$,
independent of~$x$, such that 
 $$ \begin{matrix}
 \mbox{for all $L' > L$ and $x \in \hat X$, at least 98\%
of} \\
 \mbox{the initial segment $\bigl\{ \phi_t(x) \}_{t = 0}^{L'}$ is in~$K$,}
 \end{matrix} $$
 and this holds only on a set of measure $1 - \epsilon$, but let us ignore
this detail.)
 Given $x,y \in \hat X$, with $x \approx y$, flow along the orbits until
$d(x_t,y_t) = 1$. Unfortunately, it may not be the case that 
 \label{xtytNotInK}
 $x_t$ and~$y_t$ are in~$K$, but, because 99\% of each orbit is in~$K$, we
may choose a nearby value~$t'$ (say, $t \le t' \le 1.1t$), such that 
 $$ \mbox{$x_{t'} \in K$ and $y_{t'} \in K$.} $$
 By \term[polynomial!divergence]{polynomial divergence}, we know that the
$y$-orbit drifts \emph{slowly} ahead of the $x$-orbit, so 
 $$ \mbox{$y_{t'} \approx \eta_{1 + \delta}(x_{t'})$ for some
small~$\delta$.} $$
 Thus, combining the above argument with a strengthened version
of~\prefp{AssNoFP} \seeex{Time1Dist} shows that $\psi^{-1}(\omega) \cap
\hat X$ has no \term[point!accumulation]{accumulation points} (hence, is
countable). This completes the proof.
 \end{proof}

The following application of the \index{shearing property}Shearing Property
is a better illustration of how it is actually used in the proof of
\term[Ratner's Theorems!Measure Classification]{Ratner's Theorem}.

\begin{defn}
 A \defit{self-joining}\index{joining|indsee{self-joining}} of~$(\eta_t,X)$
is a probability measure~$\hat\mu$ on $X \times X$, such that
 \begin{enumerate}
 \item $\hat\mu$ is \term[measure!invariant]{invariant} under the diagonal
flow $\eta_t \times \eta_t$,
 and
 \item $\hat\mu$ projects to~$\mu_G$ on each factor of the product;
that is, $\hat\mu(A \times Y) = \mu_G(A)$ and $\hat\mu(Y \times B) =
\mu_G(B)$.
 \end{enumerate}
 \end{defn}

\begin{eg} \label{JoinEgs} \ 
 \begin{enumerate}
 \item The \term[measure!product]{product measure} $\hat\mu = \mu_G
\times \mu_G$ is a \term{self-joining}.
 \item \label{JoinEgs-diag}
 There is a natural \term{diagonal embedding} $x \mapsto (x,x)$ of $X$ in
$X \times X$. This is clearly \term{equivariant}, so $\mu_G$ pushes to
an $(\eta_t \times \eta_t)$-\term[measure!invariant]{invariant measure} on
$X \times X$. It is a self-joining, called the
\defit[self-joining!diagonal]{diagonal self-joining}.
 \item \label{JoinEgs-conj}
 Replacing the identity map $x \mapsto x$ with
\index{map!covering}covering maps yields a generalization of
\pref{JoinEgs-diag}: For some $g \in G$, let $\Gamma' = \Gamma \cap
(g^{-1} \Gamma g)$, and assume $\Gamma'$ has finite index in~$\Gamma$.
There are two natural \term[map!covering]{covering maps} from $X'$ to~$X$:
 \begin{itemize}
 \item $\psi_1(\Gamma' x) = \Gamma x$,
 and
 \item $\psi_2(\Gamma' x) = \Gamma g x$
 \end{itemize}
 \seeex{JoinEgsConjCont}.
 Define $\psi \colon X' \to X \times X$ by
 $$ \psi(x) = \bigl( \psi_1(x), \psi_2(x) \bigr) .$$
 Then 
 \begin{itemize}
 \item $\psi$ is \term{equivariant} (because $\psi_1$ and~$\psi_2$ are
equivariant), so the $G$-\term[measure!invariant]{invariant} measure
$\mu'_G$ on~$X'$ \term[push-forward (of a measure)]{pushes}
 \nindex{$\psi_* \mu$ = push-forward of~$\mu$ by the map~$\psi$}
 to an \term[measure!invariant]{invariant measure} $\hat\mu = \psi_*
\mu'_G$ on $X \times X$, defined by 
 $$(\psi_* \mu'_G)(A) = \mu'_G \bigl( \psi^{-1}(A) \bigr), $$
 and
 \item $\hat\mu$ is a \term{self-joining} (because $\psi_1$ and~$\psi_2$
are measure preserving).
 \end{itemize}
 This is called a \defit[self-joining!finite cover]{finite-cover
self-joining}.
 \end{enumerate}
 \end{eg}

 \begin{figure}
 \begin{center}
 \includegraphics[scale=0.44035]{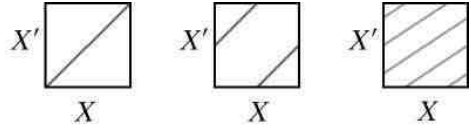}
 \caption{The diagonal self-joining and some other finite-cover
self-joinings.}
 \label{FiniteCoverFigs}
 \end{center}
 \end{figure}

For \term[unipotent!flow]{unipotent flows} on $\Gamma \backslash
\SL(2,\real)$, Ratner showed that these are the only
\term[self-joining]{product self-joinings}.

\begin{cor}[(Ratner's Joinings Theorem)] \label{RatnerJoin}
 Any \term[ergodic!measure]{ergodic} \term{self-joining} of a
\term[horocycle!flow]{horocycle flow} must be either
 \begin{enumerate}
 \item a \term[self-joining!finite cover]{finite cover},
 or
 \item the \term[self-joining]{product self-joining}.
 \end{enumerate}
 \end{cor}

This follows quite easily from \term[Ratner's Theorems!Measure Classification]{Ratner's Theorem} \pref{RatnerMeas}
\seeex{RatnerJoiningsEx}, but we give a direct proof of the following
weaker statement. (Note that if the self-joining $\hat\mu$ is a
\term[self-joining!finite cover]{finite cover}, then $\hat\mu$ has
\term[finite fibers|(]{finite fibers}; that is, $\mu $ is \term[support!of
a measure]{supported} on a set with only finitely many points from each
horizontal or vertical line \seeex{Cover->FinFib}).
Corollary~\ref{FinFib->FinCov} will complete the proof
of~\pref{RatnerJoin}.

\begin{cor} \label{JoiningsThm} \index{Ratner!Joinings Theorem}
 If $\hat\mu$ is an \term[ergodic!measure]{ergodic} \term{self-joining} of
\term[horocycle!flow]{$\eta_t$}, then either
 \begin{enumerate}
 \item  $\hat\mu$ is the \term[self-joining!product]{product joining}, 
 or
 \item $\hat\mu$ has \term[self-joining!finite fibers]{finite fibers}.
 \end{enumerate}
 \end{cor}

\begin{proof} We omit some details (see Exer.~\ref{JoiningsThmRigor} and
Rem.~\ref{TechPtJoinings}).

 Consider two points $(x,a)$ and $(x,b)$ in the same
\term[fiber!vertical]{vertical fiber}. If the fiber is infinite (and $X$ is
compact), then we may assume $a \approx b$. By the \term[shearing
property]{Shearing Property} \pref{ShearingSL2R(1step)}, there is some~$t$
with $a_t \approx \eta_1(b_t)$. Let
 \nindex{$\xi_t$ = vertical flow on $X \times X$}
 $\xi_t$ be the
\defit[flow!vertical]{vertical flow} on~$X \times X$, defined by
 $$ \xi_t(x,y) = \bigl( x, \eta_t(y) \bigr) .$$ 
 Then
 $$(x,a)_t = (x_t,a_t) \approx \bigl(
x_t,\eta_1(b_t) \bigr) = \xi_1\bigl( (x,b)_t \bigr) .$$
 We now consider two cases.

\setcounter{case}{0}

\begin{case} \label{JoiningsThmVertCase}
 Assume $\hat\mu$ is $\xi_t$-\term[measure!invariant]{invariant}.
 \end{case}
 Then the \term[ergodic!flow]{ergodicity} of~$\eta_t$ implies that
$\hat\mu$ is the \term[self-joining!product]{product joining}
\seeex{VertInvJoining}.

\begin{case}
 Assume $\hat\mu$ is \emph{not}
$\xi_t$-\term[measure!invariant]{invariant}.
 \end{case}
 In other words, we have $(\xi_1)_*(\hat\mu) \neq \hat\mu$. On the other
hand, $(\xi_1)_*(\hat\mu)$ is $\eta_t$-\term[measure!invariant]{invariant}
(because $\xi_1$ commutes with~$\eta_t$ \seeex{CommuteInvariant}). It
is a general fact that any two \term[ergodic!measure]{ergodic measures} for
the same flow must be mutually \term[measure!singular]{singular}
\seeex{ErgSingular}, so $(\xi_1)_*(\hat\mu) \perp \hat\mu$; thus, there is
a conull subset~$\hat X$ of $X \times X$, such that $\hat\mu\bigl( \xi_1(X)
\bigr) = 0$. From this, it is not difficult to see that there is a compact
subset~$K$ of $X \times X$, such that
 \begin{equation} \label{JoinKdist>0}
  \mbox{$\hat\mu(K) > 0.99$ and $d \bigl( K, \xi_1(K)
\bigr) > 0$}
 \end{equation}
 \seeex{JoiningShrinkCpct}.

 To complete the proof, we show:

\begin{claim} \index{Claim}
 Any \term[generic!set]{generic set} for~$K$ intersects each
\term[fiber!vertical]{vertical fiber} $\{x\} \times X$ in a countable set.
 \end{claim}
 Suppose not. (This will lead to a contradiction.) Because the fiber
is uncountable, there exist $(x,a)$ and $(x,b)$ in the
\term[generic!set]{generic set}, with $a \approx b$. Flow along the orbits
until 
 $$a_t \approx \eta_1(b_t) ,$$ and
assume (as is true 98\% of the time) that $(x,a)_t$ and $(x,b)_t$ belong
to~$K$. Then
 $$ K \ni (x,a)_t = (x_t, a_t)
 \approx (x_t, \eta_1(b_t))
 = \xi_1 \bigl( (x,b)_t \bigr) \in \xi_1(K), $$
 so $d \bigl(K, \xi_1(K) \bigr) = 0$. This contradicts \pref{JoinKdist>0}.
 \end{proof}

 \begin{figure}
 \begin{center}
 \includegraphics[scale=0.44035]{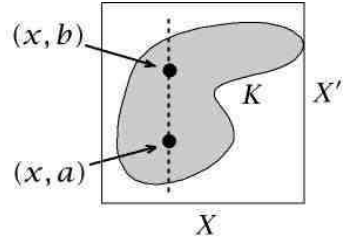}
 \caption{Two points $(x,a)$ and $(x,b)$ on the same vertical fiber.}
 \label{joinK}
 \end{center}
 \end{figure}

\begin{rem} \label{TechPtJoinings}
 The above proof ignores an important technical point that also arose on
p.~\pageref{xtytNotInK} of the proof of Cor.~\ref{QuotFinFib}: at the
precise time~$t$ when $a_t \approx \eta_1(b_t)$, it may not be the case
that $(x,a)_t$ and $(x,b)_t$ belong to~$K$. We choose a nearby
value~$t'$ (say, $t \le t' \le 1.1t$), such that 
 $(x,a)_{t'}$ and $(x,b)_{t'}$ belong to~$K$. By polynomial divergence, we
know that $a_{t'} \approx \eta_{1+\delta}(b_{t'})$ for some small~$\delta$.

Hence, the final stage of the proof requires $\xi_{1+\delta}(K)$ to be
disjoint from~$K$. Since $K$ must be chosen before we know the value
of~$\delta$, this is a serious issue. Fortunately, a single set~$K$ can
be chosen that works for all values of~$\delta$ \cf{KnotKc}.
 \end{rem}

\begin{exercises}

\item \label{ConjLat->IsoFlows}
 Suppose $\Gamma$ and $\Gamma'$ are \term[lattice]{lattices} in $G =
\SL(2,\real)$. Show that if \term[subgroup!conjugate]{$\Gamma' = g^{-1}
\Gamma g$}, for some $g \in G$, then the map $\psi \colon \Gamma
\backslash \SL(2,\real) \to \Gamma' \backslash \SL(2,\real)$, defined by
$\psi(\Gamma x) = \Gamma' g x$,
 \begin{enumerate}
 \item is well defined,
 and
 \item is \term{equivariant}; that is, $\psi \circ \eta_t = \eta'_t \circ
\psi$.
 \end{enumerate}

\item 
 A nonempty, closed subset~$C$ of $X \times X$ is
\defit[set!minimal]{minimal} for $\eta_t \times \eta_t$ if the orbit of
every point in~$C$ is dense in~$C$. Show that if $C$ is a compact minimal
set for $\eta_t \times \eta_t$, then $C$ has \term{finite fibers}.
 \hint{Use the proof of \pref{JoiningsThm}.}

\item \label{DiagPushInvMeas}
 Suppose 
 \begin{itemize}
 \item $(X,\mu)$ and $(X',\mu')$ are Borel measure spaces,
 \item $\varphi_t$ and $\varphi'_t$ are \term[flow!measurable]{(measurable)
flows} on~$X$ and~$X'$, respectively,
 \item $\psi  \colon X \to X'$ is a \term{measure-preserving} map, such
that $\psi \circ \varphi_t = \varphi'_t \circ \psi$ (a.e.),
 and 
 \item $\mu_{\times}$ is the Borel measure on $X \times X'$ that is defined
by 
 $$\mu_{\times}(\Omega) = \mu \bigl\{ \, x \in X \mid \bigl( x, \psi(x)
\bigr) \in \Omega \, \bigr\} .$$
 \end{itemize}
 Show:
 \begin{enumerate}
 \item $\mu_{\times}$ is  $\varphi_t \times
\varphi'_t$-\term[measure!invariant]{invariant}.
 \item If $\mu$ is \term[ergodic!measure]{ergodic} (for $\varphi_t$), then
$\mu_{\times}$ is \term[ergodic!measure]{ergodic} (for $\varphi_t \times
\varphi'_t$).
 \end{enumerate}

\item \label{UnipInSL2xSL2}
 The product $\SL(2,\real) \times \SL(2,\real)$ has a natural embedding in
\index{SL(4,R)*$\SL(4,\real)$}$\SL(4,\real)$ (as block diagonal matrices).
Show that if $u$ and~$v$ are \term[unipotent!matrix]{unipotent matrices}
in $\SL(2,\real)$, then the image of $(u,v)$ is a
\term[unipotent!matrix]{unipotent matrix} in $\SL(4,\real)$.

\item \label{GrfGrp->Homo}
 Suppose $\psi$ is a function from a group~$G$ to a group~$H$. Show
$\psi$ is a \index{homomorphism!of Lie groups}homomorphism if and only if
the \index{graph}graph of~$\psi$ is a \term[subgroup!of $G \times
H$]{subgroup of $G \times H$}.

\item \label{GrfOrbit->Aff}
 Suppose
 \begin{itemize}
 \item $G_1$ and $G_2$ are groups,
 \item $\Gamma_1$ and $\Gamma_2$ are subgroups of $G_1$ and~$G_2$,
respectively,
 \item $\psi$ is a function from $\Gamma_1 \backslash G_1$ to $\Gamma_2
\backslash G_2$,
 \item $S$ is a subgroup of $G_1 \times G_2$,
 and
 \item the \term{graph} of~$\psi$ is equal to $x S$, for some $x \in
\Gamma_1 \backslash G_1 \times \Gamma_2 \backslash G_2$.
 \end{itemize}
 Show:
 \begin{enumerate}
 \item \label{GrfOrbit->Aff-Aff}
 If $S \cap (e \times G_2)$ is trivial, then $\psi$ is an
\term[map!affine]{affine map}.
 \item If $\psi$ is surjective, and $\Gamma_2$ does not contain any
nontrivial \index{subgroup!normal}normal subgroup of~$G_2$, then $S \cap
(e \times G_2)$ is trivial.
 \end{enumerate}
 \hint{\pref{GrfOrbit->Aff-Aff}~$S$~is the \index{graph}graph of a
\index{homomorphism!of Lie groups}homomorphism from~$G_1$ to~$G_2$
\seeex{GrfGrp->Homo}.}

\item \label{AffMap->ConjLatt}
 Suppose $\Gamma$ and $\Gamma'$ are \term[lattice]{lattices} in a (simply connected)
\term{Lie group}~$G$. 
 \begin{enumerate}
 \item Show that if there is a bijective \term[map!affine]{affine map} from
$\Gamma \backslash G$ to $\Gamma' \backslash G$, then there is an
\term{automorphism}~$\alpha$ of~$G$, such that $\alpha(\Gamma) = \Gamma'$.
 \item Show that if $\alpha$ is an automorphism of $\SL(2,\real)$, such
that $\alpha(u)$ is conjugate to~$u$, then $\alpha$ is an
\defit[automorphism!inner]{inner automorphism}; that is, there is some $g
\in \SL(2,\real)$, such that $\alpha(x) = g^{-1} x g$, for all $x \in
\SL(2,\real)$.
 \item Show that if there is a bijective affine map $\psi \colon \Gamma
\backslash G \to \Gamma' \backslash G$, such that $\psi(x u) =
\psi(x) u$, for all $x \in \Gamma \backslash G$, then
$\Gamma$ is \term[subgroup!conjugate]{conjugate} to~$\Gamma'$.
 \end{enumerate}

\item \label{quot->surj}
 Show that if $\psi \colon X \to \Omega$ is a \term{measure-preserving}
map, then $\psi(X)$ is a conull subset of~$\Omega$.

\item \label{QuotExLatts}
 Verify Eg.~\fullref{QuotEgs}{latts}.

\item \label{QuotExTori}
 Verify Eg.~\fullref{QuotEgs}{tori}.

\item \label{RatnerQuotientsEx}
 Give a short proof of \term[Ratner!Quotients
Theorem]{Cor.~\ref{RatnerQuotients}}, by using \term[Ratner's Theorems!Measure Classification]{Ratner's Measure Classification Theorem}.

\item \label{LattFinInd}
 Suppose $\Gamma$ is a \term{lattice} in a \term{Lie group}~$G$. Show that
a subgroup $\Gamma'$ of~$\Gamma$ is a \term{lattice} if and only if $\Gamma'$ has
finite index in~$\Gamma$.

\item \label{AffQuotFinFib}
 Suppose $\Gamma$ and~$\Gamma'$ are \term[lattice]{lattices} in a \term{Lie
group}~$G$, such that $\Gamma \subset \Gamma'$. Show that the natural map
$\Gamma \backslash G \to \Gamma' \backslash G$ has \term[finite
fibers|)]{finite fibers}.

\item \label{ConjByUEx}
 Verify \eqref{ConjByU}.

\item \label{ConjByAEx}
 Verify \eqref{ConjByA}.

\item \label{QuotErgodic}
 Suppose $(\varphi_t',\Omega')$ is a \term[quotient!of a flow]{quotient}
of a flow $(\varphi_t,\Omega)$. Show that if $\varphi$ is
\term[ergodic!flow]{ergodic}, then $\varphi'$ is ergodic.

\item \label{PolyProportion}
 Given any natural number~$d$, and any $\delta > 0$, show there is some
$\epsilon > 0$, such that if 
 \begin{itemize}
 \item $f(x)$ is any real \term[function!polynomial]{polynomial} of degree
$\le d$, 
 \item $C \in \real^+$,
 \item $[k, k + \ell]$ is any real interval,
 and
 \item $|f(t)| < C$ for all $t \in [k,k + \ell]$, 
 \end{itemize}
 then $|f(t)| < (1 + \delta) C$ for all $t \in [k, k + (1+\epsilon) \ell]$.

\item \label{BottomLeftLarge}
 Given positive constants $\epsilon < L$, show there exists $\epsilon_0 >
0$, such that if $|\alpha|,|\mathsf{b}|,|\mathsf{c}|,|\mathsf{d}| <
\epsilon_0$, and $N > 0$, and we have
 $$ \mbox{$| \mathsf{c}  - (\alpha -\mathsf{d})t - \mathsf{b} t^2| < L$
for all $t \in [0,N]$,} $$
 then $|\alpha + \mathsf{b} t| + |\mathsf{d} - \mathsf{b} t | < \epsilon$
for all $t \in [0,N]$.

\item \label{NoFP->d>0}
 Suppose $\psi$ is a \index{homeomorphism}homeomorphism of a compact
metric space $(X,d)$, and that $\psi$ has no \term[fixed point]{fixed
points}. Show there exists $\epsilon > 0$, such that,  for all $x \in X$,
we have $d \bigl( \psi(x), x \bigr) > \epsilon$.

\item \label{ProbMeasIsRegular}
 (Probability measures are \term[measure!regular]{regular})
  Suppose 
 \begin{itemize}
 \item $X$ is a metric space that is \term{separable} and \term{locally
compact},
 \item $\mu$ is a Borel probability measure on~$X$,
 \item $\epsilon > 0$,
 and
 \item $A$ is a measurable subset of~$X$.
 \end{itemize}
 Show:
 \begin{enumerate}
 \item \label{ProbMeasIsRegular-reg}
 there exist a compact set~$C$ and an open set~$V$,
such that $C \subset A \subset V$ and $\mu(V \smallsetminus C) < \epsilon$,
 and
 \item \label{ProbMeasIsRegular-func}
 there is a continuous function~$f$ on~$X$, such that 
 $$ \mu \{\, x \in X \mid \chi_A(x) \neq f(x)\,\} < \epsilon ,$$
 where $\chi_A$ is the characteristic function of~$A$.
 \end{enumerate}
 \hint{Recall that ``\defit{separable}" means $X$ has a
countable, dense subset, and that ``\defit{locally compact}" means every
point of~$X$ is contained in an open set with compact closure.
 \pref{ProbMeasIsRegular-reg} Show the collection~$\mathcal{A}$ of
sets~$A$ such that $C$ and~$V$ exist for every~$\epsilon$ is a
$\sigma$-algebra.
 \pref{ProbMeasIsRegular-func} Note that
 $$ \frac{d(x,X \smallsetminus V)}{d(x,X \smallsetminus V) + d(x,C)} $$
  is a continuous function that is $1$ on~$C$ and~$0$ outside of~$V$.}

\item \label{LusinThmEx}
 (\term[Theorem!Lusin]{Lusin's Theorem}) Suppose 
 \begin{itemize}
 \item $X$ is a metric space that is \term{separable} and \term{locally
compact},
 \item $\mu$ is a Borel probability measure on~$X$,
 \item $\epsilon > 0$,
 and
 \item $\psi \colon X \to \real$ is measurable.
 \end{itemize}
 Show there is a continuous function~$f$ on~$X$, such that 
 $$ \mu \{\, x \in X \mid \psi(x) \neq f(x)\,\} < \epsilon .$$
 \nolinehint{Construct step functions $\psi_n$ that \term{converge
uniformly} to~$\psi$ on a set of measure $1 - (\epsilon/2)$.
 (Recall that a \defit[function!step]{step function} is a linear
combination of  characteristic functions of sets.)
 Now~$\psi_n$ is equal to a continuous function $f_n$ on a set
of measure $1 - 2^{-n}$ \cfex{ProbMeasIsRegular}. Then $\{f_n\}$
converges to~$f$ \term[converge uniformly]{uniformly} on a set of measure
$> 1 - \epsilon$.}

\item \label{Time1Dist}
 Suppose 
 \begin{itemize}
 \item $(X,d)$ is a metric space, 
 \item $\mu$ is a probability \term[measure!invariant]{measure} on~$X$,
 \item $\varphi_t$ is an \term[ergodic!flow]{ergodic}, continuous,
\term[flow!measure-preserving]{measure-preserving} flow on~$X$,
 and
 \item $\epsilon > 0$.
 \end{itemize}
 Show that either
 \begin{enumerate}
 \item some orbit of $\varphi_t$ has measure~$1$,
 or
 \item there exist $\delta > 0$ and a compact subset $K$ of~$X$, such that
 \begin{itemize}
 \item $\mu(K) > 1 - \epsilon$
 and
 \item $d \bigl( \varphi_t(x), x \bigr) > \delta$, for all $t \in
(1-\delta,1+\delta)$.
 \end{itemize}
 \end{enumerate}

\item \label{JoinEgsConjCont}
 Show that the maps $\psi_1$ and $\psi_2$ of Eg.~\fullref{JoinEgs}{conj}
are well defined and continuous.

\item \label{RatnerJoiningsEx}
 Derive \term[Ratner!Joinings Theorem]{Cor.~\ref{RatnerJoin}} from
\term[Ratner's Theorems!Measure Classification]{Ratner's Measure
Classification Theorem}.

\item \label{Cover->FinFib}
 Show that if $\hat\mu$ is a \term[self-joining!finite cover]{finite-cover
self-joining}, then there is a $\hat\mu$-conull subset~$\Omega$ of $X
\times X$, such that $\bigl( \{x\} \times X \bigr) \cap \Omega$ and
$\bigl( X \times \{x\} \bigr) \cap \Omega$ are finite, for every $x \in X$.

\item \label{JoiningsThmRigor}
 Write a rigorous (direct) proof of \term[Ratner!Joinings
Theorem]{Cor.~\ref{JoiningsThm}}, by choosing appropriate conull subsets
of~$X$, and so forth.
 \hint{You may assume (without proof) that there is a compact subset~$K$
of $X \times X$, such that $\mu(K) > 0.99$ and $K \cap \xi_s(K) =
\emptyset$ for all $s \in \real$ with $(\xi_s)_* \hat\mu \neq \hat\mu$
\cf{TechPtJoinings}.}

\item \label{VertInvJoining}
 Verify that $\hat\mu$ must be the \term[self-joining!product]{product
joining} in Case~\ref{JoiningsThmVertCase} of the proof of
\term[Ratner!Joinings Theorem]{Cor.~\ref{JoiningsThm}}.

\item \label{CommuteInvariant}
 Suppose
 \begin{itemize}
 \item $\varphi_t$ is a \term[flow!measurable]{(measurable) flow} on a
measure space~$X$,
 \item $\mu$ is a $\varphi_t$-\term[measure!invariant]{invariant
probability measure} on~$X$,
 and
 \item $\psi \colon X \to X$ is a Borel map that commutes with $\varphi_t$.
 \end{itemize}
 Show that $\psi_* \mu$ is $\varphi_t$-\term[measure!invariant]{invariant}.

\item \label{SingPartOfMeas}
 Suppose $\mu$ and~$\nu$ are probability measures on a measure space~$X$. 
Show $\nu$ has a unique decomposition $\nu = \nu_1 + \nu_2$, where $\nu_1
\perp \mu$ and $\nu_2 = f \mu$, for some $f \in L^1(\mu)$.
 (Recall that the notation 
 \nindex{$\mu_1 \perp \mu_2$ = measures $\mu_1$ and $\mu_2$ are singular
to each other}{$\mu_1 \perp \mu_2$}
 means the measures $\mu_1$ and~$\mu_2$ are
\defit[measure!singular]{singular} to each other; that is, some
$\mu_1$-conull set is $\mu_2$-null, and vice-versa.)
 \hint{The map $\phi \mapsto \int \phi\, d\mu$ is a linear functional on
$L^2(X, \mu + \nu)$, so it is represented by integration against a function
$\psi \in L^2(X, \mu + \nu)$. Let $\nu_1$ be the restriction of~$\nu$ to
$\psi^{-1}(0)$, and let $f = (1 - \psi)/\psi$.}

\item \label{ErgSingular}
 Suppose
 \begin{itemize}
 \item $\varphi_t$ is a \term[flow!measurable]{(measurable) flow} on a
space~$X$,
 and
 \item $\mu_1$ and~$\mu_2$ are two different
\term[ergodic!measure]{ergodic},
$\varphi_t$-\term[measure!invariant]{invariant} probability measures
on~$X$. 
 \end{itemize}
 Show that $\mu_1$ and~$\mu_2$ are \defit[measure!singular]{singular} to
each other.
 \hint{Exer.~\ref{SingPartOfMeas}.}

\item \label{JoiningShrinkCpct}
 Suppose
  \begin{itemize}
 \item $X$ is a locally compact, \term{separable} metric space,
 \item $\mu$ is a probability measure on~$X$,
 and
 \item $\psi \colon X \to X$ is a Borel map, such that $\psi_* \mu$ and
$\mu$ are \term[measure!singular]{singular} to each other.
 \end{itemize}
 Show:
 \begin{enumerate}
 \item There is a conull subset~$\Omega$ of~$X$, such that
$\psi^{-1}(\Omega)$ is disjoint from~$\Omega$.
 \item For every $\epsilon > 0$, there is a compact subset~$K$ of~$X$,
such that $\mu(K) > 1 - \epsilon$ and $\psi^{-1}(K)$ is disjoint from~$K$.
 \end{enumerate}

\end{exercises}

\section{The Shearing Property for larger groups}
\label{GenShearSect}

If $G$ is \index{SL(3,R)*$\SL(3,\real)$}$\SL(3,\real)$, or some other
group larger than $\SL(2,\real)$, then the \term[shearing
property]{Shearing Property} is usually not true as stated in
\pref{ShearingSL2R} or \pref{ShearingSL2R(1step)}. This is because the
\term{centralizer} of the subgroup~$u_t$ is usually larger than $\{u_t\}$. 
 
\begin{eg}
 If $y = xq$, with $q \in C_G(u_t)$, then $u^{-t} q u^t =
q$ for all~$t$, so, contrary to \pref{ShearingSL2R(1step)}, the points~$x$
and~$y$ move together, along parallel orbits; there is no \term{relative
motion} at all.
 \end{eg}

In a case where there is \term{relative motion} (that is, when $q \notin
C_G(u^t)$), the fastest relative motion will usually not be along the
orbits of~$u^t$, but, rather, along some other direction in the
\term{centralizer} of~$u^t$. (We saw an example of this in the proof that
\term[self-joining!finite fibers]{self-joinings} have finite fibers
\seecor{JoiningsThm}: under the \term[unipotent!flow]{unipotent flow}
$\eta_t \times \eta_t$, the points $(x,a)$ and $(y,b)$ move apart in the
direction of the flow~\term[flow!vertical]{$\xi_t$}, not~$\eta_t \times
\eta_t$.)

\begin{prop}[{(\term[shearing property]{Generalized Shearing Property})}]
\label{GenShearing}
 The fastest \term{relative motion} between two \term[point!nearby]{nearby
points} is along some direction in the \term{centralizer} of~$u^t$.

 More precisely, if
 \begin{itemize}
 \item $\{u^t\}$ is a \term[unipotent!subgroup!one-parameter]{unipotent
one-parameter subgroup} of~$G$,
 and
 \item $x$ and~$y$ are \term[point!nearby]{nearby points} in $\Gamma
\backslash G$, 
 \end{itemize}
 then either
 \begin{enumerate}
 \item there exists $t > 0$ and $c \in C_G(u^t)$, such that 
 \begin{enumerate}
 \item $\|c\| = 1$,
 and
 \item $x u^t \approx y u^t c$,
 \end{enumerate}
 or
 \item there exists $c \in C_G(u^t)$, with $c \approx \Id$, such that $y =
xc$.
 \end{enumerate}
 \end{prop}

\begin{proof}[Proof \rm(\emphit{Requires some Lie theory}).]
  Write $y = xq$, with $q \approx \Id$. It is easiest to work with
\term[exponential!coordinates]{exponential coordinates} in the \term{Lie
algebra}; for $g \in G$ (with $g$ near~$\Id$), let 
 \nindex{$\lie g$ = element of~$\Lie G$ with $\exp \lie g = g$}
 $\lie g$ be the (unique) small element of~$\Lie G$ with $\exp \lie g =
g$. In particular, choose
 \begin{itemize}
 \item $\lie u \in \Lie G$ with $\exp (t \lie u) = u^t$,
 and
 \item
  $\lie q \in \Lie G$ with $\exp \lie q = q$.
 \end{itemize}
 Then
 \begin{align*}
 \lie{u^{-t} q u^t}
 &= \lie q (\Ad u^t)
 = \lie q \exp \bigl( \ad (t \lie u) \bigr) \\
 &= {\textstyle \lie q
 + \lie q (\ad \lie u) t
 + \frac{1}{2} \lie q (\ad \lie u)^2 t^2
 + \frac{1}{6} \lie q (\ad \lie u)^3 t^3
 +
\cdots}
% &= {\textstyle \lie q
% \pm t [\lie q, \lie u]
% \pm \frac{1}{2} t^2 [\lie q, \lie u, \lie u]
% \pm \frac{1}{6} t^3 [\lie q, \lie u, \lie u, \lie u]
% \pm
%\cdots}
 . \end{align*}
 For large~$t$, the \term[largest!term]{largest term} is the one with the
\term{highest power of~$t$}; that is, the last nonzero term $\lie q (\ad
\lie u)^k$. Then
 $$ [\lie q (\ad \lie u)^k, \lie u]
 = \bigl( \lie q (\ad \lie u)^k \bigr) (\ad \lie u)
 = \lie q (\ad \lie u)^{k+1}
 = 0$$
 (because the next term does not appear), so $\lie q (\ad \lie u)^k$ is in
the \term{centralizer} of~$u^t$.
 \end{proof}

The above proposition shows that the direction of \term[relative
motion!fastest]{fastest relative motion} is along the \term{centralizer}
of~$u^t$. This direction may or may not belong to~$\{u^t\}$ itself. In the
proof of \term[Ratner's Theorems]{Ratner's Theorem}, it turns out that we wish to ignore
motion \emph{along} the orbits, and consider, instead, only the
\index{component!transverse}component of the relative motion that is
\defit[transverse divergence]{transverse} (or perpendicular) to the
orbits of the flow. This direction, by definition, does not belong to
$\{u^t\}$. It may or may not belong to the \term{centralizer} of~$\{u^t\}$.
 \index{divergence, transverse|indsee{transverse~divergence}}

\begin{eg}
 Assume $G = \SL(2,\real)$, and suppose $x$ and~$y$ are two points in
$\Gamma \backslash G$ with $x \approx y$. Then, by continuity, $x_t \approx
y_t$ for a long time. Eventually, we will be able to see a difference
between $x_t$ and~$y_t$. The \term[shearing property!in
$\SL(2,\real)$]{Shearing Property} \pref{ShearingSL2R} tells us that, when
this first happens, $x_t$ will be indistinguishable from some point on the
orbit of~$y$; that is, $x_t \approx y_{t'}$ for some~$t'$. This will
continue for another long time (with $t'$ some function of~$t$), but we can
expect that $x_t$ will eventually diverge from the orbit of~$y$ --- this
is \defit[transverse divergence]{transverse divergence}. (Note that this
transverse divergence is a \term[second order effect]{second-order
effect}; it is only apparent after we mod out the relative motion along
the orbit.) Letting $y_{t'}$ be the point on the orbit of~$y$ that is
closest to~$x_t$, we write $x_t = y_{t'} g$ for some $g \in G$. Then $g -
\Id$ represents the transverse divergence. When this transverse divergence
first becomes macroscopic, we wish to understand which of the matrix
entries of $g - \Id$ are macroscopic.

 In the matrix on the RHS of Eq.~\pref{ConjByU}, we have already observed
that the \term[largest!entry]{largest entry} is in the bottom left corner,
the direction of~$\{u^t\}$. If we ignore that entry, then the two diagonal
entries are the largest entries. The diagonal corresponds to the
subgroup~$\{a^t\}$ (or, in geometric terms, to the
\term[geodesic!flow]{geodesic flow}~$\gamma_t$). Thus, the fastest
\defit[transverse divergence]{transverse} divergence is in the
direction of $\{a^t\}$. Notice that $\{a^t\}$
\term[normalizer]{normalizes}~$\{u^t\}$ \seeex{ANormsU}.
 \end{eg}

\begin{prop} \label{High1DTransN}
 The fastest \defit[transverse divergence]{transverse} motion is
along some direction in the \term{normalizer} of~$u^t$.
 \end{prop}

\begin{proof}
 In the calculations of the proof of Prop.~\ref{GenShearing}, any term
that belongs to~$\Lie U$ represents motion along~$\{u^t\}$. Thus, the
fastest transverse motion is represented by the \emph{last term} $\lie q
(\ad \lie u)^k$ that is \emph{not} in~$\Lie U$. Then $\lie q (\ad \lie
u)^{k+1} \in \Lie U$, or, in other words, 
 $$[\lie q (\ad \lie u)^k, \lie u] \in \Lie U .$$
 Therefore $\lie q (\ad \lie u)^k$ normalizes~$\Lie U$.
 \end{proof}

By combining this observation with ideas from the proof that joinings have
finite fibers\index{self-joining!finite fibers} \seecor{JoiningsThm},
we see that the fastest \term[transverse divergence]{transverse
divergence} is almost always in the direction of 
 \nindex{$\Stab_G(\mu)$ = stabilizer = $\{\mbox{elements of
$G$ that preserve $\mu$}\}$}
 \index{stabilizer!of a measure}$\Stab_G(\mu)$, the subgroup
consisting of elements of~$G$ that preserve~$\mu$. More precisely:

\begin{cor} \label{FastestInStab1D} \index{shearing property}
 There is a conull subset $X'$ of~$X$, such that, for all $x,y \in X'$,
with $x \approx y$, the fastest \term[transverse divergence]{transverse}
motion is along some direction in \term[stabilizer!of a
measure]{$\Stab_G(\mu)$}.
 \end{cor}

\begin{proof}
  Because the fastest \term[transverse divergence]{transverse}
motion is along the \term{normalizer}, we know that
 $$y u^{t'} \approx x u^t c ,$$
 for some $t,t' \in \real$ and $c \in N_G(u^t)$. 

Suppose $c \notin \Stab_G(\mu)$\index{stabilizer!of a measure}. Then, as
in the proof of \term[Ratner!Joinings Theorem]{\pref{JoiningsThm}}, we may
assume $x u^t, y u^{t'} \in K$, where $K$ is a large compact set, such
that $K \cap K c = \emptyset$. (Note that $t'$ is used, instead of~$t$, in
order to eliminate relative motion \emphit{along} the $\{u^t\}$-orbit.) We
have $d(K, Kc) > 0$, and this
contradicts the fact that $x u^t c \approx y u^{t'}$.
 \end{proof}

\begin{rem}
 We note an important difference between the preceding two results:
 \begin{enumerate}
 \item Proposition~\ref{High1DTransN} is purely algebraic, and applies to
all $x,y \in \Gamma \backslash G$ with $x \approx y$.
 \item Corollary~\ref{FastestInStab1D} depends on the measure~$\mu$ --- it
applies only on a conull subset of $\Gamma \backslash G$.
 \end{enumerate}
 \end{rem}

We have considered only the case of a
\term[unipotent!subgroup!one-parameter|(]{one-parameter
subgroup}~$\{u^t\}$, but, for the proof of \term[Ratner's Theorems!Measure
Classification]{Ratner's Theorem} in general, it is important to know that
the analogue of Prop.~\ref{High1DTransN} is also true for actions of larger
\term[unipotent!subgroup]{unipotent subgroups}~$U$: 
 \begin{equation} \label{HighUTransN}
 \begin{matrix}
 \mbox{the fastest \term[transverse divergence]{transverse} motion
is along} \\
 \mbox{some direction in the normalizer of~$U$.} 
 \end{matrix}
 \end{equation}
 \index{normalizer}
 To make a more precise statement, consider two points $x,y \in X$, with
$x \approx y$. When we apply larger and larger elements~$u$ of~$U$ to~$x$,
we will eventually reach a point where we can see that $xu \notin yU$.
When we first reach this stage, there will be an element~$c$ of the
\term{normalizer} $N_G(U)$, such that $x u c$ appears to be in~$y U$; that
is,
 \begin{equation} \label{shear:xun=yu'}
 \mbox{$x u c \approx y u'$, for some $u' \in U$.} 
 \end{equation}

This implies that the analogue of Cor.~\ref{FastestInStab1D} is also true:

\begin{cor} \label{FastestInStab}
 There is a conull subset $X'$ of~$X$, such that, for all $x,y \in X'$,
with $x \approx y$, the fastest \term[transverse divergence]{transverse}
motion to the $U$-orbits is along some direction in
\term[stabilizer!of a measure]{$\Stab_G(\mu)$}.
 \end{cor}

To illustrate the importance of these results, let us prove the following
special case of \term[Ratner's Theorems!Measure Classification]{Ratner's
Measure Classification Theorem}. It is a major step forward. It shows, for
example, that if $\mu$ is not \term[support!of a measure]{supported} on a single $u^t$-orbit,
then there must be other translations in~$G$ that preserve~$\mu$.

\begin{prop} \label{RatMeas-S=U}
 Let
 \begin{itemize}
 \item $\Gamma$ be a \term{lattice} in a \term{Lie group}~$G$,
 \item $u^t$ be a \term[unipotent!subgroup!one-parameter]{unipotent
one-parameter subgroup} of~$G$,
 and
 \item $\mu$ be an \term[ergodic!measure|(]{ergodic}
$u^t$-\term[measure!invariant|(]{invariant probability measure} on $\Gamma
\backslash G$.
 \end{itemize}
 If $U = \Stab_G(\mu)^\circ$\index{stabilizer!of a
measure} is \term[unipotent!subgroup]{unipotent}, then
$\mu$ is \term[support!of a measure]{supported} on a single $U$-orbit.
 \end{prop}

\begin{proof} This is similar to the proof that joinings have finite
fibers\index{self-joining!finite fibers} \seecor{JoiningsThm}. We ignore
some details (these may be taken to be exercises for the reader). For
example, let us ignore the distinction between
$\Stab_G(\mu)$\index{stabilizer!of a measure} and its identity
\index{component!identity}component $\Stab_G(\mu)^\circ$
\seeex{RatMeas-S=U(IdComp)}.

 By \term[ergodic!action]{ergodicity}, it suffices to find a $U$-orbit
of positive measure, so let us suppose all $U$-orbits have measure~$0$.
Actually, let us make the stronger assumption that all
\term[normalizer]{$N_G(U)$}-orbits have measure~$0$. This will lead to a
contradiction, so we can conclude that $\mu$ is \term[support!of a
measure]{supported} on an orbit of $N_G(U)$. It is easy to finish from
there \seeex{RatMeas-S=U(N(U))}.

By our assumption of the preceding paragraph, for almost every $x \in
\Gamma \backslash G$, 
 there exists $y \approx x$, such that
 \begin{itemize}
 \item $y \notin x N_G(U)$,
 and
 \item $y$ is in the \term[support!of a measure]{support} of~$\mu$.
 \end{itemize}
 Because $y \notin x \, N_G(U)$, the $U$-orbit of~$y$ has nontrivial
\term[transverse divergence]{transverse} divergence from the $U$-orbit
of~$x$ \seeex{notN(U)->transverse}, so
 $$y u' \approx xu c ,$$
 for some $u,u' \in U$ and $c \notin U$. From
Cor.~\ref{FastestInStab}, we know that $c \in
\Stab_G(\mu)$\index{stabilizer!of a measure}. This contradicts the fact
that $U = \Stab_G(\mu)$.
 \end{proof}

\begin{exercises}

 \item \label{RatMeas-S=U(IdComp)}
 The proof we gave of Prop.~\ref{RatMeas-S=U} assumes that
\term[stabilizer!of a measure]{$\Stab_G(\mu)$} is
\term[unipotent!subgroup]{unipotent}. Correct the proof to use only the
weaker assumption that $\Stab_G(\mu)^\circ$ is unipotent.

\item \label{ErgForNormSubgrp}
 Suppose 
 \begin{itemize}
 \item $\Gamma$ is a closed subgroup of a \term{Lie group}~$G$,
 \item $U$ is a \term[unipotent!subgroup]{unipotent},
\index{subgroup!normal}normal subgroup of~$G$,
 and
 \item $\mu$ is an \term[ergodic!measure]{ergodic}
$U$-\term[measure!invariant]{invariant probability measure} on $\Gamma
\backslash G$.
 \end{itemize}
 Show that $\mu$ is \term[support!of a measure]{supported} on a single
orbit of \term[stabilizer!of a measure]{$\Stab_G(\mu)$}.
 \hint{For each \term[normalizer]{$g \in N_G(U)$}, such that $g \notin
\Stab_G(\mu)$, there is a conull subset $\Omega$ of $\Gamma \backslash G$,
such that $\Omega \cap g \Omega = \emptyset$
\seeexs{ErgSingular}{JoiningShrinkCpct}. You may assume, without proof,
that this set can be chosen independent of~$g$: there is a conull subset
$\Omega$ of $\Gamma \backslash G$, such that if $g \in N_G(U)$ and $g
\notin \Stab_G(\mu)$, then $\Omega \cap g \Omega = \emptyset$. (This will
be proved in \pref{KnotKc}.)}

\item \label{RatMeas-S=U(N(U))}
 Suppose
 \begin{itemize}
 \item $\Gamma$ is a \term{lattice} in a \term{Lie group}~$G$,
 \item $\mu$ is a $U$-\term[measure!invariant]{invariant probability
measure} on $\Gamma \backslash G$,
 and
 \item $\mu$ is \term[support!of a measure]{supported} on a single
\term[normalizer]{$N_G(U)$}-orbit.
 \end{itemize}
 Show that $\mu$ is \term[support!of a measure]{supported} on a single $U$-orbit.
 \hint{Reduce to the case where $N_G(U) = G$, and use
Exer.~\ref{ErgForNormSubgrp}.}

\item \label{notN(U)->transverse}
 In the situation of Prop.~\ref{RatMeas-S=U}, show that if $x,y \in \Gamma
\backslash G$, and $y \notin x \, N_G(U)$, then the $U$-orbit of~$y$ has
nontrivial \term[transverse divergence]{transverse divergence} from
the $U$-orbit of~$x$.

 \end{exercises}

\section{Entropy and a proof for $G = \SL(2,\real)$}
\label{EntropyIntroSect}

The \term[shearing property]Shearing Property (and consequences such as
\pref{RatMeas-S=U}) are an important part of the proof of
\term{Ratner's Theorems}, but there are two additional ingredients. We
discuss the role of \term[entropy!of a dynamical system]{entropy} in this
section. The other ingredient, exploiting the direction of
\term[transverse divergence]{transverse divergence}, is the topic of the
following section.

To illustrate, let us prove \term[Ratner's Theorems!Measure Classification]{Ratner's Measure Classification Theorem}
\pref{RatnerMeas} for the case $G = \SL(2,\real)$:

\begin{thm} \label{RatMeasThmSL2R}
 If
 \begin{itemize}
 \item $G = \SL(2,\real)$,
 \item $\Gamma$ is any \term{lattice} in~$G$,
 and
 \item \term[horocycle!flow]{$\eta_t$} is the usual
\term[unipotent!flow]{unipotent flow} on $\Gamma \backslash G$,
corresponding to the \term[unipotent!subgroup!one-parameter|)]{unipotent
one-parameter subgroup}
 $u^t = \begin{bmatrix}
 1 & 0 \\
 t & 1
 \end{bmatrix}$
 \see{SL2etcNotation},
 \end{itemize}
 then every \term[ergodic!measure]{ergodic}
$\eta_t$-\term[measure!invariant]{invariant probability measure} on $\Gamma
\backslash G$ is \term[homogeneous!measure]{homogeneous}.
 \end{thm}

\begin{proof}
 Let
 \begin{itemize}
 \item $\mu$ be any \term[ergodic!measure]{ergodic}
$\eta_t$-\term[measure!invariant]{invariant probability measure} on
$\Gamma \backslash G$,
 and
 \item $S = \Stab_G(\mu)$\index{stabilizer!of a measure}.
 \end{itemize}
 We wish to show that $\mu$ is \term[support!of a measure]{supported} on a
single $S$-orbit.

Because $\mu$ is $\eta_t$-\term[measure!invariant]{invariant}, we know
that $\{u^t\} \subset S$.
 We may assume $\{u^t\} \neq S^\circ$. (Otherwise, it is obvious that
$S^\circ$ is \term[unipotent!subgroup]{unipotent}, so
Prop.~\ref{RatMeas-S=U} applies.) Therefore, $S^\circ$ contains the
\term[subgroup!one-parameter!diagonal]{diagonal one-parameter subgroup}
 $$ a^s = \begin{bmatrix} e^s & 0 \\ 0 & e^{-s} \end{bmatrix} $$
 \seeex{ContainU}.
 To complete the proof, we will show $S$ also contains the
\term[unipotent!subgroup!opposite]{opposite unipotent subgroup}
 \begin{equation} \label{vOppUnipDefn}
 v^r = \begin{bmatrix}
 1 & r \\
 0 & 1
 \end{bmatrix}
 .
 \end{equation}
 Because $\{u^t\}$, $\{a^s\}$, and~$\{v^r\}$, taken together, generate all
of~$G$, this implies $S = G$, so $\mu$ must be the
$G$-\term[measure!invariant]{invariant} \term[measure!Haar]{(Haar) measure}
on $\Gamma \backslash G$, which is obviously
\term[homogeneous!measure]{homogeneous}.

Because $\{a^s\} \subset S$, we know that $a^s$ preserves~$\mu$. Instead
of continuing to exploit dynamical properties of the unipotent subgroup
$\{u^t\}$, we complete the proof by working with~$\{a^s\}$.

Let \term[geodesic!flow]{$\gamma_s$} be the flow corresponding to~$a^s$
\seeNot{SL2etcNotation}. The map $\gamma_s$ is \emph{not} an
\term{isometry}: 
 \begin{itemize}
 \item $\gamma_s$ \term[infinitesimal!distance]{multiplies
infinitesimal distances} in $u^t$-orbits by $e^{2s}$, 
 \item $\gamma_s$ \term[infinitesimal!distance]{multiplies
infinitesimal distances} in $v^r$-orbits by $e^{-2s}$,
 and
 \item $\gamma_s$ does act as an isometry on $a^s$-orbits; it
\term[infinitesimal!distance]{multiplies infinitesimal distances} along
$a^s$-orbits by~$1$
 \end{itemize}
 \seeex{GeodCommRelns}. The map $\gamma_s$ is \term{volume preserving}
because these factors cancel exactly: $e^{2s} \cdot e^{-2s} \cdot 1 = 1$.

The fact that $\gamma_s$ preserves the usual volume form on $\Gamma
\backslash G$ led to the equation $e^{2s} \cdot e^{-2s} \cdot 1 = 1$. 
 Let us find the analogous conclusion that results from the fact that
$\gamma_s$ preserves the measure~$\mu$:
 \begin{itemize}
 \item Because $\{a^s\}$ \term[normalizer]{normalizes} $\{u^t\}$
\seeex{ANormsU}, 
 $$ B = \{\, a^s u^t \mid s,t \in \real\,\}$$
 is a subgroup of~$G$. 
 \item Choose a small (2-dimensional)\index{dimension!of a manifold}
\term{disk}~$D$ in some $B$-orbit.
 \item For some (fixed) small $\epsilon > 0$, and each $d \in D$, let 
$\mathcal{B}_d = \{\, d v^r \mid 0 \le r \le \epsilon \,\}$.
 \item Let $\mathcal{B} = \bigcup_{d \in D} \mathcal{B}_d$.
 \item Then $\mathcal{B}$ is the disjoint union of the fibers $\{
\mathcal{B}_d\}_{d \in D}$, so the restriction $\mu|_{\mathcal{B}}$ can be
decomposed as an \term[measure!direct integral]{integral of probability
measures} on the \term[fiber]{fibers}:
 $$ \mu|_{\mathcal{B}} = \int_D \mu_d \, \nu(d) ,$$
 where $\nu_d$ is a probability measure on~$\mathcal{B}_d$
\see{FiberMeasures}.
 \item The map $\gamma_s$ \term[multiplies areas]{multiplies areas} in~$D$
by $e^{2s} \cdot 1 = e^{2s}$.
 \item Then, because $\mu$ is
$\gamma_s$-\term[measure!invariant]{invariant}, the contraction along the
fibers $\mathcal{B}_d$ must exactly cancel this: for $X \subset
\mathcal{B}_d$, we have
 $$ \mu_{\gamma_s(d)} \bigl( \gamma_s(X) \bigl) = e^{-2s} \mu_d(X) .$$
 \end{itemize}
 \smallskip
 The conclusion is that the fiber measures $\mu_d$ \term[scale like
Lebesgue measure]{scale} exactly like the \term[measure!Lebesgue]{Lebesgue
measure} on $[0,\epsilon]$. This implies, for example, that $\mu_d$ cannot
be a \term[point!mass]point mass. In fact, one can use this conclusion to
show that $\mu_d$ must be precisely the \term[measure!Lebesgue]{Lebesgue
measure}. (From this, it follows immediately that $\mu$ is the
\term[measure!Haar]{Haar measure} on $\Gamma \backslash G$.) As will be
explained below, the concept of \term[entropy!of a dynamical
system]{entropy} provides a convenient means to formalize the argument.
 \end{proof}

\begin{notation}
 As will be explained in Chap.~\ref{EntropyChap}, one can define the
\defit[entropy!of a dynamical system]{entropy} of any
measure-preserving transformation on any measure space. (Roughly speaking,
it is a number that describes how quickly orbits of the transformation
diverge from each other.)
 For any $g \in G$ and any $g$-\term[measure!invariant]{invariant}
probability measure~$\mu$ on $\Gamma \backslash G$, let
 \nindex{$h_\mu(g)$ = entropy of the translation by~$g$}
 $h_\mu(g)$ denote the \term[entropy!of a dynamical system]{entropy} of
the translation by~$g$.
 \end{notation}

A general lemma relates \term[entropy!of a dynamical system]{entropy} to
the rates at which the flow \term[entropy!vs.\ stretching]{expands the
volume form} on certain transverse foliations \see{EntLem}. In the special
case of~$a^s$ in $\SL(2,\real)$, it can be stated as follows. 

\begin{lem} \label{EntropyExpandSL2}
 Suppose $\mu$ is an
\term[geodesic!flow]{$a^s$}-\term[measure!invariant]{invariant probability
measure} on\/ $\Gamma \backslash \SL(2,\real)$.

 We have \term[entropy!of a dynamical system]{$h_\mu(a^s) \le 2|s|$}, with
equality  if and only if $\mu$ is
$\{u^t\}$-\term[measure!invariant]{invariant}.
 \end{lem}

We also have the following general fact \seeex{h(T)=h(Tinv)Ex}:

\begin{lem} \label{h(T)=h(Tinv)}
 The \term[entropy!of a dynamical system]{entropy} of any invertible
\term{measure-preserving} transformation is equal to the \term[entropy!of
a dynamical system]{entropy} of its inverse.
 \end{lem}

Combining these two facts yields the following conclusion, which completes
the proof of Thm.~\ref{RatMeasThmSL2R}.

\begin{cor} \label{Entropy->SL2inv}
 Let $\mu$ be an \term[ergodic!measure]{ergodic}
$\{u^t\}$-\term[measure!invariant]{invariant probability measure} on\/
$\Gamma \backslash \SL(2,\real)$. 

If $\mu$ is
\term[geodesic!flow]{$\{a^s\}$}-\term[measure!invariant]{invariant}, then
$\mu$ is\/ $\SL(2,\real)$-invariant.
 \end{cor}

\begin{proof}
 From the equality clause of Lem.~\ref{EntropyExpandSL2}, we have
$h_\mu(a^s) = 2|s|$, so Lem.~\ref{h(T)=h(Tinv)} asserts that
$h_\mu(a^{-s}) = 2|s|$.

On the other hand, there is an automorphism of $\SL(2,\real)$ that maps
$a^s$ to $a^{-s}$, and interchanges $\{u^t\}$ with $\{v^r\}$. Thus
Lem.~\ref{EntropyExpandSL2} implies:
 $$\begin{matrix}
 \text{$h_\mu(a^{-s}) \le 2|s|$,} \\
 \text{with equality if and only if $\mu$ is $\{v^r\}$-invariant.}
 \end{matrix} $$
 Combining this with the conclusion of the preceding paragraph, we
conclude that $\mu$ is $\{v^r\}$-\term[measure!invariant]{invariant}.

Because $v^r$, $a^s$, and~$u^t$, taken together, generate the entire
$\SL(2,\real)$, we conclude that $\mu$ is
$\SL(2,\real)$-\term[measure!invariant]{invariant}. \end{proof}

\begin{exercises}

\item \label{ContainULie}
 Let $T = \begin{bmatrix} 1&\mathsf{a}&c \\ 0&1&\mathsf{b} \\ 0&0&1
\end{bmatrix}$, with $\mathsf{a},\mathsf{b} \neq 0$.
 Show that if $V$ is a vector subspace of~$\real^3$, such that $T(V)
\subset V$ and $\dim V > 1$, then $\{ (0,*,0) \} \subset V$.

\item \label{ContainU}
 [\emphit{Requires some Lie theory}] Show that if $H$ is a
\term[subgroup!connected|(]{connected subgroup} of $\SL(2,\real)$ that
contains $\{u^t\}$ as a proper subgroup, then $\{a^s\} \subset H$.
 \hint{The \term{Lie algebra} of~$H$ must be invariant under $\Ad_G u^t$.
For the appropriate basis of the Lie algebra $\LieSL(2,\real)$, the desired
conclusion follows from Exer.~\ref{ContainULie}.}

\item \label{GeodCommRelns}
 Show:
 \begin{enumerate}
 \item $\gamma_s (x u^t) = \gamma_s(x) \, u^{e^{2s} t}$,
 \item $\gamma_s (x v^t) = \gamma_s(x) \, v^{e^{-2s} t}$,
 and
 \item $\gamma_s (x a^t) = \gamma_s(x) \, a^t$.
 \end{enumerate}

 \end{exercises}

\section{Direction of divergence and a joinings proof}
\label{JoiningSect}

In \S\ref{Poly+ShearSect}, we proved only a weak form of the
\term[Ratner!Joinings Theorem]{Joinings Theorem} \pref{RatnerJoin}. To
complete the proof of \pref{RatnerJoin} and, more importantly, to
illustrate another important ingredient of Ratner's proof, we provide a
direct proof of the following fact:

\begin{cor}[(Ratner)] \label{FinFib->FinCov}
 If
 \begin{itemize}
 \item $\hat\mu$ is an \term[ergodic!measure]{ergodic} \term{self-joining}
of~$\eta_t$,
 and
 \item $\hat\mu$ has \term{finite fibers},
 \end{itemize}
 then $\hat\mu$ is a \term[self-joining!finite cover]{finite cover}.
 \end{cor}

\begin{notation} \label{SL2xSL2Notn}
 We fix some notation for the duration of this section. Let
 \begin{itemize}
 \item $\Gamma$ be a \term{lattice} in $G = \SL(2,\real)$,
 \item $X = \Gamma \backslash G$,
 \item $U = \{u^t\}$,
 \item $A = \{a^s\}$,
 \item $V = \{v^r\}$,
 \item $\cover{\ } \colon G \to G \times G$ be the natural \term{diagonal
embedding}.
 \end{itemize}
 \end{notation}

\vskip-\smallskipamount

At a certain point in the proof of \term[Ratner's Theorems!Measure Classification]{Ratner's Measure Classification Theorem},
we will know, for certain points $x$ and $y = xg$, that the direction of
fastest \term[transverse divergence]{transverse divergence} of the
orbits belongs to a certain subgroup. This leads to a restriction on~$g$.
In the setting of Cor.~\ref{FinFib->FinCov}, this crucial observation
amounts to the following lemma.

\begin{lem} \label{StildeForSL2xSL2}
 Let $x,y \in X \times X$. If
 \begin{itemize}
 \item $x \approx y$,
 \item $y \in x(V \times V)$, 
 and
 \item the direction of fastest \term[transverse divergence]{transverse
divergence} of the $\cover{U}$-orbits of $x$ and~$y$ belongs to
$\cover{A}$,
 \end{itemize}
 then $y \in x \cover{V}$.
 \end{lem}

\begin{proof}
 We have $y = x v$ for some $v \in V \times V$. Write $x = (x_1,x_2)$, $y
= (y_1,y_2)$ and $v = (v_1,v_2) = (v^{r_1}, v^{r_2})$. To determine the
direction of fastest \term[transverse divergence]{transverse
divergence}, we calculate 
 \begin{align*}
  \cover{u^{-t}} v \cover{u^t} - (\Id,\Id) 
 &= ( u^{-t} v_1 u^t
- \Id, u^{-t} v_2 u^t - \Id) \\
 &\approx 
 \left(
 \begin{bmatrix}
 r_1 t & 0 \\
 -r_1 t^2 & -r_1 t
 \end{bmatrix}
 , 
 \begin{bmatrix}
 r_2 t & 0 \\
 -r_2 t^2 & -r_2 t
 \end{bmatrix}
 \right)
 \end{align*}
 \cf{ConjByU}.
 By assumption, the \term[largest!term]{largest terms} of the two
\index{component!of an ordered pair}components must be (essentially) equal,
so we conclude that $r_1 = r_2$. Therefore $v \in \cover{V}$, as desired.
 \end{proof}

Also, as in the preceding section, the proof of Cor.~\ref{FinFib->FinCov}
relies on the relation of \term[entropy!vs.\ stretching]{entropy} to the
rates at which a flow expands the volume form on transverse foliations.
For the case of interest to us here, the general lemma \pref{EntLem} can
be stated as follows.

\begin{lem} \label{EntJoining}
 Let
 \begin{itemize}
 \item $\hat\mu$ be an $\cover{a^s}$-\term[measure!invariant]{invariant
probability measure} on $X \times X$,
 and
 \item $\hat{V}$ be a \term[subgroup!connected]{connected subgroup} of $V
\times V$.
 \end{itemize}
 Then:
 \begin{enumerate}
 \item  \label{EntJoining-invt}
 If $\hat\mu$ is $\hat{V}$-\term[measure!invariant]{invariant}, then
$h_{\hat\mu}(\cover{a^s}) \ge 2 |s| \dim \hat{V}$.

 \item  \label{EntJoining-fol}
 If there is a conull, Borel subset~$\Omega$ of $X \times X$, such that
$\Omega \cap x (V \times V) \subset x \hat{V}$, for every $x \in \Omega$,
then $h_{\hat\mu}(\cover{a^s}) \le 2 |s| \dim \hat{V}$. 

 \item  \label{EntJoining-equal}
 If the hypotheses of \pref{EntLem-GminInW} are satisfied, and
equality holds in its conclusion, then $\hat\mu$~is
$\hat{V}$-\term[measure!invariant]{invariant}.

 \end{enumerate}
 \end{lem}

\begin{proof}[Proof of Cor.~\ref{FinFib->FinCov}]
 We will show that $\hat\mu$ is
$\cover{G}$-\term[measure!invariant]{invariant}, and is \term[support!of a measure]{supported}
on a single $\cover{G}$-orbit. (Actually, we will first replace
$\cover{G}$ by a conjugate subgroup.) Then it is easy to see that
$\hat\mu$ is a finite cover \seeex{OrbIsFiniteCover}.

It is obvious that $\hat\mu$ is not \term[support!of a measure]{supported} on a
single $\cover{U}$-orbit (because $\hat\mu$ must project to the
\term[measure!Haar]{Haar measure} on each factor of $X \times X$), so, by
combining \pref{HighUTransN} with \pref{FastestInStab} (and
Exer.~\ref{RatMeas-S=U(N(U))}), we see that \term[stabilizer!of a
measure]{$\Stab_{G \times G}(\hat\mu)$} must contain a
\term[subgroup!connected]{connected subgroup} of \term[normalizer]{$N_{G
\times G}(\cover{U})$} that is not contained in~$\cover{U}$. (Note that
$N_{G \times G}(\cover{U}) = \cover{A} \ltimes ( U \times U )$
\seeex{N(diagU)}.) Using the fact that $\hat\mu$ has \term{finite fibers},
we conclude that $\Stab_{G \times G}(\hat\mu)$ contains a
\term[subgroup!conjugate]{conjugate} of~$\cover{A}$
\seeexs{FinFib->NoFactor}{NotInDiagU}. Let us assume, without loss of
generality, that $\cover{A} \subset \Stab_{G \times G}(\hat\mu)$
\seeex{Stab(translate)}; then 
  \begin{equation} \label{FinFib->FinCovPf-N(U)}
  N_{G \times G}(\cover{U}) \cap \Stab_{G \times G}(\hat\mu) = \cover{AU} 
 \end{equation}
 \seeex{NoFactorsSL2xSL2}.
 Combining  \pref{HighUTransN}, \pref{FastestInStab},
and~\pref{FinFib->FinCovPf-N(U)} yields a conull subset $(X \times X)'$ of
$X \times X$, such that if $x,y \in (X \times X)'$ (with $x \approx y$),
then the direction of fastest \term[transverse divergence]{transverse
divergence} between the $\cover{U}$-orbits of $x$ and~$y$ is an element of
$\cover{AU}$.
 Thus, Lem.~\ref{StildeForSL2xSL2} implies that $(X \times X)' \cap x (V
\times V) \subset x \cover{V}$, so an \term[entropy!of a dynamical
system]{entropy} argument, based on Lem.~\ref{EntJoining}, shows that 
 \begin{equation} \label{FinCovPf-Vinvt}
 \mbox{$\hat\mu$ is $\cover{V}$-\term[measure!invariant]{invariant}}
 \end{equation}
 \seeex{EntCalcForFinCovPf}.

Because $\cover{U}$, $\cover{A}$, and $\cover{V}$, taken together,
generate $\cover{G}$, we conclude that $\hat\mu$ is
$\cover{G}$-\term[measure!invariant]{invariant}. Then, because $\hat\mu$
has \term{finite fibers} (and is \term[ergodic!measure]{ergodic}), it is
easy to see that $\hat\mu$ is \term[support!of a measure]{supported} on a single
$\cover{G}$-orbit \seeex{GinvFinCov->Orbit}.
 \end{proof}

\begin{exercises}

\item \label{QuotSmooth}
 Obtain \term[Ratner!Quotients Theorem]{Cor.~\ref{RatnerQuotients}} by
combining Cors.~\ref{QuotFinFib} and~\ref{FinFib->FinCov}.

\item \label{OrbIsFiniteCover}
 In the notation of \pref{FinFib->FinCov} and \pref{SL2xSL2Notn}, show
that if $\hat\mu$ is $\cover{G}$-\term[measure!invariant]{invariant}, and
is \term[support!of a measure]{supported} on a single $\cover{G}$-orbit in $X \times X$, then
$\hat\mu$ is a \term[self-joining!finite cover]{finite-cover joining}.
 \hint{The $\cover{G}$-orbit \term[support!of a measure]{support}ing~$\hat\mu$ can be identified
with $\Gamma' \backslash G$, for some \term{lattice}~$\Gamma'$ in~$G$.}

\item \label{N(diagU)}
 In the notation of \pref{SL2xSL2Notn}, show that
\term[normalizer]{$N_{G\times G}(\cover{U}) = \cover{A} \ltimes (U \times
U)$}.

\item \label{FinFib->NoFactor}
 In the notation of \pref{FinFib->FinCov}, show that if $\hat\mu$ has
\term{finite fibers}, then \index{stabilizer!of a measure}$\Stab_{G \times
G}(\mu) \cap \bigl( G \times \{e\} \bigr)$ is trivial.

\item \label{NotInDiagU}
 In the notation of \pref{SL2xSL2Notn}, show that if $H$ is a
\term[subgroup!connected]{connected subgroup} of $\cover{A} \ltimes (U
\times U)$, such that
 \begin{itemize}
 \item $H \not\subset U \times U$,
 and
 \item $H \cap \bigl( G \times \{e\} \bigr)$ and $H \cap \bigl( \{e\}
\times G \bigr)$ are trivial,
 \end{itemize}
 then $H$ contains a \term[subgroup!conjugate]{conjugate} of~$\cover{A}$.

\item \label{Stab(translate)}
 Suppose
 \begin{itemize}
 \item $\Gamma$ is a \term{lattice} in a \term{Lie group}~$G$,
 \item $\mu$ is a measure on $\Gamma \backslash G$,
 and
 \item $g \in G$.
 \end{itemize}
 Show \index{stabilizer!of a
measure}\index{subgroup!conjugate}$\Stab_G(g_* \mu) = g^{-1} \,
\Stab_G(\mu) \, g$.

\item \label{NoFactorsSL2xSL2}
 In the notation of \pref{SL2xSL2Notn}, show that if $H$ is a subgroup of
$(A \times A) \ltimes (U \times U)$, such that
 \begin{itemize}
 \item $\cover{AU} \subset H$,
 and 
 \item $H \cap \bigl( G \times \{e\} \bigr)$ and $H \cap \bigl( \{e\}
\times G \bigr)$ are trivial,
 \end{itemize}
 then $H = \cover{AU}$.

\item \label{EntCalcForFinCovPf} 
 Establish \pref{FinCovPf-Vinvt}.

\item \label{GinvFinCov->Orbit}
 In the notation of \pref{FinFib->FinCov} and \pref{SL2xSL2Notn}, show
that if $\hat\mu$ is $\cover{G}$-\term[measure!invariant]{invariant}, then
$\hat\mu$ is \term[support!of a measure]{supported} on a single $\cover{G}$-orbit.
 \hint{First show that $\hat\mu$ is \term[support!of a measure]{supported} on a finite union of
$\cover{G}$-orbits, and then use the fact that $\hat\mu$ is
\term[ergodic!measure]{ergodic}.}

\end{exercises}

\section{From measures to orbit closures} \label{MeasuresToOrbitsSect}

In this section, we sketch the main ideas used to derive \term[Ratner's Theorems!Orbit Closure]{Ratner's Orbit Closure Theorem} \pref{RatnerOrb}
from her \term[Ratner's Theorems!Measure Classification]{Measure Classification Theorem}
\pref{RatnerMeas}. This is a generalization of \pref{UniqErg->Minl}, and
is proved along the same lines. Instead of establishing only
\pref{RatnerOrb}, the proof yields the much stronger \term[Ratner's Theorems!Equidistribution]{Equidistribution Theorem}
\pref{RatnerUnifDistPrecise}.

\begin{proof}[{Proof of the Ratner Equidistribution Theorem.}] 
 To simplify matters, let us
 \begin{enumerate} \renewcommand{\theenumi}{\Alph{enumi}}
 \item \label{RUDPPf-AssCpct}
 assume that $\Gamma \backslash G$ is compact, 
 and
 \item \label{RUDPPf-AssErg}
 ignore the fact that not all measures are \term[ergodic!measure]{ergodic}.
 \setcounter{saveenumi}{\arabic{enumi}}
 \end{enumerate}
 Remarks \ref{RatnerUnifDistNonCpctPf} and~\ref{RatnerUnifDistNonErgPf}
indicate how to modify the proof to eliminate these assumptions.

Fix $x \in G$. By passing to a subgroup of~$G$, we may assume
 \begin{enumerate} \renewcommand{\theenumi}{\Alph{enumi}}
 \setcounter{enumi}{\arabic{saveenumi}}
 \item \label{RUDPPf-noS}
 there does not exist any \term[subgroup!connected]{connected}, closed,
proper subgroup~$S$ of~$G$, such that 
 \begin{enumerate}
  \item $\{u^t\}_{t \in \real} \subset S$,
 \item the image $[x S]$ of $x S$ in $\Gamma \backslash G$ is closed, and
has \term[finite volume]{finite $S$-invariant volume}.
 \end{enumerate}
 \end{enumerate}
 We wish to show that the $u^t$-orbit of~$[x]$ is \term{uniformly
distributed} in all of $\Gamma \backslash G$, with respect to the
\term[finite volume]{$G$-invariant volume} on $\Gamma \backslash G$. That
is, letting
 \begin{itemize}
 \item $x_t = x u^t$
 and
 \item $\displaystyle \mu_L(f) = \frac{1}{L} \int_0^L f \bigl( [x_t] \bigr)
\, dt$,
 \end{itemize}
 we wish to show that the measures~$\mu_L$ converge to $\vol_{\Gamma
\backslash G}$, as $L \to \infty$. 

Assume, for simplicity, that $\Gamma \backslash G$ is compact
\see{RUDPPf-AssCpct}. Then the space of probability measures on $\Gamma
\backslash G$ is compact (in an appropriate weak$^*$ topology), so it
suffices to show that 
 \index{measure!limit}
 $$\mbox{if 
 $\mu_{L_n}$ is any convergent sequence,
 then the limit~$\mu_\infty$ is $\vol_{\Gamma \backslash G}$.}$$
 It is easy to see that $\mu_\infty$ is
$u^t$-\term[measure!invariant]{invariant}. Assume for simplicity, that it
is also \term[ergodic!measure]{ergodic} \see{RUDPPf-AssErg}. Then
\term[Ratner's Theorems!Measure Classification]{Ratner's
Measure Classification Theorem} \pref{RatnerMeas} implies that  there
is a \term[subgroup!connected]{connected}, closed subgroup~$S$ of~$G$, and
some point~$x'$ of~$G$, such that
 \begin{enumerate}
 \item $\{u^t\}_{t \in \real} \subset S$,
 \item the image $[x' S]$ of $x' S$ in $\Gamma \backslash G$ is closed, and
has \term[finite volume]{finite $S$-invariant volume},
 and
 \item $\mu_\infty = \vol_{[x' S]}$.
 \end{enumerate}
 It suffices to show that $[x] \in [x' S]$, for then \pref{RUDPPf-noS}
implies that $S = G$, so 
 $$\mu_\infty = \vol_{[x' S]} = \vol_{[x' G]} = \vol_{\Gamma \backslash
G}, $$
 as desired.

To simplify the remaining details, let us assume, for the moment, that $S$
is trivial, so $\mu_\infty$ is the \term[point!mass]{point mass} at the
point~$[x']$. (Actually, this is not possible, because $\{u^t\} \subset
S$, but let us ignore this inconsistency.) This means, for any
neighborhood~$\neigh$ of~$[x']$, no matter how small, that the orbit
of~$[x]$ spends more than 99\% of its life in~$\neigh$. By
\term[polynomial!divergence]{Polynomial Divergence of Orbits}
\cf{PolyOrbsStayClose}, this implies that if we enlarge~$\neigh$ slightly,
then the orbit is \emph{always} in~$\neigh$. Let $\widetilde{\neigh}$ be
the inverse image of~$\neigh$ in~$G$. Then, for some connected
\index{component!connected}component~$\widetilde{\neigh}^\circ$
of~$\widetilde{\neigh}$, we have $x u^t \in \widetilde{\neigh}^\circ$, for
all~$t$. But $\widetilde{\neigh}^\circ$ is a small set (it has the same
diameter as its image~$\neigh$ in $\Gamma \backslash G$), so this implies
that $x u^t$ is a bounded function of~$t$. A
\term[function!polynomial!bounded]{bounded polynomial} is constant, so we
conclude that 
 $$ \mbox{$x u^t = x$ for all $t \in \real$}. $$
 Because $[x']$ is in the closure of the orbit of~$[x]$,  this implies
that $[x] = [x'] \in [x' S]$, as desired.

To complete the proof, we point out that a similar argument applies even if
$S$ is not trivial. We are ignoring some technicalities, but the idea is
simply that the orbit of~$[x]$ must spend more than 99\% of its life very
close to~$[x' S]$. By \term[polynomial!divergence]{Polynomial Divergence of
Orbits}, this implies that the orbit spends all of its life fairly close
to~$[x' S]$. Because the distance to $[x' S]$ is a
\term[function!polynomial]{polynomial} function, we conclude that it is a
constant, and that this constant must be~$0$. So $[x] \in [x' S]$, as
desired.
 \end{proof}

The following two remarks indicate how to eliminate the assumptions
\pref{RUDPPf-AssCpct} and \pref{RUDPPf-AssErg} from the proof of
\pref{RatnerUnifDistPrecise}.

\begin{rem} \label{RatnerUnifDistNonCpctPf}
 If $\Gamma \backslash G$ is not compact, we consider its \term{one-point
compactification} $X = (\Gamma \backslash G) \cup \{\infty\}$. Then
 \begin{itemize}
 \item the set
 \nindex{$\Prob(X)$ = $\{\,\mbox{probability measures on~$X$}\,\}$}
 $\Prob(X)$ of \term[measure!probability]{probability measures} on~$X$ is
compact,
 and
 \item $\Prob(\Gamma \backslash G) = \bigset{ \mu \in \Prob(X) }{
\mu\bigl( \{\infty\} \bigr) = 0 }$.
 \end{itemize}
 Thus, we need only show that the \term[measure!limit]{limit measure}
$\mu_\infty$ gives measure~$0$ to the point~$\infty$. In spirit, this is a
consequence of the \term[polynomial!divergence]{Polynomial Divergence of
Orbits}, much as in the above proof of \pref{RatnerUnifDistPrecise},
putting $\infty$ in the role of~$x'$. It takes considerable ingenuity to
make the idea work, but it is indeed possible. A formal statement of the
result is given in the following theorem.
 \end{rem}

\begin{thm}[(Dani-Margulis)] \label{UnipNotNearInfty}
 Suppose
 \begin{itemize}
 \item $\Gamma$ is a \term{lattice} in a \term{Lie group}~$G$,
 \item $u^t$ is a \term[unipotent!subgroup!one-parameter]{unipotent
one-parameter subgroup} of~$G$,
 \item $x \in \Gamma \backslash G$,
 \item $\epsilon > 0$,
 and
 \item $\lambda$ is the \term[measure!Lebesgue]{Lebesgue measure}
on~$\real$.
 \end{itemize}
 Then there is a compact subset~$K$ of\/ $\Gamma \backslash G$,
such that
 $$ \limsup_{L \to \infty}  \frac{
 \lambda \{\, t \in [0,L] \mid x u^t \notin K \,\}
 }{ L} < \epsilon .$$
 \end{thm}

\begin{rem} \label{RatnerUnifDistNonErgPf}
 Even if the \term[measure!limit]{limit measure}~$\mu_\infty$ is not
\term[ergodic!measure|)]{ergodic}, \term[Ratner's Theorems!Measure
Classification]{Ratner's Measure Classification Theorem} tells us that
each of its \term[ergodic!component]{ergodic components} is
\term[measure!homogeneous]{homogeneous}. That is, for each ergodic
component~$\mu_z$, there exist
 \begin{itemize}
 \item a point $x_z \in G$, and
 \item a closed, \term[subgroup!connected]{connected subgroup} $S_z$
of~$G$,
 \end{itemize}
  such that
 \begin{enumerate}
 \item $\{u^t\}_{t \in \real} \subset S_z$,
 \item the image $[x_z S_z]$ of $x_z S_z$ in $\Gamma \backslash G$ is
closed, and has \term[finite volume]{finite $S_z$-invariant volume},
 and
 \item $\mu_z = \vol_{[x S_z]}$.
 \end{enumerate}
 Arguments from algebra, based on the \term[Theorem!Borel Density]{Borel
Density Theorem}, tell us that:
 \begin{enumerate} \renewcommand{\theenumi}{\alph{enumi}}
 \item up to \term[subgroup!conjugate]{conjugacy}, there are only countable
many possibilities for the subgroups~$S_z$ \seeex{CtblyManyS},
 and
 \item for each subgroup~$S_z$, the point~$x_z$ must belong to a countable
collection of orbits of the \term{normalizer} $N_G(S_z)$.
 \end{enumerate}
 The \defit{singular set} $\singular(u^t)$ corresponding to~$u^t$ is the
union of all of these countably many $N_G(S_z)$-orbits for all of the
possible subgroups~$S_z$. Thus:
 \begin{enumerate}
 \item $\singular(u^t)$ is a countable union of
lower-dimensional\index{dimension!of a manifold}
\term[submanifold]{submanifolds} of $\Gamma \backslash G$,
 and
 \item if $\mu'$ is any $u^t$-\term[measure!invariant]{invariant}
probability measure on $\Gamma \backslash G$, such that $\mu' \bigl(
\singular(u^t) \bigr) = 0$, then $\mu'$ is the
\term[measure!Lebesgue]{Lebesgue measure}.
 \end{enumerate}
 So we simply wish to show that $\mu_\infty\bigl( \singular(u^t) \bigr) =
0$. 

This conclusion follows from the \term[speed!polynomial]{polynomial speed}
of \term[unipotent!flow]{unipotent flows}. Indeed, for every $\epsilon >
0$, because $x \notin \singular(u^t)$, one can show there is an open
neighborhood~$\neigh$ of~$\singular(u^t)$, such that
 \begin{equation} \label{SingSet<Epsilon}
 \frac{\lambda\{\, t \in [0,L] \mid x u^t \in \neigh \,\}}{L} <
\epsilon \text{\quad  for every $L > 0$}, 
 \end{equation}
 where $\lambda$ is the Lebesgue measure on~$\real$.
 \end{rem}

For many applications, it is useful to have the following stronger
(``uniform") version of the Equidistribution Theorem \pref{RatnerUnifDist}:

\begin{thm} \label{DaniMargUnifDist}
 Suppose
 \begin{itemize}
 \item $\Gamma$ is a \term{lattice} in a
\index{subgroup!connected}connected \term{Lie group}~$G$,
 \item $\mu$ is the $G$-\term[measure!invariant|)]{invariant} probability
measure on $\Gamma \backslash G$,
 \item $\{u^t_n\}$ is a sequence of one-parameter
\term[unipotent!subgroup!one-parameter]{unipotent subgroups} of~$G$,
converging to a one-parameter subgroup~$u^t$ {\rm(}that is, $u^t_n \to
u^t$ for all~$t$\/{\rm)},
 \item $\{x_n\}$ is a convergent sequence of points in $\Gamma \backslash
G$, such that $\lim_{n \to \infty} x_n \notin \singular(u^t)$,
 \item $\{L_n\}$ is a sequence of real numbers tending to~$\infty$,
 and
 \item $f$ is any bounded, continuous function on $\Gamma \backslash G$.
 \end{itemize}
 Then
 $$ \lim_{n \to \infty} \frac{1}{L_n} \int_0^{L_n} f (x_n u^t_n) \, dt
 = \int_{\Gamma \backslash G} f \, d\mu .$$
 \end{thm}

\begin{exercises}

\item Reversing the logical order, prove that Thm.~\ref{UnipNotNearInfty}
is a corollary of the \term[Ratner's Theorems!Equidistribution]{Equidistribution Theorem} \pref{RatnerUnifDist}.

\item \label{NormalFP}
 Suppose $S$ is a subgroup of~$G$, and $H$ is a subgroup of~$S$. Show, for
all \term[normalizer]{$g \in N_G(S)$} and all $h \in H$, that $S g h = S
g$.

 \end{exercises}

\section*{Brief history of Ratner's Theorems}
 \index{Ratner's Theorems}
 \index{Ratner's Theorems!Orbit Closure}
 \index{Ratner's Theorems!Measure Classification}
 \index{Ratner's Theorems!Equidistribution}

\begingroup \parskip = 1pt plus 0.2pt

In the 1930's, G.~Hedlund \cite{Hedlund1, Hedlund2, Hedlund3} proved that
if $G = \SL(2,\real)$ and $\Gamma \backslash G$ is compact, then
\index{horocycle!flow}\term[unipotent!flow]{unipotent flows} on $\Gamma
\backslash G$ are \term[ergodic!flow]{ergodic} and
\term[flow!minimal]{minimal}. 

It was not until 1970 that H.~Furstenberg \cite{FurstenbergUnique} proved
these flows are \term[ergodic!uniquely]{uniquely ergodic}, thus
establishing the Measure Classification Theorem for this case. At about
the same time, W.~Parry \cite{ParryProp, ParryIso}  proved an Orbit
Closure Theorem, Measure Classification Theorem, and Equidistribution
Theorem for the case where $G$ is \term[group!nilpotent]{nilpotent}, and
G.~A.~Margulis \cite{MargulisRecur}  used the
\term[speed!polynomial]{polynomial speed} of
\term[unipotent!flow|(]{unipotent flows} to prove the important fact that
unipotent orbits cannot go off to infinity.

Inspired by these and other early results, M.~S.~Raghunathan conjectured a
version of the Orbit Closure Theorem, and showed that it would imply the
\term{Oppenheim Conjecture}. Apparently, he did not publish this
conjecture, but it appeared in a paper of S.~G.~Dani \cite{DaniHoro} in
1981. In this paper, Dani conjectured a version of the
Measure Classification Theorem. 

Dani \cite{DaniHoroRk1} also generalized Furstenberg's Theorem to the case
where $\Gamma \backslash \SL(2,\real)$ is not compact. 
 Publications of  R.~Bowen \cite{Bowen}, S.~G.~Dani \cite{DaniHoroRk1,
DaniHoro}, R.~Ellis and W.~Perrizo \cite{EllisPerrizo}, and  W.~Veech
\cite{Veech} proved further generalizations for the case where the
\term[unipotent!subgroup]{unipotent subgroup} $U$ is
\term[subgroup!horospherical]{horospherical} (see \ref{horospherical}
for the definition). (Results for
\term[subgroup!horospherical]{horosphericals} also follow from a method in
the thesis of G.~A.~Margulis \cite[Lem.~5.2]{MargulisThesis}
\cfex{HoroUniqErg}.)

M.~Ratner began her work on the subject at about this time, proving
her Rigidity Theorem \pref{RatnerRigThm}, Quotients Theorem
\pref{RatnerQuotients}, Joinings Theorem \pref{RatnerJoin}, and other
fundamental results in the early 1980's \cite{RatnerRig, RatnerQuot,
RatnerJoin}. (See \cite{RatnerHyp} for an overview of this early work.) 
 Using Ratner's methods, D.~Witte \cite{WitteRig, WitteZero} generalized
her rigidity theorem to all~$G$.

S.~G.~Dani and J.~Smillie \cite{DaniSmillie} proved the Equidistribution
Theorem when $G = \SL(2,\real)$. S.~G.~Dani \cite{DaniRecur} showed that
unipotent orbits spend only a negligible fraction of their life near
infinity.
 A.~Starkov \cite{StarkovSolv} proved an \term[Ratner's Theorems!Orbit
Closure]{orbit closure theorem} for the case where $G$ is
\term[group!solvable]{solvable}.

Using \term[unipotent!flow]{unipotent} flows, G.~A.~Margulis'
\cite{MargulisOppenheim} proved the \term{Oppenheim Conjecture}
\pref{Oppenheim} on values of \term[quadratic form]{quadratic forms}. He
and S.~G.~Dani \cite{DaniMargulis1, DaniMargulis2, DaniMargulis3} then
proved a number of results, including the first example of an
\term[Ratner's Theorems!Orbit Closure]{orbit closure theorem} for actions
of non-\term[subgroup!horospherical]{horospherical} unipotent subgroups of
a \term[group!semisimple]{semisimple} Lie group --- namely, for ``generic"
\index{unipotent!subgroup!one-parameter}one-parameter unipotent subgroups
of $\SL(3,\real)$. (G.~A~Margulis \cite[\S3.8, top of
p.~319]{MargulisOppenheimSurvey} has pointed out that the methods could
yield a proof of the general case of the \index{Ratner's Theorems!Orbit
Closure}Orbit Closure Theorem.)

Then M.~Ratner \cite{RatnerSolvable, RatnerSS, RatnerMeasure, RatnerOrbit}
proved her amazing theorems (largely independently of the work of
others),  by expanding the ideas from her earlier study of
\index{horocycle!flow}horocycle flows.
   (In the meantime, N.~Shah \cite{ShahGeneric} showed that the
\index{Ratner's Theorems!Measure Classification}Measure Classification
Theorem implied an \term[Ratner's Theorems!Equidistribution]{Equidistribution Theorem} for many cases when $G
= \SL(3,\real)$\index{SL(3,R)*$\SL(3,\real)$}.)

\index{Ratner's Theorems}Ratner's Theorems were soon generalized to
\term[p-adic*$p$-adic]{$p$-adic groups}, by M.~Ratner \cite{RatnerPadic}
and, independently, by G.~A.~Margulis and G.~Tomanov
\cite{MargulisTomanov-Ratner, MargulisTomanovAlmost}. N.~Shah
\cite{ShahUnipGen} generalized the results to subgroups \term{generated by
unipotent elements} \pref{RatOrbH}. (For
\term[subgroup!connected|)]{\emph{connected} subgroups} generated by
unipotent elements, this was proved in Ratner's original papers.)

\endgroup % parskip

\begin{notes}

\notesect{WhatisRatnerSect}

See \cite{BekkaMayer} for an excellent introduction to the general
area of flows on \term[homogeneous!space]{homogeneous spaces}. Surveys at
an advanced level are given in \cite{Dani-review, Dani-Encyclopedia,
KleinbockShahStarkov, MargulisICM, StarkovBook}. Discussions of
\term{Ratner's Theorems} can be found in \cite{BekkaMayer, Dani-review,
Ghys(Bourbaki), RatnerICM, RatnerTata, StarkovBook}. 

Raghunathan's book \cite{RaghunathanBook} is the standard reference for
basic properties of \term{lattice} subgroups.

The dynamical behavior of the \term[geodesic!flow]{geodesic flow} can be
studied by associating a \term{continued fraction} to each point of $\Gamma
\backslash G$. (See \cite{Arnoux} for an elementary explanation.) In this
representation, the fact that some orbit closures are \term{fractal} sets
\pref{FractalOrbit} is an easy observation.

See \cite[Thms.~1.12 and 1.13, pp.~22--23]{RaghunathanBook} for solutions
of Exers.~\ref{DivInG/Gamma} and~\ref{Lattice->ClosedOrbit}.

\notesect{MargOppSect}

\term[Theorem!Margulis (on quadratic forms)]{Margulis' Theorem on values of
quadratic forms} \pref{Oppenheim} was proved in \cite{MargulisOppenheim},
by using \index{unipotent!flow}unipotent flows. For a discussion and
history of this theorem, and its previous life as the \term{Oppenheim
Conjecture}, see \cite{MargulisOppenheimSurvey}. An elementary proof is
given in \cite{DaniMargulis3}, \cite[Chap.~6]{BekkaMayer}, and
\cite{Dani-Oppenheim}.

 % Should find a reference for the connected components of $\SO(Q)$?

\notesect{MeasVersSect}

M.~Ratner proved her \term[Ratner's Theorems!Measure Classification]{Measure Classification Theorem}
\pref{RatnerMeas} in \cite{RatnerSolvable, RatnerSS, RatnerMeasure}. She
then derived her \term[Ratner's Theorems!Equidistribution]{Equidistribution
Theorem} \pref{RatnerUnifDist} and her \term[Ratner's Theorems!Orbit
Closure]{Orbit Closure Theorem} \pref{RatnerOrb} in \cite{RatnerOrbit}. A
derivation also appears in  \cite{DaniMargulis-LimDist}, and an outline
can be found in \cite[\S11]{MargulisTomanov-Ratner}.

In her original proof of the \term[Ratner's Theorems!Measure
Classification]{Measure Classification Theorem}, Ratner only assumed that
$\Gamma$~is \defit[subgroup!discrete]{discrete}, not that it is a
\term{lattice}.
 D.~Witte \cite[\S3]{WitteQuot} observed that discreteness is also not
necessary (Rem.~\fullref{RatMeasThmRems}{notlatt}).

See \cite[\S12]{Phelps-Choquet} for a discussion of
\term[Theorem!Choquet]{Choquet's Theorem}, including a solution of
Exer.~\ref{ChoquetEx}.

\notesect{ApplSect}

The quantitative version \pref{QuantOppThm} of
\term[Theorem!Margulis (on quadratic forms)]{Margulis' Theorem} on values
of quadratic forms is due to A.~Eskin, G.~A.~Margulis, and S.~Mozes
\cite{EskinMargulisMozes}. See \cite{MargulisOppenheimSurvey} for more
discussion of the proof, and the partial results that preceded it.

See \cite{Lindenstrauss-QUE} for a discussion of \index{Quantum!Unique
Ergodicity}Quantum Unique Ergodicity and related results.
Conjecture~\ref{QUEConj} (in a more general form) is due to Z.~Rudnick and
P.~Sarnack. Theorem~\ref{AQUE} was proved by E.~Lindenstrauss
\cite{Lindenstrauss-QUE}. The crucial fact that $h_{\hat\mu}(a_t) \neq 0$
was proved by J.~Bourgain and E.~Lindenstrauss
\cite{BourgainLindenstrauss}.

The results of \S\ref{ApplOh} are due to H.~Oh \cite{OhOppUnip}.

\notesect{Poly+ShearSect}

The \term[Ratner!Rigidity Theorem]{Ratner Rigidity Theorem}
\pref{RatnerRigThm} was proved in \cite{RatnerRig}.
Remark~\fullref{RatRigRem}{Witte} was proved in \cite{WitteRig, WitteZero}.

Flows by diagonal subgroups were proved to be \term[Bernoulli
shift]{Bernoulli} \fullsee{RatRigRem}{Bernoulli} by S.~G.~Dani
\cite{Dani-Bernoulli}, using methods of D.~Ornstein and B.~Weiss
\cite{OrnsteinWeiss}.

The \term[Ratner!Quotients Theorem]{Ratner Quotients Theorem}
\pref{RatnerQuotients} was proved in
\cite{RatnerQuot}, together with Cor.~\ref{QuotFinFib}.

The crucial property \pref{PolyOrbsStayClose} of
\term[polynomial!divergence]{polynomial divergence} was introduced by
M.~Ratner \cite[\S2]{RatnerRig} for \term[unipotent!flow|)]{unipotent
flows} on homogeneous spaces of $\SL(2,\real)$. Similar ideas had
previously been used by Margulis in \cite{MargulisRecur} for more general
unipotent flows.

The Shearing Property (\ref{ShearingSL2R} and \ref{ShearingSL2R(1step)})
was introduced by M.~Ratner \cite[Lem.~2.1]{RatnerQuot} in the proof of her
\term[Ratner!Quotients Theorem]{Quotients Theorem} \pref{RatnerQuotients},
and was also a crucial ingredient in the proof \cite{RatnerJoin} of her
Joinings Theorem \pref{RatnerJoin}. She \cite[Defn.~1]{RatnerJoin} named a
certain precise version of this the
``\index{H-property|indsee{property, H-}}\term[property!H-]{H-property},"
in honor of the \term[horocycle!flow]{horocycle flow}. 

An introduction to \term{Nonstandard Analysis} (the rigorous theory of
\term[infinitesimal]{infinitesimals}) can be found in \cite{Robinson-NonStand} or
\cite{StroyanLuxemburg}.

Lusin's Theorem (Exer.~\ref{LusinThmEx}) and the decomposition of a measure
into a \term[measure!singular]{singular part} and an
\term[measure!absolutely continuous]{absolutely continuous part}
\seeex{SingPartOfMeas} appear in many graduate analysis texts, such as
\cite[Thms.~2.23 and 6.9]{Rudin-RealCplx}.

 \index{SL(2,R)*$\SL(2,\real)$|)}

\notesect{GenShearSect}

The generalization \pref{GenShearing} of the \term[shearing
property]{Shearing Property} to other \term[Lie group]{Lie groups} appears
in \cite[\S6]{WitteRig}, and was called the ``\term[Ratner!property]{Ratner
property}\index{property!Ratner|indsee{Ratner~property}}." The important
extension \pref{HighUTransN} to \term{transverse divergence} of actions of
\index{dimension!of a Lie group}higher-dimensional
\index{subgroup!unipotent}unipotent subgroups is implicit in the
``\index{property!R-}R-property" introduced by M.~Ratner
\cite[Thm.~3.1]{RatnerSolvable}. In fact, the \index{property!R-}R-property
combines \pref{HighUTransN} with \index{polynomial!divergence}polynomial
divergence. It played an essential role in Ratner's proof of the
\index{Ratner's Theorems!Measure Classification}Measure Classification
Theorem.

The arguments used in the proofs of \pref{FastestInStab1D} and
\pref{RatMeas-S=U} appear in \cite[Thm.~4.1]{RatnerSL2}.

\notesect{EntropyIntroSect}

Theorem~\ref{RatMeasThmSL2R} was proved by S.~G.~Dani \cite{DaniHoroRk1},
using methods of H.~Furstenberg \cite{FurstenbergUnique}. Elementary
proofs based on Ratner's ideas (without using \term[entropy!of a
dynamical system]{entropy}) can be found in \cite{RatnerSL2},
\cite[\S4.3]{BekkaMayer}, \cite{Ghys(Bourbaki)}, and
\cite[\S16]{StarkovBook}. 

The \term[entropy!of a dynamical system]{entropy} estimates
\pref{EntropyExpandSL2} and \pref{EntJoining} are special cases of a
result of G.~A.~Margulis and G.~Tomanov
\cite[Thm.~9.7]{MargulisTomanov-Ratner}. (Margulis and Tomanov were
influenced by a theorem of F.~Ledrappier and L.-S.~Young
\cite{LedrappierYoung}.) We discuss the Margulis-Tomanov result in
\S\ref{h(g)Sect}, and give a sketch of the proof in
\S\ref{EntropyEstimateSect}.

\notesect{JoiningSect}

The subgroup~$\cover{V}$ will be called
$\widetilde{S}_-$ in \S\ref{DefineStildeSect}. The proof of
Lem.~\ref{StildeForSL2xSL2} essentially amounts to a verification of
Eg.~\fullref{StildeEg}{SL2xSL2}.

The key point \pref{FinFib->FinCovPf-N(U)} in the proof of
Cor.~\ref{FinFib->FinCov} is a special case of Prop.~\ref{Stilde=S}.

\notesect{MeasuresToOrbitsSect}

G.~A.~Margulis \cite{MargulisRecur} proved a weak version of
Thm.~\ref{UnipNotNearInfty} in 1971. Namely, for each $x \in \Gamma
\backslash G$, he showed there is a compact subset~$K$ of $\Gamma
\backslash G$, such that
  \begin{equation} \label{WhenUinCpct}
 \{\, t \in [0,\infty) \mid [x u^t] \in K \,\}
 \end{equation}
 is unbounded. The argument is elementary, but ingenious.
 A very nice version of the proof appears in
\cite[Appendix]{DaniMargulis3} (and \cite[\S V.5]{BekkaMayer}).

 Fifteen years later, S.~G.~Dani \cite{DaniRecur} refined Margulis' proof,
and obtained \pref{UnipNotNearInfty}, by showing that the set
\pref{WhenUinCpct} not only is unbounded, but has density $> 1 - \epsilon$.
 The special case of Thm.~\ref{UnipNotNearInfty} in which $G =
\SL(2,\real)$\index{SL(2,R)*$\SL(2,\real)$} and $\Gamma =
\SL(2,\integer)$\index{SL(2,Z)*$\SL(2,\integer)$} can be proved easily (see
\cite[Thm.~3.1]{RatnerSL2} or \cite[Thm.~12.2, p.~96]{StarkovBook}).

The uniform version \pref{DaniMargUnifDist} of the Equidistribution
Theorem was proved by S.~G.~Dani and G.~A.~Margulis
\cite{DaniMargulis-LimDist}. The crucial inequality \pref{SingSet<Epsilon}
is obtained from the Dani-Margulis \defit{linearization method} introduced
in \cite[\S3]{DaniMargulis-LimDist}.

 \end{notes}

%% file: RatnerEntropy.tex
\mychapter{Introduction to Entropy} \label{EntropyChap}

 \index{entropy|(}

The entropy of a dynamical system can be intuitively described as a
number that expresses the amount of ``unpredictability" in the system.
Before beginning the formal development, let us illustrate the idea by
contrasting two examples.

\section{Two dynamical systems} \label{TwoEgsSect}

\begin{defn}
 In classical \term[ergodic!theory]{ergodic theory}, a \defit{dynamical
system} (with \term{discrete time}) is an action of~$\integer$ on a measure
space, or, in other words, a measurable transformation $T \colon \Omega
\to \Omega$ on a measure space~$\Omega$. (The intuition is that the points
of~$\Omega$ are in motion. A particle that is at a certain point $\omega
\in \Omega$ will move to a point $T(\omega)$ after a unit of time.) We
assume:
 \begin{enumerate}
 \item $T$ has a (measurable) inverse $T^{-1} \colon \Omega \to \Omega$, 
 and
 \item there is a $T$-\term[measure!invariant]{invariant probability
measure}~$\mu$ on~$\Omega$.
 \end{enumerate}
 The assumption that $\mu$ is \defit[measure!invariant]{$T$-invariant}
means $\mu \bigl( T(A) \bigr) = \mu(A)$, for every measurable subset~$A$
of~$\Omega$. (This generalizes the notion of an \emphit{incompressible
fluid} in fluid dynamics.)
 \end{defn}

\begin{eg}[{(\term[irrational!rotation]{Irrational rotation
of the circle})}]
 Let $\torus = \real / \integer$ be the circle group;  for any $\beta \in
\real$, we have a measurable transformation $T_\beta \colon \torus \to
\torus$ defined by 
 \nindex{$T_\beta$ = irrational rotation}
 $$ T_\beta(t) = t + \beta .$$
 (The usual arc-length \term[measure!Lebesgue]{Lebesgue measure} is
$T_\beta$-\term[measure!invariant]{invariant}.) In physical terms, we have
a circular hoop of circumference~$1$ that is rotating at a speed of~$\beta$
\seefig{hoop}. Note that we are taking the \emph{circumference,} not
the radius, to be~$1$. 

\begin{figure}
 \includegraphics[scale=0.44035]{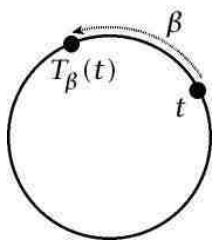}
 \caption{$T_\beta$ moves each point on the the circle a distance of
$\beta$. In other words, $T_\beta$ rotates the circle $360 \beta$
degrees.}
 \label{hoop}
 \end{figure}

 If $\beta$ is \term[irrational!number]{irrational}, it is well known that
every orbit of this dynamical system is uniformly distributed in~$\torus$
\seeex{UnifDistRot} (so the dynamical system is
\defit[ergodic!uniquely]{uniquely ergodic}).
 \end{eg}

\begin{eg}[(\term{Bernoulli shift})]
 Our other basic example comes from the study of \term{coin tossing}.
Assuming we have a fair coin, this is modeled by a two event universe
$\COIN = \{\head, \tail\}$,  in which each event has probability $1/2$. 
The probability space for tossing two coins (independently) is the product
space  $\COIN \times \COIN$, with the product measure (each of the four
possible events has probability $(1/2)^2 = 1/4$). For  $n$~tosses, we take
the product measure on  $\COIN^n$ (each of the $2^n$ possible events has
probability $(1/2)^n$).  Now consider tossing a coin once each day for all
eternity (this is a doubly infinite sequence of coin tosses).  The
probability space is an infinite cartesian product
 $$ \COIN^\infty = \{\, f \colon \integer \to \COIN \,\} $$
 with the \defit[measure!product]{product measure}:  for any two disjoint
finite subsets $\headset$  and~$\tailset$  of~$\integer$, the probability
that  
 \begin{itemize}
 \item $f(n) = \head$ for all $n \in  \headset$,
 and 
 \item $f(n) = \tail$ for all $n \in \tailset$
 \end{itemize}
  is exactly  $(1/2)^{|\headset| + |\tailset|}$. 

 A particular \term[coin tossing!history]{coin-tossing history} is
represented by a single element $f \in  \COIN^\infty$. Specifically,
$f(0)$  is the result of today's coin toss, $f(n)$ is the result of the
toss $n$~days from now (assuming $n > 0$), and $f(-n)$ is the result of
the toss $n$~days ago.  Tomorrow, the history will be represented by an
element $g \in  \COIN^\infty$  that is closely related to~$f$, namely,
$f(n+1) = g(n)$  for every~$n$. (Saying \emph{today} that ``I will toss a
head $n+1$ days from now" is equivalent to saying \emph{tomorrow} that ``I
will toss a head $n$~days from now.") This suggests a dynamical system
(called a \defit{Bernoulli shift}) that consists of translating each
sequence of  $\head$'s  and~$\tail$'s  one notch to the left:
 \nindex{$\Bern$ = Bernoulli shift}
 $$ \mbox{$\Bern \colon \COIN^\infty \to \COIN^\infty$ is defined by 
$(\Bern f)(n) = f(n+1)$} .$$
% (One can construct a Bernoulli shift with any other probability space in
%place of~$\COIN$.)
 It is well known that almost every coin-tossing history
consists (in the limit) of half heads and half tails. More generally, it
can be shown that almost every orbit in this dynamical system is
\term{uniformly distributed}, so it is
``\defit[ergodic!dynamical system]{ergodic}." 
 \end{eg}

It is also helpful to see the \term{Bernoulli shift} from a more concrete
point of view:

\begin{eg}[(\term{Baker's Transformation})]
 A baker starts with a lump of \term{dough} in the shape of the unit square
$[0,1]^2$. She cuts the dough in half along the line $x = 1/2$
\seefig{BakerLoafFig}, and places the right half $[1/2,1] \times [0,1]$
above the left half (to become $[0,1/2] \times [1,2]$). The dough is now
$2$~units tall (and $1/2$~unit wide). She then pushes down on the dough,
reducing its height, so it widens to retain the same area, until the
dough is again square. (The pushing applies the linear map $(x,y)
\mapsto (2x,y/2)$ to the dough.) More formally, 
 \nindex{$\Baker$ = Baker's Transformation}
 \begin{equation} \label{BakerFormula}
 \Baker(x,y) = 
 \begin{cases}
 (2x,y/2) & \mbox{if $x \le 1/2$} \\
 \bigl( 2x-1,(y+1)/2 \bigr) & \mbox{if $x \ge 1/2$} 
 .
 \end{cases} 
 \end{equation}
 (This is not well defined on the set $\{ x = 1/2\}$ of measure zero.)

\begin{figure}
 \includegraphics[scale=0.41667]{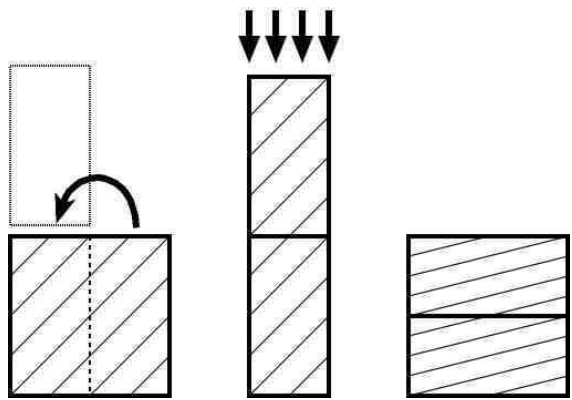}
 \caption{Baker's Transformation: the right half of the loaf is placed on
top of the left, and then the pile is pressed down to its original
height.}
 \label{BakerLoafFig}
 \end{figure}

 Any point $(x,y)$ in $[0,1]^2$
can be represented, in binary, by two strings of $0$'s and~$1$'s:
 $$ (x,y) = (0.x_0 x_1 x_2 \ldots, 0.y_1 y_2 y_3 \ldots ) ,$$
 and we have
 $$ \Baker (0.x_0 x_1 x_2 \ldots, 0.y_1 y_2 y_3 \ldots )
 = (0.x_1 x_2 \ldots, 0.x_0 y_1 y_2 y_3 \ldots ) ,$$
 so we see that
 \begin{equation} \label{Bern=Baker}
 \begin{matrix}
 \mbox{$\Bern$ is \term[isomorphic, measurably]{isomorphic}
to $\Baker$} \hfill \\
 \mbox{(modulo a set of measure zero),}
 \end{matrix}
 \end{equation}
 by identifying $f \in \COIN^\infty$ with
 $$ \bigl( 0 . \widehat{f(0)} \, \widehat{f(1)} \, \widehat{f(2)} \ldots,
\, 0 . \widehat{f(-1)} \, \widehat{f(-2)} \, \widehat{f(-3)} \ldots \bigr)
,$$
 where $\widehat{\head} = 0$ and $\widehat{\tail} = 1$ (or vice-versa).
 \end{eg}

\begin{exercises}

\item Show $\beta$ is rational if and only if there is some positive
integer~$k$, such that $(T_\beta)^k(x) = x$ for all $x \in \torus$.

\item \label{UnifDistRot}
 Show that if $\beta$ is \term[irrational!number]{irrational}, then every
orbit of~$T_\beta$ is \term{uniformly distributed} on the circle; that is,
if $I$ is any arc of the circle, and $x$~is any point in~$\torus$, show
that
 $$ \lim_{N \to \infty} \frac{
 \# \bigl\{\, k \in \{1,2,\ldots,N\} \mid
T_\beta^k(x) \in I \, \bigr\}
 }{N}
  = \operatorname{length}(I) .$$
 \hint{Exer.~\ref{T2UnifDistOpen}.}

 \end{exercises}

\section{Unpredictability}
 \index{unpredictable}
 \label{UnpreditableSect}

There is a fundamental difference between our two examples:  the
\term{Bernoulli shift} (or \term{Baker's Transformation}) is much more
``random" or ``\term{unpredictable}" than an
\term[irrational!rotation]{irrational rotation}. (Formally, this will be
reflected in the fact that the \term[entropy!of a dynamical
system|(]{entropy} of a \term{Bernoulli shift} is nonzero, whereas that of
an \term[irrational!rotation]{irrational rotation} is zero.) Both of these
dynamical systems are deterministic, so, from a certain point of view,
neither is unpredictable.  But the issue here is to predict behavior of
the dynamical system from imperfect knowledge, and these two examples look
fundamentally different from this point of view.

\begin{eg} \index{irrational!rotation}
 Suppose we have know the \term{location} of some point~$x$
of~$\torus$ to within a distance of less than 0.1; that is, we have a
point $x_0 \in \torus$, and we know that $d(x,  x_0) < 0.1$.  Then, for
every~$n$, we also know the location of $T_\beta^n(x)$ to within a
distance of 0.1. Namely, $d \bigl( T_\beta^n(x), T_\beta^n(x_0) \bigr) <
0.1$, because $T_\beta$ is an \term{isometry} of the circle. Thus, we can
predict the location of $T_\beta^n(x)$ fairly
accurately.
 \end{eg}

The \term{Baker's Transformation} is \emph{not}
\term[unpredictable]{predictable} in this sense:

\begin{eg} \label{Unpredict(Bake)Eg} \index{Baker's Transformation}
 Suppose there is an \term{impurity} in the baker's bread \term{dough}, and
we know its \term{location} to within a distance of less than 0.1. After
the dough has been kneaded once, our uncertainty in the horizontal
direction has doubled \seefig{baker-oval}. Kneading again doubles the
horizontal uncertainty, and perhaps adds a second possible vertical
location (if the cut $\{x = 1/2\}$ goes through the region where the
impurity might be). As the baker continues to knead, our hold on the
impurity continues to deteriorate (very quickly --- at an exponential
rate!).  After, say, 20~kneadings, the impurity could be almost anywhere
in the loaf --- every point in the dough is at a distance less than 0.001
from a point where the impurity could be \seeex{BreadAnywhere}. In
particular, we now have no idea whether the impurity is in the left half
of the loaf or the right half.

\begin{figure}
 \includegraphics[scale=0.44035]{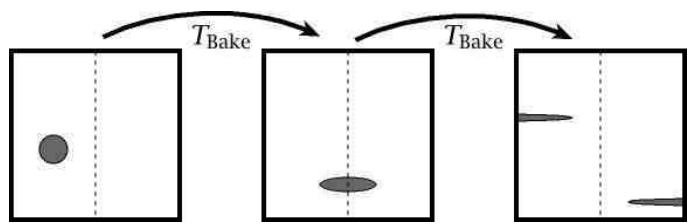}
 \caption{Kneading the dough stretches a circular disk to a
narrow, horizontal ellipse.}
 \label{baker-oval}
 \end{figure}

The upshot is that a small change in an initial position can quickly lead
to a large change in the subsequent location of a particle. Thus, errors
in measurement make it impossible to predict the future course of the
system with any reasonable precision. (Many scientists believe that
weather forecasting suffers from this difficulty --- it is said that a
\term{butterfly flapping its wings} in Africa may affect the next month's
rainfall in Chicago.)
 \end{eg}

To understand \term[entropy!of a dynamical
system]{entropy}, it is important to look at
\term[unpredictable]{unpredictability} from a different point of view,
which can easily be illustrated by the \term{Bernoulli shift}.

\begin{eg} \label{Unpredict(Bern)Eg} \index{Bernoulli shift}
 Suppose we have tossed a fair coin 1000 times. Knowing the
results of these tosses does not help us predict the next toss: there
is still a 50\% chance of heads and a 50\% chance of tails. 

More formally, define a function $\chi \colon \COIN^\infty \to
\{\head,\tail\}$ by $\chi(f) = f(0)$. Then the values 
 $$ \chi \bigl( T^{-1000}(f) \bigr), \chi \bigl( T^{-999}(f) \bigr),
\chi \bigl( T^{-998}(f) \bigr), \ldots, \chi \bigl( T^{-1}(f) \bigr) $$
 give no \term{information} about the value of~$\chi(f)$. Thus, $\Bern$ is
quite unpredictable, even if we have a lot of past history to go on. As we
will see, this means that the \term[entropy!of a dynamical
system]{entropy} of~$\Bern$ is not zero.
 \end{eg}

\begin{eg} \label{Predict(rot)Eg}
 In contrast, consider an \term[irrational!rotation]{irrational
rotation}~$T_\beta$. For concreteness, let us take $\beta = \sqrt{3}/100 =
.01732\ldots$, and, for convenience, let us identify $\torus$ with the
half-open interval $[0,1)$ (in the natural way). Let $\chi \colon [0,1)
\to \{0,1\}$ be the characteristic function of $[1/2,1)$. Then, because
$\beta \approx 1/60$, the sequence
 $$ \ldots, \chi \bigl( T_\beta^{-2}(x) \bigr), \chi \bigl(
T_\beta^{-1}(x) \bigr), \chi \bigl( T_\beta^{0}(x) \bigr), \chi \bigl(
T_\beta^{1}(x) \bigr), \chi \bigl( T_\beta^{2}(x) \bigr), \ldots $$ 
 consists of alternating strings of about thirty~$0$'s and about
thirty~$1$'s. Thus,
 $$ \begin{matrix}
 \mbox{if $\chi \bigl( T_\beta^{-1}(x) \bigr) = 0$ and $\chi \bigl(
T_\beta^{0}(x) \bigr) = 1$, then}
 \hfill\vphantom{\Bigl(}\\
 \mbox{we know that $\chi \bigl( T_\beta^{k}(x) \bigr) = 1$ for $k =
1,\ldots, 25$.}
 \end{matrix} $$
 So the results are somewhat \term{predictable}. (In contrast, consider the
fortune that could be won by predicting 25~\term[coin tossing]{coin tosses}
on a fairly regular basis!)

But that is only the beginning. The values of
  $$ \chi \bigl( T_\beta^{-1000}(x) \bigr), \chi \bigl( T_\beta^{-999}(x)
\bigr), \chi \bigl( T_\beta^{-998}(x) \bigr), \ldots, \chi \bigl(
T_\beta^{-1}(x) \bigr) $$
 can be used to determine the position of~$x$ fairly accurately.
Using this more subtle \term{information}, one can predict far more than
just 25 values of~$\chi$ --- it is only when $T_\beta^n(x)$ happens to
land very close to $0$ or~$1/2$ that the value of $\chi \bigl(
T_\beta^{n}(x) \bigr)$ provides any new \term{information}. Indeed,
knowing more and more values of $\chi \bigl( T_\beta^k(x) \bigr)$ allows
us to make longer and longer strings of predictions. In the limit, the
amount of \term[unpredictable]{unpredictability} goes to~$0$, so it turns
out that the \term[entropy!of a dynamical
system]{entropy} of~$T_\beta$ is~$0$. 
 \end{eg}

We remark that the relationship between \term[entropy!of a dynamical
system]{entropy} and
\term[unpredictable]{unpredictability} can be formalized as follows
\seeex{h-zero}.

\begin{prop} \label{past=future->h(T)=0}
 The \term[entropy!of a dynamical
system]{entropy} of a transformation is~$0$ if and only if the
\term{past determines the future} {\upshape(}almost surely{\upshape)}.

More precisely, the \term[entropy!of a dynamical
system]{entropy} of~$T$ is~$0$ if and only if, for each
partition~$\partA$ of~$\Omega$ into finitely many measurable sets, there
is a conull subset~$\Omega'$, such that, for $x,y \in \Omega'$,
 $$ \mbox{if $T^k(x) \sim T^k(y)$, for all $k < 0$, then 
 $T^k(x) \sim T^k(y)$, for all $k$} ,$$
 where $\sim$ is the \term{equivalence relation} corresponding to the
partition~$\partA$: namely, $x \sim y$ if there exists $A \in
\partA$ with $\{x,y\} \subset A$.
 \end{prop}

\begin{figure}
 \includegraphics[scale=0.44035]{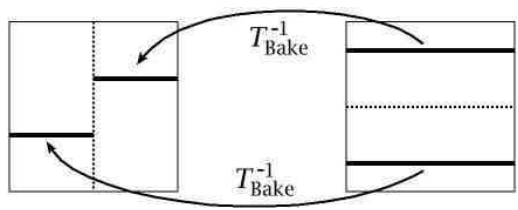}
 \caption{The inverse image of horizontal line segments.}
 \label{TBakeInv}
 \end{figure}

\begin{eg} \label{PastNotFutureEg} \ 
 \begin{enumerate}
 \item Knowing the entire past history of a fair coin does not tell us
what the next toss will be, so, in accord with
Eg.~\ref{Unpredict(Bern)Eg}, Prop.~\ref{past=future->h(T)=0} implies
the \term[entropy!of a dynamical
system]{entropy} of~\term[Bernoulli shift]{$\Bern$} is \emph{not}~$0$.

 \item \label{PastNotFutureEg-Baker}
 For the \term{Baker's Transformation}, let 
 $$ \partA = \bigl\{ [0,1/2) \times [0,1], [1/2,1] \times [0,1] \bigr\} $$
 be the partition of the unit square into a left half and a right half.
The inverse image of any horizontal line segment lies entirely in one of
these halves of the square (and is horizontal) \seefig{TBakeInv}, so (by
induction), if $x$ and~$y$ lie on the same horizontal line segment, then
$\Baker^k(x) \sim \Baker^k(y)$ for all $k < 0$. On the other hand, it is
(obviously) easy to find two points $x$ and~$y$ on a horizontal line
segment, such that $x$ and~$y$ are in opposite halves of the partition.
So the past does not determine the present (or the future). This is an
illustration of the fact that the \term[entropy!of a dynamical
system]{entropy}~$h(\Baker)$ of~$\Baker$ is
\emph{not}~$0$ \see{Unpredict(Bake)Eg}. 

\item Let $\partA = \{I, \torus \smallsetminus I\}$, for some (proper)
arc~$I$ of~$\torus$. If $\beta$ is \term[irrational!number]{irrational},
then, for any $x \in \torus$, the set $\{\, T_\beta^k (x) \mid k < 0 \,\}$
is \term[dense orbit]{dense} in~$\torus$. From this observation, it is not
difficult to see that if $T_\beta^k(x) \sim T_\beta^k(y)$ for all $k < 0$,
then $x = y$. Hence, for an \term[irrational!rotation]{irrational
rotation}, the past does determine the future. This is a manifestation of
the fact that the \term[entropy!of a dynamical
system]{entropy} $h(T_\beta)$ of~$T_\beta$ is~$0$
\see{Predict(rot)Eg}.

 \end{enumerate}
 \end{eg}

\onlyoneexercise
 \begin{exercises}

 \item \label{BreadAnywhere} \index{Baker's Transformation}
 Show, for any $x$ and~$y$ in the unit square, that there exists $x'$,
such that $d(x,x') < 0.1$ and $d \bigl( \Baker^{20}(x),y \bigr) < 0.01$.

\end{exercises}
 \endgroup

\section{Definition of entropy}
 \index{entropy!of a dynamical system}
 \label{DefnEntopySect}

The fundamental difference between the behavior of the above two examples
will be formalized by the notion of the \term[entropy!of a dynamical
system]{entropy} of a dynamical system, but, first, we define the
\term[entropy!of a partition]{entropy} of a partition.

\begin{rem}
 Let us give some motivation for the following definition. Suppose we are
interested in the location of some point~$\omega$ in some
probability space~$\Omega$.
 \begin{itemize}
 \item If $\Omega$ has been divided into $2$~sets of equal measure,
and we know which of these sets the point~$\omega$ belongs to, then we
have $1$~\term[bit (of information)]{bit of information}\index{information}
about the location of~$\omega$. (The two halves of~$\Omega$ can be
labelled `0' and~`1'.)
 \item If $\Omega$ has been divided into $4$~sets of equal measure,
and we know which of these sets the point~$\omega$ belongs to, then we
have $2$~bits of information about the location of~$\omega$. (The four
quarters of~$\Omega$ can be labelled `00', `01', `10', and~`11'.)
 \item More generally, if $\Omega$ has been divided into $2^k$~sets of
equal measure, and we know which of these sets the point~$\omega$ belongs
to, then we have exactly~$k$ bits of information about the location
of~$\omega$.
 \item The preceding observation suggests that if $\Omega$ has been
divided into $n$~sets of equal measure, and we know which of these sets
the point~$\omega$ belongs to, then we have exactly $\log_2 n$ bits of
information about the location of~$\omega$.
 \item But there is no need to actually divide $\Omega$ into pieces: if
we have a certain subset~$A$ of~$\Omega$, with $\mu(A) = 1/n$, and we
know that $\omega$ belongs to~$A$, then we can say that we have exactly
$\log_2 n$ bits of information about the location of~$\omega$.
 \item More generally, if we know that $\omega$ belongs to a certain
subset~$A$ of~$\Omega$, and $\mu(A) = p$, then we should say that we have
exactly $\log_2 (1/p)$ bits of information about the location of~$\omega$.
 \item Now, suppose $\Omega$ has been partitioned into finitely many
subsets $A_1,\ldots,A_m$, of measure $p_1,\ldots,p_m$, and that we will
be told which of these sets a random point~$\omega$ belongs to. Then the
right-hand side of \eqref{H(partition)eqn} is the amount of information
that we can expect (in the sense of probability theory) to receive about
the location of~$\omega$.
 \end{itemize}
 \end{rem}

\begin{defn} \label{H(partition)}
 Suppose 
 \nindex{$\partA$, $\partB$ = partition}
 $\partA = \{A_1,\ldots,A_m\}$ is a \term{partition} of a 
probability space $(\Omega, \mu)$ into finitely many measurable sets  of
measure  $p_1,\ldots,p_m$ respectively. 
 (Each set~$A_i$ is called an \defit[atom (of a partition)]{atom}
of~$\partA$.) The \defit[entropy!of a partition]{entropy} of this
partition is
 \nindex{$H(\partA)$ = entropy of the partition~$\partA$}
 \begin{equation} \label{H(partition)eqn}
 H(\partA)
 = H(p_1,\ldots,p_m)
 = \sum_{i=1}^m p_i \, \log \frac{1}{p_i} .
 \end{equation}
 If $p_i = 0$, then $p_i \, \log(1/p_i) = 0$, by convention.
 Different authors may use different bases for the logarithm, usually
either~$e$ or~$2$, so this definition can be varied by a scalar multiple.
Note that \term[entropy!of a partition]{entropy} is never negative.
 \end{defn}

\begin{rem} \label{ErgPartMotiv}
 Let us motivate the definition in another way. Think of the
partition~$\partA$ as representing an \term{experiment} with  $m$~(mutually
exclusive) possible outcomes. (The probability of the
$i^{\text{th}}$~outcome is~$p_i$.) We wish $H(\partA)$ to represent the
amount of \term{information} one can expect to gain by performing the
experiment. Alternatively, it can be thought of as the amount of
uncertainty regarding the outcome of the experiment.  For example, 
$H(\{\Omega\}) = 0$, because we gain no new \term{information} by
performing an experiment whose outcome is known in advance.  

Let us list some properties of~$H$ that one would expect, if it is to fit
our intuitive understanding of it as the \term{information} gained from an
\term{experiment}. 
 \begin{enumerate}

 \item \label{ErgPartMotiv-prob}
  The \term[entropy!of a partition]{entropy} does not depend on the
particular subsets chosen for the partition, but only on their
probabilities.  Thus, for $p_1,\ldots,p_n \ge 0$  with $\sum_i p_i = 1$,
we have a real number $H(p_1,\ldots,p_n) \ge 0$. Furthermore, permuting
the probabilities $p_1,\ldots,p_n$ does not change the value of the entropy
$H(p_1,\ldots,p_n)$.

\item \label{ErgPartMotiv-certain}
 An experiment yields no \term{information} if and only if we can predict
its outcome with certainty.  Therefore, we have $H(p_1,\ldots,p_n) = 0$ if
and only if $p_i = 1$ for some~$i$.

\item \label{ErgPartMotiv-impossible}
 If a certain outcome of an experiment is impossible, then there
is no harm in eliminating it from the description of the experiment, so
$H(p_1,\ldots,p_n,0) = H(p_1,\ldots,p_n)$.

 \item  \label{ErgPartMotiv-max}
 The least predictable experiment is one in which all outcomes are
equally likely, so, for given~$n$, the function  $H(p_1, \ldots, p_n)$  is
maximized at  $(1/n,\ldots,1/n)$.

\item \label{ErgPartMotiv-continuous}
 $H(p_1, \ldots, p_n)$ is a continuous function of its arguments.

\item \label{ErgPartMotiv-conditional}
 Our final property is somewhat more
sophisticated, but still intuitive.  Suppose we have two finite partitions
$\partA$  and~$\partB$ (not necessarily independent), and
let 
 \nindex{$\partA \vee \partB$ = join of the partitions}
 $$ \partA \vee \partB
 = \{\, A \cap B \mid A \in \partA, \ B \in \partB \,\} $$
  be their \defit[join (of partitions)]{join}. The join corresponds to
performing \emph{both} experiments.  We would expect to get the same amount
of \term{information} by performing the two experiments in succession, so
 $$ H(\partA \vee \partB) = H(\partA ) +
H(\partB \mid \partA) ,$$
 where 
 $H(\partB \mid \partA)$, the \defit[entropy!expected]{expected}
(or \defit[entropy!conditional]{conditional}) \emph{entropy} of~$\partB$,
given~$\partA$, is the amount of \term{information} expected from
performing experiment~$\partB$, given that experiment~$\partA$  has
already been performed. 

More precisely, if experiment~$\partA$ has been performed, then some 
event~$A$  has been observed. The amount of
\term{information}~$H_A(\partB)$ now expected from
experiment~$\partB$ is given by the \term[entropy!of a partition]{entropy}
of the partition $\{B_j \cap A \}_{j=1}^n$  of~$A$, so
  $$ H_A(\partB) = H \left( \frac{\mu(B_1 \cap A)}{\mu(A)}, 
 \frac{\mu(B_2 \cap A)}{\mu(A)},
 \ldots,
 \frac{\mu(B_n \cap A)}{\mu(A)} \right) .$$
 This should be weighted by the probability
of~$A$, so
 \nindex{$H(\partB \mid \partA)$ = conditional entropy}
 $$ H(\partB \mid \partA) = \sum_{A \in \partA} 
H_A(\partB) \, \mu(A) .$$
 \end{enumerate}
 An elementary (but certainly nontrivial) argument shows that any
function~$H$ satisfying all of these conditions must be as described in
Defn.~\ref{H(partition)} (for some choice of the base of the logarithm).
 \end{rem}

The \term[entropy!of a partition]{entropy} of a partition leads directly
to the definition of the \term[entropy!of a dynamical
system]{entropy} of a dynamical system. To motivate this, think about
repeating the same \term{experiment} every day.  The first day we
presumably obtain some \term{information} from the outcome of the
experiment.  The second day may yield some additional information (the
result of an experiment --- such as recording the time of sunrise --- may
change from day to day).  And so on. If the dynamical system is
``predictable," then later experiments do not yield much new information. 
On the other hand, in a truly \term{unpredictable} system, such as a
\term[coin tossing]{coin toss}, we learn something new every day --- the
expected total amount of information is directly proportional to the
number of times the experiment has been repeated.  The total amount of
\term{information} expected to be obtained after $k$~daily experiments
(starting today) is
 \nindex{$E^k(T, \partA)$ = information expected from $k$ experiments}
 $$ E^k(T, \partA) = H \bigl( \partA \vee T^{-1}(\partA) \vee
T^{-2}(\partA) \vee \cdots \vee T^{-(k-1)}(\partA) \bigr) , $$
 where 
 \nindex{$T^\ell(\partA)$ = transform of~$\partA$ by~$T^\ell$}
 $$ T^\ell(\partA) = \{\, T^\ell(A) \mid A \in \partA \,\} $$
 \seeex{SeveralDays}.
 It is not difficult to see that the limit
 \nindex{$h(T,\partA)$ = entropy of~$T$ w.r.t.\ $\partA$}
 $$ h(T,\partA) = \lim_{k \to \infty} \frac{ E^k(T, \partA)
}{k} $$
 exists \seeex{lim(Ek/k)exists}. The \term[entropy!of a dynamical
system]{entropy}
of~$T$ is the value of this limit for the most effective experiment:

\begin{defn} \label{h(T)defn}
 The \defit{entropy}
 \nindex{$h(T) = h_\mu(T)$ = entropy of~$T$}
 $h(T)$ is the supremum of $h(T,\partA)$ over all partitions~$\partA$
of~$\Omega$ into finitely many measurable sets.
 \end{defn}

\begin{notation}
 The \term[entropy!of a dynamical
system]{entropy} of~$T$ may depend on the choice of the
\term[measure!invariant]{invariant measure}~$\mu$, so, to avoid confusion,
it may sometimes be denoted $h_\mu(T)$.
 \end{notation}

\begin{rem}
 The \defit[entropy!of a flow]{entropy} of a flow is defined to be the
entropy of its time-one map; that is, $h \bigl( \{\varphi_t\} \bigr) = h(
\varphi_1)$.
 \end{rem}

\begin{rem} \label{TopEntropy}
 Although we make no use of it in these lectures, we mention that there is
also a notion of \defit[entropy!topological]{entropy} that is purely
topological. Note that $E^k(T,\partA)$ is large if the partition
 $$ \partA \vee T^{-1}(\partA) \vee T^{-2}(\partA)
\vee \cdots  \vee T^{-(k-1)}(\partA)$$
 consists of very many small sets. In pure topology, without a measure,
one cannot say whether or not the sets in a collection are ``small," but
one can say whether or not there are very many of them, and that is the
basis of the definition. However, the topological definition uses open
covers of the space, instead of measurable partitions of the space.

Specifically, suppose $T$ is a \index{homeomorphism}homeomorphism of a
compact metric space~$X$.
 \begin{enumerate}
 \item For each open cover~$\partA$ of~$X$, let $N(\partA)$ be the minimal
cardinality of a subcover.
 \item If $\partA$ and~$\partB$ are open covers, let
 $$ \partA \vee \partB = \{\, A \cap B \mid A \in \partA, \ B \in
\partB \,\} .$$
 \item Define %the \defit[entropy!topological]{topological entropy}
 \nindex{$h_{\text{top}}(T)$ = topological entropy of~$T$}
 $$ h_{\text{top}}(T) = \sup_{\partA} \lim_{k \to \infty} 
 \frac{ \log N \bigl( \partA \vee T^{-1}(\partA) %\vee T^{-2}(\partA)
\vee \cdots  \vee T^{-(k-1)}(\partA) \bigr)}
 {k} .$$
 \end{enumerate}
 It can be shown, for every $T$-\term[measure!invariant]{invariant
probability measure}~$\mu$ on~$X$, that $h_\mu(T) \le h_{\text{top}}(T)$.
 \end{rem}

\begin{exercises}
 \index{entropy!of a partition}

\item Show the function~$H$ defined in Eq.~\pref{H(partition)eqn}
satisfies the formulas in:
 \begin{enumerate}
 \item Rem.~\fullref{ErgPartMotiv}{certain},
 \item Rem.~\fullref{ErgPartMotiv}{impossible},
 \item Rem.~\fullref{ErgPartMotiv}{max},
 \item Rem.~\fullref{ErgPartMotiv}{continuous},
 and
 \item Rem.~\fullref{ErgPartMotiv}{conditional}.
 \end{enumerate}

 \item \label{h-H(AinB)}
 Intuitively, it is clear that altering an \term{experiment} to produce
more refined outcomes will not reduce the amount of \term{information}
provided by the experiment. Formally, show that  if $\partA \subset
\partB$, then $H(\partA) \le H(\partB)$.
 (We write $\partA \subset \partB$ if each \term[atom (of a
partition)]{atom} of~$\partA$ is a union of atoms of~$\partB$ (up to
measure zero).)

 \item \label{h-wedge}
  It is easy to calculate that the \term[entropy!of a partition]{entropy}
of a combination of
\term[experiment!statistically independent]{experiments} is precisely the
sum of the entropies of the individual experiments. Intuitively, it is
reasonable to expect that independent experiments provide the most
\term{information} (because they have no redundancy). Formally, show 
 $$H(\partA_1 \vee \cdots \vee \partA_n) \le \sum_{i=1}^n H(\partA_i) .$$
 \hint{Assume $n = 2$ and use \term{Lagrange Multipliers}.}

\item \label{h-A|B}
 Show $H(\partA \mid \partB) \le H(\partA)$.
 \hint{Exer.~\ref{h-wedge}.}

 \item \label{h-=}
 Show $H(T^\ell \partA) = H(\partA)$, for all $\ell \in \integer$.
 \hint{$T^\ell$ is measure preserving.}

\item \label{SeveralDays}
 For $x,y \in \Omega$, show that $x$ and~$y$ are in the same \term[atom (of a
partition)]{atom} of 
 $\bigvee_{\ell=0}^{k-1}  T^{-\ell}(\partA)$
 if and only if, for all $\ell \in \{ 0,\ldots,k-1 \}$, the two
points $T^\ell(x)$ and $T^\ell(y)$ are in the same atom of~$\partA$.

\item \label{h(T)=h(Tinv)Ex}
 \index{entropy}
 Show $h(T) = h(T^{-1})$.
 \hint{Exer.~\ref{h-=} implies $E^k(T,\partA) = E^k(T^{-1},\partA)$.}

\item Show $h(T^\ell) = |\ell| \, h(T)$, for all $\ell \in \integer$.
 \hint{For $\ell > 0$, consider $E^k \bigl( T^\ell, \bigvee_{i=0}^{\ell-1}
T^{-i} \partA \bigr)$.}

\item \label{EntropyIso}
 Show that \term[entropy!of a dynamical system]{entropy} is an
\term[invariant!under isomorphism]{isomorphism invariant}. That is, if
 \begin{itemize}
 \item $\psi \colon (\Omega,\mu) \to (\Omega',\mu')$ is a
measure-preserving map,
 such that 
 \item $\psi \bigl( T(\omega) \bigr) = T' \bigl( \psi(\omega) \bigr)$ a.e.,
 \end{itemize}
 then $h_\mu(T) = h_{\mu'}(T')$.

 \item \label{h-T<A}
 Show $h(T,\partA) \le H(\partA)$.
 \hint{Exers.~\ref{h-=} and \ref{h-wedge}.}

 \item \label{h-h(AinB)}
 Show that if $\partA \subset \partB$, then $h(T,\partA) \le h(T,\partB)$.

 \item \label{h-dist(AB)}
 Show
 $|h(T,\partA) - h(T,\partB)| \le H(\partA \mid \partB) + H(\partB
\mid \partA)$.
 \hint{Reduce to the case where $\partA \subset
\partB$, by using the fact that $\partA \vee \partB$ contains
both~$\partA$ and~$\partB$.}

\item \label{Eksubadditive}
 Show that the sequence $E^k(T,\partA)$ is \defit[subadditive
sequence]{subadditive}; that is, $E^{k+\ell}(T,\partA) \le E^k(T,\partA) +
E^\ell(T,\partA)$.

\item \label{lim(Ek/k)exists}
 Show that $\lim_{k \to \infty} \frac{1}{k} E^k(T,\partA)$ exists.
 \hint{Exer.~\ref{Eksubadditive}.}

\item \label{htop(isom)=0}
 Show that if $T$ is an \term{isometry} of a compact metric space, then
\term[entropy!topological]{$h_{\text{top}}(T) = 0$}.
 \hint{If $\partB_\epsilon$ is the open cover by balls of
radius~$\epsilon$, then $T^\ell(\partB_\epsilon) = \partB_\epsilon$.
Choose $\epsilon$ to be a \defit[Lebesgue number (of an open
cover)]{Lebesgue number} of the open cover~$\partA$; that is, every ball of
radius~$\epsilon$ is contained in some element of~$\partA$.}

\end{exercises}

\section{How to calculate entropy}
 \index{entropy!of a dynamical system}
 \label{CalcEntropySect}

Definition~\ref{h(T)defn} is difficult to apply in practice, because of
the supremum over all possible finite partitions. The following theorem
eliminates this difficulty, by allowing us to consider only a single
partition. (See Exer.~\ref{h-gen} for the proof.)

\begin{thm} \label{genpart}
 If $\partA$ is any finite \defit[partition!generating]{generating
partition} for~$T$, that is, if
 $$ \bigcup_{k=-\infty}^\infty T^k (\partA) $$
 generates the $\sigma$-algebra of all measurable sets {\upshape(}up to
measure~$0${\upshape)}, then
 \index{entropy}
 $$ h(T) = h(T , \partA) .$$
 \end{thm}

\begin{cor} \label{h(Tbeta)}
 \index{irrational!rotation}
 For any $\beta \in \real$, $h(T_\beta) = 0$.
 \end{cor}

\begin{proof}
 Let us assume $\beta$ is \term[irrational!number]{irrational}. (The other
case is easy; see Exer.~\ref{h(ratrot)}.) Let 
 \begin{itemize}
 \item $I$ be any (nonempty) proper arc of~$\torus$,
 and
 \item $\partA = \{I, \torus \smallsetminus I\}$.
 \end{itemize}
 It is easy to see that if $\partB$ is any finite partition of~$\torus$
into connected sets, then $\#(\partA \vee \partB) \le 2 + \#\partB$.
Hence 
 \begin{equation} \label{RotNumSets}
 \# \bigl( \partA \vee T_\beta^{-1}(\partA) \vee
T_\beta^{-2}(\partA) \vee \cdots \vee T_\beta^{-(k-1)}(\partA) \bigr)
 \le 2k 
 \end{equation}
 \seeex{RotNumSetsEx}, so, using \fullref{ErgPartMotiv}{max}, we see that
 \begin{align*}
 E^k(T_\beta,\partA)
 &= H \bigl(  \partA \vee T_\beta^{-1}(\partA) \vee
T_\beta^{-2}(\partA) \vee \cdots \vee T_\beta^{-(k-1)}(\partA) \bigr) \\
 & \le H \left( \frac{1}{2k}, \frac{1}{2k}, \ldots, \frac{1}{2k}
\right) \\
 &= \sum_{i=1}^{2k}  \frac{1}{2k} \log(2k) \\
 &= \log(2k) .
 \end{align*}
 Therefore
 $$ h(T_\beta, \partA) = \lim_{k \to \infty} \frac{E^k(T_\beta,\partA)}{k}
 \le \lim_{k \to \infty} \frac{\log(2k)}{k}
 = 0 .$$
 One can show that $\partA$ is a \defit[partition!generating]{generating
partition} for~$T_\beta$
\seeex{RotArcGen}, so Thm.~\ref{genpart} implies that
$h(T_\beta) = 0$.
 \end{proof}

\begin{cor} \label{h(Bern)} \label{h(Baker)}
 \index{Bernoulli shift} \index{Baker's Transformation}
 $h(\Bern) = h(\Baker) = \log 2$.
 \end{cor}

\begin{proof}
 Because $\Bern$ and~$\Baker$ are \term[isomorphic, measurably]{isomorphic}
\see{Bern=Baker}, they have the same \term[entropy!of a dynamical
system]{entropy} \seeex{EntropyIso}. Thus, we need only calculate
$h(\Bern)$.

 Let
 $$ A = \{\, f \in \COIN^\infty \mid f(0) = 1 \,\}
 \mbox{ \qquad and \qquad }
 \partA = \{ A, \COIN^\infty \smallsetminus A \} .$$
 Then 
 $$ \partA \vee T^{-1}(\partA) \vee
T^{-2}(\partA) \vee \cdots \vee T^{-(k-1)}(\partA) $$
 consists of the $2^k$ sets of the form
 $$ \COIN_{\epsilon_0,\epsilon_1,\ldots,\epsilon_{k-1}}
 = \{\, f \in \COIN^\infty \mid \mbox{$f(\ell) = \epsilon_\ell$, for
$\ell = 0,1,2,\ldots,k-1$} \,\} ,$$
 each of which has measure $1/2^k$. Therefore,
 \begin{align*}
 h(\Bern, \partA) &= \lim_{k \to \infty}
\frac{E^k(\Bern,\partA)}{k}
 = \lim_{k \to \infty} \frac{ 2^k \cdot \left[ \frac{1}{2^k} \log 2^k
\right]}{k} \\
 &= \lim_{k \to \infty} \frac{ k \, \log 2}{k}
 = \log 2
 .
 \end{align*}
  One can show that $\partA$ is a \defit[partition!generating]{generating
partition} for~$\Bern$
\seeex{BernGen}, so Thm.~\ref{genpart} implies that $h(\Bern) = \log
2$.
 \end{proof}

\begin{rem} \label{CtblPartn}
 One need not restrict to finite partitions; $h(T)$ is the
supremum of $h(T,\partA)$, not only over all finite measurable
partitions, but over all \term[partition!countable]{countable
partitions}~$\partA$, such that $H(\partA) < \infty$. When considering
countable partitions (of finite \term[entropy!of a partition]{entropy}),
some sums have infinitely many terms, but, because the terms are positive,
they can be rearranged at will.  Thus, essentially the same proofs can be
applied.
 \end{rem}

We noted, in Prop.~\ref{past=future->h(T)=0}, that if $h(T) = 0$, then the
past determines the present (and the future). 
 \index{past determines the future}
 That is, if  we know the results of all past experiments, then we
can predict the result of today's experiment. Thus, $0$~is the amount of
\term{information} we can expect to get by performing today's
\term{experiment}. More generally, Thm.~\ref{h(T)FromPast} below shows that
$h(T)$ is always the amount of \term{information} we expect to obtain by
performing today's experiment.

\begin{eg}
 \index{Baker's Transformation}
 As in Eg.~\fullref{PastNotFutureEg}{Baker}, let $\partA$ be the partition
of the unit square into a left half and a right half. It is not difficult
to see that
 \begin{align*}
 \mbox{$x$ and~$y$ lie on the same horizontal line segment} \\
 \Leftrightarrow 
 \quad
 \mbox{$\Baker^k(x) \sim \Baker^k(y)$ for all $k < 0$} . 
 \end{align*}
 (We ignore points for which one of the coordinates is a dyadic rational
--- they are a set of measure zero.) Thus, the results of past
\term[experiment]{experiments} tell us which horizontal line segment
contains the point~$\omega$ (and provide no other \term{information}). The
partition~$\partA$ cuts this line segment precisely in half, so the two
possible results are equally likely in today's experiment; the expected
amount of \term{information} is $\log 2$, which, as we know, is the
\term[entropy!of a dynamical system]{entropy} of~$\Baker$
\see{h(Bern)}.
 \end{eg}

\begin{notation} \label{ConditionOnPast} \ 
 \begin{itemize}
 \item Let
 \nindex{$\past{\partA} = \bigvee_{\ell=1}^\infty T^{\ell}
\partA$}
 $\past{\partA} = \bigvee_{\ell=1}^\infty T^{\ell} \partA$. 
 Thus, $\past{\partA}$ is the partition that corresponds to knowing the
results of all past \term[experiment]{experiments} \seeex{PastPartition}. 
 \nindex{$H(\partB \mid \past{\partA})$ = conditional entropy}
 \item Let $H(\partB \mid \past{\partA})$ denote the
\term[entropy!conditional]{conditional entropy} of a partition~$\partB$
with respect to~$\past{\partA}$. More precisely:
 \begin{itemize}
 \item the
measure~$\mu$ has a \defit[measure!conditional]{conditional
measure}~$\mu_A$ on each \term[atom (of a
partition)]{atom}~$A$ of~$\past{\partA}$;
 \item the partition~$\partB$ induces a partition $\partB_A$ of each
atom~$A$ of~$\past{\partA}$;
 \item we have the \term[entropy!of a partition]{entropy} $H(\partB_A)$
(with respect to the probability measure $\mu_A$);
 and
 \item $H(\partB \mid \past{\partA})$ is the integral of $H(\partB_A)$ over
all of~$\Omega$.
 \end{itemize} 
 \end{itemize}
 Thus, $H(\partB \mid \past{\partA})$ represents the amount of
\term{information} we expect to obtain by performing experiment~$\partB$,
given that we know all previous results of experiment~$\partA$.
 \end{notation}

See Exer.~\ref{h-past} for the proof of the following theorem.

\begin{thm} \label{h(T)FromPast}
 If $\partA$ is any finite \defit[partition!generating]{generating
partition} for~$T$, then
 $h(T) = H(\partA \mid \past{\partA})$.
 \end{thm}

Because $T^{-1}\past{\partA} = \partA \vee \past{\partA}$, the following
corollary is immediate.

\begin{cor} \label{h(T)=H(TA|A)}
 If $\partA$ is any finite \defit[partition!generating]{generating
partition} for~$T$, then
 $h(T) = H(T^{-1} \past{\partA} \mid \past{\partA})$.
 \end{cor}

\begin{exercises}
 \index{entropy}

\item \label{h(ratrot)}
 \index{irrational!rotation}
 Show, directly from the definition of \term[entropy!of a dynamical
system]{entropy}, that if $\beta$ is
rational, then $h(T_\beta) = 0$.

\item \label{RotNumComps}
 Show that if
 \begin{itemize}
 \item $\partA = \{\torus \smallsetminus I, I\}$, where $I$~is a proper
arc of~$\torus$,
 and
 \item $\partB$ is a finite partition of~$\torus$ into connected sets,
 \end{itemize}
 then \index{component!connected}
 $$ \sum_{C \in \partA \vee \partB}
 \text{(\# components of~$C$)}
 \le 2 + \#\partB .$$

\item \label{RotNumSetsEx}
 Prove \eqref{RotNumSets}.

\item \label{RotArcGen}
 Show that if
 \begin{itemize}
 \item $\partA = \{\torus \smallsetminus I, I\}$, where $I$~is a
nonempty, proper arc of~$\torus$,
 and
 \item $\beta$ is \term[irrational!number]{irrational},
 \end{itemize}
 then $\partA$ is a \term[partition!generating]{generating partition}
for~$T_\beta$.
 \hint{If $n$ is large, then $\bigvee_{k = 0}^n T^k \partA$ is a
partition of~$\torus$ into small intervals. Thus, any open interval is a
countable union of sets in $\bigcup_{k = 0}^\infty T^k \partA$.}

\item \label{BernGen}
 Show that the partition~$\partA$ in the proof of Cor.~\ref{h(Tbeta)} is
a \term[partition!generating]{generating partition}
for~\term[Bernoulli shift]{$\Bern$}.

\item \label{OtherBern}
 The construction of a  \defit{Bernoulli shift} can be generalized, by
using any probability space in place of~$\COIN$. Show that if a Bernoulli
shift~$T$ is constructed from a measure space with probabilities
$\{p_1,\ldots,p_n\}$, then $h(T) = H(p_1,\ldots,p_n)$.

 \item \label{h-finite}
 Show that if $S$ is any finite, nonempty set of integers, then
 $h(T, \bigvee_{k \in S} T^k \partA) = h(T,\partA)$.

\item \label{GenPart->everyset}
 Suppose
 \begin{itemize}
 \item $\partA$ is a finite \defit[partition!generating]{generating
partition} for~$T$, 
 \item $\partB$ is a finite partition of~$\Omega$,
 and
 \item  $\epsilon > 0$.
 \end{itemize}
 Show that there is a finite set~$S$ of integers, such that $H(\partB
\mid \partA_S) < \epsilon$, where $\partA_S = \bigvee_{\ell \in S} T^\ell
\partA$.

 \item \label{h-gen}
 Show that if $\partA$ is a finite \defit[partition!generating]{generating
partition} for~$T$, then
$h(T) = h(T,\partA)$.
 \hint{Exer.~\ref{GenPart->everyset}.}

\item \label{PastPartition}
 For $x,y \in \Omega$, show that $x$ and~$y$ belong to the same \term[atom (of a
partition)]{atom}
of~$\past{\partA}$ if and only if $T^k(x)$ and $T^k(y)$ belong
to the same atom of~$\partA$, for every $k < 0$.

 \item \label{h-past}
 Show:
 \begin{enumerate}
 \item $E^k(T,\partA) = H(\partA) + 
 \sum_{\ell=1}^{k-1} H \bigl( \partA \mid \bigvee_{j=1}^{\ell} T^{j}
\partA \bigr)$.
 \item $H \bigl( \partA \mid \bigvee_{\ell=1}^{k} T^{\ell} \partA \bigr)$
is a decreasing sequence.
 \item $h(T,\partA) = H \bigl( \partA \mathrel{\big|} \past{\partA}
\bigr)$.
 \item $h(T,\partA) = H(T^{-1}\past{\partA}
\mid \past{\partA})$.
 \end{enumerate}
 \hint{You may assume, without proof, that 
 $$ \textstyle H \bigl( \partA \mathrel{\big|}
\bigvee_{\ell=1}^\infty T^{\ell} \partA \bigr) = \lim_{k\to\infty} H
\bigl( \partA \mid \bigvee_{\ell=1}^{k} T^{\ell}\partA \bigr) ,$$
 if the limit exists.}

 \item \label{h-zero}
 Show $h(T,\partA) = 0$ if and only if $\partA \subset \past{\partA}$
{\upshape(}up to measure zero{\upshape)}.
 \hint{Exer.~\ref{h-past}.}

 \end{exercises}

\section{Stretching and the entropy of a translation}
 \index{entropy!of a dynamical system}
 \label{h(g)Sect}

\begin{rem} \label{RemStretch} \ 
 \begin{enumerate}
 \item \label{RemStretch-isom}
 Note that \term[irrational!rotation]{$T_\beta$} is an \term{isometry}
of~$\torus$. Hence, Cor.~\ref{h(Tbeta)} is a particular case of the
general fact that if $T$ is an \term{isometry} (and $\Omega$ is a compact
metric space), then $h(T) = 0$ \seeexs{h(Isom)=0}{htop(isom)=0}.
 \item \label{RemStretch-Baker}
 Note that \term[Baker's Transformation]{$\Baker$} is \emph{far} from being
an \term{isometry} of the unit square. Indeed:
 \begin{itemize}
 \item the unit square can be foliated into horizontal line segments, and
$\Baker$ stretches (local) distances on the leaves of
this foliation by a factor of~$2$ \cf{BakerFormula}, so horizontal
distances grow exponentially fast (by a factor of~$2^n$) under iterates
of~$\Baker$;
 and
  \item distances in the complementary (vertical) direction are
contracted  exponentially fast (by a factor of~$1/2^n$).
 \end{itemize}
 It is not a coincidence that $\log 2$, the logarithm of the
\term{stretching} factor, is the \term[entropy!of a dynamical
system]{entropy} of $\Baker$. 
 \end{enumerate}
 \end{rem}

The following theorem states a precise relationship between stretching and
\term[entropy!of a dynamical system]{entropy}. (It can be stated in more
general versions that apply to non-smooth maps, such as~$\Baker$.) Roughly
speaking, \term[entropy!of a dynamical system]{entropy} is calculated by
adding contributions from all of the independent directions that are
stretched at exponential rates (and ignoring directions that are
contracted).

\begin{thm} \label{h(T)=stretch}
 Suppose
 \begin{itemize}
 \item $\Omega = M$ is a smooth, compact \index{manifold}manifold,
 \item $\vol$ is a volume form on~$M$,
 \item $T$ is a volume-preserving \index{diffeomorphism}diffeomorphism,
 \item $\tau_1,\ldots,\tau_k$ are {\upshape(}positive\/{\upshape)}
real numbers,
 and
 \item the tangent bundle $\T M$ is a direct sum of
$T$-\term[invariant!subbundle]{invariant subbundles} $\E_1,\ldots,\E_k$,
such that $\| dT (\xi) \| = \tau_i \| \xi \|, $ for each tangent
vector~$\xi \in \E_i$,
 \end{itemize}
 then 
 $$ h_{\vol}(T) = \sum_{\tau_i > 1} (\dim \E_i) \, \log \tau_i .$$
 \end{thm}

\begin{eg} \label{h(SmoothIsom)}
 If $T$ is an \term{isometry} of~$M$, let $\tau_1 = 1$ and $\E_1 = \T M$.
Then the theorem asserts that $h(T) = 0$. This establishes
Rem.~\fullref{RemStretch}{isom} in the special case where $\Omega$ is a
smooth \index{manifold}manifold and $T$~is a
\index{diffeomorphism}diffeomorphism.
 \end{eg}

\begin{eg}
 For \term[Baker's Transformation]{$\Baker$}, let 
 $\tau_1 = 2$, $\tau_2 = 1/2$, 
 $$ \mbox{$\E_1 =
\{\mbox{horizontal vectors}\}$ and $\E_2 = \{\mbox{vertical vectors}\}$.}
$$
 Then, if we ignore technical problems arising from the
nondifferentiability of~$\Baker$ and the boundary of the unit square, the
theorem confirms our calculation that \term[entropy]{$h(\Baker) = \log 2$}
\see{h(Baker)}.
 \end{eg}

For the special case where $T$ is the translation by an element of~$G$,
Thm.~\ref{h(T)=stretch} can be rephrased in the following form.

\begin{notation}
 Suppose $g$ is an element of~$G$
 \begin{itemize}
 \item Let
 \nindex{$G_+$ = horospherical subgroup}
 $G_+ = \{\, u \in G \mid \lim_{k \to -\infty} g^{-k} u g^k = e
\,\}$. (Note that $k$ tends to \emph{negative} infinity.) Then $G_+$ is a
closed, \index{subgroup!connected|(}connected subgroup of~$G$
\seeex{HoroSubgrpLie}. 
 \item Let
 $$J(g,G_+) = \left| \det \bigl( (\Ad g)|_{\Lie G_+} \bigr) \right|$$
 be the \index{Jacobian}Jacobian of~$g$ acting by conjugation on~$G_+$.
 \end{itemize}
 \end{notation}

\begin{rem} \label{horospherical}
 $G_+$ is called the (expanding)
\defit[subgroup!horospherical]{horospherical subgroup} corresponding
to~$g$. (Although this is not reflected in the notation, one should keep
in mind that the subgroup~$G_+$ depends on the choice of~$g$.) Conjugation
by~$g$ \emphit{expands} the elements of~$G_+$ because, by definition,
conjugation by~$g^{-1}$ contracts them. 

 There is also a \defit[subgroup!horospherical]{contracting} horospherical
subgroup~$G_-$, consisting of the elements that are contracted by~$g$. It
is defined by
 $$ G_- = \{\, u \in G \mid \lim_{k \to \infty} g^{-k} u g^k = e \,\} ,$$
 the only difference being that the limit is now taken as $k$ tends to
\emph{positive} infinity. Thus, $G_-$ is the expanding horospherical
subgroup corresponding to~$g^{-1}$.
 \end{rem}

\begin{cor} \label{h(G/Gamma)}
 Let
 \begin{itemize}
 \item $G$ be a \index{subgroup!connected}connected \term{Lie group},
 \item $\Gamma$ be a \term{lattice} in~$G$,
 \item $\vol$ be a $G$-\term[measure!invariant|(]{invariant volume} on\/
$\Gamma \backslash G$,
 and
 \item $g \in G$.
 \end{itemize}
 Then
 $$ h_{\vol}(g) = \log J(g,G_+) ,$$
 \index{Jacobian}
 where, abusing notation, we write $h_{\vol}(g)$ for the \term[entropy!of
a dynamical system]{entropy} of the translation $T_g \colon \Gamma
\backslash G \to \Gamma \backslash G$, defined by $T_g(x) = x g$.
 \end{cor}

\begin{cor} \label{h(unip)=0}
 \index{unipotent!flow}
 If $u \in G$ is \term[unipotent!element]{unipotent}, then
$h_{\vol}(u) = 0$.
 \end{cor}

\begin{cor} \label{h(diagSL2)}
 \index{geodesic!flow}
 If 
 \begin{itemize}
 \item $G = \SL(2,\real)$\index{SL(2,R)*$\SL(2,\real)$},
 and
 \item $a^s = \begin{bmatrix} e^s & 0 \\ 0 & e^{-s} \end{bmatrix}$,
 \end{itemize}
 then $h_{\vol}(a^s) = 2|s|$.
 \end{cor}

Corollary~\ref{h(G/Gamma)} calculates the \term[entropy!of a dynamical
system]{entropy} of~$g$ with respect to the natural volume form on $\Gamma
\backslash G$. The following generalization provides an estimate (not
always an exact calculation) for other \term[measure!invariant]{invariant
measures}. If one accepts that \term[entropy!of a dynamical
system]{entropy} is determined by the amount of \term{stretching}, in the
spirit of Thm.~\ref{h(T)=stretch} and Cor.~\ref{h(G/Gamma)}, then the
first two parts of the following proposition are fairly obvious at an
intuitive level. Namely:\index{entropy!vs.\ stretching}
 \begin{enumerate}
 \item[\pref{EntropyLemma-Inv->}] The hypothesis of
\fullref{EntropyLemma}{Inv->} implies that stretching along any direction
in~$W$ will contribute to the calculation of $h_\mu(a)$. This yields only
an inequality, because there may be other directions, not along~$W$, that
also contribute to $h_\mu(a)$; that is, there may well be other directions
that are being stretched by~$a$ and belong to the \term[support!of a
measure]{support} of~$\mu$.
 \item[\pref{EntropyLemma-GminInW}] Roughly speaking, the hypothesis
of~\fullref{EntropyLemma}{GminInW} states that any direction stretched
by~$a$ and belonging to the \term[support!of a measure]{support} of~$\mu$
must lie in~$W$. Thus, only directions in~$W$ contribute to the
calculation of $h_\mu(a)$. This yields only an inequality, because some
directions in~$W$ may not belong to the \term[support!of a
measure]{support} of~$\mu$.
 \end{enumerate}

\begin{notation} \label{aWNotn}
 Suppose 
 \begin{itemize}
 \item $g$ is an element of~$G$, with corresponding
\term[subgroup!horospherical]{horospherical subgroup}~$G_+$,
 and
 \item $W$ is a \index{subgroup!connected}connected \term{Lie subgroup}
of~$G_+$ that is \term[normalizer]{normalized} by~$g$.
 \end{itemize}
 Let 
 \nindex{$J(g, W)$ = Jacobian of~$g$ on~$W$}
 \index{Jacobian}
 $$J(g, W) = \left| \det \bigl( (\Ad g)|_{\Lie W} \bigr) \right|$$
 be the \index{Jacobian}Jacobian of~$g$ on~$W$. Thus, 
 $$ \log J(g, W) = \sum_\lambda \log |\lambda| ,$$
 where the sum is over all \term[eigenvalue]{eigenvalues} of $(\Ad
g)|_{\Lie W}$, counted with multiplicity, and $\Lie W$~is the \term{Lie
algebra} of~$W$.
 \end{notation}

\begin{prop} \label{EntropyLemma}
 Suppose 
 \begin{itemize}
 \item $g$ is an element of~$G$,
 \item $\mu$ is an \term[ergodic!measure]{measure}
$g$-\term[measure!invariant]{invariant probability measure} on $\Gamma
\backslash G$,
 and
 \item $W$ is a \index{subgroup!connected}connected \term{Lie subgroup}
of~$G_+$ that is \term[normalizer]{normalized} by~$g$.
 \end{itemize}
 Then:
 \begin{enumerate}

 \item \label{EntropyLemma-Inv->}
 If $\mu$ is $W$-\term[measure!invariant]{invariant}, then $h_\mu(g) \ge
\log J(g, W)$\index{Jacobian}.

 \item \label{EntropyLemma-GminInW}
 If there is a conull, Borel subset~$\Omega$ of $\Gamma \backslash
G$, such that $\Omega \cap x G_- \subset x W$, for every $x \in \Omega$,
then $h_\mu(g) \le \log J(g, W)$\index{Jacobian}. 

 \item \label{EntropyLemma-equal}
 If the hypotheses of \pref{EntLem-GminInW} are satisfied, and
equality holds in its conclusion, then $\mu$~is
$W$-\term[measure!invariant]{invariant}.

 \end{enumerate}
 \end{prop}

See \S\ref{EntropyEstimateSect} for a sketch of the proof.

\medskip

Although we have no need for it in these lectures, let us state the
following vast generalization of Thm.~\ref{h(T)=stretch} that calculates
the \term[entropy!of a dynamical
system]{entropy} of any \index{diffeomorphism}diffeomorphism.

\begin{notation}
 Suppose $T$ is a \index{diffeomorphism}diffeomorphism of a smooth
\index{manifold}manifold~$M$.
 \begin{enumerate}
 \item for each $x \in M$ and $\lambda \ge 0$,  we let
 $$ E_\lambda(x) = \bigset{\, v \in \T_x M }{ \limsup_{n \to \infty}
\frac{ \log \|d(T^{-n})_x(v) \|}{n} \le -\lambda } .$$
 Note that 
 \begin{itemize}
 \item $E_\lambda(x)$ is a vector subspace of~$\T_x M$, for each $x$
and~$\lambda$,
 and
 \item we have $E_{\lambda_1}(x) \subset E_{\lambda_2}(x)$ if $\lambda_1
\le \lambda_2$.
 \end{itemize}
 \item For each $\lambda > 0$, the \defit[Lyapunov exponent!multiplicity
of]{multiplicity} of~$\lambda$ at~$x$ is 
 \nindex{$m_x(\lambda)$ = multiplicity of Lyapunov exponent}
 $$m_x(\lambda) = \min_{\mu < \lambda} \bigl( \dim E_\lambda(x) - \dim
E_\mu(x) \bigr) .$$
 By convention, $m_x(0) = \dim E_0(x)$.
 \item We use
 $$\Lyap(T,x) = \{\, \lambda \ge 0 \mid m_x(\lambda) \neq 0 \,\} $$
 to denote the set of \defit[Lyapunov exponent]{Lyapunov exponents} of~$T$
at~$x$. Note that $\sum_{\lambda \in \Lyap(T,x)} m_x(\lambda) = \dim M$,
so $\Lyap(T,x)$ is a finite set, for each~$x$.
 \end{enumerate}
 \end{notation}

\begin{thm}[(\term{Pesin's Entropy Formula})] \label{Pesin}
 Suppose
 \begin{itemize}
 \item $\Omega = M$ is a smooth, compact \index{manifold}manifold,
 \item $\vol$ is a volume form on~$M$,
 and
 \item $T$ is a volume-preserving \index{diffeomorphism}diffeomorphism.
 \end{itemize}
 Then
 $$ h_{\vol}(T) = \int_M \left( \sum_{\lambda \in \Lyap(T,x)} m_x(\lambda)
\, \lambda \right) \, d \vol(x) . $$
 \end{thm}

\begin{exercises}
 \index{entropy}

\item \label{h(Isom)=0}
 Suppose
 \begin{itemize}
 \item $T$ is an \term{isometry} of a compact metric space~$\Omega$,
 \item $\mu$ is a $T$-\term[measure!invariant]{invariant probability
measure} on~$\Omega$,
 and
 \item $\{T^{-k}x\}_{k=1}^\infty$ is dense in the \term[support!of a
measure]{support} of~$\mu$, for a.e.\ $x \in \Omega$.
 \end{itemize}
 Use Thm.~\ref{h(T)FromPast} (and Rem.~\ref{CtblPartn}) to show that
$h_\mu(T) = 0$.
 \hint{Choose a \term[point!of density]{point of density}~$x_0$ for~$\mu$,
and let $\partA$ be a \term[partition!countable]{countable partition}
of~$\Omega$, such that $H(\partA) < \infty$ and $\lim_{x \to x_0} \diam
\bigl( \partA(x) \bigr) = 0$, where $\partA(x)$ denotes the \term[atom (of
a partition)]{atom} of~$\partA$ that contains~$x$. Show, for a conull
subset of~$\Omega$, that each atom of $\partA^+$ is a single point.}

\item Let $G = \SL(2,\real)$\index{SL(2,R)*$\SL(2,\real)$}, and define
$u^t$ and~$a^s$ as usual \see{SL2etcNotation}. Show that if $s > 0$, then
$\{u^t\}$ is the (expanding) \term[subgroup!horospherical]{horospherical
subgroup} corresponding to~$a^s$.

\item \label{HoroSubgrp}
 For each $g \in G$, show that the corresponding
\term[subgroup!horospherical]{horospherical subgroup}~$G_+$ is indeed a
\emph{subgroup} of~$G$.

\item \label{HoroSubalg}
 Given $g \in G$, let
   \begin{equation} \label{HoroLieDefn}
 \textstyle
 \Lie G_+ = \{\, v \in \Lie G \mid \lim_{k \to -\infty}  v (\Ad_G g)^k = 0
\,\} .
 \end{equation}
 \vskip-\smallskipamount
 \begin{enumerate}
 \item Show that $\Lie G_+$ is a Lie subalgebra of~$\Lie G$. 
 \item Show that if $\Ad_G g$ is \index{diagonalizable!over
R*over~$\real$}diagonalizable over~$\real$, then $\Lie G_+$ is the sum of
the \term[eigenspace]{eigenspaces} corresponding to
\term[eigenvalue]{eigenvalues} of $\Ad_G g$ whose absolute value is
(strictly) greater than~$1$.
 \end{enumerate}

\item \label{HoroSubgrpLie}
 Given $g \in G$, show that the corresponding
\term[subgroup!horospherical]{horospherical subgroup}~$G_+$ is the
\index{subgroup!connected|)}connected subgroup of~$G$ whose \term{Lie
algebra} is~$\Lie G_+$.
 \hint{For $u \in G_+$, there is some $k \in \integer$ and some $v \in
\Lie G$, such that $\exp v = g^{-k} u g^k$.}

\item Derive Cor.~\ref{h(G/Gamma)} from Thm.~\ref{h(T)=stretch}, under the
additional assumptions that:
 \begin{enumerate}
 \item $\Gamma \backslash G$ is compact, and
 \item \label{h(G/Gamma)-gss}
 $\Ad_G g$ is \index{diagonalizable!over C*over~$\complex$}diagonalizable
(over~$\complex$). 
 \end{enumerate}

\item \label{h(unip)Ex}
 Derive Cor.~\ref{h(unip)=0} from Cor.~\ref{h(G/Gamma)}.

\item \label{h(diagSL2)Ex}
 Derive Cor.~\ref{h(diagSL2)} from Cor.~\ref{h(G/Gamma)}.

\item \label{h(unip)=0Ex}
 \index{unipotent!flow}
 Give a direct proof of Cor.~\ref{h(unip)=0}.
 \hint{Fix 
 \begin{itemize} \itemsep0pt
 \item a small set~$\Omega$ of positive measure, 
 \item $n \colon \Omega \to \integer^+$ with $x u^{n(x)} \in \Omega$ for
a.e.\ $x \in \Omega$ and $\int_\Omega n \, d\mu < \infty$,
 \item $\lambda > 1$, such that $d(x u^1, y u^1) < \lambda \, d(x,y)$
for $x,y \in \Gamma \backslash G$,
 \item a partition $\partA$ of~$\Omega$, such that $\diam \bigl( \partA(x)
\bigr) \le \epsilon \, \lambda^{-n(x)}$ \seeex{L1Partition}.
 \end{itemize}
 Use the argument of Exer.~\ref{h(Isom)=0}.}

\item \label{L1Partition}
 Suppose
 \begin{itemize}
 \item $\Omega$ is a precompact subset of a \index{manifold}manifold~$M$,
 \item $\mu$ is a probability measure on~$\Omega$,
 and
 \item $\rho \in L^1(\Omega,\mu)$.
 \end{itemize}
 Show there is a \term[partition!countable]{countable partition} $\partA$
of~$\Omega$, such that
 \begin{enumerate}
 \item $\partA$ has finite \term[entropy!of a partition]{entropy},
 and
 \item for a.e.\ $x \in \Omega$, we have $\diam \partA(x) \le e^{\rho(x)}$,
where $\partA(x)$\label{A(x)Defn} denotes the \term[atom (of a
partition)]{atom} of~$\partA$ containing~$x$.
 \end{enumerate}
 \hint{For each~$n$, there is a partition~$\partB_n$ of~$\Omega$ into sets
of diameter less than $e^{-(n+1)}$, such that $\# \partB_n \le C e^{n \dim
M}$. Let $\partA$ be the partition into the sets of the form $B_n \cap
R_n$, where $B_n \in \partB_n$ and $R_n = \rho^{-1}[n,n+1)$.
 Then 
 $ H(\partA) \le H(\partR) + \sum_{n=0}^\infty \mu(R_n) \, \log \#\partB_n
< \infty$.}

\item Derive Cor.~\ref{h(G/Gamma)} from Prop.~\ref{EntropyLemma}.

\item Use Prop.~\ref{EntropyLemma} to prove
 \begin{enumerate}
 \item Lem.~\ref{EntropyExpandSL2},
 and
 \item Lem.~\ref{EntJoining}.
 \end{enumerate}

\item Derive Thm.~\ref{h(T)=stretch} from Thm.~\ref{Pesin}.

 \end{exercises}

\section{Proof of the entropy estimate}
 \index{entropy!of a dynamical system}
 \label{EntropyEstimateSect}

For simplicity, we prove only the special case where $g$~is a diagonal
matrix in $G = \SL(2,\real)$\index{SL(2,R)*$\SL(2,\real)$}. The same
method applies in general, if  ideas from the solution of
Exer.~\ref{h(unip)=0Ex} are added.

\begin{thmref}{EntropyExpandSL2}
 \begin{prop} \label{EstimateSL2}
 Let $G = \SL(2,\real)$\index{SL(2,R)*$\SL(2,\real)$}, and suppose $\mu$
is an $a^s$-\term[measure!invariant]{invariant probability measure} on
$\Gamma \backslash G$.
 \begin{enumerate}
 \item \label{EstimateSL2-Uinv}
 If $\mu$ is $\{u^t\}$-\term[measure!invariant]{invariant}, then
$h(a^s,\mu) = 2|s|$.
 \item \label{EstimateSL2-always}
 We have $h_\mu(a^s) \le 2|s|$.
 \item \label{EstimateSL2-=}
 If $h_\mu(a^s) = 2|s|$ {\rm(}and $s \neq 0${\rm)}, then $\mu$ is
$\{u^t\}$-\term[measure!invariant]{invariant}.
 \end{enumerate}
 \end{prop}
 \end{thmref}

\begin{notation} \label{EntropyProofNotn} \ 
 \begin{itemize}
 \item Let $U = \{u^t\}$.
 \item Let $v^r$ be the opposite
\term[unipotent!subgroup!one-parameter]{unipotent one-parameter subgroup}
\see{vOppUnipDefn}.
 \item Let $a = a^s$, where $s > 0$ is sufficiently large that $e^{-s} <
1/10$, say. Note that 
 \begin{equation} \label{EntPfAContract}
 \lim_{k \to -\infty} a^{-k} u^t a^k = e
 \qquad \mbox{and} \qquad
 \lim_{k \to \infty} a^{-k} v^t a^k = e .
 \end{equation}
 \item Let $x_0$ be a point in the \term[support!of a measure]{support} of~$\mu$.
 \item Choose some small $\epsilon > 0$.
 \item Let
 \begin{itemize}
 \item $U_\epsilon = \{\, u^t \mid -\epsilon < t < \epsilon \,\}$,
 and
 \item $D$ be a small $2$-disk through~$x_0$ that is transverse to the
$U$-orbits,
 \end{itemize}
 so $D U_\epsilon$ is a neighborhood of~$x_0$ that is naturally
homeomorphic to $D \times U_\epsilon$. 
 \item For any subset~$A$ of~$D U_\epsilon$, and any $x \in D$, the
intersection $A \cap x U_\epsilon$ is called a \defit{plaque} of~$A$.
 \end{itemize}
 \end{notation}

\begin{lem} \label{AhasGoodPlaques}
 There is an open neighborhood~$A$ of~$x_0$ in $D U_\epsilon$, such that,
for any \term{plaque}~$F$ of~$A$, and any $k \in \integer^+$,
 \begin{equation} \label{NestedPlaques}
 \mbox{if $F \cap A a^k \neq \emptyset$, then $F \subset A a^k$.}
 \end{equation}
 \end{lem}

\begin{proof}
 We may restate \pref{NestedPlaques} to say:
 $$  \mbox{if $F a^{-k} \cap A \neq \emptyset$, then $F a^{-k} \subset A$.}
$$
 Let $A_0$ be any very small neighborhood of~$x_0$. If $F a^{-k}$
intersects~$A_0$, then we need to add it to~$A$. Thus, we need to add
 $$ A_1 = \bigcup 
 \bigset{ F a^{-k} }{
 \begin{matrix}
 \mbox{$F$ is a \term{plaque} of~$A_0$,} \\
 \mbox{$F a^{-k} \cap A_0 \neq \emptyset$,} \\
 k > 0
 \end{matrix}
 }
 .$$
 This does not complete the proof, because it may be the case that, for
some \term{plaque}~$F$ of~$A_0$, a translate $F a^{-k}$ intersects~$A_1$, but does
not intersect~$A_0$. Thus, we need to add more \term[plaque]{plaques} to~$A$, and continue
inductively:
 \begin{itemize}
 \item Define $\displaystyle A_{n+1} = \bigcup
 \bigset{ F a^{-k} }{
 \begin{matrix}
 \mbox{$F$ is a \term{plaque} of~$A_0$,} \\
 \mbox{$F a^{-k} \cap A_n \neq \emptyset$,} \\
 k > n
 \end{matrix}
 } $.
 \item Let $A = \cup_{n=0}^\infty A_n$.
 \end{itemize}
 It is crucial to note that we may restrict to $k > n$ in the definition
of~$A_{n+1}$ \seeex{PlaqueIntersect}. Because conjugation by~$a^{-k}$
contracts distances along~$U$ exponentially \see{EntPfAContract}, this
implies that $\diam A$ is bounded by a geometric series that converges
rapidly. By keeping $\diam A$ sufficiently small, we guarantee that $A
\subset D U_\epsilon$.
 \end{proof}

\begin{notation} \ 
 \begin{itemize}
 \item Let 
 $$\partA = \{ A, (\Gamma \backslash G) \smallsetminus A \} ,$$
 where $A$ was constructed in Lem.~\ref{AhasGoodPlaques}.
 (Technically, this is not quite correct --- the proof of
Lem.~\ref{subordinate} shows that we should take a similar, but more
complicated, partition of~$\Gamma \backslash G$.)
 \item Let $\past{\partA} = \bigvee_{k=1}^\infty \partA a^k$
\cf{ConditionOnPast}.
 \item Let
 \nindex{$\partA(x)$ or $\past{\partA}(x)$ = atom of a partition,
\pageref{A(x)Defn}}
 $\past{\partA}(x)$ be the \term[atom (of a
partition)]{atom} of~$\past{\partA}$ containing~$x$, for
each $x \in \Gamma \backslash G$.
 \end{itemize}
 \end{notation}

\begin{assump}
 Let us assume that the measure~$\mu$ is \term[ergodic!measure]{measure}
for~$a$. (The general case can be obtained from this by considering the
\term[ergodic!decomposition]{ergodic decomposition} of~$\mu$.)
 \end{assump}

\begin{lem} \label{subordinate}
 The partition $\past{\partA}$ is
\defit[partition!subordinate]{subordinate} to~$U$. That is, for a.e.\ $x
\in \Gamma \backslash G$,
 \begin{enumerate}
 \item \label{subordinate-inU}
 $\past{\partA}(x) \subset xU$, 
 and
 \item \label{subordinate-neigh}
 more precisely, $\past{\partA}(x)$ is a relatively compact, open
neighborhood of~$x$ {\rm(}with respect to the orbit topology of~$xU${\rm)}.
 \end{enumerate}
 \end{lem}

\begin{proof}
 For a.e.\ $x \in A$, we will show that $\past{\partA}(x)$ is simply the
\term{plaque} of~$A$ that contains~$x$. Thus, \pref{subordinate-inU}
and~\pref{subordinate-neigh} hold for a.e.\ $x \in A$. By
\term[ergodic!action]{ergodicity}, it immediately follows that the
conditions hold for a.e.\ $x \in \Gamma \backslash G$.

If $F$ is any \term{plaque} of~$A$, then \pref{NestedPlaques} implies, for each
$k > 0$, that $F$ is contained in a single \term[atom (of a
partition)]{atom} of~$\partA a^k$.
Therefore, $F$ is contained in a single \term[atom (of a partition)]{atom}
of~$\past{\partA}$. The problem is to show that $\past{\partA}(x)$
contains only a single \term{plaque}.

 Let $V_\epsilon = \{v^r\}_{r=-\epsilon}^\epsilon$, and pretend, for the
moment, that $x_0 U_\epsilon V_\epsilon$ is a neighborhood of~$x_0$.
(Thus, we are we are ignoring $\{a^s\}$, and pretending that $G$ is
2-dimensional\index{dimension!of a Lie group}.) For $k > 0$, we know that
conjugation by~$a^k$ contracts~$\{v^r\}$, so $A a^k$ is very thin in the
$\{v^r\}$-direction (and correspondingly long in the $\{u^t\}$-direction).
In the limit, we conclude that the \term[atom (of a
partition)]{atoms} of $\past{\partA}$ are \index{infinitely!thin}infinitely
thin in the $\{v^r\}$-direction. The union of any two
\term[plaque]{plaques} has a nonzero length in the $\{v^r\}$-direction, so
we conclude that an \term[atom (of a partition)]{atom} of $\past{\partA}$
contains only one \term{plaque}, as desired.

To complete the proof, we need to deal with the $\{a^s\}$-direction.
(Unfortunately, this direction is not contracted by~$a^k$, so the argument
of the preceding paragraph does not apply.) To do this, we alter the
definition of~$\partA$. 
 \begin{itemize}
 \item Let $\partB$ be a \term[partition!countable]{countable partition}
of~$D$, such that $H(\partB) < \infty$ and $\lim_{x \to x_0} \diam \bigl(
\partB(x) \bigr) = 0$.
 \item Let $\hatpartA$ be the corresponding partition of~$A$:
 $$ \hatpartA = \{\, (B U_\epsilon) \cap A \mid B \in \partB\,\} .$$
 \item Let $\partA = \hatpartA \cup \{ (\Gamma \backslash G)
\smallsetminus A \}$.
 \end{itemize}
 Then $\partA$ is a countable partition of $\Gamma \backslash G$ with
$H(\partA) < \infty$, so $h_\mu(a) = H(\past{\partA} a^{-1} \mid
\past{\partA})$ \seeand{CtblPartn}{h(T)=H(TA|A)}. Ergodicity implies that
$x a^k$ is close to~$x_0$ for some values of~$k$. From the choice
of~$\partB$, this implies that $\past{\partA}(x)$ has small length in the
$\{a^s\}$-direction. In the limit, $\past{\partA}(x)$ must be
\index{infinitely!thin}infinitely thin in the $\{a^s\}$-direction.
 \end{proof}
 \index{infinitely!small|indsee{infinitesimal}}

\begin{proof}[Proof of \fullref{EstimateSL2}{Uinv}]
 We wish to show $H(\past{\partA} a^{-1} \mid \past{\partA}) = 2s$
\see{h(T)=H(TA|A)}.

For any $x \in \Gamma \backslash G$, let
 \begin{itemize}
 \item $\mu_x$ be the \term[measure!conditional]{conditional
measure} induced by~$\mu$ on~$\past{\partA}(x)$,
 and
 \item $\lambda$ be the \term[measure!Haar]{Haar measure} (that is, the
\term[measure!Lebesgue]{Lebesgue measure}) on~$xU$.
 \end{itemize}
 Because $\past{\partA}(x) \subset x U$, and $\mu$~is
$U$-\term[measure!invariant]{invariant}, we know that 
 \begin{equation} \label{EstMu=Lambda}
 \mbox{$\mu_x$ is the restriction of~$\lambda$ to~$\past{\partA}(x)$ (up
to a scalar multiple).} 
 \end{equation}
 Now $\past{\partA} \subset \past{\partA} a^{-1}$, so $\past{\partA}
a^{-1}$ induces a partition $\partA_x$ of~$\past{\partA}(x)$. By
definition, we have
 $$ H (\past{\partA} a^{-1} \mid \past{\partA}) 
 = -\int_{\Gamma \backslash G} \log f \, d\mu ,$$
 where
 \begin{equation} \label{EstFmu=Flambda}
 f(x) = \mu_x \bigl( \partA_x(x) \bigr) = \frac{\lambda \bigl( \partA_x(x)
\bigr)}{\lambda \bigl( \past{\partA}(x) \bigr)} .
 \end{equation}
 Note that, because translating by~$a$ transforms $\past{\partA} a^{-1}$
to~$\past{\partA}$,  we have
 $$ x \, \{\, u \in U \mid x u \in \partA_x(x) \,\} \, a
 = x a \, \{\, u \in U \mid x a u \in \past{\partA}(x) \,\} .$$
 Conjugating by~$a$ expands $\lambda$ by a factor of
$e^{2s}$, so this implies
 $$ \frac{ \lambda \bigl( \past{\partA}(x a) \bigr) }{ \lambda \bigl(
\past{\partA}(x) \bigr)}
 =  \frac{ e^{2s} \lambda \bigl( \partA_x(x) \bigr) }{ \lambda \bigl(
\past{\partA}(x) \bigr)}
 = e^{2s} f(x)
 .$$
 Because $0 \le f(x) \le 1$ and $e^{2s}$ is a constant, we conclude that 
 $$ \Bigl( \log \lambda \bigl( \past{\partA}(x a) \bigr) - \log \lambda
\bigl( \past{\partA}(x) \bigr) \Bigr)^+ \in L^1(\Gamma \backslash G, \mu)
,$$
 so Lem.~\ref{Int(f-f)=0} below implies
 $$ - \int_{\Gamma \backslash G} \log f \, d\mu = \log e^{2s} = 2s ,$$
 as desired.
 \end{proof}

The following observation is obvious (from the
\term[measure!invariant]{invariance} of~$\mu$) if $\psi \in L^1(\Gamma
\backslash G, \mu)$. The general case is proved in
Exer.~\ref{Int(f-f)=0Ex}.

\begin{lem} \label{Int(f-f)=0}
 Suppose 
 \begin{itemize}
 \item $\mu$ is an $a$-\term[measure!invariant]{invariant probability
measure} on $\Gamma \backslash G$,
 and
 \item $\psi$ is a real-valued, measurable function on
$\Gamma \backslash G$, 
 \end{itemize}
 such that 
 $$ \bigl( \psi(x a) - \psi(x) \bigr)^+ \in L^1(\Gamma \backslash G,
\mu) ,$$
 where $(\alpha)^+ = \max(\alpha,0)$. Then
 $$ \int_{\Gamma \backslash G} \bigl( \psi(x a) - \psi(x) \bigr) \,
d\mu(x) = 0 .$$
 \end{lem}

\begin{proof}[Proof of \fullref{EstimateSL2}{always}]
 This is similar to the proof of \pref{EstimateSL2-Uinv}. 
 Let $\lambda_x = \lambda/\lambda \bigl( \past{\partA}(x) \bigr)$ be the
normalization of~$\lambda$ to a probability measure on $\past{\partA}(x)$. (In
the proof of \pref{EstimateSL2-Uinv}, we had $\lambda_x = \mu_x$
\see{EstMu=Lambda}, so we did not bother to define~$\lambda_x$.)
 Also, define
 $$ f_\mu(x) = \mu_x \bigl( \partA_x(x) \bigr)
 \qquad \mbox{and} \qquad 
 f_\lambda(x) =  \lambda_x \bigl( \partA_x(x) \bigr) .$$
 (In the proof of \pref{EstimateSL2-Uinv}, we had $f_\mu = f_\lambda$
\see{EstFmu=Flambda}; we simply called the function~$f$.)

 We have
 $$ h_\mu(a) = H (\past{\partA} a^{-1} \mid \past{\partA}) 
 = -\int_{\Gamma \backslash G} \log f_\mu \, d\mu ,$$
 and the proof of \pref{EstimateSL2-Uinv} shows that
 $$ - \int_{\Gamma \backslash G} \log f_\lambda \, d\mu = 2s ,$$
 so it suffices to show
 $$ \int_{\Gamma \backslash G} \log f_\lambda \, d\mu \le 
 \int_{\Gamma \backslash G} \log f_\mu \, d\mu .$$
 Thus, we need only show, for a.e.\ $x \in \Gamma \backslash G$, that
 \begin{equation} \label{EstimateSL2AlwaysPf-IntA<intA}
 \int_{\past{\partA}(x)} \log f_\lambda \, d\mu_x
 \le \int_{\past{\partA}(x)} \log f_\mu \, d\mu_x .
 \end{equation}
 
 Write $\partA_x = \{A_1,\ldots,A_n\}$. For $y \in A_i$, we have
 $$ f_\lambda(y) = \lambda_x(A_i) 
 \qquad \mbox{and} \qquad
 f_\mu(y) = \mu_x(A_i),$$
 so
 $$  \int_{\past{\partA}(x)} (\log f_\lambda - \log f_\mu) \, d\mu_x
 = \sum_{i=1}^n \log \frac{\lambda_x(A_i)}{\mu_x(A_i)} \, \mu_x(A_i) .$$
 Because
 $$ \sum_{i=1}^n \frac{\lambda_x(A_i)}{\mu_x(A_i)} \mu_x(A_i)
 = \sum_{i=1}^n \lambda_x(A_i) = \lambda_x \bigl( \past{\partA}(x) \bigr) =
1,$$
 the concavity of the $\log$ function implies 
 \begin{align} \label{EstimateSL2AlwaysPf-sum<0}
 \sum_{i=1}^n \log \frac{\lambda_x(A_i)}{\mu_x(A_i)} \mu_x(A_i)
 &\le \log \sum_{i=1}^n \frac{\lambda_x(A_i)}{\mu_x(A_i)} \mu_x(A_i)
 \\&= \log 1 \notag
 \\&= 0 . \notag
 \end{align}
 This completes the proof.
 \end{proof}

\begin{proof}[Proof of \fullref{EstimateSL2}{=}]
 Let $\mu_U$ be the \term[measure!conditional]{conditional measure}
induced by~$\mu$ on the orbit~$xU$. To show that $\mu$ is
$U$-\term[measure!invariant|)]{invariant}, we wish to show that $\mu$ is
equal to~$\lambda$ (up to a scalar multiple).

 We must have equality in the proof of~\fullref{EstimateSL2}{always}.
Specifically, for a.e.\ $x \in \Gamma \backslash G$, we must have equality
in \pref{EstimateSL2AlwaysPf-IntA<intA}, so we must have equality in
\pref{EstimateSL2AlwaysPf-sum<0}. Because the log function is strictly
concave, we conclude that 
 $$ \frac{\lambda_x(A_i)}{\mu_x(A_i)} = \frac{\lambda_x(A_j)}{\mu_x(A_j)}
$$ for all~$i,j$. Since
 $$ \sum_{i=1}^n \lambda_x(A_i) = \lambda_x(\partA_x) = 1
 = \mu_x(\partA_x) = \sum_{i=1}^n \mu_x(A_i) ,$$
 we conclude that $\lambda_x(A_i) = \mu_x(A_i)$. This means that $\mu_U(A)
= \lambda(A)$ for all \term[atom (of a
partition)]{atoms} of $\partA_x$. By applying the same argument
with $a^k$ in the place of~$a$ (for all $k \in \integer^+$), we conclude
that $\mu_U(A) = \lambda(A)$ for all~$A$ in a collection that generates
the $\sigma$-algebra of all measurable sets in~$xU$. Therefore $\mu_U =
\lambda$.
 \end{proof}

 \onlyoneexercise
 \begin{exercises}

 \item \label{PlaqueIntersect}
 In the proof of Lem.~\ref{AhasGoodPlaques}, show that if 
 \begin{itemize}
 \item $F$ is a \term{plaque} of~$A_0$,
 \item $1 \le k \le n$,
 and
 \item $F a^{-k} \cap A_n \neq \emptyset$,
 \end{itemize}
 then $F a^{-k} \cap (A_0 \cup \cdots \cup A_{n-1}) \neq \emptyset$.
 \hint{Induction on~$k$.}

 \end{exercises}
 \endgroup

\begin{notes}
 \index{entropy}

The \term[entropy!of a dynamical system]{entropy} of dynamical systems is
a standard topic that is discussed in many textbooks, including
\cite[\S4.3--\S4.5]{KatokHasselblatt-intro} and \cite[Chap.~4]{Walters}.

\notesect{TwoEgsSect}

\term[irrational!rotation]{Irrational rotations} $T_\beta$ and
\term[Bernoulli shift]{Bernoulli shifts} $\Bern$ are standard examples. The
\term{Baker's Transformation}~$\Baker$ is less common, but it appears in
\cite[p.~22]{BekkaMayer}, for example.

\notesect{UnpreditableSect}

Proposition~\ref{past=future->h(T)=0}\index{past determines the future}
appears in standard texts, including \cite[Cor.~4.14.1]{Walters}.

\notesect{DefnEntopySect}

This material is standard (including the properties of \term{entropy}
developed in the exercises).

Our treatment of the \term[entropy!of a partition]{entropy} of a partition
is based on \cite{Khinchin}. The elementary argument mentioned at the end
of Rem.~\fullref{ErgPartMotiv}{conditional} appears in
\cite[pp.~9--13]{Khinchin}.

It is said that the \term[entropy!of a dynamical
system]{entropy of a dynamical system} was first defined
by A.~N.~Kolmogorov \cite{Kolmogorov1, Kolmogorov2}, and that much of the
basic theory is due to Ya.~Sinai \cite{Sinai1, Sinai2}. 

\term[entropy!topological]{Topological entropy} was defined
by R.~L.~Adler, A.~G.~Konheim, and M.~H.~McAndrew \cite{AdlerTopEnt}.  Our
discussion in Rem.~\ref{TopEntropy} is taken from
\cite{HasselblattKatok-principal}. 

L.~W.~Goodwyn \cite{Goodwyn} proved the inequality $h_\mu(T) \le
h_{\text{top}}(T)$. A simple proof of a stronger result appears in
\cite{Misiurewicz-Variational}.

\notesect{CalcEntropySect}

This material is standard.

Theorem~\ref{genpart} is due to Ya.~Sinai.

\notesect{h(g)Sect}

Corollary~\ref{h(G/Gamma)} was proved by R.~Bowen \cite{Bowen} when
$\Gamma \backslash G$ is compact. The general case (when $\Gamma
\backslash G$ has finite volume) was apparently already known to
dynamicists in the Soviet Union. For example, it follows from the argument
that proves \cite[(8.35), p.~68]{MargulisThesis}.

A complete proof of the crucial \term[entropy!of a dynamical
system]{entropy} estimate \pref{EntropyLemma}
appears in \cite[Thm.~9.7]{MargulisTomanov-Ratner}.
 It is based on ideas from \cite{LedrappierYoung}.

Pesin's Formula \pref{Pesin} was proved in \cite{PesinFormula}. Another
proof appears in \cite{Mane}.

Exercise~\ref{L1Partition} is \cite[Lem.~2]{Mane}.

\notesect{EntropyEstimateSect}

This section is based on \cite[\S9]{MargulisTomanov-Ratner}.

Lemma~\ref{Int(f-f)=0} is proved in \cite[Prop.~2.2]{LedrappierStrelcyn}.

\end{notes}

%% file: RatnerErgodic.tex
\mychapter{Facts from Ergodic Theory} \label{ErgodicChap}

This chapter simply gathers some necessary background results, mostly
without proof.

\section{Pointwise Ergodic Theorem} \label{PtwiseErgSect}

In the proof of \term[Ratner's Theorems]{Ratner's Theorem} (and in many other situations), one
wants to know that the orbits of a flow are \term{uniformly distributed}.
It is rarely the case that \emph{every} orbit is uniformly distributed
(that is what it means to say the flow is \defit[ergodic!uniquely]{uniquely
ergodic}), but the \term[Theorem!Pointwise Ergodic]{Pointwise Ergodic
Theorem} \pref{PtwiseErgThmFlow} shows that if the flow is
``\term[ergodic!flow|(]{ergodic}," a much weaker condition, then
\emph{almost every} orbit is uniformly distributed. (See the exercises for
a proof.)

\begin{defn}
 A \term[flow!measure-preserving]{measure-preserving flow} $\varphi_t$ on a
probability space $(X,\mu)$ is \defit[ergodic!flow]{ergodic} if, for each
$\varphi_t$-\term[set!invariant]{invariant subset}~$A$ of~$X$, we have
either $\mu(A) = 0$ or $\mu(A) = 1$.
 \end{defn}

\begin{eg} For $G = \SL(2,\real)$\index{SL(2,R)*$\SL(2,\real)$} and
$\Gamma = \SL(2,\integer)$\index{SL(2,Z)*$\SL(2,\integer)$}, the
\term[horocycle!flow]{horocycle flow} $\eta_t$ and the
\term[geodesic!flow]{geodesic flow} $\gamma_t$ are
\term[ergodic!flow]{ergodic} on $\Gamma \backslash G$ (with respect to the
\term[measure!Haar]{Haar measure} on $\Gamma \backslash G$)
\seeand{HoroFlowErg}{GeodFlowErg}. These are special cases of the 
\term[Theorem!Moore Ergodicity]{Moore Ergodicity Theorem} \pref{MooreErg},
which implies that most flows of
\term[subgroup!one-parameter]{one-parameter subgroups} on $\Gamma
\backslash \SL(n,\real)$\index{SL(l,R)*$\SL(\ell,\real)$} are ergodic, but
the ergodicity of $\gamma_t$ can easily be proved from scratch
\seeex{GeodFlowErgEx}.
 \end{eg}

\begin{thm}[{(\term[Theorem!Pointwise Ergodic]{Pointwise Ergodic
Theorem})}] \label{PtwiseErgThmFlow}
 Suppose 
 \begin{itemize}
 \item $\mu$ is a probability measure on a \term{locally compact},
\term{separable} metric space~$X$,
 \item $\varphi_t$ is an \term[ergodic!flow]{ergodic},
\term[flow!measure-preserving]{measure-preserving flow} on~$X$,
 and
 \item $f \in L^1(X,\mu)$.
 \end{itemize}
 \nindex{$L^1(X,\mu)$ = Banach space of real-valued $L^1$ functions on~$X$}
 Then
 \begin{equation} \label{PtwiseErgThmFlow-limit}
 \frac{1}{T} \int_0^T f\bigl( \varphi_t(x) \bigr) \, dt
 \to \int_X f \, d \mu
 ,
 \end{equation}
 for a.e.\ $x \in X$.
 \end{thm}

\begin{defn} \label{GenericDefn}
 A point $x \in X$ is \defit[generic!point]{generic} for~$\mu$ if
\pref{PtwiseErgThmFlow-limit} holds for every
\term[function!continuous!uniformly]{uniformly continuous, bounded
function} on~$X$. In other words, a point is generic for~$\mu$ if its
orbit is \term{uniformly distributed} in~$X$.
 \end{defn}

\begin{cor} \label{GenericPtae}
 If $\varphi_t$ is \term[ergodic!flow]{ergodic}, then almost every point
of~$X$ is \term[generic!point]{generic} for~$\mu$.
 \end{cor}

The converse of this corollary is true \seeex{Generic->Erg}.

\begin{exercises}

\item Prove Cor.~\ref{GenericPtae} from Thm.~\ref{PtwiseErgThmFlow}.

\item \label{Generic->Erg}
 Let $\varphi_t$ be a \term[flow!measure-preserving]{measure-preserving
flow} on $(X,\mu)$. Show that if $\varphi_t$ is \emph{not}
\term[ergodic!flow]{ergodic}, then almost \emph{no} point of~$X$ is
\term[generic!point]{generic} for~$\mu$.

\item \label{UnbddHitsAE}
 Let
 \begin{itemize}
 \item $\varphi_t$ be an \term[ergodic!flow]{ergodic}
\term[flow!measure-preserving]{measure-preserving flow} on $(X,\mu)$,
 and
 \item $\Omega$ be a non-null subset of~$X$.
 \end{itemize}
 Show, for a.e.\ $x \in X$, that
 $$ \{\, t \in \real^+ \mid \varphi_t(x) \in \Omega \,\}$$
 is unbounded.
 \hint{Use the Pointwise Ergodic Theorem.}

\item \label{PtwiseErgZMaxl}
 Suppose 
 \begin{itemize}
 \item $\phi \colon X \to X$ is a measurable bijection of~$X$,
 \item $\mu$~is a $\phi$-\term[measure!invariant|(]{invariant} probability
measure on~$X$,
 \item $f \in L^1(X,\mu)$,
 and
 \item $S_n(x) =  f(x) + f \bigl( \phi(x) \bigr) +
\cdots + f \bigl( \phi^{n-1}(x) \bigr)$.
 \end{itemize}
  Prove the \defit[Theorem!Maximal Ergodic]{Maximal Ergodic
Theorem}: for every $\alpha \in \real$, if we let
 $$ E = \bigset{ x \in X }{ \sup_n \frac{S_n(x)}{n} > \alpha} ,$$
 then $\int_E f \, d\mu \ge \alpha \, \mu(E)$.
 \hint{Assume $\alpha = 0$. Let $S^+_n(x) = \max_{0\le k\le n} S_k(x)$,
and $E_n = \{\, x \mid S^+_n > 0\,\}$, so $E = \cup_n E_n$. For $x \in
E_n$, we have $f(x) \ge S^+_n(x) - S^+_n \bigl( \phi(x) \bigr)$, so
$\int_{E_n} f\, d\mu \ge 0$.}

\item \label{PtwiseErgZ}
 Prove the \defit[Theorem!Pointwise Ergodic]{Pointwise Ergodic Theorem}
for $\phi$, $\mu$, $f$, and~$S_n$ as in Exer.~\ref{PtwiseErgZMaxl}. That
is, if $\phi$ is \term[ergodic!dynamical system]{ergodic}, show, for
a.e.~$x$, that 
 $$ \lim_{n \to \infty} \frac{S_n(x)}{n} = \int_X f\, d\mu .$$
 \hint{If $\{\, x \mid \limsup S_n(x)/n >
\alpha \,\}$ is not null, then it must be conull, by ergodicity. So the
Maximal Ergodic Theorem (Exer.~\ref{PtwiseErgZMaxl}) implies $\int_X f \,
d\mu \ge \alpha$.}

\item \label{PtwiseZNotErg}
 For $\phi$, $\mu$, $f$,
and~$S_n$ as in Exer.~\ref{PtwiseErgZMaxl}, show there is a function
$f^* \in L^1(X,\mu)$, such that:
 \begin{enumerate}
 \item for a.e.~$x$, we have
 $ \lim_{n \to \infty} {S_n(x)}/{n} = f^*(x)$,
 \item for a.e.~$x$, we have
 $ f^* \bigl( \phi(x) \bigr) = f^*(x)$,
 and
 \item $\int_X f^* \, d\mu = \int_X f \, d\mu$.
 \end{enumerate}
 \smallskip
 (This generalizes Exer.~\ref{PtwiseErgZ}, because we do \emph{not} assume
$\phi$ is \term[ergodic!dynamical system]{ergodic}.)
 \hint{For $\alpha < \beta$, replacing $X$ with the
$\phi$-\term[set!invariant]{invariant} set
 $$ X_\alpha^\beta = \{\, x \mid \liminf {S_n(x)}/{n} < \alpha < \beta <
\limsup {S_n(x)}/{n} \,\} $$
 and applying Exer.~\ref{PtwiseErgZMaxl} yields
 $\int_{X_\alpha^\beta} f \, d\mu \le \alpha \, \mu(X_\alpha^\beta)$ and
 $\int_{X_\alpha^\beta} f \, d\mu \ge \beta \, \mu(X_\alpha^\beta)$.}

\item \label{PtwiseErgThmFlowPfEx}
 Prove Thm.~\ref{PtwiseErgThmFlow}.
 \hint{Assume $f \ge 0$ and apply Exer.~\ref{PtwiseErgZ} to the function
$\overline{f}(x) = \int_0^1 f \bigl( \varphi_t(x) \bigr) \, dt$.}

\item \label{Int(f-f)=0Ex}
 Prove Lem.~\ref{Int(f-f)=0}.
 \hint{The \term[Theorem!Pointwise Ergodic]{Pointwise Ergodic Theorem}
\pref{PtwiseZNotErg} remains valid if $f = f_+ - f_-$, with $f_+ \ge 0$,
$f_- \le 0$, and $f_+ \in L^1(X,\mu)$, but the limit $f^*$ can be 
$-\infty$ on a set of positive measure. Applying this to 
 $f(x) = \psi(x a) - \psi(x)$,
 we conclude that 
 $$f^*(x) =
 \lim_{n \to \infty}  \frac{\psi(x a^n) - \psi(x)}{n}$$
 exists a.e.\ (but may be~$-\infty$). Furthermore,
 $\int_{\Gamma \backslash G} f^* \, d\mu = \int_{\Gamma \backslash G} f \,
dx$. Since $\psi(x a^n)/n \to 0$ in measure, there is a sequence $n_k \to
\infty$, such that $\psi(x a^{n_k})/n_k \to 0$~a.e. So $f^*(x) = 0$ a.e.}

\item \label{UnifErg<>UnifConv}
 Suppose 
 \begin{itemize}
 \item $X$ is a compact metric space,
 \item $\phi \colon X \to X$ is a \index{homeomorphism}homeomorphism,
 and
 \item $\mu$ is a $\phi$-\term[measure!invariant]{invariant} probability
measure on~$X$.
 \end{itemize}
 Show $\phi$ is \term[ergodic!uniquely]{uniquely ergodic} if and only if,
for every continuous function~$f$ on~$X$, there is a constant~$C$,
depending on~$f$, such that
 $$ \lim_{n \to \infty} \frac{1}{n} \sum_{k=1}^n f \bigl( \phi^k(x) \bigr)
= C ,$$
 uniformly over $x \in X$.

 \end{exercises}

\section{Mautner Phenomenon} \label{MautnerSect}

We prove that the \term[geodesic!flow]{geodesic flow} is
\term[ergodic!flow]{ergodic} \seecor{GeodFlowErg}. The same methods apply
to many other flows on $\Gamma \backslash G$.

\begin{defn}
 Suppose $\varphi_t$ is a flow on a measure space $(X,\mu)$, and $f$~is a
measurable function on~$X$.
 \begin{itemize}
 \item $f$ is \defit[function!invariant!essentially]{essentially
invariant} if, for each $t \in \real$, we have $f \bigl( \varphi_t(x)
\bigr) = f(x)$ for a.e.\ $x \in X$.
 \item $f$ is \defit[function!constant (essentially)]{essentially constant}
if $f(x) = f(y)$ for a.e.\ $x,y \in X$.
 \end{itemize}
 \end{defn}

\begin{rem} \label{EssConst<>Ergodic}
 It is obvious that any \term[function!constant (essentially)]{essentially
constant} function is \term[function!invariant!essentially]{essentially
invariant}. The converse holds if and only if $\varphi_t$ is
\term[ergodic!flow]{ergodic} \seeex{EssConst->Ergodic}. 
 \end{rem} 

See Exers.~\ref{MautnerSL2PfEx} and~\ref{GeodFlowErgEx} for the proof of
the following proposition and the first corollary.

\begin{prop}[(Mautner Phenomenon)] \label{MautnerSL2}
 Suppose 
 \begin{itemize}
 \item $\mu$ is a probability measure on $\Gamma \backslash G$,
 \item $f \in L^2(\Gamma \backslash G, \mu)$,
 and
 \item $u^t$ and $a^s$ are
\term[subgroup!one-parameter!hyperbolic]{one-parameter subgroups} 
 \index{unipotent!subgroup!one-parameter}
 of~$G$,
 \end{itemize}
 such that
 \begin{itemize}
 \item $a^{-s} u^t a^s = u^{e^s t}$,
 \item $\mu$ is \term[measure!invariant]{invariant} under both $u^t$
and~$a^s$,
 and
 \item $f$ is essentially
$a^s$-\term[function!invariant!essentially]{invariant}.
 \end{itemize}
 Then $f$ is essentially
$u^t$-\term[function!invariant!essentially]{invariant}.
 \end{prop}

\begin{cor} \label{GeodFlowErg}
 The \term[geodesic!flow]{geodesic flow} $\gamma_t$ is
\term[ergodic!flow]{ergodic} on $\Gamma \backslash
\SL(2,\real)$\index{SL(2,R)*$\SL(2,\real)$}.
 \end{cor}

The following corollary is obtained by combining Prop.~\ref{MautnerSL2}
with Rem.~\ref{EssConst<>Ergodic}.

\begin{cor} \label{MautnerErgSL2}
 Suppose 
 \begin{itemize}
 \item $\mu$ is a probability measure on $\Gamma \backslash G$,
 and
 \item $u^t$ and $a^s$ are
\term[subgroup!one-parameter!hyperbolic]{one-parameter subgroups}
 \index{unipotent!subgroup!one-parameter}
 of~$G$,
 \end{itemize}
 such that
 \begin{itemize}
 \item $a^{-s} u^t a^s = u^{e^s t}$,
 \item $\mu$ is \term[measure!invariant]{invariant} under both $u^t$
and~$a^s$,
 and
 \item $\mu$ is \term[ergodic!measure]{ergodic} for~$u^t$.
 \end{itemize}
 Then $\mu$ is ergodic for~$a^s$.
 \end{cor}

The following result shows that flows on $\Gamma \backslash G$ are often
\term[ergodic!flow]{ergodic}. It is a vast generalization of the fact
that the \term[horocycle!flow]{horocycle flow}~$\eta_t$ and the
\term[geodesic!flow]{geodesic flow}~$\gamma_t$ are
\term[ergodic!flow|)]{ergodic} on $\Gamma \backslash
\SL(2,\real)$\index{SL(2,R)*$\SL(2,\real)$}.
 \index{flow!ergodic|indsee{ergodic~flow}}

\begin{thm}[{(\term[Theorem!Moore Ergodicity]{Moore Ergodicity Theorem})}]
\label{MooreErg}
 Suppose
 \begin{itemize}
 \item $G$ is a \index{subgroup!connected}connected, \index{group!simple
(or almost)}simple \term{Lie group} with finite center,
 \item $\Gamma$ is a \term{lattice} in~$G$,
 and
 \item $g^t$ is a \term[subgroup!one-parameter]{one-parameter subgroup}
of~$G$, such that its closure $\closure{\{g^t\}}$ is not compact.
 \end{itemize}
 Then $g^t$ is \term[ergodic!flow]{ergodic} on $\Gamma \backslash G$
{\rm(}w.r.t.\ the \term[measure!Haar]{Haar measure} on $\Gamma \backslash
G${\rm)}.
 \end{thm}

\begin{cor} \label{HoroFlowErg}
 The \term[horocycle!flow]{horocycle flow} $\eta_t$ is
\term[ergodic!flow]{ergodic} on $\Gamma \backslash
\SL(2,\real)$\index{SL(2,R)*$\SL(2,\real)$}.
 \end{cor}

\begin{rem} \label{MixingRem}
 The conclusion of Thm.~\ref{MooreErg} can be strengthened: not only is
$g^t$ \term[ergodic!flow]{ergodic} on $\Gamma \backslash G$, but it is
\defit{mixing}. That is, if
 \begin{itemize}
 \item $A$ and~$B$ are any two measurable subsets of $\Gamma \backslash G$,
 and
 \item $\mu$ is the $G$-\term[measure!invariant]{invariant} probability
measure on $\Gamma \backslash G$,
 \end{itemize}
 then $(A g^t) \cap B \to \mu(A) \, \mu(B)$ as $t \to \infty$. 
 \end{rem}

The following theorem is a restatement of this remark in terms of
functions. 

\begin{thm} \label{matcoeffs->0}
 If
 \begin{itemize}
 \item $G$ is a \index{subgroup!connected}connected, \index{group!simple
(or almost)}simple \term{Lie group} with finite center,
 \item $\Gamma$ is a \term{lattice} in~$G$, 
 \item $\mu$ is the $G$-\term[measure!invariant]{invariant} probability
measure on $\Gamma \backslash G$, and
 \item $g^t$ is a \term[subgroup!one-parameter]{one-parameter subgroup}
of~$G$, such that $\closure{\{g_t\}}$ is not compact,
 \end{itemize}
 then 
 $$\lim_{t \to \infty} \int_{\Gamma \backslash G} \phi( x g^t) \, \psi(x)
\, d\mu =  \|\phi\|_2 \, \|\psi\|_2 ,$$
 for every $\phi,\psi \in L^2(\Gamma \backslash G, \mu)$.
 \end{thm}

\begin{rem}
 For an elementary (but very instructive) case of the following proof,
assume
 $G = \SL(2,\real)$\index{SL(2,R)*$\SL(2,\real)$}, and $g^t = a^t$ is
diagonal. Then only Case~\ref{matcoeffs->0-A} is needed, and we have
 $$ U = \begin{bmatrix} 1 & 0 \\ * & 1 \end{bmatrix}
 \text{ \qquad and \qquad}
 V = \begin{bmatrix} 1 & * \\ 0 & 1 \end{bmatrix}
 .$$
 (Note that $\langle U, V \rangle = G$; that is, $U$ and~$V$, taken
together, generate~$G$.)
 \end{rem}

\begin{proof} (\emphit{Requires some Lie theory and Functional Analysis})
  \begin{itemize}
 \item Let $\hilbert = 1^\perp$ be the (closed) subspace of $L^2(\Gamma
\backslash G, \mu)$ consisting of the functions of integral~$0$. Because
the desired conclusion is obvious if $\phi$ or~$\psi$ is constant, we may
assume $\phi,\psi \in \hilbert$.

 \item For each $g \in G$, define the \term{unitary operator} $g^\rho$ on
$\hilbert$ by
 $ (\phi g^\rho)(x) = \phi(x g^{-1})$.

\item Define $\langle \phi \mid \psi \rangle = \int_{\Gamma \backslash G}
\phi \psi \, d\mu$.

\item Instead of taking the limit only along a
\term[subgroup!one-parameter]{one-parameter subgroup}~$g^t$, we allow a
more general limit along any sequence $g_j$, such that $g_j \to \infty$
in~$G$; that is, $\{g_j\}$ has no convergent subsequences.

 \end{itemize}

 \setcounter{case}{0}

\begin{case} \label{matcoeffs->0-A}
 Assume $\{g_j\}$ is contained in a \term[torus!hyperbolic]{hyperbolic
torus}~$A$.
 \end{case} 
 By passing to a subsequence, we may assume $g_j^\rho$ \term[weak
convergence]{converges weakly}, to some operator~$E$; that is,
 $$ \langle \phi g_j^\rho \mid \psi \rangle
 \to \langle \phi E \mid \psi \rangle
 \mbox{ \ for every $\phi,\psi \in \hilbert$} .$$
 Let 
 \begin{align*}
 U &= \{\, u \in G \mid g_j u g_j^{-1} \to e \,\} \\
 \noalign{and}
 V &= \{\, v \in G \mid g_j^{-1} v g_j \to e \,\} 
 . \end{align*}
 For $v \in V$, we have
 \begin{align*}
 \langle \phi v^\rho E \mid \psi \rangle
 &= \lim_{j \to \infty} \langle \phi v^\rho g_j^\rho \mid \psi \rangle
 \\&= \lim_{j \to \infty} \langle \phi g_j^\rho (g_j^{-1} v g_j)^\rho \mid
\psi \rangle
 \\&= \lim_{j \to \infty} \langle \phi g_j^\rho \mid
\psi \rangle
 \\&= \langle \phi E \mid \psi \rangle ,
 \end{align*}
 so $v^\rho E = E$. Therefore, $E$ annihilates the image of $v^\rho -
\Id$, for every $v \in V$. Now, these images span a dense subspace of the
orthogonal complement $(\hilbert^V)^\perp$ of the subspace $\hilbert^V$
of elements of~$\hilbert$ that are \index{fixed point}fixed by every
element of~$V$. Hence, $E$ annihilates $(\hilbert^V)^\perp$.

Using $*$ to denote the \term[adjoint!of an operator]{adjoint}, we have
 $$ \langle \phi E^* \mid \psi \rangle
 =  \langle \phi \mid \psi E \rangle
 = \lim_{j \to \infty} \langle \phi \mid \psi g_j^\rho \rangle
 = \lim_{j \to \infty} \langle  \phi (g_j^{-1})^\rho \mid  \psi \rangle ,$$
 so the same argument, with $E^*$ in the place of~$E$ and $g_j^{-1}$ in the
place of~$g_j$, shows that $E^*$ annihilates $(\hilbert^U)^\perp$.

 Because $g^\rho$ is \term[unitary operator]{unitary}, it is
\term[normal!operator]{normal} (that is, commutes with its adjoint); thus,
the limit~$E$ is also normal: we have $E^* E = E E^*$. Therefore
 \begin{align*}
 \lVert \phi E \rVert^2
 &= \langle \phi E \mid \phi E \rangle
 = \langle \phi (E E^*) \mid \phi \rangle 
 \\&= \langle \phi (E^* E) \mid \phi \rangle
 = \langle \phi E^* \mid \phi E^* \rangle
 = \lVert \phi E^* \rVert^2 ,
 \end{align*}
 so $\ker E = \ker E^* $. 
 Hence \begin{align*}
 \ker E
 &= \ker E + \ker E^*
 \supset (\hilbert^{V})^\perp + (\hilbert^U)^\perp
 \\&= (\hilbert^{V} \cap \hilbert^U)^\perp
 = (\hilbert^{\langle U, V \rangle})^\perp
 . \end{align*}
 By passing to a subsequence (so $\{g_j\}$ is contained in a single
\term{Weyl chamber}), we may assume $\langle U, V \rangle = G$. Then
 $\hilbert^{\langle U, V \rangle} = \hilbert^G = 0$,
 so $\ker E
 \supset 0^\perp
 = \hilbert$.
 Hence, for all $\phi,\psi \in \hilbert$, we have
 $$ \lim \langle \phi g_j^\rho \mid \psi \rangle
 = \langle \phi E \mid \psi \rangle
 = \langle 0 \mid \psi \rangle
 = 0 ,$$
 as desired.

\begin{case}
 The general case.
 \end{case}
 From the \term{Cartan Decomposition} $G = KAK$, we may write $g_j = c'_j
a_j c_j$, with $c'_j,c_j \in K$ and $a_j \in A$. Because $K$ is compact,
we may assume, by passing to a subsequence, that $\{c'_j\}$ and $\{c_j\}$
converge: say, $c'_j \to c'$ and $c_j \to c$. Then
 \begin{align*}
 \lim_{j \to \infty} \langle \phi g_j^\rho \mid \psi \rangle
 &= \lim_{j \to \infty} \langle \phi  (c'_j a_j c_j)^\rho \mid \psi \rangle \\
 &= \lim_{j \to \infty} \langle  \phi (c'_j)^\rho a_j^\rho \mid \psi (c_j^{-1})^\rho
\rangle \\
 &= \lim_{j \to \infty} \bigl\langle  \phi (c')^\rho a_j^\rho 
 \mathrel{\big|}  \psi (c^{-1})^\rho \bigr\rangle \\
 &= 0 ,
 \end{align*}
 by Case~\ref{matcoeffs->0-A}.
 \end{proof}

\begin{rem} \label{matcoeffs->0(infty)}
 \ 
 \begin{enumerate}
 \item \label{matcoeffs->0(infty)-discrete}
 If
 \begin{itemize}
 \item $G$ and $\{g^t\}$ are as in Thm.~\ref{matcoeffs->0},
 \item $\Gamma$ is any \term[subgroup!discrete]{discrete subgroup} of~$G$,
that is \emph{not} a \term{lattice},
 and
 \item $\mu$ is the (infinite) $G$-\term[measure!invariant]{invariant}
measure on $\Gamma \backslash G$,
 \end{itemize}
 then the above proof (with $\hilbert = L^2(H \backslash G, \mu)$) shows
that
 $$ \lim_{t \to \infty} \int_{H \backslash G} \phi( x g^t) \, \psi(x) \,
d\mu = 0 ,$$
 for every $\phi,\psi \in L^2(\Gamma \backslash G, \mu)$.
 \item Furthermore, the \term[subgroup!discrete]{discrete
subgroup}~$\Gamma$ can be replaced with any closed subgroup~$H$ of~$G$,
such that $H \backslash G$ has a $G$-\term[measure!invariant]{invariant}
measure~$\mu$ that is finite on compact sets.
 \begin{itemize}
 \item If the measure of $H \backslash G$ is \term[finite volume]{finite},
then the conclusion is as in Thm.~\ref{matcoeffs->0}. 
 \item If the measure is infinite, then the conclusion is as
in~\pref{matcoeffs->0(infty)-discrete}.
 \end{itemize}
 \end{enumerate}
 \end{rem}

\begin{exercises}

\item \label{EssConst->Ergodic}
 Suppose $\varphi_t$ is a flow on a measure space $(X,\mu)$. Show that
$\varphi_t$ is \term[ergodic!flow]{ergodic} if and only if every
essentially \term[function!invariant!essentially]{invariant} measurable
function is \term[function!constant (essentially)]{essentially constant}.

\item \label{MautnerSL2PfEx}
 Prove Prop.~\ref{MautnerSL2} (without quoting other theorems of the
text).
 \hint{We have
 $$f(x u) = f(x u a^s) = f
\bigl( (x a^s) (a^{-s} u a^s) \bigr) \approx f(x a^s) = f(x) ,$$
 because $a^{-s} u a^s \approx e$.}

\item \label{GeodFlowErgEx}
 Derive Cor.~\ref{GeodFlowErg} from Prop.~\ref{MautnerSL2}.
 \hint{If $f$ is essentially
$a^s$-\term[function!invariant!essentially]{invariant}, then the Mautner
Phenomenon implies that it is also essentially $u^t$-invariant and
essentially $v^r$-invariant.}

\item Show that any \term{mixing} flow on $\Gamma \backslash G$ is
\term[ergodic!flow]{ergodic}.
 \hint{Let $A = B$ be a $g^t$-\term[set!invariant]{invariant} subset
of $\Gamma \backslash G$.}

\item Derive Rem.~\ref{MixingRem} from Thm.~\ref{matcoeffs->0}.

\item \label{MautnerSubgrpEx}
 Derive Thm.~\ref{matcoeffs->0} from Rem.~\ref{MixingRem}.
 \hint{Any $L^2$ function can be approximated by \index{function!step}step
functions.}

\item Suppose
 \begin{itemize}
 \item $G$ and $\{g^t\}$ are as in Thm.~\ref{matcoeffs->0}, 
 \item $\Gamma$ is any \term[subgroup!discrete]{discrete subgroup} of~$G$, 
 \item $\mu$~is the $G$-\term[measure!invariant]{invariant} measure on
$\Gamma \backslash G$,
 \item $\phi \in L^p(\Gamma \backslash G, \mu)$, for some $p < \infty$,
 and
 \item $\phi$ is essentially
$g^t$-\term[function!invariant!essentially]{invariant}.
 \end{itemize}
 Show that $\phi$ is essentially
$G$-\term[function!invariant!essentially]{invariant}.
 \hint{Some power of~$\phi$ is in~$L^2$. Use Thm.~\ref{matcoeffs->0} and
Rem.~\ref{matcoeffs->0(infty)}.}

\item \label{ParabUniqErg}
 Let
 \begin{itemize}
 \item $\Gamma$ be a \term{lattice} in $G =
\SL(2,\real)$\index{SL(2,R)*$\SL(2,\real)$},
 and
 \item $\mu$ be a probability measure on $\Gamma \backslash G$.
 \end{itemize}
 Show that if $\mu$ is \term[measure!invariant]{invariant} under both
$a^s$ and~$u^t$, then $\mu$ is the \term[measure!Haar]{Haar measure}.
 \hint{Let $\lambda$ be the Haar measure on $\Gamma \backslash G$, let
$U_\epsilon = \{\, u^t \mid 0 \le t \le \epsilon\}$, and define
$A_\epsilon$ and $V_\epsilon$ to be similar small intervals in $\{a^s\}$,
and $\{v^r\}$, respectively. If $f$ is continuous with compact
\term[support!of a function]{support}, then
 $$ \lim_{s \to \infty} \int_{y U_\epsilon A_\epsilon V_\epsilon} f( x
a^s) \, d\lambda(x) = 
 \lambda(y U_\epsilon A_\epsilon V_\epsilon) \int_{\Gamma \backslash G} f
\, d \lambda ,$$
 for all $y \in \Gamma \backslash G$ \seethm{matcoeffs->0}.
 Because $f$ is \term[function!continuous!uniformly]{uniformly continuous},
we see that 
 $$ \int_{y U_\epsilon A_\epsilon V_\epsilon} f( x a^s) \, d\lambda(x)
 = \int_{y U_\epsilon A_\epsilon V_\epsilon a^s} f \, d\lambda
 $$
 is approximately 
 $$ \frac{\lambda(y U_\epsilon A_\epsilon V_\epsilon a^s)}{e^{2s \epsilon}}
\int_0^{e^{2s} \epsilon} f(y u^t) \, dt .$$
 By choosing $y$ and $\{s_k\}$ such that $y a^{s_k} \to y$ and applying the
\term[Theorem!Pointwise Ergodic]{Pointwise Ergodic Theorem}, conclude that
$\lambda = \mu$.}

\end{exercises}

\section{Ergodic decomposition}
 \index{ergodic!decomposition}
 \label{ErgDecompSect}

Every \term[flow!measure-preserving]{measure-preserving flow} can be
decomposed into a union of \term[ergodic!flow]{ergodic} flows.

\begin{eg} \label{ErgDecOfT3Eg}
 Let 
 \begin{itemize}
 \item $v = (\alpha,1,0) \in \real^3$, for some
irrational~$\alpha$,
 \item $\varphi_t$ be the corresponding flow on $\torus^3 = \real^3 /
\integer^3$,
 and
 \item $\mu$ be the \term[measure!Lebesgue]{Lebesgue measure}
on~$\torus^3$.
 \end{itemize}
 Then $\varphi_t$ is \emph{not} \term[ergodic!flow]{ergodic}, because sets
of the form $A \times \torus^2$ are \term[set!invariant]{invariant}.

However, the flow decomposes into a union of \term[ergodic!flow]{ergodic}
flows: for each $z \in \torus$, let
 \begin{itemize}
 \item $T_z = \{z\} \times \torus^2$,
 and
 \item $\mu_z$ be the \term[measure!Lebesgue]{Lebesgue measure} on the
torus~$T_z$.
 \end{itemize}
 Then:
 \begin{enumerate}
 \item $\torus^3$ is the disjoint union $\bigcup_z T_z$,
 \item the restriction of~$\varphi_t$ to each subtorus~$T_z$ is
\term[ergodic!flow]{ergodic} (with respect to~$\mu_z$),
 and
 \item the measure~$\mu$ is the integral of the measures~$\mu_z$ (by
\term[Theorem!Fubini]{Fubini's Theorem}).
 \end{enumerate}
 \end{eg}

The following proposition shows that every measure~$\mu$ can be decomposed
into \term[ergodic!measure]{ergodic measures}. Each ergodic measure~$\mu_z$
is called an
 \index{component!ergodic|indsee{ergodic~component}}
 \defit[ergodic!component]{ergodic component} of~$\mu$.

\begin{prop} \label{Mu=IntegralErgodic}
 If $\mu$ is any $\varphi_t$-\term[measure!invariant]{invariant}
probability measure on~$X$, then there exist
 \begin{itemize}
 \item a measure~$\nu$ on a space~$Z$,
 and
 \item a {\rm(}measurable\/{\rm)} family $\{\mu_z\}_{z \in Z}$ of
\term[ergodic!measure]{ergodic measures} on~$X$,
 \end{itemize}
 such that $\mu = \int_Z \mu_z \,d\nu$; that is,
 $\int_X f \, d\mu = \int_Z \int_X f \, d\mu_z \, d\nu(z)$, 
 for every $f \in L^1(X,\mu)$.
 \end{prop}

\begin{proof}[Proof {\normalfont\itshape{\rm(}requires some
\term{Functional Analysis}\/{\rm)}}.]
 Let $\measures$ be the set of
$\varphi_t$-\term[measure!invariant]{invariant} probability measures
on~$X$. This is a \term{weak$^\ast$-compact}, \term[convex!set]{convex}
subset of the \term[dual (of a Banach space)]{dual} of a certain
\term{Banach space}, the continuous functions on~$X$ that \term{vanish
at~$\infty$}. So \term[Theorem!Choquet]{Choquet's Theorem} asserts that
any point in~$\measures$ is a \term[convex!combination]{convex
combination} of \term[point!extreme]{extreme points} of~$\measures$. That
is, if we let $Z$ be the set of extreme points, then there is a
probability measure~$\nu$ on~$Z$, such that $\mu = \int_Z z \, d\nu(z)$.
Simply letting $\mu_z = z$, and noting that the extreme points
of~$\measures$ are precisely the \term[ergodic!measure]{ergodic measures}
\seeex{Erg<>Extreme} yields the desired conclusion.
 \end{proof}

The above proposition yields a decomposition of the measure~$\mu$, but,
unlike Eg.~\ref{ErgDecOfT3Eg}, it does not provide a decomposition of the
space~$X$. However, any two \term[ergodic!measure]{ergodic measures} must
be mutually singular \seeex{Erg->Sing}, so a little more work yields the
following geometric version of the \term[ergodic!decomposition]{ergodic
decomposition}.  This often allows one to reduce a general question to the
case where the flow is \term[ergodic!flow]{ergodic}.

\begin{thm}[(Ergodic decomposition)] \label{ErgDecomp}
 If $\mu$ is a $\varphi_t$-\term[measure!invariant]{invariant}
probability measure on~$X$, then there exist
 \begin{itemize}
 \item a {\rm(}measurable\/{\rm)} family $\{\mu_z\}_{z \in Z}$ of
\term[ergodic!measure]{ergodic measures} on~$X$,
 \item a measure~$\nu$ on~$Z$,
 and
 \item a measurable function $\psi \colon X \to Z$,
 \end{itemize}
 such that
 \begin{enumerate}
 \item \label{ErgDecomp-integral}
 $\mu = \int_Z \mu_z \,d\nu$,
 and
 \item \label{ErgDecomp-support}
 $\mu_z$ is \term[support!of a measure]{supported} on $\psi^{-1}(z)$, for a.e.\ $z \in Z$.
 \end{enumerate}
 \end{thm}

\begin{proof}[Sketch of Proof.]
 Let $\Func \subset L^1(X,\mu)$ be the collection of $\{0,1\}$-valued
functions that are essentially
$\varphi_t$-\term[function!invariant!essentially]{invariant}. Because the
Banach space $L^1(X,\mu)$ is \term{separable}, we may choose a countable
dense subset $\Func_0 = \{\psi_n\}$ of~$\Func$. This defines a Borel
function $\psi \colon X \to \{0,1\}^\infty$.
 (By changing each of the functions in~$\Func_0$ on a set of measure~$0$,
we may assume $\psi$ is
$\varphi_t$-\term[function!invariant]{invariant}, not merely
essentially invariant.)
 Let $Z = \{0,1\}^\infty$ and $\nu = \psi_* \mu$.
 Proposition~\ref{FiberMeasures} below yields a (measurable) family
$\{\mu_z\}_{z \in Z}$ of probability measures on~$X$, such that
\pref{ErgDecomp-integral} and \pref{ErgDecomp-support} hold.

All that remains is to show that $\mu_z$ is \term[ergodic!measure]{ergodic}
for a.e.\ $z \in Z$. Thus, let us suppose that 
 $$Z_{\text{bad}} = \{\, z \in Z \mid \text{$\mu_z$ is not ergodic} \,\} $$
 is not a null set. For each $z \in Z_{\text{bad}}$, there is a
$\{0,1\}$-valued function $f_z \in L^1(X, \mu_z)$ that is essentially
$\varphi_t$-\term[function!invariant!essentially]{invariant}, but not
\term[function!constant (essentially)]{essentially constant}. The
functions~$f_z$ can be chosen to depend measurably on~$z$ (this is a
consequence of the \term[Theorem!Von Neumann Selection]{Von Neumann
Selection Theorem}); thus, there is a single measurable function~$f$
on~$Z$, such that 
 \begin{itemize}
 \item $f = f_z$ a.e.[$\mu_z$] for $z \in Z_{\text{bad}}$,
 and
 \item $f = 0$ on $Z \smallsetminus Z_{\text{bad}}$.
 \end{itemize}
 Because each $f_z$ is essentially $\varphi_t$-invariant, we know that $f$
is essentially $\varphi_t$-invariant; thus, $f \in \Func$. On the other
hand, $f$ is not \term[function!constant (essentially)]{essentially
constant} on the fibers of~$\psi$, so $f$ is not in the closure
of~$\Func$. This is a contradiction.
 \end{proof}

The above proof relies on the following very useful generalization of
\term[Theorem!Fubini]{Fubini's Theorem}.

\begin{prop} \label{FiberMeasures}
 Let
 \begin{itemize}
 \item $X$ and $Y$ be complete, \term{separable} metric spaces,
 \item $\mu$ and $\nu$ be probability measures on $X$ and~$Y$,
respectively,
 and 
 \item $\psi \colon X \to Y$ be a
\term[map!measure-preserving]{measure-preserving map} Borel map.
 \end{itemize}
 Then there is a Borel map $\lambda \colon Y \to
\Prob(X)$, such that
 \begin{enumerate}
 \item $\mu = \int_Y \lambda_y \, d\nu(y)$,
 and
 \item $\lambda_y \bigl( \psi^{-1}(y) \bigr) = 1$, for all $y \in Y$.
 \end{enumerate}
 Furthermore, $\lambda$ is unique {\rm(}up to measure zero{\rm)}.
 \end{prop}

\begin{exercises}

\item \label{Erg<>Extreme}
 In the notation of the proof of Prop.~\ref{Mu=IntegralErgodic}, show that
a point~$\mu$ of~$\measures$ is \term[ergodic!measure]{ergodic} if and only
if it is an extreme point of~$\measures$. (A point~$\mu$ of~$\measures$ is
an \defit[point!extreme]{extreme point} if it is \emph{not} a 
\term[convex!combination]{convex combination} of two other points
of~$\measures$; that is, if there do not exist $\mu_1,\mu_2 \in
\measures$, and $t \in (0,1)$, such that $\mu = t\mu_1 + (1-t)\mu_2$ and
$\mu_1 \neq \mu_2$.)

\item \label{Erg->Sing}
 Suppose $\mu_1$ and~$\mu_2$ are \term[ergodic!measure]{ergodic},
$\varphi_t$-\term[measure!invariant]{invariant} probability measures
on~$X$. Show that if $\mu_1 \neq \mu_2$, then there exist subsets
$\Omega_1$ and~$\Omega_2$ of~$X$, such that, for $i,j \in \{1,2\}$, we
have 
 $$\mu_i(\Omega_j) = 
 \begin{cases}
 1 & \mbox{if $i = j$,} \\
 0 & \mbox{if $i \neq j$}.
 \end{cases}
 $$

\item \label{CtblFibers}
  Let
 \begin{itemize}
 \item $X$ and $X'$ be complete, \term{separable} metric spaces,
 \item $\mu$ and~$\mu'$ be probability measures on $X$ and~$X'$,
respectively,
 \item $\varphi_t$ and $\varphi'_t$ be \term[ergodic!flow]{ergodic},
\term[flow!measure-preserving]{measure-preserving flows} on $X$ and~$X'$,
respectively,
 \item $\psi \colon X \to Y$ be a
\term[map!measure-preserving]{measure-preserving map},
\term[equivariant]{equivariant Borel map},
 \index{map!equivariant|indsee{equivariant}}
 and
 \item $\Omega$ be a conull subset of~$X$, such that $\psi^{-1}(y) \cap
\Omega$ is countable, for a.e.\ $y \in Y$.
 \end{itemize}
 Show there is a conull subset $\Omega'$ of~$X$, such that $\psi^{-1}(y)
\cap \Omega$ is finite, for a.e.\ $y \in Y$.
 \hint{The function $f(x) = \lambda_{\psi(x)} \bigl( \{x\} \bigr)$ is
essentially $\varphi_t$-\term[function!invariant!essentially]{invariant},
so it must be essentially constant. A probability measure with all atoms of
the same weight must have only finitely many atoms.} 

 \end{exercises}

\section{Averaging sets} \label{AvgSetSect}

The proof of \term[Ratner's Theorems]{Ratner's Theorem} uses a version of the
\term[Theorem!Pointwise Ergodic]{Pointwise Ergodic Theorem} that applies to
\term[unipotent!subgroup]{(unipotent) groups} that are not just
one-dimensional\index{dimension!of a Lie group}. The classical version
\pref{PtwiseErgThmFlow} asserts that averaging a function over larger and
larger intervals of almost any orbit will converge to the integral of the
function. Note that the average is over intervals, not over arbitrary
large subsets of the orbit. In the setting of higher-dimensional groups,
we will average over ``\term[averaging!set]{averaging sets}."

\begin{defn}
 Suppose
 \begin{itemize}
 \item $U$ is a \index{subgroup!connected}connected,
\term[unipotent!subgroup]{unipotent subgroup} of~$G$,
 \item $a$ is a \term[hyperbolic!element]{hyperbolic element} of~$G$ that
\term[normalizer]{normalizes}~$U$
 \item $a^{-n} u a^n \to e$ as $n \to -\infty$ (note that this is
$-\infty$, not~$\infty$!),
 and 
 \item $E$ is a ball in~$U$ (or, more generally, $E$ is any bounded,
non-null, Borel subset of~$U$).
 \end{itemize}
 Then:
 \begin{enumerate}
 \item we say that $a$ is an \defit[automorphism!expanding]{expanding
automorphism} of~$U$,
 \item for each $n \ge 0$, we call $E_n = a^{-n} E a^n$ an
\defit[averaging!set]{averaging set},
 \nindex{$E_n$ = averaging set $a^{-n} E a^n$}
 and
 \item we call $\{E_n\}_{n=0}^\infty$ an
\defit[averaging!sequence]{averaging sequence}.
 \end{enumerate}
 \end{defn}

\begin{rem} \label{AvgSetRem} \ 
 \begin{enumerate}
 \item  \label{AvgSetRem-expand}
 By assumption, conjugating by $a^n$ contracts~$U$
when $n < 0$. Conversely, conjugating by $a^n$ expands~$U$ when $n > 0$.
Thus, $E_1$, $E_2$,$\dots$ are larger and larger subsets of~$U$. (This
justifies calling~$a$ an ``expanding" automorphism.)
 \item Typically, one takes $E$ to be a nice set (perhaps a ball) that
contains~$e$, with $E \subset a^{-1} E a$. In this case,
$\{E_n\}_{n=0}^\infty$ is an increasing \index{Folner sequence*F\o lner
sequence}{F\o lner
sequence} \seeex{FolnerEx}, but, for technical reasons, we will employ a
more general choice of~$E$ at one point in our argument (namely, in
\ref{BasLemPf}, the proof of Prop.~\ref{BasLem}).
 \end{enumerate}
 \end{rem}

\begin{thm}[{(\term[Theorem!Pointwise Ergodic]{Pointwise Ergodic
Theorem})}] \label{PtwiseErg}
 If
 \begin{itemize}
 \item $U$ is a \index{subgroup!connected}connected,
\term[unipotent!subgroup]{unipotent subgroup} of~$G$,
 \item $a$ is an \defit[automorphism!expanding]{expanding
automorphism} of~$U$,
 \item $\nu_U$ is the \term[measure!Haar]{Haar measure} on~$U$,
 \item $\mu$ is an \term[ergodic!measure]{ergodic}
$U$-\term[measure!invariant]{invariant} probability measure on\/ $\Gamma
\backslash G$,
 and
 \item $f$ is a continuous function on\/ $\Gamma \backslash G$ with compact
\term[support!of a function]{support},
 \end{itemize}
 then there exists a $U$-\term[set!invariant]{invariant}
subset~$\Omega$ of\/~$\Gamma \backslash G$ with $\mu(\Omega) = 1$, such
that 
  $$  \frac{1}{\nu_U(E_n)}
 \int_{E_n} f(x u) \, d \nu_U(u)
 \to \int_{\Gamma \backslash G} f(y) \, d \mu(y)
 \mbox{\quad as $n \to \infty$}
 .
 $$
 for every $x \in \Omega$
 and every \term[averaging!sequence]{averaging sequence} $\{E_n\}$ in~$U$.
 \end{thm}

To overcome some technical difficulties, we will also use the following
uniform approximate version \seeex{UnifErgEx}. It is ``uniform," because
the same number~$N$ works for all points $x \in \Omega_\epsilon$, and
the same set $\Omega_\epsilon$ works for all functions~$f$.

\begin{cor}[{(\term[Theorem!Pointwise Ergodic!Uniform]{Uniform Pointwise
Ergodic Theorem})}] \label{UnifErg}
 If
 \begin{itemize}
 \item $U$ is a \index{subgroup!connected}connected,
\term[unipotent!subgroup]{unipotent subgroup} of~$G$,
 \item $a$~is an \defit[automorphism!expanding]{expanding
automorphism} of~$U$,
 \item $\nu_U$ is the \term[measure!Haar]{Haar measure} on~$U$,
 \item $\mu$ is an \term[ergodic!measure]{ergodic}
$U$-\term[measure!invariant]{invariant} probability measure on $\Gamma
\backslash G$,
 and
 \item $\epsilon > 0$,
 \end{itemize}
 then
 there exists a subset~$\Omega_\epsilon$ of\/~$\Gamma
\backslash G$ with $\mu(\Omega_\epsilon) > 1 - \epsilon$, such that for
 \begin{itemize}
 \item every continuous function~$f$ on\/ $\Gamma \backslash G$
with compact \term[support!of a function]{support},
 \item every \term[averaging!sequence]{averaging sequence} $\{E_n\}$
in~$U$,
 and
 \item every $\delta > 0$,
 \end{itemize}
 there is some $N \in \natural$, such that
 $$ \left| \frac{1}{\nu_U(E_n)}
 \int_{E_n} f(x u) \, d \nu_U(u)
 - \int_{\Gamma \backslash G} f(y) \, d \mu(y) \right|
 < \delta, $$
 for all $x \in \Omega_\epsilon$ and all $n \ge N$.
 \end{cor}

\begin{rem} \label{FolnerRems} \ 
 \begin{enumerate}
 \item \label{FolnerRems-amen}
 A \term{Lie group}~$G$ said to be \defit[Lie group!amenable]{amenable} if it
has a \index{Folner sequence*F\o lner sequence}{F\o lner sequence}. 
 \item \label{FolnerRems-solv+cpct}
 It is known that a \emph{\index{subgroup!connected}connected} \term{Lie
group}~$G$ is \term[Lie group!amenable]{amenable} if and only if there are
closed, \index{subgroup!connected}connected, \index{subgroup!normal}normal
subgroups~$U$ and~$R$ of~$G$, such that 
 \begin{itemize}
 \item $U$ is \term[unipotent!subgroup]{unipotent},
 \item $U \subset R$, 
 \item $R/U$ is \term[group!abelian]{abelian},
 and 
 \item $G/R$ is \index{Lie group!compact}compact.
 \end{itemize}
 \item \label{FolnerRems-Elon}
 There are examples to show that not every \index{Folner sequence*F\o lner
sequence}{F\o lner sequence} $\{E_n\}$ can be used as an
\term[averaging!sequence]{averaging sequence}, but it is always the case
that some subsequence of $\{E_n\}$ can be used as the averaging sequence
for a \term[Theorem!Pointwise Ergodic]{pointwise ergodic theorem}. 
 \end{enumerate}
 \end{rem}

\begin{exercises}

\item \label{FolnerEx}
 Suppose
 \begin{itemize}
 \item $U$ is a \index{subgroup!connected}connected,
\term[unipotent!subgroup]{unipotent subgroup} of~$G$,
 \item $a$ is an \defit[automorphism!expanding]{expanding
automorphism} of~$U$, 
 \item $\nu_U$ is the \term[measure!Haar]{Haar measure} on~$U$,
 and
 \item $E$ is a precompact, open subset of~$U$, such that $a^{-1} E a
\subset E$.
 \end{itemize}
 Show that the \term[averaging!sequence]{averaging sequence} $E_n$ is an
increasing \index{Folner sequence*F\o lner sequence}\emph{F\o lner
sequence}; that is,  
 \begin{enumerate}
 \item for each nonempty compact subset~$C$ of~$U$, we have $\nu_U
\bigl( (C E_n) \symmdiff E_n \bigr) /\nu_U(E_n) \to 0$ as $n \to \infty$,
 and
 \item $E_n \subset E_{n+1}$, for each~$n$.
 \end{enumerate}

\item \label{Folner->InvtMeas}
 Show that if~$G$ is \term[Lie group!amenable]{amenable}, then there is an
\term[measure!invariant|)]{invariant} probability measure for any action
of~$G$ on a compact metric space. More precisely, suppose 
 \begin{itemize}
 \item $\{E_n\}$ is a \index{Folner sequence*F\o lner sequence}{F\o lner
sequence} in a \term{Lie group}~$G$,
 \item $X$ is a compact metric space,
 and
 \item $G$ acts continuously on~$X$.
 \end{itemize}
 Show there is a $G$-\term[measure!invariant]{invariant} probability
measure on~$X$.
 \hint{\term[measure!Haar]{Haar measure} restricts to a measure $\nu_n$
on~$E_n$. \term[push-forward (of a measure)]{Pushing} this to~$X$ (and
normalizing) yields a probability measure~$\mu_n$ on~$X$.
 Any \term{weak$^\ast$-limit} of~$\{\mu_n\}$ is
$G$-\term[measure!invariant]{invariant}.}

\item \label{UnifErgEx}
 Derive Cor.~\ref{UnifErg} from Thm.~\ref{PtwiseErg}.

 \end{exercises}

\begin{notes}

A few of the many introductory books on \term[ergodic!theory]{Ergodic
Theory} are \cite{Halmos, KatokHasselblatt-intro, Walters}.

\notesect{PtwiseErgSect}

This material is standard.

The  \term[Theorem!Pointwise Ergodic]{Pointwise Ergodic Theorem} is due to
G.~D.~Birkhoff \cite{Birkhoff-ErgThm}. There are now many different proofs,
such as \cite{Kamae, KatznelsonWeiss-SimplePtwise}. (See also
\cite[Thm.~I.2.5, p.~17]{BekkaMayer}). The hints for
Exers.~\ref{PtwiseErgZMaxl} and \ref{PtwiseErgZ} are adapted from
\cite[pp.~19--24]{FriedmanIntro}. 

Exercise~\ref{Int(f-f)=0Ex} is \cite[Prop.~2.2]{LedrappierStrlcyn}.

A solution to Exer.~\ref{UnifErg<>UnifConv} appears in \cite[Thm.~I.3.8,
p.~33]{BekkaMayer}.

\notesect{MautnerSect}

The \term[Theorem!Moore Ergodicity]{Moore Ergodicity Theorem}
\pref{MooreErg} was first proved by C.~C.~Moore \cite{Moore-ErgThm}.
Later, he \cite{Moore-Mautner} extended this to a very general version of
the \term{Mautner Phenomenon} \pref{MautnerSL2}.

\term[mixing]{Mixing} is a standard topic (see, e.g.,
\cite{KatokHasselblatt-intro, Walters} and \cite[pp.~21--28]{BekkaMayer}.)
Our proof of Thm.~\ref{matcoeffs->0} is taken from \cite{EllisNerurkar}.
Proofs can also be found in \cite[Chap.~3]{BekkaMayer}, \cite[\S
II.3]{MargulisBook} and \cite[Chap.~2]{ZimmerBook}.

A solution to Exer.~\ref{EssConst->Ergodic} appears in \cite[Thm.~I.1.3,
p.~3]{BekkaMayer}.

The hint to Exer.~\ref{ParabUniqErg} is adapted from \cite[Lem.~5.2,
p.~31]{MargulisThesis}.

\notesect{ErgDecompSect}

This material is standard.

A complete proof of Prop.~\ref{Mu=IntegralErgodic} from
\term[Theorem!Choquet]{Choquet's Theorem} appears in
\cite[\S12]{Phelps-Choquet}.

 See \cite[\S8]{Oxtoby} for a brief history (and proof) of the
\term[ergodic!decomposition]{ergodic decomposition} \pref{ErgDecomp}.

 Proposition~\ref{FiberMeasures} appears in \cite[\S3]{Rohlin}.

Exer.~\ref{Erg<>Extreme} is solved in \cite[Prop.~3.1, p.~30]{BekkaMayer}.

\notesect{AvgSetSect}

For any \index{Lie group!amenable}amenable Lie group, a theorem of
A.~Tempelman \cite{Tempelman-GenErg}, generalized by W.~R.~Emerson
\cite{Emerson}, states that certain F\o lner sequences can be used as
\term[averaging!sequence]{averaging sequences} in a
\term[Theorem!Pointwise Ergodic]{pointwise ergodic theorem}. (A proof also
appears in \cite[Cor.~6.3.2, p.~218]{TempelmanBook}.) The
\term[Theorem!Pointwise Ergodic!Uniform]{Uniform Pointwise Ergodic
Theorem} \pref{UnifErg} is deduced from this in \cite[\S7.2 and
\S7.3]{MargulisTomanov}.

The book of Greenleaf \cite{Greenleaf-Amenable} is the classic source for
information on \term[Lie group!amenable]{amenable groups}.

The converse of Exer.~\ref{Folner->InvtMeas} is true
\cite[Thm.~3.6.2]{Greenleaf-Amenable}. Indeed, the existence of
\term[measure!invariant]{invariant} measures is often taken as the
definition of \term[Lie group!amenable]{amenability}. See
\cite[\S4.1]{ZimmerBook} for a discussion of
\term[Lie group!amenable]{amenable} groups from this point of view, including
the characterization mentioned in Rem.~\fullref{FolnerRems}{solv+cpct}. 

Remark~\fullref{FolnerRems}{Elon} is a theorem of E.~Lindenstrauss 
\cite{Lindenstrauss-PtwiseErg}.

\end{notes}

%% file: RatnerAlgGrps.tex
\mychapter{Facts about Algebraic Groups} \label{AlgGrpsChap}

\index{function!polynomial|(}

In the theory of \term[Lie group]{Lie groups}, all
\term[homomorphism!of Lie groups]{homomorphisms} (and other maps) are
generally assumed to be $C^\infty$ functions \seeSect{LieGrpSect}. The
theory of \term[algebraic group!theory of]{algebraic groups} describes the
conclusions that can be obtained from the stronger assumption that the
maps are \term[homomorphism!polynomial]{polynomial} functions (or, at
least, \index{function!rational}rational functions). Because the
polynomial nature of \index{unipotent!flow}unipotent flows plays
such an important role in the arguments of Chapter~\ref{IntroChap} (see,
for example, Prop.~\ref{PolyDivSL2}), it is natural to expect that a good
understanding of polynomials will be essential at some points in the more
complete proof presented in Chapter~\ref{ProofChap}. However, the reader
may wish to skip over this chapter, and refer back when necessary.

\section{Algebraic groups} \label{AlgicGrpsSect}

\begin{defn} \label{AlgicGrpDefn}
 \ 
 \begin{itemize}
 \item We use 
 \nindex{$\real[x_{1,1}, \ldots, x_{\ell,\ell}]$ = real polynomials in $\{
x_{i,j} \}$}
 $\real[x_{1,1}, \ldots, x_{\ell,\ell}]$ to
denote the set of real polynomials in the
$\ell^2$ variables $\{\, x_{i,j} \mid 1 \le i,j \le \ell\, \}$. 
 \item For any $Q \in \real[x_{1,1}, \ldots, x_{\ell,\ell}]$, and any $n
\times n$ matrix~$g$, we use $Q(g)$ to denote the value obtained by
substituting the matrix entries~$g_{i,j}$ into the variables~$x_{i,j}$. For
example:
 \nindex{$g_{i,j}$ = entries of matrix~$g$}
 \begin{itemize}
 \item If $Q= x_{1,1} + x_{2,2} + \cdots + x_{\ell,\ell}$, then $Q(g)$ is
the \index{trace (of a matrix)}trace of~$g$.
 \item If $Q= x_{1,1} x_{2,2} - x_{1,2} x_{2,1}$, then
$Q(g)$ is the \index{determinant}determinant of the first principal $2
\times 2$ minor of~$g$.
 \end{itemize}
 \item For any subset~$\mathcal{Q}$ of
$\real[x_{1,1}, \ldots, x_{\ell,\ell}]$, let
 \nindex{$\mathcal{Q}$ = some subset of $\real[x_{1,1}, \ldots,
x_{\ell,\ell}]$}
 \nindex{$\Var(\mathcal{Q})$ = variety associated to set $\mathcal{Q}$ of
polynomials}
  $$\Var(\mathcal{Q}) = \{\, g \in \SL(\ell,\real) \mid Q(g) = 0, \
\forall Q \in \mathcal{Q} \,\} .$$
 \index{SL(l,R)*$\SL(\ell,\real)$|(}
 This is the \defit{variety} associated to~$\mathcal{Q}$.
 \item A subset~$H$ of $\SL(\ell,\real)$ is \defit[Zariski!closed]{Zariski
closed}\index{Zariski!closed|(} if there is a subset~$\mathcal{Q}$ of
$\real[x_{1,1}, \ldots, x_{\ell,\ell}]$, such that $H =
\Var(\mathcal{Q})$. (In the special case where $H$ is a sub\emphit{group}
of $\SL(\ell,\real)$, we may also say that $H$ is a \defit[algebraic
group!real]{real algebraic group}\index{algebraic group!real|(} or an
\defit[algebraic group!over~$\real$]{algebraic group that is defined
over~$\real$}.)
 \end{itemize}
 \end{defn}

\begin{eg} \label{EgsOfAlgGrps}
 Each of the following is a real \index{algebraic group!real}algebraic
group \seeex{EgsOfAlgGrpsExer}:
 \begin{enumerate}

 \item \label{EgsOfAlgGrps-SL(l,R)}
 $\SL(\ell,\real)$.

%\item 
% The unipotent group $\{u^t\} = \left\{ \begin{bmatrix} 1 & 0 \\ t & 1
%\end{bmatrix} \right\}$.

\item \label{EgsOfAlgGrps-D}
  The group 
 \nindex{$\mathbb{D}_\ell$ = $\{ \text{diagonal matrices in
$\SL(\ell,\real)$} \}$}
 $$\mathbb{D}_\ell  = \begin{bmatrix}
 * &  & \vbox to 0pt{\vss\hbox to 0 pt{\hss\Huge $0$}\vss\vss} \\
 & \ddots \\
 \vbox to 0pt{\vss\hbox to 0 pt{\Huge $0$\hss}}& & * \\
 \end{bmatrix}
 \subset \SL(\ell,\real) $$
 of diagonal matrices in $\SL(\ell,\real)$.

\item  \label{EgsOfAlgGrps-U}
 The group
 \nindex{$\mathbb{U}_\ell$ = $\{\text{unipotent lower-triangular matrices
in $\SL(\ell,\real)$}\}$}
 $$ \mathbb{U}_\ell = \begin{bmatrix}
 1 &  & \vbox to 0pt{\vss\hbox to 0 pt{\hss\Huge $0$}\vss\vss} \\
 & \ddots \\
 \vbox to 0pt{\vss\hbox to 0 pt{\Huge $*$\hss}}& & 1 \\
 \end{bmatrix}
 \subset \SL(\ell,\real) $$
 of \index{lower-triangular matrices!unipotent}lower-triangular matrices
with $1$'s on the diagonal.
 \index{unipotent!subgroup}
% In other words, $\mathbb{U}_\ell$ is the subset of $\SL(\ell,\real)$
%consisting of all matrices~$u$, such that
% $$ u_{ij} =
% \begin{cases}
% 1 & \mbox{if $i = j$} \\
% 0 & \mbox{if $i < j$} .
% \end{cases} $$

\item  \label{EgsOfAlgGrps-tri}
 The group
 $$ \mathbb{D}_\ell \mathbb{U}_\ell = \begin{bmatrix}
 * &  & \vbox to 0pt{\vss\hbox to 0 pt{\hss\Huge $0$}\vss\vss} \\
 & \ddots \\
 \vbox to 0pt{\vss\hbox to 0 pt{\Huge $*$\hss}}& & * \\
 \end{bmatrix}
 \subset \SL(\ell,\real) $$
 of \index{lower-triangular matrices}lower-triangular matrices in
$\SL(\ell,\real)$.

 \item \label{EgsOfAlgGrps-SL(n,R)}
 The copy of $\SL(n,\real)$ in the top left corner of $\SL(\ell,\real)$
(if $n < \ell$). 

 \item \label{EgsOfAlgGrps-Stab(v)}
 The \defit[stabilizer!of a vector]{stabilizer}
 \par\medskip\noindent \text{\hskip1cm}%
 \nindex{$\Stab_{\SL(\ell,\real)}(v)$ = stabilizer of vector~$v$}%
 $\displaystyle \Stab_{\SL(\ell,\real)}(v) = \{\, g \in \SL(\ell,\real)
\mid v g = v \,\}$
 \par\medskip\noindent
 of any vector $v \in \real^\ell$.

 \item \label{EgsOfAlgGrps-Stab(V)}
 The \defit[stabilizer!of a subspace]{stabilizer}
 \par\medskip\noindent \text{\hskip1cm}%
 \nindex{$\Stab_{\SL(\ell,\real)}(V)$ = stabilizer of subspace~$V$}%
  $ \displaystyle \Stab_{\SL(\ell,\real)}(V) = \{\, g \in \SL(\ell,\real)
\mid \forall v \in V, \ v g \in V \,\}$
 \par\medskip\noindent
 of any linear subspace $V$ of~$\real^\ell$.

 \item \label{EgsOfAlgGrps-SO(Q)}
 The special orthogonal group 
 \index{SO(Q)*$\SO(Q)$}$\SO(Q)$ of a quadratic form~$Q$ on~$\real^\ell$
\seeDefn{QuadFormDefn}.

 \end{enumerate}
 \end{eg}

It is important to realize that most closed subsets of $\SL(\ell,\real)$
are \emph{not} Zariski closed. In particular, the following important
theorem tells us that an infinite, \index{set!discrete}discrete subset can
never be Zariski closed. (It is a generalization of the fact that any
nontrivial polynomial function on~$\real$ has
only finitely many zeroes.) We omit the proof.

\begin{thm}[(Whitney)] \label{Zar->AlmConn}
 Any Zariski closed subset of\/ $\SL(\ell,\real)$ has only
finitely many \index{component!connected}components {\rm(}with respect to
the usual topology of\/ $\SL(\ell,\real)$ as a \term{Lie group}{\rm)}.
 \end{thm}

\begin{eg}
 From Thm.~\ref{Zar->AlmConn}, we know that the
\index{subgroup!discrete}discrete group $\SL(\ell,\integer)$ is \emph{not}
Zariski closed. In fact, we will see that $\SL(\ell,\integer)$ is not
contained in any Zariski closed, proper subgroup of $\SL(\ell,\real)$
\seeex{SLZarDensinSLREx}.
 \end{eg}
 \index{group!discrete|indsee{subgroup,~discrete}}

\begin{rem} \label{ZarSingularSet}
 Zariski closed sets need not be \index{submanifold}submanifolds of
$\SL(\ell,\real)$. This follows from Exer.~\ref{UnionZarClosed}, for
example, because the union of two submanifolds that intersect is usually
not a submanifold --- the intersection is a \term{singularity}. 

 Exercise~\ref{ZarDim} defines the \defit[dimension!of a Zariski closed
set]{dimension} of any Zariski closed set~$Z$.
  Although we do not prove this, it can be shown that (if $Z$ is
nonempty), there is a unique smallest Zariski closed subset~$S$ of~$Z$,
such that 
 \begin{itemize}
 \item $\dim S < \dim Z$,
 \item $Z \smallsetminus S$ is a $C^\infty$ \index{submanifold}submanifold
of $\SL(\ell,\real)$,
 and
 \item $\dim Z$ (as defined below) is equal to the
\index{dimension!of a manifold}dimension of $Z \smallsetminus S$ as a
\index{manifold}manifold.
 \end{itemize}
 The set~$S$ is the \emph{singular set} of~$Z$. From the uniqueness of~$S$,
it follows that any Zariski closed sub\emph{group} of $\SL(\ell,\real)$ is
a $C^\infty$ \index{submanifold}submanifold of $\SL(\ell,\real)$
\seeex{AlgGrp->LieGrp}; 
 \end{rem}

\begin{exercises}

\item \label{Zar->Closed}
 Show that every Zariski closed subset of $\SL(\ell,\real)$ is closed
(in the usual topology of $\SL(\ell,\real)$ as a \term{Lie group}).

\item \label{EgsOfAlgGrpsExer}
 Verify that each of the groups in Eg.~\ref{EgsOfAlgGrps} is Zariski
closed.
 \hint{\pref{EgsOfAlgGrps-SL(l,R)}~Let $\mathcal{Q} = \emptyset$.
 \pref{EgsOfAlgGrps-D}~Let $\mathcal{Q} = \{\, x_{i,j} \mid i \neq j
\,\}$.
 \pref{EgsOfAlgGrps-Stab(v)}~Let 
 $\mathcal{Q} = \{\, v_1 x_{1,j} + \cdots + v_\ell x_{\ell,j} - v_j \mid 1
\le j \le \ell \,\}$,
 where $v = (v_1,\ldots,v_\ell)$.%
% \pref{EgsOfAlgGrps-U}~Let $\mathcal{Q} = \{\, x_{i,i} - 1 \,\} \cup \{\,
%x_{i,j} \mid i < j \,\}$.
% \pref{EgsOfAlgGrps-SO(l)}~Let
% $\mathcal{Q} = \bigset{ \delta_i^j - \sum_{k=1}^\ell x_{i,k} x_{j,k} }{1
%\le i \le j \le \ell}$.
% \pref{EgsOfAlgGrps-SL(n,R)}~Let $\mathcal{Q} = \{\, x_{i,j} - \delta_i^j
%\mid \max\{i,j\} > n \,\}$.
 }

\item Show that if $Z$ is a \term[Zariski!closed]{Zariski closed} subset
of $\SL(\ell,\real)$, and $g \in \SL(\ell,\real)$, then $Z g$ is Zariski
closed.

\item \label{UnionZarClosed}
 Suppose $Z_1$ and~$Z_2$ are \term[Zariski!closed]{Zariski closed}
subsets of $\SL(\ell,\real)$. Show that the union $Z_1 \cup Z_2$ is
\term[Zariski!closed]{Zariski closed}.

\item \label{AlgGrp->LieGrp}
 Show that if $G$ is a Zariski closed subgroup of $\SL(\ell,\real)$, then
$G$ is a $C^\infty$ \index{submanifold}submanifold of $\SL(\ell,\real)$
(so $G$ is a \term{Lie group}).
 \hint{Uniqueness of the \term{singular set}~$S$ \see{ZarSingularSet}
implies $Sg = S$ for all $g \in G$, so $S = \emptyset$.}

 \medskip
 \noindent {\itshape The remaining exercises present some\/ {\rm(}more
technical\/{\rm)} information about Zariski closed sets, including the
notion of \index{dimension!of a Zariski closed set}dimension.}
 \medskip

\item For any subset $Z$ of $\SL(\ell,\real)$, let 
 $\ideal(Z)$ be the collection of polynomials that vanish on~$Z$; that is,%
 \nindex{$\ideal(Z)$ = ideal of polynomials vanishing on~$Z$}
 $$ \ideal(Z) = \{\, Q \in \real[x_{1,1}, \ldots, x_{\ell,\ell}] \mid
\forall z \in Z, \ Q(z) = 0 \,\} .$$
 \begin{enumerate}
 \item Show $Z$ is Zariski closed if and only if \index{variety}$Z = \Var
\bigl( \ideal(Z) \bigr)$.
 \item Show that $\ideal(Z)$ is an ideal; that is,
 \begin{enumerate}
 \item $0 \in \ideal(Z)$,
 \item for all $Q_1,Q_2 \in \ideal(Z)$, we have $Q_1 + Q_2 \in \ideal(Z)$,
 and
 \item for all $Q_1 \in \ideal(Z)$ and $Q_2 \in \real[x_{1,1}, \ldots,
x_{\ell,\ell}]$, we have $Q_1 Q_2 \in \ideal(Z)$.
 \end{enumerate}
 \end{enumerate}

\item \label{Noetherian}
 Recall that a ring~$R$ is \defit[Noetherian ring]{Noetherian} if it has
the \term[chain condition!ascending]{ascending chain condition} on ideas;
this means that if $I_1 \subset I_2 \subset \cdots$ is any increasing
chain of ideals, then we have $I_n = I_{n+1} = \cdots$ for some~$n$.
 \begin{enumerate}
 \item Show that a commutative ring~$R$ is Noetherian if and only if all of
its ideals are \term[finitely generated ideal]{finitely generated}; that
is, for each ideal~$I$ of~$R$, there is a \emph{finite} subset~$F$
of~$I$, such that $I$ is the smallest ideal of~$R$ that contains~$F$.
 \item \label{Noetherian-R[x]}
 Show that $\real[x_{1,1}, \ldots, x_{\ell,\ell}]$ is Noetherian.
 \item \label{Noetherian-FiniteVariety}
 Show that if $Z$ is a \term[Zariski!closed]{Zariski closed} subset of
$\SL(\ell,\real)$, then there is a \emph{finite} subset~$Q$ of
$\real[x_{1,1}, \ldots, x_{\ell,\ell}]$, such that $Z =
\Var(Q)$\index{variety}.
 \item Prove that the collection of \term[Zariski!closed]{Zariski closed}
subsets of $\SL(\ell,\real)$ has the \term[chain
condition!descending]{descending chain condition}: if $Z_1 \supset Z_2
\supset \cdots$ is a decreasing chain of Zariski closed sets, then we have
$Z_n = Z_{n+1} = \cdots$ for some~$n$.
 \end{enumerate}
 \hint{\pref{Noetherian-R[x]}~Show that if $R$ is Noetherian, then the
polynomial ring $R[x]$ is Noetherian: If $I$ is an ideal in $R[x]$, let
 $$I_n = \{\, r \in R \mid \exists Q \in R[x], \ \mbox{$r x^n + Q \in I$
and $\deg Q < n$}\,\} .$$
 Then $I_n \subset I_{n+1} \subset \cdots$ is an increasing chain of
ideals.}

\item \label{ZarIrred}
  A \term[Zariski!closed]{Zariski closed} subset of $\SL(\ell,\real)$ is
\defit[irreducible (Zariski closed set)]{irreducible} if it is \emph{not}
the union of two Zariski closed \emph{proper} subsets.

 Let $Z$ be a \term[Zariski!closed]{Zariski closed} subset of
$\SL(\ell,\real)$.
 \begin{enumerate}
 \item \label{ZarIrred-Union}
 Show that $Z$ is the union of finitely many irreducible
\term[Zariski!closed]{Zariski closed} subsets.
 \item \label{ZarIrred-Comp}
 An \defit[component!irreducible]{irreducible component} of~$Z$ is an
irreducible Zari\-ski closed subset of~$Z$ that is not not properly
contained in any irreducible Zariski closed subset of~$Z$.
 \begin{enumerate}
 \item Show that $Z$ is the union of its irreducible components.
 \item Show that $Z$ has only finitely many irreducible components.
 \end{enumerate}
 \end{enumerate}
 \hint{\pref{ZarIrred-Union}~Proof by contradiction: use the \term[chain
condition!descending]{descending chain condition}.
 \pref{ZarIrred-Comp}~Use~\pref{ZarIrred-Union}.}

\item Suppose $G$ is a Zariski closed subgroup of $\SL(\ell,\real)$. 
 \begin{enumerate}
 \item Show that the irreducible \index{component!irreducible}components
of~$G$ are disjoint.
 \item Show that the irreducible components of~$G$ are cosets of a Zariski
closed subgroup of~$G$.
 \end{enumerate}

 \item \label{ZarDim}
 The \defit[dimension!of a Zariski closed set]{dimension} of a Zariski
closed set~$Z$ is the largest~$r$, such that there is a chain $Z_0 \subset
Z_1 \subset \cdots \subset Z_r$ of nonempty, irreducible Zariski closed
subsets of~$Z$. 

It can be shown (and you may assume) that $\dim Z$ is the
largest~$r$, for which there is a linear map $T \colon \Mat_{\ell\times
\ell}(\real) \to \real^r$, such that $T(Z)$ contains a nonempty open subset
of~$\real^r$.
 \begin{enumerate}
 \item Show $\dim Z = 0$ if and only if $Z$ is finite and nonempty.
 \item Show $\dim Z_1 \le \dim Z_2$ if $Z_1 \subset Z_2$.
 \item Show $\dim (Z_1 \cup Z_2) = \max\{ \dim Z_1, \dim Z_2 \}$ if $Z_1$
and~$Z_2$ are Zariski closed.
 \item Show $\dim \SL(\ell,\real) = \ell^2 - 1$.
 \item Show that the collection of \term[irreducible (Zariski closed
set)]{irreducible} Zariski closed subsets of $\SL(\ell,\real)$ has the
\term[chain condition!ascending]{ascending chain condition}: if $Z_1
\subset Z_2 \subset \cdots$ is an increasing chain of irreducible Zariski
closed sets, then we have $Z_n = Z_{n+1} = \cdots$ for some~$n$.
  \end{enumerate}

\item \label{SubvarSmallerDim}
 Suppose $V$ and~$W$ are Zariski closed sets in $\SL(\ell,\real)$. Show
that if
 \begin{itemize}
 \item $V \subset W$,
 \item $W$ is irreducible,
 and
 \item $\dim V = \dim W$,
 \end{itemize}
 then $V = W$.

 \end{exercises}

\section{Zariski closure} \label{ZarClosureSect}
 \index{Zariski!closure|(}

 \begin{defn}
 The \defit[Zariski!closure]{Zariski closure} of a subset~$H$ of
$\SL(\ell,\real)$ is the (unique) smallest Zariski closed subset of
$\SL(\ell,\real)$ that contains~$H$ \seeex{ZarClosUniq}. 
 We use \nindex{$\Zar{H}$ = Zariski closure of~$H$}{$\Zar{H}$} to denote
the Zariski closure of~$H$. 
 \end{defn}

\begin{rem} \ 
 \begin{enumerate}
 \item Obviously, $H$ is Zariski closed if and only if $\Zar{H} = H$. 
 \item One can show that if $H$ is a subgroup of $\SL(\ell,\real)$, then
$\Zar{H}$ is also a subgroup of $\SL(\ell,\real)$ \seeex{ZarIsGrp}.
 \end{enumerate}
 \end{rem} 

Every Zariski closed subgroup of $\SL(\ell,\real)$ is closed
\seeex{Zar->Closed} and has only finitely many connected
\index{component!connected}components \see{Zar->AlmConn}. The converse is
false:

\begin{eg} \label{DiagNotZar}
 Let
 $$ A = \bigset{ \begin{bmatrix} t & 0 \\ 0 & 1/t \end{bmatrix} }{ t \in
\real^+ } \subset \SL(2,\real) .$$
 \index{SL(2,R)*$\SL(2,\real)$}
 Then
 \begin{enumerate}
 \item $A$  is closed,
 \item $A$ is \index{subgroup!connected}connected (so it has only one
connected \index{component!connected}component),
 and
 \item \label{DiagNotZar-Zar}
 $\Zar{A} = \bigset{ \begin{bmatrix} t & 0 \\ 0 & 1/t \end{bmatrix}
}{ t \in \real \smallsetminus \{0\} }$
 \seeex{DiagNotZar-ZarEx}.
 \end{enumerate}
 So $\Zar{A} \iso \real \smallsetminus \{0\}$ has two connected
\index{component!connected}components. Since $\Zar{A} \neq A$, we know
that $A$ is not Zariski closed.
 \end{eg}

Although $A$ is not exactly equal to $\Zar{A}$ in Eg.~\ref{DiagNotZar},
there is very little difference: $A$ has finite index in~$\Zar{A}$. For
most purposes, a finite group can be ignored, so we make the following
definition.

 \begin{defn}
 A subgroup~$H$ of $\SL(\ell,\real)$ is
\defit[Zariski!closed!almost]{almost Zariski closed} if $H$ is a
finite-index subgroup of~$\Zar{H}$. 
 \index{Zariski!closed!almost|(}
 \end{defn}

\begin{rem}
 Any finite-index subgroup of a \term{Lie group} is closed
\seeex{FinInd->Closed}, so any subgroup of $\SL(\ell,\real)$ that is
\term[Zariski!closed!almost]{almost Zariski closed} must be closed.
 \end{rem}

The reader may find it helpful to have some alternative characterizations
\seeex{AlmZarAlternateEx}:

\begin{rem} \label{AlmZarAlternate} \ 
 \begin{enumerate}
 \item A \index{subgroup!connected}connected subgroup~$H$  of
$\SL(\ell,\real)$ is \term[Zariski!closed!almost]{almost Zariski closed}
if and only if it is the identity \index{component!identity}component of a
subgroup that is Zariski closed.
 \item A subgroup~$H$ of $\SL(\ell,\real)$ is
\term[Zariski!closed!almost]{almost Zariski closed} if and only if it is
the union of (finitely many) \index{component!connected}components of a
Zariski closed group.
 \item Suppose $H$ has only finitely many connected
\index{component!connected}components. Then $H$ is
\term[Zariski!closed!almost]{almost Zariski closed} if and only if its
identity \index{component!identity}component $H^\circ$ is almost Zariski
closed.
 \item Suppose $H$ is a Lie subgroup of $\SL(\ell,\real)$. Then $H$ is
\term[Zariski!closed!almost]{almost Zariski closed} if and only if $\dim
\Zar{H} = \dim H$\index{dimension!of a Lie group}.
 \end{enumerate}
 \end{rem}

Note that if $H$ is \term[Zariski!closed!almost]{almost Zariski closed},
then it is closed, and has only finitely many connected
\index{component!connected}components. 
 Here are two examples to show that the converse is false. (Both examples
are closed and \index{subgroup!connected}connected.)
Corollary~\ref{GTAlmZar} below implies that all examples of this phenomenon
must be based on similar constructions.

\begin{eg} \  
 \begin{enumerate}
 \item For any \index{irrational!number}irrational number~$\alpha$, let 
 $$T = \bigset{
 \begin{bmatrix}
 t^\alpha & 0 & 0 \\
 0 & t & 0 \\
 0 & 0 & 1/t^{\alpha + 1}
 \end{bmatrix}
 }{ t \in \real^+}
 \subset \SL(3,\real) . $$
 \index{SL(3,R)*$\SL(3,\real)$}
 Then
 $$\Zar{T} = \bigset{
 \begin{bmatrix}
 s & 0 & 0 \\
 0 & t & 0 \\
 0 & 0 & 1/(st)
 \end{bmatrix}
 }{ s,t \in \real \smallsetminus \{0\}} . $$
 Since $\dim T = 1 \neq 2 =
\dim \Zar{T}$, we conclude that $T$ is not almost Zariski closed.

The calculation of~$\Zar{T}$ follows easily from Cor.~\ref{SplitTori}
below. Intuitively, the idea is simply that, for elements~$g$ of~$T$, the
relation between $g_{1,1}$ and $g_{2,2}$ is \term{transcendental}, not
algebraic, so it cannot be captured by a polynomial. Thus, as far as
polynomials are concerned, there is no relation at all between $g_{1,1}$
and $g_{2,2}$ --- they can vary independently of one another. This
independence is reflected in the Zariski closure.

 \item Let 
 $$ H = \bigset{
  \begin{bmatrix}
 e^t & 0 & 0 & 0 \\
 0 & e^{-t} & 0 & 0 \\
 0 & 0 & 1 & t \\
 0 & 0 & 0 & 1
 \end{bmatrix}
 }{ t \in \real}
 \subset \SL(4,\real) . $$
 \index{SL(4,R)*$\SL(4,\real)$}
 Then 
 $$ \Zar{H} = \bigset{
  \begin{bmatrix}
 e^s & 0 & 0 & 0 \\
 0 & e^{-s} & 0 & 0 \\
 0 & 0 & 1 & t \\
 0 & 0 & 0 & 1
 \end{bmatrix}
 }{ s, t \in \real}
 \subset \SL(4,\real) . $$
 \index{SL(4,R)*$\SL(4,\real)$}
 Since $\dim H = 1 \neq 2 = \dim \Zar{H}$, we conclude that $H$ is not
almost Zariski closed.

Formally, the fact that $H$ is not almost-Zariski closed follows from
Thm.~\ref{SplitTori} below. Intuitively, the \term{transcendental}
relation between $g_{1,1}$ and $g_{3,4}$ is lost in the Zariski closure.
 \end{enumerate}

 \end{eg}

\begin{exercises}

\item \label{ZarClosUniq}
 For each subset~$H$ of $\SL(\ell,\real)$, show there is a unique
Zariski closed subset $\Zar{H}$ of $\SL(\ell,\real)$ containing~$H$, such
that if $C$ is any Zariski closed subset $\Zar{H}$ of $\SL(\ell,\real)$
that contains~$H$, then $\Zar{H} \subset C$.
 \hint{Any intersection of Zariski closed sets is Zariski closed.}

\item \label{CentClosed}
 Show that if $Z$ is any subset of an \index{algebraic group!real}algebraic
group~$G$, then the \term{centralizer} $C_G(Z)$ is Zariski closed.

\item \label{DiagNotZar-ZarEx} Verify \fullref{DiagNotZar}{Zar}. 
 \hint{Let $\mathcal{Q} = \{\, x_{1,2}, x_{2,1} \,\}$. 
 If $Q(x_{1,1},x_{2.2})$ is a polynomial, such that $Q(t,1/t) = 0$ for all
$t > 0$, then $Q(t,1/t) = 0$ for all $t \in \real$.}

\item Show that if $H$ is a \index{subgroup!connected}\emph{connected}
subgroup of $\SL(\ell,\real)$, then $\Zar{H}$ is \term[irreducible
(Zariski closed set)]{irreducible}.

\item \label{FinInd->Closed}
 (\emphit{Requires some Lie theory})
 Suppose $H$ is a finite-index subgroup of a \term{Lie group}~$G$. Show
that $H$ is an open subgroup of~$G$. (So $H$ is closed.)
 \hint{There exists $n \in \integer^+$, such that $g^n \in H$ for all $g
\in G$. Therefore $\exp(x) = \exp \bigl( (1/n) x \bigr)^n \in H$ for every
element~$x$ of the \index{Lie algebra}Lie algebra of~$G$.}

\item \label{AlmZarAlternateEx}
 Verify each part of Rem.~\ref{AlmZarAlternate}.

 \end{exercises}

\section{Real Jordan decomposition} \label{JordanDecSect}
 \index{Jordan!decomposition!real|(}

The real Jordan decomposition writes any matrix as a combination of
matrices of three basic types.

\begin{defn}
 Let $g \in \SL(\ell,\real)$.

 \begin{itemize}

 \item $g$ is \defit[unipotent!matrix]{unipotent} if $1$ is the
only \index{eigenvalue}eigenvalue of~$g$ (over~$\complex$); in other words,
$(g-1)^\ell = 0$ \see{UnipMatDefn}.
 \index{unipotent!element|(}

 \item $g$ is \defit[hyperbolic!element]{hyperbolic} (or
\defit[R-split*$\real$-split!element]{$\real$-split}) if it is
\index{diagonalizable!over R*over~$\real$}diagonalizable over~$\real$, and
all of its \index{eigenvalue}eigenvalues are positive; that is, if $h^{-1}
g h$ is a diagonal matrix with no negative entries, for some $h \in
\SL(\ell,\real)$.
 \index{hyperbolic!element|(}

\item $g$ is \defit[elliptic element]{elliptic} if it is
\index{diagonalizable!over C*over~$\complex$}diagonalizable
over~$\complex$, and all of its \index{eigenvalue}eigenvalues are of
absolute value~$1$.
 \index{elliptic element|(}

 \end{itemize}
 \end{defn}

\begin{eg}
 For all $t \in \real$:
 \begin{enumerate}
 \item $\begin{bmatrix} 1 & 0 \\ t & 1 \end{bmatrix}$ is unipotent,
 \item $\begin{bmatrix} e^t & 0 \\ 0 & e^{-t} \end{bmatrix}$ is hyperbolic,
 \item $\begin{bmatrix} \cos t & \sin t \\ -\sin t & \cos t \end{bmatrix}$
is elliptic \seeex{RotElliptic}.
 \end{enumerate}
 See Exer.~\ref{SL2JordanTypes} for an easy way to tell whether an element
of $\SL(2,\real)$\index{SL(2,R)*$\SL(2,\real)$} is unipotent, hyperbolic,
or elliptic.
 \end{eg}

\begin{prop}[{(\term[Jordan!decomposition!real]{Real Jordan
decomposition})}] \label{RealJordanDecomp}
 For any $g \in \SL(\ell,\real)$, there exist unique $g_u, g_h,g_e \in
\SL(\ell,\real)$, such that
 \begin{enumerate}
 \item $g = g_u g_h g_e$,
 \item $g_u$ is unipotent,\
 \item $g_h$ is hyperbolic,
 \item $g_e$ is elliptic,
  and
 \item $g_u$, $g_h$, and $g_e$ all commute with each other.
 \end{enumerate}
 \end{prop}
 \nindex{$g_u$, $g_h$, $g_e$ = Jordan components of matrix~$g$}

\begin{proof}
 (Existence)
 The usual \term[Jordan!decomposition]{Jordan decomposition} of Linear
Algebra (also known as ``\term[Jordan!Canonical Form]{Jordan Canonical
Form}") implies there exist $h \in \SL(\ell,\complex)$, a
\index{matrix!nilpotent}nilpotent matrix~$N$, and a diagonal matrix~$D$,
such that $h^{-1} g h = N + D$, and $N$~commutes with~$D$.  This is an
additive decomposition. By factoring out~$D$, we obtain a multiplicative
decomposition:
 $$h^{-1} g h = (ND^{-1} + \Id) D = u D ,$$
 where $u =  ND^{-1} + \Id$ is unipotent (because $u - \Id = ND^{-1}$ is
nilpotent, since $N$ commutes with $D^{-1}$).

 Now, because any complex number~$z$ has a (unique)
\term{polar form} $z = r e^{i \theta}$, we may write $D = D_h D_e$, where
$D_h$~is hyperbolic, $D_e$~is elliptic, and both matrices are diagonal, so
they commute with each other (and, from the structure of the
\index{Jordan!Canonical Form}Jordan Canonical Form, they both commute
with~$N$). Conjugating by~$h^{-1}$, we obtain
 $$ g = h (u D_h D_e) h^{-1} = g_u g_h g_e ,$$
 where $g_u = h u h^{-1}$, $g_h = h D_h h^{-1}$, and $g_e = h D_e h^{-1}$.
This is the desired decomposition.

\smallskip
 (Uniqueness)
 The uniqueness of the decomposition is, perhaps, not so interesting to the
reader, so we relegate it to the exercises
\seeexs{JorDecCommute}{JorDecUnique}. Uniqueness is, however, often of
vital importance. For example, it can be used to address a technical
difficulty that was ignored in the above proof: from our construction, it
appears that the matrices $g_u$, $g_h$, and $g_e$ may have complex
entries, not real. However, using an overline to denote complex
conjugation, we have $\overline{g} = \overline{g_u} \, \overline{g_h} \,
\overline{g_e}$. Since $\overline{g} = g$, the uniqueness of the
decomposition implies $\overline{g_u} = g_u$, $\overline{g_h} = g_h$, and
$\overline{g_e} = g_e$. Therefore, $g_u, g_h, g_e \in \SL(\ell,\real)$, as
desired.
 \end{proof}

The \term{uniqueness} of the \term[Jordan!decomposition]{Jordan
decomposition} implies, for $g,h \in \SL(\ell,\real)$, that if $g$
commutes with~$h$, then the Jordan
 \index{Jordan!component|indsee{component,~Jordan}}
 \index{component!Jordan}components
$g_u$, $g_h$, and~$g_e$ commute with~$h$ (see also
Exer.~\ref{JorDecCommute}). In other words, if the \term{centralizer}
$C_{\SL(\ell,\real)}(h)$ contains~$g$, then it must also contain the
Jordan \index{component!Jordan}components of~$g$. Because the
centralizer is Zariski closed \seeex{CentClosed}, this is a special case
of the following important result.

\begin{thm} \label{JorDecinG}
 If 
 \begin{itemize}
 \item $G$ is a Zariski closed subgroup of $\SL(\ell,\real)$,
 and
 \item $g \in G$,
 \end{itemize}
 then $g_u,g_h,g_e \in G$.
 \end{thm}

We postpone the proof to \S\ref{ChevalleySect}.

As mentioned at the start of the chapter, we should assume that
\term[homomorphism!polynomial]{homomorphisms} are polynomial functions.
(But some other types of functions will be allowed to be more general
rational functions, which are not defined when the denominator is~$0$.)

\begin{defn} \label{RegularFuncDefn}
 Let $H$ be a subset of $\SL(\ell,\real)$.
 \begin{enumerate}
 \item \label{RegularFuncDefn-poly}
 A function $\phi \colon H \to \real$ is a
\defit[function!polynomial]{polynomial} (or is
\defit[function!regular]{regular}) if there exists $Q \in \real[x_{1,1},
\ldots, x_{\ell,\ell}]$, such that $\phi(h) = Q(h)$ for all $h \in H$.
 \item A real-valued function $\psi$ defined a subset of~$H$ is
\defit[function!rational]{rational}
 if there exist polynomials $\phi_1,\phi_2
 \colon H \to \real$, such that
 \begin{enumerate}
 \item the domain of~$\psi$ is $\{\, h \in H \mid \phi_2(h) \neq 0 \,\}$,
 and
 \item $\psi(h) = \phi_1(h)/\phi_2(h)$ for all~$h$ in the domain
of~$\psi$.
 \end{enumerate}
 \item A function $\phi \colon H \to \SL(n,\real)$ is a
\defit[function!polynomial]{polynomial} if, for each $1 \le i,j \le n$,
the \term[matrix!entry]{matrix entry} $\phi(h)_{i,j}$ is a polynomial
function of $h \in H$. Similarly, $\psi$ is
\defit[function!rational]{rational} if each $\psi(h)_{i,j}$ is a rational
function of $h \in H$.
 \end{enumerate}
 \end{defn}

We now show that any 
\term[homomorphism!polynomial]{polynomial homomorphism} respects the real
Jordan decomposition; that is, $\rho(g_u) = \rho(g)_u$, $\rho(g_h) =
\rho(g)_h$, and $\rho(g_e) = \rho(g)_e$.

\begin{cor} \label{rho(unip/hyp)}
 Suppose 
 \begin{itemize}
 \item $G$ is a real \index{algebraic group!real}algebraic group,
 and
 \item $\rho \colon G \to \SL(m,\real)$ is a
\term[homomorphism!polynomial]{polynomial homomorphism}.
 \end{itemize}
 Then:
 \begin{enumerate}
 \item \label{rho(unip/hyp)-unip}
 If $u$ is any unipotent element of~$G$, then $\rho(u)$ is a
unipotent element of\/ $\SL(m,\real)$.
 \item If $a$ is any hyperbolic element of~$G$, then $\rho(a)$ is a
hyperbolic element of\/ $\SL(m,\real)$.
 \item If $k$ is any elliptic element of~$G$, then $\rho(k)$ is an
elliptic element of\/ $\SL(m,\real)$.
 \end{enumerate}
 \end{cor}

\begin{proof}
 Note that the \index{graph}graph of~$\rho$ is a Zariski closed subgroup of
$G \times H$ \seeex{GrfZarClosed}. 

 We prove only \pref{rho(unip/hyp)-unip}; the others are similar. Since
$u$ is unipotent, we have $u_u = u$, $u_h = e$, and $u_e = e$. Therefore,
the real Jordan decomposition of $\bigl( u, \rho(u) \bigr)$ is
 $$ \bigl( u, \rho(u) \bigr) =  \bigl( u, \rho(u)_u \bigr) \bigl( e,
\rho(u)_h \bigr) \bigl( e, \rho(u)_e \bigr) .$$
 Since $\bigl( u, \rho(u) \bigr) \in \graph \rho$, Thm.~\ref{JorDecinG}
implies 
 $$ \bigl( u, \rho(u)_u \bigr) = \bigl( u, \rho(u) \bigr)_u \in \graph
\rho .$$
 Let $y = \rho(u)_u$. Since $(u,y) \in \graph \rho$\index{graph}, we have
$\rho(u) = y$. Hence $\rho(u) = \rho(u)_u$ is unipotent.
 \end{proof}

\begin{exercises}

\item \label{UlUnip}
 Show that every element of $\mathbb{U}_\ell$ is unipotent.

\item \label{RotElliptic}
 Show $\begin{bmatrix} \cos t & \sin t \\ -\sin t & \cos t \end{bmatrix}$
is an elliptic element of $\SL(2,\real)$\index{SL(2,R)*$\SL(2,\real)$},
for every $t \in \real$.

\item \label{SL2JordanTypes}
 Let $g \in \SL(2,\real)$\index{SL(2,R)*$\SL(2,\real)$}. Recall that
 \index{trace (of a matrix)}\nindex{$\trace g$ = trace of the
matrix~$g$}{$\trace g$}
 is the sum of the diagonal entries of~$g$. Show:
 \begin{enumerate}
 \item $g$ is unipotent if and only if $\trace g = 2$.
 \item $g$ is hyperbolic if and only if $\trace g > 2$.
 \item $g$ is elliptic if and only if $-2 < \trace g < 2$.
 \item $g$ is neither unipotent, hyperbolic, nor elliptic if and
only if $\trace g \le -2$.
 \end{enumerate}

\item \label{CommuteJordan}
 Suppose $g$ and~$h$ are elements of $\SL(\ell,\real)$, such that $gh =
hg$. Show:
 \begin{enumerate}
 \item If $g$ and~$h$ are unipotent, then $gh$ is unipotent.
 \item If $g$ and~$h$ are hyperbolic, then $gh$ is hyperbolic.
 \item If $g$ and~$h$ are elliptic, then $gh$ is elliptic.
 \end{enumerate}

\item \label{JorDecCommute}
 Suppose $g, g_u, g_h, g_e \in
\SL(\ell,\complex)$\index{SL(l,C)*$\SL(\ell,\complex)$}, and these matrices
are as described in the conclusion of Prop.~\ref{RealJordanDecomp}.
 Show (without using the \term[Jordan!decomposition]{Jordan
decomposition} or any of its properties) that if $x \in
\SL(\ell,\complex)$, and $x$~commutes with~$g$, then $x$ also commutes
with each of $g_u$, $g_h$, and~$g_e$.
 \hint{Passing to a conjugate, assume $g_h$ and~$g_e$ are diagonal. We have
 $g_h^{-n} x g_h^n = (g_u g_e)^n x (g_u g_e)^{-n}$. Since each matrix
entry of the LHS is an \index{exponential!function}exponential function
of~$n$, but each matrix entry on the RHS grows at most polynomially, we
see that the LHS must be constant. So $x$ commutes with~$g_h$. Then
$g_u^{-n} x g_u^n = g_e^n x g_e^{-n}$. Since a bounded polynomial must be
constant, we see that $x$ commutes with $g_u$ and~$g_e$.}

\item \label{JorDecUnique} 
 Show that the real Jordan decomposition is unique.
 \index{diagonalizable!over C*over~$\complex$}
 \hint{If $g = g_u g_h g_e = g_u' g_h' g_e'$, then $g_u^{-1} g_u' = g_h
g_e (g_h' g_e')^{-1}$ is both unipotent and diagonalizable
over~$\complex$ (this requires Exer.~\ref{JorDecCommute}). Therefore $g_u
= g_u'$. Similarly, $g_h = g_h'$ and $g_e = g_e'$.}

\item \label{JordanEigenvec}
 Suppose $g \in \SL(\ell,\real)$, $v \in \real^\ell$, and $v$~is an
\index{eigenvector}eigenvector for~$g$. Show that $v$ is also an
eigenvector for $g_u$, $g_h$, and~$g_e$.
 \hint{Let $W$ be the \index{eigenspace}eigenspace corresponding to the
\index{eigenvalue}eigenvalue~$\lambda$ associated to~$v$. Because $g_u$,
$g_h$, and~$g_e$ commute with~$g$, they preserve~$W$. The
\term[Jordan!decomposition]{Jordan decomposition} of $g|_W$, the
restriction of~$g$ to~$W$, is
 $(g|_W)_u (g|_W)_h (g|_W)_e$.}

\item \label{SimulDiag}
 Show that any commuting set of \index{diagonalizable!over
R*over~$\real$}diagonalizable matrices can be diagonalized simultaneously.
More precisely, suppose
 \begin{itemize}
 \item $S \subset \SL(\ell,\real)$,
 \item each $s \in S$ is hyperbolic,
 and
 \item the elements of~$S$ all commute with each other.
 \end{itemize}
 Show there exists $h \in \SL(\ell,\real)$, such that every element of
$h^{-1} S h$ is diagonal. 

\item \label{MultIsPoly}
 Suppose $G$ is an subgroup of $\SL(\ell,\real)$ that is almost Zariski
closed.
 \begin{enumerate}
 \item For $i(g) = g^{-1}$, show that $i$ is a polynomial function
from~$G$ to~$G$.
 \item For $m(g,h) = g h$, show that $m$ is a polynomial function
from~$G \times G$ to~$G$. (Note that $G \times G$ can naturally be
realized as a subgroup of $\SL(2\ell,\real)$ that is almost Zariski
closed.)
 \end{enumerate}
 \hint{Cramer's Rule provides a polynomial formula for the inverse of a
matrix of \index{determinant}determinant one. The usual formula for the
product of two matrices is a polynomial.}

\item \label{InvImgAlmZar}
 Show that if 
 \begin{itemize}
 \item $f \colon \SL(\ell,\real) \to \SL(m,\real)$ is a polynomial,
 and
 \item $H$ is a Zariski closed subgroup of $\SL(m,\real)$,
 \end{itemize}
 then $f^{-1}(H)$ is \term[Zariski!closed]{Zariski closed}.

\item \label{ZarIsGrp}
 Show that if $H$ is any subgroup of $\SL(\ell,\real)$, then $\Zar{H}$ is
also a subgroup of $\SL(\ell,\real)$.
 \hint{Exercises \ref{MultIsPoly} and \ref{InvImgAlmZar}.}

\item  \label{norm-closed}
 Show that if $H$ is a \index{subgroup!connected|(}connected Lie subgroup
of $\SL(\ell,\real)$, then the \term{normalizer} $N_{\SL(\ell,\real)}(H)$
is \term[Zariski!closed]{Zariski closed}.
 \hint{The \index{homomorphism!of Lie groups}homomorphism $\Ad \colon
\SL(\ell,\real) \to \SL \bigl( \LieSL(\ell,\real)
\bigr)$ is a polynomial.}
 \index{SL(V)*$\SL(V)$}

\item 
 Show that if $G$ is any \index{subgroup!connected}connected subgroup of
$\SL(\ell,\real)$, then $G$ is a \term[subgroup!normal]{\emph{normal}}
subgroup of~$\Zar{G}$.
 \index{normal!subgroup|indsee{subgroup,~normal}}

\item \label{GxHZarClosed}
 There is a natural embedding of $\SL(\ell,\real) \times
\SL(m,\real)$ in $\SL(\ell+m,\real)$. Show that if $G$ and~$H$ are Zariski
closed subgroups of $\SL(\ell,\real)$ and $\SL(m,\real)$, respectively,
then $G \times H$ is Zariski closed in $\SL(\ell+m,\real)$.

\item \label{GrfZarClosed}
 Suppose $G$ is a Zariski closed subgroup of $\SL(\ell,\real)$, and $\rho
\colon G \to \SL(m,\real)$ is a \term[homomorphism!polynomial]{polynomial
homomorphism}. There is a natural embedding of the \index{graph}graph
of~$\rho$ in $\SL(\ell + m, \real)$ \cfex{GxHZarClosed}. Show that the
graph of $\rho$ is Zariski closed.

 \end{exercises}
 \index{Jordan!decomposition!real|)}
  \index{elliptic element|)}

\section{Structure of almost-Zariski closed groups}
\label{StructureAlgGrpSect}

The main result of this section is that any \index{algebraic
group!real}algebraic group can be decomposed into subgroups of three basic
types:  unipotent, torus, and \index{Lie group!semisimple}semisimple
\seethm{AlgGrpTLU}.

\begin{defn} \label{UnipTorusSSDefn} \
 \begin{itemize}

\item  A subgroup~$U$ of $\SL(\ell,\real)$ is
\defit[unipotent!group]{unipotent} if and only if it is conjugate to a
subgroup of~$\mathbb{U}_\ell$.

\item A subgroup~$T$ of $\SL(\ell,\real)$ is a
\defit[torus!algebraic]{torus} if 
 \begin{itemize}
 \item $T$ is conjugate (over~$\complex$) to a group of diagonal matrices;
that is, $h^{-1} T h$ consists entirely of diagonal matrices, for some $h
\in \SL(\ell,\complex)$\index{SL(l,C)*$\SL(\ell,\complex)$}), 
 \item $T$ is \index{subgroup!connected}connected,
 and 
 \item $T$ is almost Zariski closed.
 \end{itemize}
  (We have required tori to be \index{subgroup!connected}connected, but
this requirement should be relaxed slightly; any subgroup of $\Zar{T}$
that contains~$T$ may also be called a torus.)

 \item A closed subgroup~$L$ of $\SL(\ell,\real)$ is
\defit[Lie group!semisimple]{semisimple} if its identity
\index{component!identity}component $L^\circ$ has no nontrivial,
\index{subgroup!connected}connected, \index{group!abelian}abelian,
\index{subgroup!normal}normal subgroups.

 \end{itemize}
 \end{defn}

\begin{rem} \label{Unips+Tori}
 Here are alternative characterizations of unipotent groups and tori:

\begin{enumerate}

 \item \label{Unips+Tori-unip} (Engel's Theorem)
 A subgroup~$U$ of $\SL(\ell,\real)$ is
\defit[unipotent!group]{unipotent} if and only if every element of~$U$ is
\term[unipotent!matrix]{unipotent} \seeex{Unip->Lower}.

\item \label{Unips+Tori-tori}
 A \index{subgroup!connected}connected subgroup~$T$ of $\SL(\ell,\real)$
is a \index{torus!algebraic}torus if and only if
 \begin{itemize}
 \item $T$ is \index{group!abelian}abelian,
 \item each individual element of~$T$ is
 \index{diagonalizable!over C*over~$\complex$}diagonalizable
(over~$\complex$),
 and
 \item $T$ is almost Zariski closed
 \end{itemize}
 \seeex{SimulDiag}.
 
 \end{enumerate}
 \end{rem}

Unipotent groups and tori are fairly elementary, but the
\term[Lie group!semisimple]{semisimple groups} are more difficult to
understand. The following fundamental theorem of Lie theory reduces
their study to simple groups (which justifies their name).

\begin{defn}
 A group~$G$ is \defit[group!simple (or almost)]{almost simple} if it has
no infinite, proper, \index{subgroup!normal}normal subgroups.
 \end{defn}

\begin{thm} \label{SemiSimp=Simples}
 Let $L$ be a \index{subgroup!connected}connected, \index{Lie
group!semisimple}semisimple subgroup of\/ $\SL(\ell,\real)$. Then, for
some~$n$, there are closed, \index{subgroup!connected}connected subgroups
$S_1,\ldots,S_n$ of~$L$, such that
 \begin{enumerate}
 \item each $S_i$ is almost \index{group!simple (or almost)}simple,
 and
 \item $L$ is isomorphic to {\rm(}a finite cover of\/{\rm)} $S_1 \times
\cdots \times S_n$.
 \end{enumerate}
 \end{thm}

The almost-simple \index{group!simple (or almost)}groups have been
classified by using the theory known as ``\term{roots and weights}." We
merely provide some typical examples, without proof.

\begin{eg} \label{SimpleGrpEgs} \ 
 \begin{enumerate}
 \item \label{SimpleGrpEgs-SL}
 $\SL(\ell,\real)$ is almost \index{group!simple (or almost)}simple (if
$\ell \ge 2$).
 \item If $Q$ is a quadratic form on~$\real^\ell$ that is \term[quadratic
form!nondegenerate]{nondegenerate} \seeDefn{QuadFormDefn}, and $\ell \ge
3$, then \index{SO(Q)*$\SO(Q)$}$\SO(Q)$ is \index{Lie
group!semisimple}semisimple (and it is almost simple if, in addition, $n
\neq 4$). (For $\ell = 2$, the groups $\SO(2)$ and $\SO(1,1)$ are tori,
not semisimple \seeex{SO2abel}.)
 \end{enumerate}
 \end{eg}

From the above almost-simple groups, it is easy to construct numerous
\index{Lie group!semisimple}semisimple groups. One example
is\index{SO(Q)*$\SO(Q)$}\index{SL(l,R)*$\SL(\ell,\real)$}
 $$\SL(3,\real) \times \SL(7,\real) \times \SO(6) \times \SO(4,7) .$$

The following structure theorem is one of the major results in the theory
of \index{algebraic group!real}algebraic groups. 

\begin{defn}
 Recall that a \term{Lie group} $G$ is a \defit{semidirect product} of
closed subgroups $A$ and~$B$ (denoted \nindex{$A \ltimes B$ = semidirect
product}{$G = A \ltimes B$}) if
 \begin{enumerate}
 \item $G = AB$,
 \item $B$ is a \index{subgroup!normal}normal subgroup of~$G$,
 and
 \item $A \cap B = \{e\}$.
 \end{enumerate}
 (In this case, the map $(a,b) \mapsto ab$ is a
\index{diffeomorphism}diffeomorphism from $A \times B$ onto~$G$. However,
it is not a group isomorphism (or even a
\index{homomorphism!of Lie groups}homomorphism) unless every element of~$A$
commutes with every element of~$B$.)
 \end{defn}

\begin{thm} \label{AlgGrpTLU}
 Let $G$ be a \index{subgroup!connected}connected subgroup of\/
$\SL(\ell,\real)$ that is almost Zariski closed. Then there exist:
 \begin{itemize}
 \item a \index{Lie group!semisimple}semisimple subgroup $L$ of~$G$,
 \item a \index{torus!algebraic}torus $T$ in~$G$,
 and
 \item a \index{unipotent!subgroup}unipotent subgroup~$U$ of~$G$,
 \end{itemize}
 such that 
 \begin{enumerate}
 \item $G = (L T) \ltimes U$,
 \item $L$, $T$, and $U$ are almost Zariski closed,
 and
 \item $L$ and $T$ \index{centralizer}centralize each other, and have
finite intersection.
 \end{enumerate}
 \end{thm}

\begin{proof}[{Sketch of proof \rm (\emphit{requires some Lie theory}).}]
 Let $R$ be the \term[radical of a Lie group]{radical} of~$G$, and let $L$
be a \term[subgroup!Levi]{Levi subgroup} of~$G$; thus, $R$~is
\term[group!solvable]{solvable}, $L$~is \term[Lie
group!semisimple]{semisimple}, $LR = G$, and $L \cap R$ is
\index{subgroup!discrete}discrete \see{G=LRad}. From the
\term[Theorem!Lie-Kolchin]{Lie-Kolchin Theorem} \pref{LieKolchin}, we know
that $R$ is conjugate (over~$\complex$) to a group of
\index{lower-triangular matrices}lower-triangular matrices. By working in
$\SL(\ell,\complex)$\index{SL(l,C)*$\SL(\ell,\complex)$}, let us assume,
for simplicity, that $R$~itself is lower triangular. That is, $R \subset
\mathbb{D}_\ell \ltimes \mathbb{U}_\ell$.

Let $\pi \colon \mathbb{D}_\ell \ltimes \mathbb{U}_\ell \to
\mathbb{D}_\ell$ be the natural projection. It is not difficult to see
that there exists $r \in R$, such that $\pi(R) \subset \Zar{\langle \pi(r)
\rangle}$ (by using \pref{T=split+cpct} and \pref{SplitTori}). Let 
 $$ \text{$T = \Zar{\langle r_s r_e \rangle}^\circ$
 and
 $U = R \cap \mathbb{U}_\ell$} .$$
 Because $\pi(r_s r_e) = \pi(r)$, we have $\pi(R) \subset \pi(T)$, so, for
any $g \in R$, there exists $t \in T$, such that $\pi(t) = \pi(g)$. Then
$\pi(t^{-1} g) = e$, so $t^{-1} g \in U$. Therefore $g \in t U \subset T
\ltimes U$. Since $g \in R$ is arbitrary, we conclude that 
 $$ R = T \ltimes U .$$
 This yields the desired decomposition $G = (L T) \ltimes U$.
 \end{proof}

\begin{rem}
 The subgroup~$U$ of \pref{AlgGrpTLU} is the \emph{unique} maximal
\index{unipotent!subgroup}unipotent \index{subgroup!normal}normal subgroup
of~$G$. It is called the \defit[unipotent!radical]{unipotent radical}
of~$G$.
 \end{rem}

It is obvious (from the \term[Jordan!decomposition]{Jordan decomposition})
that every element of a compact real \index{algebraic
group!compact}algebraic group is \index{elliptic element}elliptic. We
conclude this section by recording (without proof) the fact that this
characterizes the compact real \index{algebraic group!compact}algebraic
groups.

\begin{thm} \label{Gcpct}
 An almost-Zariski closed subgroup of\/ $\SL(\ell,\real)$ is
\index{algebraic group!compact}compact if and only if all of its elements
are elliptic.
 \end{thm}

\begin{cor} \label{AlgGrpCpct} \ 
 \begin{enumerate}

 \item \label{AlgGrpCpct-unip}
 A nontrivial \index{unipotent!subgroup}unipotent subgroup~$U$ of\/
$\SL(\ell,\real)$ is never \index{algebraic group!compact}compact.

 \item  \label{AlgGrpCpct-torus}
 A \index{torus!compact}torus $T$ in $\SL(\ell,\real)$ is \index{algebraic
group!compact}compact if and only if none of its nontrivial elements are
hyperbolic.

 \item \label{AlgGrpCpct-semi} % \label{SemiNoUnip}
 A \index{subgroup!connected}connected, \index{Lie
group!semisimple}semisimple subgroup~$L$ of\/ $\SL(\ell,\real)$ is
\index{algebraic group!compact}compact if and only if it has no nontrivial
unipotent elements {\rm(}also, if and only if it has no nontrivial
hyperbolic elements{\rm)}.

 \end{enumerate}
 \end{cor}

We conclude this section with two basic results about tori.

 \begin{defn} 
 A torus~$T$ is
\defit[hyperbolic!torus|indsee{torus,~hyperbolic}]{hyperbolic} (or
\defit[R-split*$\real$-split!torus]{$\real$-split}) if every element
of~$T$ is hyperbolic.
 \end{defn}

\begin{cor} \label{T=split+cpct}
 Any \index{subgroup!connected}connected \index{torus!algebraic}torus~$T$
has a unique decomposition into a direct product $T = T_h \times T_c$,
where
 \begin{enumerate}
 \item $T_h$ is a \index{torus!hyperbolic}hyperbolic torus,
 and
 \item $T_c$ is a \index{torus!compact}compact torus.
 \end{enumerate}
 \end{cor}

\begin{proof}
 Let 
 $$ T_h = \{\, g \in T \mid \text{$g$ is hyperbolic} \,\} $$
 and
 $$ T_c = \{\, g \in T \mid \text{$g$ is elliptic} \,\}. $$
 Because $T$ is \index{group!abelian}abelian, it is easy to see that $T_h$
and $T_c$ are subgroups of~$T$ \seeex{CommuteJordan}. It is immediate from
the real Jordan decomposition that $T = T_h \times T_c$.

All that remains is to show that $T_h$ and $T_c$ are almost Zariski
closed. 
 \begin{enumerate}
 \item[($T_h$)] Since $T_h$ is a set of commuting matrices that are
\index{diagonalizable!over R*over~$\real$}diagonalizable over~$\real$,
there exists $h \in \SL(\ell,\real)$, such that $h^{-1} T_h h \subset
\mathbb{D}_\ell$ \seeex{SimulDiag}. Hence, $T_h = T \cap (h
\mathbb{D}_\ell h^{-1})$ is almost Zariski closed.
 \item[($T_e$)] Let 
 $$ \text{$\mathbb{D}_\ell^\complex$ be the group of diagonal matrices in
$\SL(\ell,\complex)$,} $$
 and \index{eigenvalue}
 $$C = \bigset{ g \in \mathbb{D}_\ell^\complex }{
 \begin{matrix}
  \text{every eigenvalue of~$g$} \\
 \text{has absolute value~$1$}
 \end{matrix}
 } .$$
 Because $T$ is a \index{torus!algebraic}torus, there exists $h \in
\SL(\ell,\complex)$, such that $h^{-1} T h \subset
\mathbb{D}_\ell^\complex$. Then 
 $T_c = T \cap h C h^{-1}$ is \index{torus!compact}compact. So it is
Zariski closed (see Prop.~\ref{Cpct->ZarClosed} below). \qedhere
 \end{enumerate}
 \end{proof}

A (real) \defit{representation} of a group is a homomorphism into
$\SL(m,\real)$, for some~$m$.
 The following result provides an explicit description of the
representations of any hyperbolic torus.

\begin{cor} \label{TorusRepn}
 Suppose
 \begin{itemize}
 \item $T$ is a {\rm(}hyperbolic{\rm)} \index{torus!hyperbolic}torus that
consists of diagonal matrices in $\SL(\ell,\real)$,
 and 
 \item $\rho \colon T \to \SL(m,\real)$ is any
\term[homomorphism!polynomial]{polynomial homomorphism}.
 \end{itemize}
 Then there exists $h \in \SL(n,\real)$, such that, letting 
 $$ \text{$\rho'(t) = h^{-1} \, \rho(t) \, h$ \quad for $t \in T$,}$$
 we have:
 \begin{enumerate}
 \item \label{TorusRepn-diag}
 $\rho'(T) \subset \mathbb{D}_m$,
 and
 \item \label{TorusRepn-powers}
 For each $j$ with $1 \le j \le m$, there are integers
$n_1,\ldots,n_\ell$, such that
 $$ \rho'(t)_{j,j} = t_{1,1}^{n_1} t_{2,2}^{n_2} \cdots
t_{\ell,\ell}^{n_\ell} $$
 for all $t \in T$.
 \end{enumerate}
 \end{cor}

\begin{proof}
 \pref{TorusRepn-diag} Since $\rho(T)$ is a set of commuting matrices that
are \index{diagonalizable!over R*over~$\real$}diagonalizable over~$\real$,
there exists $h \in \SL(m,\real)$, such that $h^{-1} \rho(T) h \subset
\mathbb{D}_m$ \seeex{SimulDiag}.

\pref{TorusRepn-powers} For each~$j$, $\rho'(t)_{j,j}$ defines a
\term[homomorphism!polynomial]{polynomial homomorphism} from~$T$
to~$\real^+$. With the help of Lie theory, it is not difficult to see that
any such homomorphism is of the given form \seeex{RegularCharacter}.
 \end{proof}

\begin{exercises}

\item \label{SO2abel} 
 Show:
 \begin{enumerate}
 \item \index{SO(Q)*$\SO(Q)$}$\SO(2)$ is a compact
\index{torus!compact}torus,
 and
 \item \index{SO(Q)*$\SO(Q)$}$\SO(1,1)^\circ$ is a hyperbolic
\index{torus!hyperbolic}torus.
 \end{enumerate}
 \hint{We have
 $$ \mbox{$\SO(2) = \left\{\begin{bmatrix} \cos\theta & \sin\theta \\
-\sin\theta & \cos\theta \end{bmatrix} \right\}$ 
 and 
 $\SO(1,1) = \left\{\pm\begin{bmatrix} \cosh t & \sinh t \\ \sinh t &
\cosh t \end{bmatrix} \right\}$,}$$
 where $\cosh t = (e^t + e^{-t})/2$ and $\sinh t = (e^t - e^{-t})/2$.}

\item \label{UnipsClosed}
 Show:
 \begin{enumerate}
 \item The set of unipotent elements of $\SL(\ell,\real)$ is Zariski
closed.
 \item If $U$ is a \index{unipotent!subgroup}unipotent subgroup of
$\SL(\ell,\real)$, then $\Zar{U}$ is also unipotent.
 \end{enumerate}

\item Prove the easy direction ($\Rightarrow$) of Thm.~\ref{Gcpct}.

\item Assume that Thm.~\ref{Gcpct} has been proved for \index{Lie
group!semisimple}semisimple groups. Prove the general case.
 \hint{Use Thm.~\ref{AlgGrpTLU}.}

\item \label{Unip->Lower}
 (\emphit{Advanced})
 Prove \term[Theorem!Engel]{Engel's Theorem} \fullref{Unips+Tori}{unip}.
 \hint{($\Leftarrow$)~It suffices to show that $U$ fixes some nonzero
vector~$v$. (For then we may consider the action of~$U$ on
$\real^\ell/\real v$, and complete the proof by induction on~$\ell$.)
There is no harm in working over~$\complex$, rather than~$\real$, and we
may assume there are no $U$-invariant\index{subspace, invariant} subspaces
of~$\complex^\ell$. Then a theorem of \term[Theorem!Burnside]{Burnside}
states that every $\ell \times \ell$ matrix~$M$ is a linear combination of
elements of~$U$. Hence, for any $u \in U$, $\trace(u M) = \trace M$. Since
$M$ is arbitrary, we conclude that $u = \Id$.}

 \end{exercises}

\section{Chevalley's Theorem and applications}
 \label{ChevalleySect}
 \index{Theorem!Chevalley}

\begin{notation}
 For a map $\rho \colon G \to Z$ and $g \in G$, we often write 
 \nindex{$g^\rho$ = $\rho(g)$}{$g^\rho$}
 for the image of~$g$ under~$\rho$. That is, $g^\rho$ is another notation
for $\rho(g)$.
 \end{notation}

\begin{prop}[(Chevalley's Theorem)] \label{ChevalleyFixLine}
\index{Theorem!Chevalley}
 A subgroup~$H$ of a real \index{algebraic group!real}algebraic group~$G$
is Zariski closed if and only if, for some~$m$, there exist 
 \begin{itemize}
 \item a \term[homomorphism!polynomial]{polynomial homomorphism} $\rho
\colon G \to \SL(m,\real)$,
 and
 \item a vector $v \in \real^m$,
 \end{itemize}
 such that 
 $H = \{\, h \in G \mid v h^\rho \in \real v \,\}$.
 \index{stabilizer!of a subspace}
 \end{prop}

\begin{proof}
 ($\Leftarrow$) This follows easily from
Eg.~\fullref{EgsOfAlgGrps}{Stab(V)} and Exer.~\ref{InvImgAlmZar}.

 ($\Rightarrow$) There is no harm in assuming $G = \SL(\ell,\real)$.
There is a \emph{finite} subset~$\mathcal{Q}$ of $\real[x_{1,1}, \ldots,
x_{\ell,\ell}]$, such that $H = \Var(\mathcal{Q})$\index{variety}
\seeex{Noetherian-FiniteVariety}. Choose $d \in \integer^+$, such that
$\deg Q < d$ for all $Q \in \mathcal{Q}$, and let
 \begin{itemize}
 \item $V = \{\, Q \in \real[x_{1,1}, \ldots,
x_{\ell,\ell}]  \mid \deg Q < d \,\}$
 and
 \item $W = \{\, Q \in V \mid \text{$Q(h_{i,j}) = 0$ for all $h \in H$}
\,\}$.
 \end{itemize}
 Thus, we have $H = \bigcap_{Q \in W} \Var(\{Q\})$\index{variety}.

There is a natural \index{homomorphism!of Lie groups}homomorphism~$\rho$
from $\SL(\ell,\real)$ to the group $\SL(V)$\index{SL(V)*$\SL(V)$} of
(special) linear transformations on~$V$, defined by%
 \begin{equation}
  (Q g^\rho) (x_{i,j}) = Q \bigl( (g x)_{i,j} \bigr)
 \label{ChevalleyFixLinePf-rho}
 \end{equation}
 \seeex{ChevalleyRho-homo}. Note that we have
$\Stab_{\SL(\ell,\real)}(W) = H$\index{stabilizer!of a subspace}
\seeex{ChevalleyRho-stab}. By taking a basis for~$V$, we may think
of~$\rho$ as a \term[homomorphism!polynomial]{polynomial homomorphism}
into $\SL(\dim V, \real)$ \seeex{ChevalleyRho-poly}. Then this is almost
exactly what we want; the only problem is that, instead of a
$1$-dimensional space $\real v$, we have the space~$W$ of (possibly)
larger \index{dimension!of a vector space}dimension.

To complete the proof, we convert $W$ into a $1$-dimensional space, by
using a standard trick of \term{multilinear algebra}. For $k = \dim W$, we
let 
 $$ \mbox{$V' = \bigwedge^k V$ and $W' =
\bigwedge^k W \subset V'$,} $$
  where $\bigwedge^k V$ denotes the $k$th \term{exterior power} of~$V$.
 Now $\rho$ naturally induces a \term[homomorphism!polynomial]{polynomial
homomorphism} $\rho' \colon \SL(\ell,\real) \to
\SL(V')$\index{stabilizer!of a subspace}, and, for this action, $H =
\Stab_{\SL(\ell,\real)}(W)$ \seeex{Stab(wedge)}. By choosing a basis
for~$V'$, we can think of~$\rho'$ as a homomorphism into $\SL\left({\dim V
\choose k}, \real\right)$\index{SL(V)*$\SL(V)$}. Since $\dim W' = {\dim W
\choose k} = 1$, we obtain the desired conclusion (with~$\rho'$ in the
place of~$\rho$) by letting~$v$ be any nonzero vector in~$W'$.
 \end{proof}

\begin{proof}[Proof of Thm.~\ref{JorDecinG}.]
 From \index{Theorem!Chevalley}Chevalley's Theorem
\pref{ChevalleyFixLine}, we know there exist
 \begin{itemize}
 \item a polynomial
\index{homomorphism!polynomial}homomorphism $\rho \colon \SL(\ell,\real)
\to \SL(m,\real)$, for some~$m$,
 and
 \item a vector $v \in \real^m$,
 \end{itemize}
 such that 
 $G = \{\, g \in \SL(\ell,\real) \mid v g^\rho \in \real v \,\}$.
 Furthermore, from the explicit description of~$\rho$ in the proof
of \index{Theorem!Chevalley}Prop.~\ref{ChevalleyFixLine}, we see that it
satisfies the conclusions of Cor.~\ref{rho(unip/hyp)} with
$\SL(\ell,\real)$ in the place of~$G$ \cfex{ChevalleyRhoJordan}. Thus, for
any $g \in \SL(\ell,\real)$, we have
 $$ \text{$(g_u)^\rho = (g^\rho)_u$, $(g_h)^\rho = (g^\rho)_h$, and
$(g_e)^\rho = (g^\rho)_e$}.$$

For any $g \in G$, we have $v g^\rho \in \real v$. In other words, $v$~is
an \index{eigenvector}eigenvector for~$g^\rho$. Then $v$ is also an
\term{eigenvector} for $(g^\rho)_u$ \seeex{JordanEigenvec}. Since
$(g_u)^\rho = (g^\rho)_u$, this implies
 $  v (g_u)^\rho \in \real v$,
 so $g_u \in G$. By the same argument, $g_h \in G$ and $g_e \in G$.
 \end{proof}

\index{Theorem!Chevalley}Chevalley's Theorem yields an explicit
description of the \index{torus!hyperbolic}hyperbolic tori.

\begin{cor} \label{SplitTori}
 Suppose $T$ is a \index{subgroup!connected}connected group of diagonal
matrices in $\SL(\ell,\real)$, and let $d = \dim T$. Then $T$ is almost
Zariski closed if and only if there are \index{linear!functional}linear
functionals $\lambda_1,\ldots,\lambda_\ell \colon \real^d \to \real$, such
that
 \begin{enumerate}
 \item $T = \bigset{
 \begin{bmatrix}
 e^{\lambda_1(x)} \\
 & e^{\lambda_2(x)} \\
 && \ddots \\
 &&&e^{\lambda_\ell(x)} 
 \end{bmatrix}
 }{ x \in \real^d} $,
 and
 \item for each~$i$, there are \emph{integers} $n_1,\ldots,n_d$, such that
 $$\text{$\lambda_i(x_1,\ldots,x_d) = n_1 x_1 + \cdots + n_d x_d$ for all
$x \in \real^d$}. $$
 \end{enumerate}
 \end{cor}

\begin{proof}
 Combine \index{Theorem!Chevalley}Prop.~\ref{ChevalleyFixLine} with
Cor.~\ref{TorusRepn} \seeex{SplitToriEx}.
 \end{proof}

\begin{exercises}

\item \label{TranslateSameDegree}
 Suppose $Q \in \real[x_{1,1}, \ldots, x_{\ell,\ell}]$ and $g  \in
\SL(\ell,\real)$. 
 \begin{itemize}
 \item Let $\phi \colon \SL(\ell,\real) \to \real$ be the polynomial
function corresponding to~$Q$, and 
 \item define $\phi' \colon \SL(\ell,\real) \to \real$ by $\phi'(x) =
\phi(gx)$.
 \end{itemize}
 Show there exists $Q' \in \real[x_{1,1}, \ldots, x_{\ell,\ell}]$, with
$\deg Q' = \deg Q$, such that $\phi'$ is the polynomial function
corresponding to~$Q'$.
 \hint{For fixed~$g$, the \index{matrix!entry}matrix entries of $gh$ are
linear functions of~$h$.}

\item \label{ChevalleyRho}
 Define $\rho \colon \SL(\ell,\real) \to \SL(V)$\index{SL(V)*$\SL(V)$} as
in \eqref{ChevalleyFixLinePf-rho}.\index{Theorem!Chevalley}
 \begin{enumerate}
 \item \label{ChevalleyRho-homo}
 Show $\rho$ is a group \index{homomorphism!of Lie groups}homomorphism.
 \item \label{ChevalleyRho-stab} 
 For the subspace~$W$ defined in the proof of Prop.~\ref{ChevalleyFixLine},
show $H = \Stab_{\SL(\ell,\real)}(W)$\index{stabilizer!of a subspace}.
 \item \label{ChevalleyRho-poly}
 By taking a basis for~$V$, we may think of $\rho$ as a map into $\SL(\dim
V, \real)$. Show $\rho$ is a polynomial.
 \end{enumerate}
 \hint{\pref{ChevalleyRho-stab} We have $Q \subset W$.}

\item \label{Stab(wedge)}
 Suppose 
 \begin{itemize}
 \item $W$ is a subspace of a real vector space~$V$,
 \item $g$ is an invertible linear transformation on~$V$,
 and
 \item $k = \dim W$.
 \end{itemize}
 Show $\bigwedge^k (Wg) = \bigwedge^k W$\index{exterior power} if and only
if $Wg = W$.

\item \label{ChevalleyRhoJordan}
  Define $\rho \colon \SL(\ell,\real) \to
\SL(V)$\index{SL(V)*$\SL(V)$} as in \eqref{ChevalleyFixLinePf-rho}.
 \begin{enumerate}
 \item \label{ChevalleyRhoJordan-hyper}
 Show that if $g$ is hyperbolic, then $\rho(g)$ is hyperbolic.
 \item \label{ChevalleyRhoJordan-ell}
 Show that if $g$ is elliptic, then $\rho(g)$ is elliptic.
 \item \label{ChevalleyRhoJordan-unip}
 Show that if $g$ is unipotent, then $\rho(g)$ is unipotent.
 \end{enumerate}
 \hint{(\ref{ChevalleyRhoJordan-hyper},\ref{ChevalleyRhoJordan-ell})~If
$g$ is diagonal, then any \term{monomial} is an
\index{eigenvector}eigenvector of $g^\rho$.}

\item \label{SplitToriEx}
 Prove Cor.~\ref{SplitTori}.

 \end{exercises}

\section{Subgroups that are almost Zariski closed} \label{AlmZarSect}

We begin the section with some results that guarantee certain types of
groups are almost Zariski closed.

\begin{prop} \label{Cpct->ZarClosed}
 Any \index{algebraic group!compact}\index{Lie group!compact}compact
subgroup of\/ $\SL(\ell,\real)$ is Zariski closed.
 \end{prop}

\begin{proof}
 Suppose $C$ is a \index{Lie group!compact}compact subgroup of
$\SL(\ell,\real)$, and $g$ is an element of $\SL(\ell,\real)
\smallsetminus C$. It suffices to find a polynomial $\phi$ on
$\SL(\ell,\real)$, such that $\phi(C) = 0$, but $\phi(g) \neq 0$.

The sets $C$ and $Cg$ are compact and disjoint, so, for any $\epsilon >
0$, the \term[Theorem!Stone-Weierstrass]{Stone-Weierstrass Theorem}
implies there is a polynomial $\phi_0$, such that $\phi_0(c) < \epsilon$
and $\phi_0(cg) > 1-\epsilon$ for all $c \in C$. (For our purposes, we may
choose any $\epsilon < 1/2$.) For each $c \in C$, let $\phi_c(x) =
\phi(cx)$, so $\phi_c$ is a polynomial of the same degree as~$\phi_0$
\seeex{TranslateSameDegree}. Define $\overline{\phi} \colon
\SL(\ell,\real) \to \real$ by averaging over $c \in C$:
 $$ \overline{\phi}(x) = \int_C \phi_c(x) \, dc ,$$
 where $dc$ is the Haar measure on~$C$, normalized to be a probability
measure.
 Then
 \begin{enumerate}
 \item $\overline{\phi}(c) < \epsilon$ for $c \in C$,
 \item $\overline{\phi}(g) > 1 - \epsilon$,
 \item $\overline{\phi}$ is constant on~$C$ (because Haar measure is
invariant\index{measure!invariant}),
 and
 \item $\overline{\phi}$ is a polynomial function (each of its coefficients
is the average of the corresponding coefficients of the $\phi_c$'s).
 \end{enumerate}
 Now let $\phi(x) = \overline{\phi}(x) - \overline{\phi}(c)$ for any $c
\in C$.
 \end{proof}

\begin{prop} \label{Unip->Closed}
 If $U$ is a \index{subgroup!connected}connected,
\index{unipotent!subgroup}unipotent subgroup of\/ $\SL(\ell,\real)$, then
$U$ is Zariski closed.
 \end{prop}

\begin{proof}[Proof {\rm(\emphit {requires some Lie theory})}.]
  By passing to a conjugate, we may assume $U \subset
\mathbb{U}_\ell$. The \index{Lie algebra}Lie algebra $\AMSfrak{U}_\ell$
of~$\mathbb{U}_\ell$ is the space of strictly \index{lower-triangular
matrices!strictly}lower-triangular matrices \seeex{LieAlgEg-U}. Because
$A^\ell = 0$ for $A \in \AMSfrak{U}_\ell$\index{matrix!nilpotent}, the
\index{exponential!map}exponential map
 $$ \exp(A) = \Id + A + \frac{1}{2} A^2 + \cdots + \frac{1}{(\ell-1)!}
A^{\ell-1} $$
 is a polynomial function on~$\AMSfrak{U}_\ell$, and its inverse, the
\term{logarithm} map
 $$ \log (\Id + N) = N - \frac{1}{2} N^2 + \frac{1}{3} N^3  \pm \cdots \pm
\frac{1}{\ell-1} N^{\ell-1}, $$
 is a polynomial function on~$\mathbb{U}_\ell$. 

Therefore $\exp$ is a bijection from $\AMSfrak{U}_\ell$
onto~$\mathbb{U}_\ell$, so $U = \exp \Lie U$, where $\Lie U$ is the
\index{Lie algebra}Lie algebra of~$U$. This means 
 $$ U = \{\, u \in \mathbb{U}_\ell \mid \log u \in \Lie U \,\} .$$
 Since $\log$ is a polynomial function (and $\Lie U$, being a linear
subspace, is defined by polynomial equations --- in fact, linear
equations), this implies that $U$ is defined by polynomial equations.
Therefore, $U$ is Zariski closed.
 \end{proof}

The following result is somewhat more difficult; we omit the proof.

\begin{thm} \label{SemiSimp->ZarClosed}
 If $L$ is any \index{subgroup!connected}connected, \index{Lie
group!semisimple}semisimple subgroup of\/ $\SL(\ell,\real)$, then $L$ is
almost Zariski closed.
 \end{thm} 

The following three results show that being almost Zariski closed is
preserved by certain natural operations. We state the first without proof.

\begin{prop} \label{AB-closed}
 If $A$~and~$B$ are almost-Zariski closed subgroups of\/ $\SL(\ell,\real)$,
such that $AB$ is a subgroup, then $AB$ is almost Zariski closed. 
 \end{prop}

\begin{cor} \label{homo-img}
 If $G$ and~$H$ are almost Zariski closed, and $\rho$ is a
\term[homomorphism!polynomial]{polynomial homomorphism} from~$G$ to~$H$,
then the image $\rho(G)$ is an almost-Zariski closed subgroup of~$H$.
 \end{cor}

\begin{proof}
 By passing to a finite-index subgroup, we may assume $G$ is
\index{subgroup!connected}connected. Write $G = (T L) \ltimes U$, as in
Thm.~\ref{AlgGrpTLU}. From Prop.~\ref{AB-closed}, it suffices to show that
$\rho(U)$, $\rho(L)$, and $\rho(T)$ are almost Zariski closed. The
subgroups $\rho(U)$ and $\rho(L)$ are handled by Prop.~\ref{Unip->Closed}
and Thm.~\ref{SemiSimp->ZarClosed}. 

Write $T = T_h \times T_c$, where $T_h$ is
\index{torus!hyperbolic}hyperbolic and $T_c$ is
\index{torus!compact}compact \seecor{T=split+cpct}. Then $\rho(T_c)$, being
compact, is Zariski closed \seeprop{Cpct->ZarClosed}. The subgroup
$\rho(T_h)$ is handled easily by combining Cors.~\ref{SplitTori}
and~\ref{TorusRepn} \seeex{rho(torus)Ex}.
 \end{proof}

\begin{cor} \label{G'almZar}
 If $G$ is any \index{subgroup!connected}connected subgroup of\/
$\SL(\ell,\real)$, then the \index{subgroup!commutator}commutator subgroup
$[G,G]$ is almost Zariski closed.
 \end{cor}

\begin{proof}
 Write $G = (L T) \ltimes U$, as in Thm.~\ref{AlgGrpTLU}. Because $T$ is
\index{group!abelian}abelian and $[L,L] = L$, we see that $[G,G]$ is a
(\index{subgroup!connected}connected subgroup of $L \ltimes U$ that
contains~$L$. Hence $[G,G] = L \ltimes \check U$, where $\check U = [G,G]
\cap U$ \seeex{SemiProdWithFactorIn}. Furthermore, since $[G,G]$ is
\index{subgroup!connected}connected, we know $\check U$ is connected, so
$\check U$ is \term[Zariski!closed]{Zariski closed} \seeprop{Unip->Closed}.
Since $L$ is almost Zariski closed \seethm{SemiSimp->ZarClosed}, this
implies $[G,G] = L \check U$ is almost Zariski closed \seeprop{AB-closed},
as desired.
 \end{proof}

\begin{cor} \label{ZarG'=G'}
 If $G$ is any \index{subgroup!connected}connected subgroup of\/
$\SL(\ell,\real)$, then $[\Zar{G},\Zar{G}] =
[G,G]$\index{subgroup!commutator}, so $\Zar{G}/G$ is
\index{group!abelian}abelian.
 \end{cor}

\begin{proof}
 Define $c \colon \Zar{G} \times \Zar{G} \to \Zar{G}$ by $c(g,h) = g^{-1}
h^{-1} g h = [g,h]$. Then $c$ is a polynomial \seeex{MultIsPoly}. Since
$c(G \times G) \subset [G,G]$ and $[G,G]$ is almost Zariski closed, we
conclude immediately that $[\Zar{G},\Zar{G}]^\circ \subset [G,G]$
\cfex{InvImgAlmZar}. This is almost what we want, but some additional
theory (which we omit) is required in order to show that
$[\Zar{G},\Zar{G}]$ is \index{subgroup!connected}connected, rather than
having finitely many \index{component!connected}components.

Because $[\Zar{G},\Zar{G}] \subset G$, it is immediate that $\Zar{G}/G$ is
\index{group!abelian}abelian.
 \end{proof}

For connected groups, we now show that tori present the only obstruction
to being almost Zariski closed.

\begin{cor} \label{GTAlmZar}
 If $G$ is any \index{subgroup!connected}connected subgroup of\/
$\SL(\ell,\real)$, then there is a connected, almost-Zariski closed
\index{torus!algebraic}torus~$T$ of\/~$\Zar{G}$, such that $G T$ is almost
Zariski closed.
 \end{cor}

\begin{proof}
 Write $\Zar{G}^\circ = (T L) \ltimes U$, with $T,L,U$ as in
Thm.~\ref{AlgGrpTLU}. Because $L = [L,L] \subset [\Zar{G}, \Zar{G}]$, we
know $L \subset G$ \seecor{ZarG'=G'}. Furthermore, because $T$
\index{normalizer}normalizes~$G$ \seeex{norm-closed}, we may assume $T
\subset G$, by replacing $G$ with $GT$. 

Therefore $G = (TL) \ltimes (U \cap G)$ \seeex{SemiProdWithFactorIn}.
Furthermore, since $G$ is \index{subgroup!connected}connected, we know
that $U \cap G$ is connected, so $U \cap G$ is Zariski closed
\seeprop{Unip->Closed}. Then Prop.~\ref{AB-closed} implies that $G = (TL)
\ltimes (U \cap G)$ is almost Zariski closed.
 \end{proof}

We will make use of the following technical result:

\begin{lem} \label{poly+U=algic}
 Show that if 
 \begin{itemize}
 \item $G$ is an almost-Zariski closed subgroup of $\SL(\ell,\real)$,
 \item $H$ and~$V$ are \index{subgroup!connected}connected subgroups
of~$G$ that are almost Zariski closed,
 and
 \item $f \colon V \to G$ is a \index{function!rational}rational function
{\rm(}not necessarily a \index{homomorphism!polynomial}homomorphism{\rm)},
with $f(e) = e$,
 \end{itemize}
 then the subgroup $\langle H, f(V) \rangle$ is almost Zariski closed.
 \end{lem}

\begin{proof}[\emphit{Plausibility argument.}]
 There is no harm in assuming that $G = \Zar{\langle f(H), H
\rangle}^\circ$, so we wish to show that $H$ and $f(V)$, taken together,
generate~$G$. Since $[G,G] H$ is
 \begin{itemize}
 \item almost Zariski closed \seeprop{AB-closed},
 \item contained in $\langle H, f(V) \rangle$ \seecor{ZarG'=G'},
 and
 \item \index{subgroup!normal}normal in~$G$ (because it contains $[G,G]$),
 \end{itemize}
 there is no harm in modding it out. Thus, 
 we may assume that $G$ is \index{group!abelian}abelian and that $H =
\{e\}$.

Now, using the fact that $G$ is \index{group!abelian}abelian, we have $G =
A \times C \times U$, where $A$~is a hyperbolic
\index{torus!hyperbolic}torus, $C$~is a compact
\index{torus!compact}torus, and $U$~is unipotent
(see Thm.~\ref{AlgGrpTLU} and Cor.~\ref{T=split+cpct}). Because these are
three completely different types of groups, it is not difficult to believe
that there are subgroups $A_V$, $C_V$, and~$U_V$ of $A$, $C$, and~$V$,
respectively, such that $\langle f(V) \rangle = A_V \times C_V \times U_V$
\cfex{f(R)algicEgEx}.

Now $U_V$, being \index{subgroup!connected}connected and
\index{unipotent!subgroup}unipotent, is Zariski closed
\seeprop{Unip->Closed}. The other two require some argument.
 \end{proof}

\begin{exercises}

\item  \label{unip-sc}
 Show that every  \index{algebraic
group!unipotent|indsee{unipotent~subgroup}} unipotent real
\index{unipotent!subgroup}algebraic group is
\index{subgroup!connected}connected and \term{simply connected}.
 \hint{See proof of \pref{Unip->Closed}.}

\item \label{rho(torus)Ex}
 Complete the proof of Cor.~\ref{homo-img}, by showing that if $T$ is a
hyperbolic \index{torus!hyperbolic}torus, and $\rho \colon T \to
\SL(m,\real)$ is a \term[homomorphism!polynomial]{polynomial
homomorphism}, then $\rho(T)$ is almost Zariski closed.
 \hint{Use Cors.~\ref{SplitTori} and~\ref{TorusRepn}.}

\item \label{SemiProdWithFactorIn}
 Show that if $G$ is a subgroup of a \term{semidirect product} $A \ltimes
B$, and $A \subset G$, then $G = A \ltimes (G \cap B)$. If, in addition,
$G$ is \index{subgroup!connected|)}connected, show that $G \cap B$ is
connected.

%\item \label{ZarG'=G'PfEx}
% Show that if $G$ is connected, then $[\Zar{G},\Zar{G}]$ is connected.
% \hint{You may assume (without proof) that if $L$ is any connected,
%semisimple subgroup of $\SL(\ell,\real)$, then $\Zar{L}/L$ is
%\index{group!abelian}abelian. Use
%the structure theorem \pref{AlgGrpTLU}.}

\item \label{f(R)algicEgEx}
 Suppose $Q \colon \real \to \real$ is any nonconstant polynomial with
$Q(0) = 0$, and define $f \colon \real \to \mathbb{D}_2 \times
\mathbb{U}_2 \subset \SL(4,\real)$\index{SL(4,R)*$\SL(4,\real)$} by
 $$ f(t) = \begin{bmatrix}
 1 + t^2 & 0 & 0 & 0 \\
 0 & 1/(1 + t^2) & 0 & 0 \\
 0 & 0 & 1 & 0 \\
 0 & 0 & Q(t) & 1
 \end{bmatrix} .$$
 Show $\langle f(\real) \rangle = \mathbb{D}_2 \times \mathbb{U}_2$.

 \end{exercises}

\section{Borel Density Theorem} \label{BDTSect}

The \term[Theorem!Borel Density]{Borel Density Theorem} \pref{BDT} is a
generalization of the important fact that if $\Gamma =
\SL(\ell,\integer)$\index{SL(l,Z)*$\SL(\ell,\integer)$}, then $\Zar{\Gamma}
= \SL(\ell,\real)$ \seeex{SLZarDensinSLREx}. Because the Zariski closure
of~$\Gamma$ is all of $\SL(\ell,\real)$, we may say that $\Gamma$ is
\defit[Zariski!dense]{Zariski dense} in $\SL(\ell,\real)$. That is why
this is known as a ``density" theorem.

\begin{prop}[{(\term[Theorem!Borel Density]{Borel Density Theorem})}]
\label{BDT}
 If\/ $\Gamma$ is any \term[lattice|(]{lattice} in any closed subgroup~$G$
of\/ $\SL(\ell,\real)$, then the Zariski closure~$\Zar{\Gamma}$
of\/~$\Gamma$ contains
 \begin{enumerate}
 \item every unipotent element of~$G$
 and
 \item every hyperbolic element of~$G$.
 \end{enumerate}
 \end{prop}

We precede the proof with a remark and two lemmas. 

\begin{rem} \
 \begin{enumerate}
 \item If $G$ is a \index{Lie group!compact}compact group, then the trivial
subgroup $\Gamma = \{e\}$ is a \term{lattice} in~$G$, and $\Zar{\Gamma} =
\{e\}$ does not contain any nontrivial elements of~$G$. This is consistent
with Prop.~\ref{BDT}, because nontrivial elements of a compact group are
neither unipotent nor hyperbolic \seecor{AlgGrpCpct}.
 \item Although we do not prove this, $\Zar{\Gamma}$ actually contains
every unipotent or hyperbolic element of~$\Zar{G}$, not only those of~$G$.
 \end{enumerate}
 \end{rem}

\begin{lem}[{(\term[Theorem!Poincar\'e Recurrence]{Poincar\'e Recurrence
Theorem})}] \label{Poincare}
 Let
 \begin{itemize}
 \item $(\Omega,d)$ be a metric space;
 \item $T \colon \Omega \to \Omega$ be a
\index{homeomorphism}homeomorphism; and
 \item $\mu$ be a $T$-invariant\index{measure!invariant} probability
measure on~$A$.
 \end{itemize}
 Then, for almost every $a \in \Omega$, there is a
sequence $n_k \to \infty$, such that $T^{n_k} a \to a$.
 \end{lem}

\begin{proof}
 Let
 $$ A_\epsilon = \{\, a \in \Omega \mid \forall m>0, \
d(T^m a, a) > \epsilon \,\}.$$
 It suffices to show $\mu(A_\epsilon) = 0$ for every
$\epsilon$.

Suppose $\mu(A_{\epsilon}) > 0$. Then we may choose a
subset~$B$ of~$A_\epsilon$, such that $\mu(B) > 0$ and
$\operatorname{diam}(B) < \epsilon$. The sets 
 $B,
T^{-1}B, T^{-2} B, \ldots$ cannot all be disjoint, because they all have
the same measure and $\mu(\Omega) < \infty$. Hence, $T^{-m} B \cap T^{-n} B
\neq \emptyset$, for some $m,n \in\integer^+$ with $m > n$. By
applying $T^{n}$, we may assume $n = 0$. For $a \in T^{-m} B \cap B$, we
have $T^m a \in B$ and $a \in B$, so 
 $$d(T^m a, a) \le \operatorname{diam}(B) < \epsilon .$$
 Since $a \in B \subset A_\epsilon$, this contradicts the definition
of~$A_\epsilon$.
 \end{proof}

\begin{notation} \ 
 \begin{itemize}
 \item Recall that the \defit{projective space} 
 \nindex{$\RP{m-1}$ = real projective space}{$\RP{m-1}$} is, by
definition, the set of one-dimensional\index{dimension!of a vector space}
subspaces of~$\real^m$. Alternatively,  $\RP{m-1}$ can be viewed as the
set of equivalence classes of the \term{equivalence relation} on $\real^m
\smallsetminus \{0\}$ defined by 
 $$ v \sim w \Leftrightarrow \text{$v = \alpha w$ for some $\alpha \in
\real \smallsetminus \{0\}$}. $$
 From the alternate description, it is easy to see that $\RP{m-1}$ is an
$(m-1)$-dimensional\index{dimension!of a manifold} smooth \term{manifold}
\seeex{RPnMfld}.

\item There is a natural action of $\SL(m,\real)$ on $\RP{m-1}$, defined
by $[v] g = [vg]$, where, for each nonzero $v \in \real^m$, we let $[v] =
\real v$ be the image of~$v$ in $\RP{m-1}$.
 \nindex{$[v]$ = image of vector~$v$ in $\RP{m-1}$}
 \goodbreak
 \end{itemize}
 \end{notation}

\begin{lem} \label{mu(PV)}
 Assume
 \begin{itemize}
 \item $g$ is an element of\/ $\SL(m,\real)$ that is either unipotent or
hyperbolic,
 \item  $\mu$ is a probability measure on the projective space\/
$\RP{m-1}$,
 and
 \item $\mu$ is invariant\index{measure!invariant} under~$g$.
 \end{itemize}
 Then $\mu$ is \term[support!of a measure]{supported} on the set of
\index{fixed point}fixed points of~$g$.
 \end{lem}

\begin{proof}
 Let $v$ be any nonzero vector in~$\real^m$.
 For definiteness, let us assume $g$ is
unipotent. (See Exer.~\ref{mu(PV)Pf-HyperEx} for a replacement of this
paragraph in the case where $g$ is hyperbolic.)
 Letting $T = g - \Id$, we
know that $T$ is \index{matrix!nilpotent}nilpotent (because $g$ is
unipotent), so there is some integer $r \ge 0$, such that $v T^r \neq 0$,
but $v T^{r+1} = 0$.
 We have 
 $$v T^r g
 = (v T^r)(\Id + T)
 = v T^r + vT^{r+1}
 = v T^r + 0
 = v T^r ,$$
 so $[v T^r] \in \RP{m-1}$ is a \index{fixed point}fixed point for~$g$.
Also, for $n \in \natural$, we have
 $$ [v] g^n 
  = \left[ \sum_{k=0}^r \binom{n}{k} v T^k \right]
 = \left[ \binom{n}{r}^{-1}
 \sum_{k=0}^r \binom{n}{k} v T^k \right]
 \to [v T^r] $$
 (because, for $k < r$, we have $\binom{n}{k}/\binom{n}{r}
\to 0$ as $n \to \infty$).
 Thus, $[v]g^n$ converges to a \index{fixed point}fixed point of~$g$, as
$n \to \infty$.

The \index{Theorem!Poincar\'e Recurrence}Poincar\'e Recurrence Theorem
\pref{Poincare} implies, for $\mu$-almost every $[v] \in \RP{m-1}$, that
there is a sequence $n_k \to \infty$, such that
$[v] g^{n_k} \to [v]$. On the other hand, we know,
from the preceding paragraph, that $[v] g^{n_k}$
converges to a \index{fixed point}fixed point of~$g$. Thus,
$\mu$-almost every element of $\RP{m-1}$ is a fixed
point of~$g$. In other words, $\mu$ is \term[support!of a measure]{supported} on
the set of \index{fixed point}fixed points of~$g$, as desired.
 \end{proof}

\begin{proof}[{Proof of the Borel Density Theorem~\pref{BDT}}]
 By \index{Theorem!Chevalley}Chevalley's Theorem \pref{ChevalleyFixLine},
there exist 
 \begin{itemize}
 \item a \term[homomorphism!polynomial]{polynomial homomorphism} $\rho
\colon \SL(\ell,\real) \to \SL(m,\real)$, for some~$m$,
 and
 \item a vector $v \in \real^m$,
 \end{itemize}
 such that 
 $\Zar{\Gamma} = \{\, g \in \SL(\ell,\real) \mid v g^\rho \in \real v
\,\}$.
 In other words, letting $[v]$ be the image of~$v$ in
$\RP{m-1}$\index{projective space}, we have 
 \begin{equation} \label{BDTPf-ZarGamma}
 \Zar{\Gamma} = \{\, g \in \SL(\ell,\real) \mid [v] g^\rho = [v] \,\}
. \end{equation}
 Since $\rho(\Gamma)$ fixes~$[v]$, the function $\rho$ induces a
well-defined map $\overline{\rho} \colon \Gamma \backslash G \to
\RP{m-1}$:
 $$ \overline{\rho}(\Gamma g) = [v] g^\rho .$$
 Because $\Gamma$ is a \term{lattice} in~$G$, there is a
$G$-invariant\index{measure!invariant} probability measure~$\mu_0$ on
$\Gamma \backslash G$. The map $\overline{\rho}$ \term[push-forward (of a
measure)]{pushes} this to a probability measure $\mu =
\overline{\rho}_*\mu_0$ on $\RP{m-1}$, defined by $\mu(A) = \mu_0 \bigl(
\overline{\rho}^{-1}(A) \bigr)$ for $A \subset \RP{m-1}$. Because $\mu_0$
is $G$-invariant and $\rho$ is a \index{homomorphism!of Lie
groups}homomorphism, it is easy to see that $\mu$ is $\rho(G)$-invariant.

Let $g$ be any element of~$G$ that is either unipotent or hyperbolic.
From the conclusion of the preceding paragraph, we know that $\mu$~is
$g^\rho$-invariant. Since $g^\rho$ is either unipotent or
hyperbolic \seecor{rho(unip/hyp)}, Lem.~\ref{mu(PV)} implies that $\mu$ is
\term[support!of a measure]{supported} on the set of \index{fixed
point}fixed points of $g^\rho$. Since $[v]$ is obviously in the
\term[support!of a measure]{support} of~$\mu$ \seeex{BDTPfSuppMuEx}, we
conclude that $[v]$ is fixed by $g^\rho$; that is, $[v] g^\rho = [v]$.
From \pref{BDTPf-ZarGamma}, we conclude that $g \in \Zar{\Gamma}$, as
desired.
 \end{proof}

\begin{exercises}

\item \label{SLZarDensinSLREx}
 Show (without using the Borel Density Theorem) that the Zariski closure
of $\SL(\ell,\integer)$\index{SL(l,Z)*$\SL(\ell,\integer)$}
 is $\SL(\ell,\real)$.
 \index{SL(l,Q)*$\SL(\ell,\rational)$}
 \hint{Let $\Gamma = \SL(\ell,\integer)$, and let $H = \Zar{\Gamma}^\circ$.
If $g \in \SL(\ell,\rational)$, then $g^{-1} \Gamma g$ contains a
finite-index subgroup of~$\Gamma$. Therefore $g$
\index{normalizer}normalizes~$H$. Because $\SL(\ell,\rational)$ is
\index{subgroup!dense}dense in $\SL(n,\real)$, this implies that $H$ is a
\index{subgroup!normal}normal subgroup of $\SL(\ell,\real)$. Now apply
Eg.~\fullref{SimpleGrpEgs}{SL}.}

\item \label{LattZarDensinSLREx}
 Use the \term[Theorem!Borel Density]{Borel Density Theorem} to show that
if $\Gamma$ is any \term{lattice} in $\SL(\ell,\real)$, then $\Zar{\Gamma}
= \SL(\ell,\real)$.
 \hint{$\SL(\ell,\real)$ is generated by its unipotent elements.}

\item \label{RPnMfld}
 Show that there is a natural \index{map!covering}covering map from the
$(m-1)$-sphere $S^{m-1}$ onto~$\RP{m-1}$\index{projective space}, so
$\RP{m-1}$ is a $C^\infty$ \index{manifold}manifold.

\item \label{mu(PV)Pf-HyperEx}
 In the notation of Lem.~\ref{mu(PV)}, show that if $g$ is
\index{hyperbolic!element}hyperbolic, and $v$ is any nonzero vector
in~$\real^m$, then $[v] g^n$ converges to a \index{fixed point}fixed point
of~$g$, as $n \to \infty$.
 \hint{Assume $g$ is diagonal. For $v = (v_1,\ldots,v_m)$, calculate $v
g^n$.}

\item \label{BDTPfSuppMuEx}
 In the notation of the proof of Prop.~\ref{BDT}, show that the
\term[support!of a measure]{support} of~$\mu$ is the closure of $[v]
G^\rho$.
 \hint{If some point of $[v] G^\rho$ is contained in an open set of
measure~$0$, then, because $\mu$ is invariant\index{measure!invariant}
under $\rho(G)$, all of $[v] G^\rho$ is contained in an open set of
measure~$0$.}

\item \label{BDTorig}
 (\emphit{The Borel Density Theorem, essentially as stated by Borel})
  Suppose
 \begin{itemize}
 \item $G$ is a \index{subgroup!connected}connected, \index{Lie
group!semisimple}semisimple subgroup of $\SL(\ell,\real)$, such that every
\index{group!simple (or almost)}simple factor of~$G$ is noncompact,
 and
 \item $\Gamma$ is a \term{lattice} in~$G$.
 \end{itemize}
 Show:
 \begin{enumerate}
 \item \label{BDTorig-Zar}
 $G \subset \Zar{\Gamma}$,
 \item \label{BDTorig-subgrp}
 $\Gamma$ is not contained in any proper, closed subgroup of~$G$
that has only finitely many connected
\index{component!connected}components,
 and
 \item \label{BDTorig-lincomb}
 if $\rho \colon G \to \SL(m,\real)$ is any continuous
\index{homomorphism!of Lie groups}homomorphism, then every element of
$\rho(G)$ is a finite linear combination {\rm(}with real
coefficients{\rm)} of elements of~$\rho(\Gamma)$.
 \end{enumerate}
 \hint{Use Prop.~\ref{BDT}. \pref{BDTorig-lincomb}~The subspace of $\Mat_{m
\times m}(\real)$ spanned by $\rho(\Gamma)$ is invariant\index{subspace,
invariant} under multiplication by~$\rho(\Gamma)$, so it must be invariant
under multiplication by $\rho(G)$.}

\item \label{CtblyManyS}
 Suppose
 \begin{itemize}
 \item $G$ is  a closed, \index{subgroup!connected}connected subgroup of
$\SL(\ell,\real)$,
 and
 \item $\Gamma$ is a \term{lattice} in~$G$.
 \end{itemize}
 Show there are only countably many closed,
\index{subgroup!connected}connected subgroups~$S$ of~$G$, such that
 \begin{enumerate}
 \item $\Gamma \cap S$ is a \term{lattice} in~$S$,
 and
 \item there is a one-parameter
\index{unipotent!subgroup!one-parameter}unipotent subgroup $u^t$ of~$S$,
such that $(\Gamma \cap S) \{u^t\}$ is dense in~$S$.
 \end{enumerate}
 \hint{You may assume, without proof, the fact that every \term{lattice}
in every \index{subgroup!connected}connected \term{Lie group} is finitely
generated.
 Show $S \subset \Zar{\Gamma \cap S}$. Conclude that $S$ is uniquely
determined by $\Gamma \cap S$.}

\item \label{NoOrbit->NoSubvar}
 Suppose
 \begin{itemize}
 \item $G$ is an almost-Zariski closed subgroup of $\SL(\ell,\real)$,
 \item $U$ is a \index{subgroup!connected}connected,
\index{unipotent!subgroup}unipotent subgroup of~$G$,
 \item $\Gamma$ is a \index{subgroup!discrete}discrete subgroup of~$G$,
 \item $\mu$ is an ergodic\index{ergodic!measure}
$U$-invariant\index{measure!invariant} probability measure on~$\Gamma
\backslash G$,
 and
 \item
  there does \emph{not} exist a subgroup~$H$ of~$G$, such that 
 \begin{itemize}
 \item $H$ is almost Zariski closed,
 \item $U \subset H$, 
	and 
 \item some $H$-orbit has full measure.
 \end{itemize}
 \end{itemize}
 Show, for all $x \in \Gamma \backslash G$ and every subset~$V$ of~$G$,
that if $\mu(x V) > 0$, then $G \subset \Zar{V}$.
 \hint{Assume $V$ is Zariski closed and irreducible. Let 
 $$ \text{$U_{xV} = \{\, u \in U \mid x V u = x V \,\}$
 and
 $U_V = \{\, u \in U \mid V u = V\,\}$.} $$
 Assuming that $V$ is minimal with $\mu(x V) > 0$, we have 
 $$ \text{$\mu( x V \cap xV u) = 0$ for $u \in U \smallsetminus U_{xV}$.}
$$
 So $U/U_{xV}$ is finite. Since $U$ is
\index{subgroup!connected}connected, then $U_{xV} = U$. Similarly (and
because $\Gamma$ is countable), $U/U_V$ is countable, so $U_V = U$.

Let $\Gamma_V = \{\, \gamma \in \Gamma \mid V \gamma = V \,\}$. Then $\mu$
defines a measure~$\mu_V$ on $\Gamma_V \backslash V$, and this
\index{push-forward (of a measure)}pushes to a measure $\overline{\mu_V}$
on $\Zar{\Gamma_V} \backslash V$. By combining
\index{Theorem!Chevalley}Chevalley's Theorem \pref{ChevalleyFixLine}, the
\index{Theorem!Borel Density}Borel Density Theorem \pref{mu(PV)}, and
the \term[ergodic!action]{ergodicity} of~$U$, conclude that
$\overline{\mu_V}$ is \term[support!of a measure]{supported} on a single
point. Letting $H = \Zar{\langle \Gamma_V, U \rangle}$, some $H$-orbit has
positive measure, and is contained in~$x V$.}

\end{exercises}

\section{Subgroups defined over~$\rational$} \label{Def/QSect}

In this section, we briefly discuss the relationship between \term{lattice}
subgroups and the integer points of a group. This material is not needed
for the proof of \term[Ratner's Theorems]{Ratner's Theorem}, but it is related, and it is
used in many applications, including \term[Theorem!Margulis (on quadratic
forms)]{Margulis' Theorem} on values of quadratic forms \pref{Oppenheim}.

\begin{defn}
 A \term[Zariski!closed]{Zariski closed} subset~$Z$ of $\SL(\ell,\real)$ is
said to be \index{algebraic group!defined over Q*defined
over $\rational$}\emph{defined over~$\rational$} if the defining
polynomials for~$Z$ can be taken to have all of their coefficients
in~$\rational$; that is, if $Z = \Var(\mathcal{Q})$\index{variety} for some
subset~$\mathcal{Q}$ of $\rational[x_{1,1},\ldots,x_{\ell,\ell}]$.
 \end{defn}
 \index{algebraic group!defined over Q*defined over $\rational$|(}

 \begin{eg} \label{EgsOverQ}
 The \index{algebraic group!defined over Q*defined
over $\rational$}algebraic groups in
(\ref{EgsOfAlgGrps-SL(l,R)}--\ref{EgsOfAlgGrps-SL(n,R)}) of
Eg.~\ref{EgsOfAlgGrps} are defined over~$\rational$. Those in
(\ref{EgsOfAlgGrps-Stab(v)}--\ref{EgsOfAlgGrps-SO(Q)}) may or may not be
defined over~$\rational$, depending on the particular choice of $v$, $V$,
or~$Q$. Namely:
 \begin{enumerate} \renewcommand{\theenumi}{\Alph{enumi}}
 \item \label{EgsOverQ-Stab(v)}
 The \index{stabilizer!of a vector}stabilizer
 of a vector~$v$ is defined over~$\rational$ if and
only if $v$ is a scalar multiple of a vector in~$\integer^\ell$
\seeex{EgsOverQ-Stab(v)Ex}. 
 \item \label{EgsOverQ-Stab(V)}
 The \index{stabilizer!of a subspace}stabilizer of a subspace~$V$
of~$\real^\ell$ is defined over~$\rational$ if and only if $V$ is spanned
by vectors in~$\integer^\ell$ \seeex{EgsOverQ-Stab(V)Ex}. 
 \item \label{EgsOverQ-SO(Q)}
 The special orthogonal group \index{SO(Q)*$\SO(Q)$}$\SO(Q)$ of a
nondegenerate quadratic form~$Q$ is defined over~$\rational$ if and only
if $Q$ is a scalar multiple of a form with integer coefficients
\seeex{EgsOverQ-SO(Q)Ex}.
 \end{enumerate}
 \end{eg}

\begin{defn}
 A polynomial function $\phi \colon H \to \SL(n,\real)$ is
\index{homomorphism!defined over Q*defined over~$\rational$}\emph{defined
over~$\rational$} if it can be obtained as in Defn.~\ref{RegularFuncDefn},
but with $\real[x_{1,1}, \ldots,x_{\ell,\ell}]$ replaced by
$\rational[x_{1,1}, \ldots,x_{\ell,\ell}]$ in
\fullref{RegularFuncDefn}{poly}. That is, only polynomials with rational
coefficients are allowed in the construction of~$\phi$.
 \end{defn}

The fact that $\integer^k$ is a \term{lattice} in~$\real^k$ has a vast
generalization:

\begin{thm}[(Borel and Harish-Chandra)] \label{GZLatt}
 Suppose 
 \begin{itemize}
 \item $G$ is a Zariski closed subgroup of\/ $\SL(\ell,\real)$,
 \item $G$ is \index{algebraic group!defined over Q*defined
over $\rational$}defined over~$\rational$,
 and
 \item no nontrivial \term[homomorphism!polynomial]{polynomial
homomorphism} from $G^\circ$ to\/ $\mathbb{D}_2$ is defined
over\/~$\rational$,
 \end{itemize}
 then $G \cap \SL(\ell,\integer)$\index{SL(l,Z)*$\SL(\ell,\integer)$} is a
\term{lattice} in~$G$.
 \end{thm}

\begin{cor} \label{SLZLatt}
 $\SL(\ell,\integer)$\index{SL(l,Z)*$\SL(\ell,\integer)$} is a
\term{lattice} in $\SL(\ell,\real)$.
 \end{cor}

\begin{eg}
 $\mathbb{D}_2 \cap \SL(2,\integer) = \{\pm
\Id\}$\index{SL(2,Z)*$\SL(2,\integer)$} is finite, so it is \emph{not}
a \term{lattice} in $\mathbb{D}_2$.
 \end{eg}

It is interesting to note that Cor.~\ref{SLZLatt} can be proved from
properties of \index{unipotent!flow}unipotent flows. (One can then use this
to obtain the general case of Thm.~\ref{GZLatt}, but this requires some
of the theory of ``\term[group!arithmetic]{arithmetic groups}"
\cfex{GZLattEx}.)

\begin{proof}[Direct proof of Cor.~\ref{SLZLatt}.]
 Let $G = \SL(\ell,\real)$ and $\Gamma =
\SL(\ell,\integer)$\index{SL(l,Z)*$\SL(\ell,\integer)$}. For
 \begin{itemize}
 \item a nontrivial, \term[unipotent!subgroup]{unipotent one-parameter
subgroup}~$u^t$,
 and
 \item a compact subset $K$ of $\Gamma \backslash G$,
 \end{itemize}
 we define $f \colon \Gamma \backslash G \to \real^{\ge 0}$ by 
 $$ f(x) = \liminf_{L \to \infty} \frac{1}{L} \int_0^L \chi_K(x u^t) \,
dt,$$
 where $\chi_K$ is the characteristic function of~$K$.

The key to the proof is that the conclusion of Thm.~\ref{UnipNotNearInfty}
can be proved by using the \index{function!polynomial}polynomial nature
of~$u^t$ --- \emph{without} knowing that $\Gamma$ is a \term{lattice}.
Furthermore, a single compact set~$K$ can be chosen to work for all~$x$ in
any compact subset of $\Gamma \backslash G$. This means that, by
choosing~$K$ appropriately, we may assume that $f > 0$ on some nonempty
open set.

Letting $\mu$ be the \index{measure!Haar}Haar measure on $\Gamma \backslash
G$, we have $\int_{\Gamma \backslash G} f \, d\mu \le \mu(K) < \infty$, so
$f \in L^1(\Gamma \backslash G, \mu)$.

It is easy to see, from the definition, that $f$ is
$u^t$-invariant\index{function!invariant}. Therefore, the
\index{Theorem!Moore Ergodicity}Moore Ergodicity Theorem implies that $f$
is essentially $G$-invariant \seeex{MautnerSubgrpEx}. So $f$ is
\index{function!constant (essentially)}essentially constant. 

If a nonzero constant is in $L^1$, then the space must have finite
measure. So  $\Gamma$ is a \term{lattice}.
 \end{proof}

\begin{exercises}

\item \label{Zpts->def/Q}
 Show that if $C$ is any subset of
$\SL(\ell,\rational)$\index{SL(l,Q)*$\SL(\ell,\rational)$}, then $\Zar{C}$
is \index{algebraic group!defined over Q*defined over $\rational$}defined
over~$\rational$.
 \hint{Suppose $\Zar{C} = \Var(\mathcal{Q})$\index{variety}, for some
$\mathcal{Q} \subset \poly^d$, where $\poly^d$ is the set of polynomials
of degree $\le d$. Because the subspace $\{\, Q \in \poly^d \mid Q(C) = 0
\,\}$ of~$\poly^d$ is defined by linear equations with rational
coefficients, it is spanned by rational vectors.}

\item \label{OverQ<>GalInvt}
 (\emphit{Requires some \term{commutative algebra}}) 
 Let $Z$ be a Zariski closed subset of $\SL(\ell,\real)$. Show that
$Z$ is \index{algebraic group!defined over Q*defined
over $\rational$}defined over~$\rational$ if and only if $\sigma(Z) = Z$,
for every Galois automorphism~$\phi$ of~$\complex$.
 \hint{($\Leftarrow$)~You may assume Hilbert's \term{Nullstellensatz},
which implies there is a subset~$\mathcal{Q}$ of
$\overline{Q}[x_{1,1},\ldots,x_{\ell,\ell}]$, such that $\Zar{C} =
\Var(\mathcal{Q})$\index{variety}, where $\overline{Q}$ is the
\term{algebraic closure} of~$Q$. Then $\overline{Q}$ may be replaced with
some finite \index{Galois!extension}Galois extension~$F$ of~$\rational$,
with \index{Galois!group}Galois group~$\Phi$. For $Q \in \mathcal{Q}$, any
\term{symmetric function} of $\{\, Q^\phi \mid \phi \in \Phi\,\}$ has
rational coefficients.}

\item \label{EgsOverQ-Stab(v)Ex}
 Verify Eg.~\fullref{EgsOverQ}{Stab(v)}.
 \hint{($\Rightarrow$)~The vector~$v$ \index{fixed point}fixed by
$\Stab_{\SL(\ell,\real)}(v)$\index{stabilizer!of a vector} is unique, up to
a scalar multiple. Thus, $v^\phi \in \real v$, for every
\index{Galois!automorphism}Galois automorphism~$\phi$ of~$\complex$.
Assuming some coordinate of~$v$ is rational (and nonzero), then all the
coordinates of~$v$ must be rational.}

\item \label{EgsOverQ-Stab(V)Ex}
 Verify Eg.~\fullref{EgsOverQ}{Stab(V)}.
 \hint{($\Rightarrow$)~Cf.~Hint to Exer.~\ref{EgsOverQ-Stab(v)Ex}. Any
nonzero vector in~$V$ with the minimal number of nonzero coordinates (and
some coordinate rational) must be \index{fixed point}fixed by each
\index{Galois!automorphism}Galois automorphism of~$\complex$. So $V$
contains a rational vector~$v$. By a similar argument, there is a rational
vector that is linearly independent from~$v$. By induction, create a basis
of rational vectors.}

\item \label{EgsOverQ-SO(Q)Ex}
 Verify Eg.~\fullref{EgsOverQ}{SO(Q)}.
 \hint{($\Rightarrow$)~Cf.~Hint to Exer.~\ref{EgsOverQ-Stab(v)Ex}. The
quadratic form~$Q$ is unique, up to a scalar multiple.}

\item \label{Q/Q->SO(Q)/Q}
 Suppose $Q$ is a \term{quadratic form} on~$\real^n$, such that
\index{SO(Q)*$\SO(Q)$}$\Zar{\SO(Q)^\circ}$ is \index{algebraic
group!defined over Q*defined over $\rational$}defined over~$\rational$.
Show that $Q$ is a scalar multiple of a form with integer coefficients.
 \hint{The invariant form corresponding to~$\SO(Q)$ is unique up to a
scalar multiple. We may assume one coefficient is~$1$, so $Q$ is fixed by
every \index{Galois!automorphism}Galois automorphism of~$\complex$.}

\item Suppose
 \begin{itemize}
 \item $G$ is a Zariski closed subgroup of $\SL(\ell,\real)$,
 \item $G^\circ$ is \index{generated by unipotent elements}generated by its
unipotent elements,
 and
 \item $G \cap \SL(\ell,\integer)$\index{SL(l,Z)*$\SL(\ell,\integer)$} is
a \term{lattice} in~$G$.
 \end{itemize}
 Show $G$ is \index{algebraic group!defined over Q*defined
over $\rational$}defined over~$\rational$.
 \hint{Use the Borel Density Theorem \pref{BDT}.}

\item \label{QHomoPts}
 Suppose $\sigma \colon G \to \SL(m,\real)$ is a polynomial
homomorphism\index{homomorphism!defined over Q*defined over~$\rational$}
that is defined over~$\rational$. Show:
 \begin{enumerate}
 \item $\sigma \bigl( G \cap \SL(\ell,\integer) \bigr) \subset
\SL(m,\rational)$
\index{SL(l,Z)*$\SL(\ell,\integer)$}
\index{SL(l,Q)*$\SL(\ell,\rational)$},
 and
 \item \label{QHomoPts-Z}
 there is a finite-index subgroup $\Gamma$ of $G \cap \SL(\ell,\integer)$,
such that $\sigma(\Gamma) \subset \SL(m,\integer)$.
 \end{enumerate}
 \hint{\pref{QHomoPts-Z}~There is a nonzero integer~$k$, such that if $g
\in G \cap \SL(\ell,\integer)$ and $g \equiv \Id \pmod{k}$, then $\sigma(g)
\in \SL(m,\integer)$.}

\item Suppose $G$ is a Zariski closed subgroup of $\SL(\ell,\real)$. Show
that if some nontrivial polynomial homomorphism from $G^\circ$ to
$\mathbb{D}_2$ is defined over~$\rational$\index{homomorphism!defined over
Q*defined over~$\rational$}, then $G \cap
\SL(\ell,\integer)$\index{SL(l,Z)*$\SL(\ell,\integer)$} is \emph{not} a
\term{lattice} in~$G$.

\item Show that if $G$ is a Zariski closed subgroup of $\SL(\ell,\real)$
that is \index{algebraic group!defined over Q*defined
over $\rational$}defined over~$\rational$, and $G^\circ$ is
\index{generated by unipotent elements}generated by its unipotent
elements, then $G \cap
\SL(\ell,\integer)$\index{SL(l,Z)*$\SL(\ell,\integer)$} is a
\term{lattice} in~$G$.

\item \label{GZLattEx}
 Suppose
 \begin{itemize}
 \item $G$ is a \index{subgroup!connected|(}connected, noncompact,
\index{group!simple (or almost)}simple subgroup of $\SL(\ell,\real)$,
 \item $\Gamma = G \cap
\SL(\ell,\integer)$\index{SL(l,Z)*$\SL(\ell,\integer)$},
 and
 \item the natural inclusion $\tau \colon \Gamma \backslash G
\hookrightarrow \SL(\ell,\integer) \backslash \SL(\ell,\real)$, defined by
$\tau(\Gamma x) = \SL(\ell,\integer) x$, is proper.
 \end{itemize}
 Show (without using Thm.~\ref{GZLatt}) that $\Gamma$ is a \term{lattice}
in~$G$.
 \hint{See the proof of Cor.~\ref{SLZLatt}.}

 \end{exercises}
 \index{unipotent!element|)}
 \index{algebraic group!defined over Q*defined over $\rational$|)}
 \index{algebraic group!real|)}
 \index{lattice|)}

\section{Appendix on Lie groups}
 \index{Lie group|(}
 \label{LieGrpSect}

In this section, we briefly recall (without proof) some facts from the
theory of Lie groups.

\begin{defn}
 A group~$G$ is a \defit{Lie group} if the underlying set is a $C^\infty$
\index{manifold}manifold, and the group operations (multiplication and
inversion) are $C^\infty$ functions.
 \end{defn}

A closed subset of a \term{Lie group} need not be a
\index{manifold}manifold (it could be a \index{fractal}Cantor set, for
example), but this phenomenon does not occur for subgroups:

\begin{thm}
 Any closed subgroup of a \term{Lie group} is a Lie group.
 \end{thm}

It is easy to see that the universal cover of a
(\index{subgroup!connected}connected) \term{Lie group} is a Lie group.

\begin{defn}
 Two \index{subgroup!connected}connected \term[Lie group]{Lie groups} $G$
and~$H$ are \defit{locally isomorphic} if their universal covers are
$C^\infty$ isomorphic.
 \end{defn}

We consider only \defit[Lie group!linear]{linear} Lie groups; that is, Lie
groups that are closed subgroups of $\SL(\ell,\real)$, for some~$\ell$.
The following classical theorem shows that, up to local isomorphism, this
results in no loss of generality.

\begin{thm}[(Ado-Iwasawa)]
 Any \index{subgroup!connected}connected \term{Lie group} is \term{locally
isomorphic} to a closed subgroup of\/ $\SL(\ell,\real)$, for some~$\ell$.
 \end{thm}

It is useful to consider subgroups that need not be closed, but may only
be immersed \index{submanifold}submanifolds:

\begin{defn}
 A subgroup~$H$ of a \term{Lie group}~$G$ is a \defit{Lie subgroup} if
there is a Lie group~$H_0$ and an injective $C^\infty$
\index{homomorphism!of Lie groups}homomorphism $\sigma \colon H_0 \to G$,
such that $H = \sigma(H_0)$. Then we consider~$H$ to be a Lie group, by
giving it a topology that makes $\sigma$ a
\index{homeomorphism}homeomorphism. (If $H$ is not closed, this is
\emph{not} the topology that $H$ acquires by being a subset of~$G$.)
 \end{defn}

\begin{defn}
 Let $G$ be a Lie subgroup of $\SL(\ell,\real)$. The \term{tangent space}
to~$G$ at the identity element~$e$ is the \defit{Lie algebra} of~$G$.
 It is, by definition, a vector subspace of the space
 $\Mat_{\ell \times \ell}(\real)$ of $\ell \times \ell$ real matrices.

The Lie algebra of a Lie group $G$, $H$, $U$, $S$, etc., is usually
denoted by the corresponding lower-case gothic letter $\Lie G$, $\Lie H$,
$\Lie U$, $\Lie S$, etc.
 \nindex{$\Mat_{\ell \times \ell}(\real)$ = all real $\ell \times \ell$
matrices}
 \nindex{$\Lie G$, $\Lie H$, $\Lie U$, $\Lie S$, $\Lie L$ = Lie algebra of
$G$, $H$, $U$, $S$, $L$}
 \index{Lie algebra|(}
 \end{defn}

\begin{eg} \label{LieAlgEg} \ 
 \begin{enumerate}
 \item \label{LieAlgEg-U}
 The Lie algebra of $\mathbb{U}_\ell$ is
 $$ \AMSfrak{U}_\ell = \begin{bmatrix}
 0 &  & \vbox to 0pt{\vss\hbox to 0 pt{\hss\Huge $0$}\vss\vss} \\
 & \ddots \\
 \vbox to 0pt{\vss\hbox to 0 pt{\Huge $*$\hss}}& & 0 \\
 \end{bmatrix}, $$
 the space of strictly \index{lower-triangular
matrices!strictly}lower-triangular matrices.
 \item \label{LieAlgEg-SL}
 Let $d \det \colon \Mat_{n
\times n}(\real) \to \real$ be the derivative of the
\index{determinant}determinant map $\det$ at the identity matrix~$\Id$.
Then $(d \det)(A) = \trace A$. Therefore the Lie algebra of
$\SL(\ell,\real)$ is 
 $$ \LieSL(\ell,\real) = \{\, A \in \Mat_{\ell\times \ell}(\real) \mid
\trace A = 0 \,\} .$$
 So the Lie algebra of any Lie subgroup of $\SL(\ell,\real)$ is contained
in $\LieSL(\ell,\real)$.
 \item \label{LieAlgEg-D} The Lie algebra of $\mathbb{D}_\ell$ is
 $$ \AMSfrak{D}_\ell = \bigset{
 \begin{bmatrix}
 a_1 &  & \vbox to 0pt{\vss\hbox to 0 pt{\hss\Huge $0$}\vss\vss} \\
 & \ddots \\
 \vbox to 0pt{\vss\hbox to 0 pt{\Huge $0$\hss}}& & a_\ell \\
 \end{bmatrix}
 }{
 a_1 + \cdots + a_\ell = 0}, $$
 the space of diagonal matrices of trace~$0$. 
 \end{enumerate}
 \end{eg}

\begin{defn} \ 
 \begin{enumerate}
 \item For $\lie x,\lie y \in \Mat_{\ell\times \ell}(\real)$, let 
 \nindex{$[\lie x,\lie y]$ = $\lie x \lie y - \lie y \lie x$ = Lie bracket
of matrices}
 $[\lie x,\lie y] = \lie x \lie y - \lie y \lie x$. This is the \defit{Lie
bracket} of~$\lie x$ and~$\lie y$.
 \item A vector subspace $\Lie H$ of $\Mat_{\ell\times\ell}(\real)$ is a
\defit{Lie subalgebra} if $[\lie x, \lie y] \in V$ for
all $\lie x, \lie y \in \Lie H$.
 \item A linear map $\tau \colon \Lie G \to \Lie H$ between Lie subalgebras
is a \defit[homomorphism!of Lie algebras]{Lie algebra homomorphism} if
$\tau \bigl( [\lie x,\lie y] \bigr) = [ \tau(\lie x), \tau(\lie y)
]$ for all $\lie x, \lie y \in \Lie G$.
 \goodbreak  %% don't want Prop alone at end of page ???
 \end{enumerate}
 \end{defn}

 \begin{prop} \ 
 \begin{enumerate}
 \item The Lie algebra of
any Lie subgroup of\/ $\SL(\ell,\real)$ is a Lie subalgebra.
 \item Any Lie subalgebra~$\Lie H$ of\/ $\LieSL(\ell,\real)$ is the Lie
algebra of a unique \emph{\index{subgroup!connected}connected} Lie
subgroup~$H$ of\/ $\SL(\ell,\real)$.
 \item The differential of a \term{Lie group}
\index{homomorphism!of Lie groups}homomorphism is a Lie algebra
\index{homomorphism!of Lie algebras}homomorphism. That is, if $\phi \colon
G \to H$ is a $C^\infty$ Lie group homomorphism, and $D\phi$ is the
derivative of~$\phi$ at~$e$, then $D\phi$ is a Lie algebra homomorphism
from~$\Lie G$ to~$\Lie H$.
 \item A \index{subgroup!connected}connected \term{Lie group} is uniquely
determined {\rm(}up to local isomorphism{\rm)} by its Lie algebra. That
is, two \index{subgroup!connected}connected Lie groups $G$ and~$H$ are
\term{locally isomorphic} if and only if their Lie algebras are isomorphic.
 \end{enumerate}
 \end{prop}

\begin{defn}
 The \defit[exponential!map]{exponential map} 
 $$\exp \colon \LieSL(\ell,\real) \to \SL(\ell,\real)$$
 is defined by the usual power series
 $$ \exp \lie x = \Id + \lie x + \frac{\lie x^2}{2!} + \frac{\lie x^3}{3!}
+ \cdots$$
 \end{defn}

\begin{eg} \label{sl2brackets}
 Let
 $$ \mbox{$\lie a = \begin{bmatrix} 1 & 0 \\ 0 & -1 \end{bmatrix}$,
 $\lie u = \begin{bmatrix} 0 & 0 \\ 1 & 0 \end{bmatrix}$,
 and
 $\lie v = \begin{bmatrix} 0 & 1 \\ 0 & 0 \end{bmatrix}$.} $$
 Then, letting
 $$ a^s = \begin{bmatrix} e^s & 0 \\ 0 & e^{-s} \end{bmatrix}
 , \ 
  u^t = \begin{bmatrix} 1 & 0 \\ t & 1 \end{bmatrix}
 , \ 
  v^r = \begin{bmatrix} 1 & r \\ 0 & 1 \end{bmatrix}
 , $$
 as usual in $\SL(2,\real)$, it is easy to see that:
 \begin{enumerate}
 \item $\exp(s \lie{a}) = a^s$,
 \item $\exp(t \lie{u}) = u^t$,
 \item $\exp(r \lie{v}) = v^r$,
 \item \label{sl2brackets-(u,a)}
 $[\lie u, \lie a] = 2 \lie u$,
 \item \label{sl2brackets-(v,a)}
 $[\lie v, \lie a] = - 2 \lie v$,
 and
 \item \label{sl2brackets-(v,u)}
 $[\lie v, \lie u] = \lie a$.
 \end{enumerate}
 \end{eg}

\begin{prop} \ 
 Let $\Lie G$ be the Lie algebra of a Lie subgroup~$G$ of\/
$\SL(\ell,\real)$. Then:
 \begin{enumerate}
 \item $\exp \Lie G \subset G$.
 \item For any $\lie g \in \Lie G$, the map $\real \to G$
defined by $g^t = \exp(t \lie g)$ is a
\term[subgroup!one-parameter]{one-parameter subgroup} of~$G$.
 \item The restriction of $\exp$ to some neighborhood of~$0$ in~$\Lie G$
is a \index{diffeomorphism}diffeomorphism onto some neighborhood of~$e$
in~$G$.
 \end{enumerate}
 \end{prop}

\begin{defn} \ 
 \begin{enumerate}
 \item A group~$G$ is \defit[group!solvable]{solvable} if there is
a chain 
 $$e = G_0 \triangleleft G_1 \triangleleft \cdots \triangleleft G_k = G$$
 of subgroups of~$G$, such that, for $1 \le i \le k$,
 \begin{enumerate}
 \item $G_{i-1}$ is a \index{subgroup!normal}normal subgroup of $G_i$,
 and
 \item the quotient group $G_i/G_{i-1}$ is \index{group!abelian}abelian.
 \end{enumerate}
 \item Any \term{Lie group} $G$ has a unique maximal closed,
\index{subgroup!connected}connected, solvable,
\index{subgroup!normal}normal subgroup. This is called the \defit[radical
of a Lie group]{radical} of~$G$, and is denoted $\Rad G$.
 \nindex{$\Rad G$ = radical of the Lie group~$G$}
 \item A \term{Lie group}~$G$ is said to be \defit[Lie
group!semisimple]{semisimple} if $\Rad G = \{e\}$.
 \end{enumerate}
 \end{defn}

\begin{rem}
 According to Defn.~\ref{UnipTorusSSDefn}, $G$ is \index{Lie
group!semisimple}semisimple if $G^\circ$ has no nontrivial,
\index{subgroup!connected}connected, \defit[group!abelian]{abelian},
\index{subgroup!normal}normal subgroups. One can show that this implies
there are also no nontrivial, \index{subgroup!connected}connected,
\defit[group!solvable]{solvable} \index{subgroup!normal}normal subgroups.
 \end{rem}

\begin{thm} \label{G=LRad}
 Any \term{Lie group}~$G$ has a closed, \index{Lie
group!semisimple}semisimple subgroup~$L$, such that 
 \begin{enumerate}
 \item $L$ is semisimple
 and
 \item $G = L \Rad G$.
 \end{enumerate}
 The subgroup~$L$ is called a \defit{Levi subgroup} of~$G$; it is usually
not unique.
 \end{thm}

\begin{warn}
 The above definition is from the theory of \term[Lie group]{Lie groups}.
In the theory of \index{algebraic group!theory of}algebraic groups, the
term \defit{Levi subgroup} is usually used to refer to a slightly different
subgroup --- namely, the subgroup $LT$ of Thm.~\ref{AlgGrpTLU}.
 \end{warn}

\begin{thm}[{(\term[Theorem!Lie-Kolchin]{Lie-Kolchin Theorem})}]
\label{LieKolchin}
 If $G$ is any \index{subgroup!connected}connected,
\index{group!solvable}solvable Lie subgroup of $\SL(\ell,\real)$, then
there exists $h \in
\SL(\ell,\complex)$\index{SL(l,C)*$\SL(\ell,\complex)$}, such that $h^{-1}
G h \subset \mathbb{D}_\ell\mathbb{U}_\ell$.
 \end{thm}
 \index{Lie group!solvable|indsee{group,~solvable}}

\begin{defn}
 Let $\Lie G$ be the Lie algebra of a Lie subgroup~$G$ of
$\SL(\ell,\real)$.
 \begin{itemize}
 \item We use $\GL(\Lie G)$ to denote the group of all
invertible linear transformations $\Lie G \to \Lie G$. This is a \term{Lie
group}, and its Lie algebra $\LieGL(\Lie G)$ consists of all (not
necessarily invertible) linear transformations $\Lie G \to \Lie G$.
 \item We define a group \index{homomorphism!of Lie groups}homomorphism
 \index{representation!adjoint}
 \index{adjoint!representation|indsee{representation,~adjoint}}
 \nindex{$\Ad_G$ = adjoint representation of~$G$ on its Lie algebra}
$\Ad_G \colon G \to \GL(\Lie G)$ by
 $$ \lie x (\Ad_G g) = g^{-1} \lie x g .$$
 Note that $\Ad_G g$ is the derivative at~$e$ of the group automorphism $x
\mapsto g^{-1} x g$, so $\Ad_G$ is a Lie algebra automorphism.
 \item We define a Lie algebra \index{homomorphism!of Lie
algebras}homomorphism $\ad_{\Lie G} \colon \Lie G \to \LieGL(\Lie G)$ by
 \nindex{$\ad_{\Lie G}$ = adjoint representation of the Lie algebra~$\Lie
G$}
 $$ \lie x (\ad_{\Lie G} \lie g) = [\lie x, \lie g] .$$
 We remark that $\ad_{\Lie g}$ is the derivative at~$e$ of~$\Ad_G$.
 \end{itemize}
 \end{defn}

\begin{rem}
 A \term{Lie group}~$G$ is \index{Lie group!unimodular}unimodular (that
is, the right \index{measure!Haar}Haar measure is also
invariant\index{measure!invariant} under left translations) if and only if
$\det \bigl( \Ad_G g \bigr) = 1$, for all $g \in G$.
 \end{rem}

\begin{prop}
 The maps $\exp$, $\Ad_G$ and $\ad_{\Lie G}$ are natural. That is, if
$\rho \colon G \to H$ is a \term{Lie group}
\index{homomorphism!of Lie groups}homomorphism, and $d\rho$ is the
derivative of~$\rho$ at~$e$, then
 \begin{enumerate}
 \item $( \exp \lie g)^\rho = \exp \bigl( d\rho(\lie g) \bigr)$,
 \item $d\rho \bigl( \lie x (\Ad_G g) \bigr) = (d\rho \lie x) \bigl( \Ad_H
g^\rho \bigr)$,
 and
 \item $d\rho \bigl( \lie x (\ad_{\Lie G} \lie g) \bigr) = (d\rho \lie x) \bigl(
\Ad_{\Lie H} d \rho(\lie g) \bigr)$.
 \end{enumerate}
 \end{prop}
 
\begin{cor}
 We have $\Ad_G (\exp \lie g) = \exp( \ad_{\Lie G} \lie g)$. That is,
 $$ \lie x \bigl( \Ad_G (\exp \lie g) \bigr)
 = \lie x + \lie x (\ad_{\Lie G} \lie g) 
 + \frac{1}{2} \lie x (\ad_{\Lie G} \lie g)^2
 + \frac{1}{3!} \lie x (\ad_{\Lie G} \lie g)^3
 + \cdots$$
 \end{cor}

The commutation relations 
 (\ref{sl2brackets-(u,a)},%
 \ref{sl2brackets-(v,a)},\ref{sl2brackets-(v,u)})
 of Eg.~\ref{sl2brackets}
 lead to a complete understanding of all $\LieSL(2,\real)$-modules:

\begin{prop} \label{sl2Rmodules}
 Suppose 
 \begin{itemize}
 \item $\module$ is a finite-dimensional\index{dimension!of a vector
space} real vector space,
 and
 \item $\rho \colon \LieSL(2,\real) \to \LieSL(\module)$ is a Lie algebra
\index{homomorphism!of Lie algebras}homomorphism.
 \end{itemize}
 Then there is a sequence $\lambda_1,\ldots,\lambda_n$ of natural
numbers, and a basis
 $$\bigset{ w_{i,j} 
 }{
 \begin{matrix}
 1 \le i \le n, \\
  0 \le j \le \lambda_i
 \end{matrix}
 }$$
 of\/ $\module$, such that, for all $i,j$, we have:
 \begin{enumerate}
 \item \label{sl2Rmodules-a}
 $w_{i,j} \lie a^\rho = (2j - \lambda_i) w_{i,j}$,
 \item \label{sl2Rmodules-u}
 $w_{i,j} \lie u^\rho = (\lambda_i - j) w_{i,j+1}$,
 and
 \item \label{sl2Rmodules-v}
 $w_{i,j} \lie v^\rho = j w_{i,j-1}$.
 \end{enumerate}
 \end{prop}

\begin{rem} The above proposition has the following immediate consequences.
 \begin{enumerate}
 \item Each $w_{i,j}$ is an \index{eigenvector}eigenvector for $\lie
a^\rho$, and all of the \index{eigenvalue}eigenvalues are integers.
(Therefore, $\lie a^\rho$ is \index{diagonalizable!over
R*over~$\real$}diagonalizable over~$\real$.)
 \item For any integer~$\lambda$, we let 
 $$\module_\lambda = \{\, w \in \module \mid w \lie a^\rho = \lambda w \,\}
.$$
 This is called the \defit{weight space} corresponding to~$\lambda$. (If
$\lambda$ is an \index{eigenvalue}eigenvalue, it is the corresponding
\term{eigenspace}; otherwise, it is $\{0\}$.) A basis of $\module_\lambda$
is given by
 $$ \bigset{ w_{i,(\lambda + \lambda_i)/2} 
 }{
 \begin{matrix}
 1 \le i \le n, \\
  \lambda_i \le |\lambda|
 \end{matrix} 
 }.$$
 \item For all $\lambda$, we have $\module_\lambda \lie u^\rho \subset
\module_{\lambda + 2}$ and $\module_\lambda \lie v^\rho \subset
\module_{\lambda - 2}$.
 \item The kernel of $\lie u^\rho$ is spanned by $\{w_{1,\lambda_1},
\ldots, w_{n, \lambda_n}\}$, and the kernel of $\lie v^\rho$ is spanned by
$\{w_{1,0}, \ldots, w_{n, 0}\}$.
 \item $\lie u^\rho$ and $\lie v^\rho$ are
\index{matrix!nilpotent}nilpotent.
 \end{enumerate}
 \end{rem}

\begin{exercises}

\item \label{BorelIsDU}
 Suppose $u^t$ is a nontrivial, one-parameter,
\index{unipotent!subgroup!one-parameter}unipotent subgroup of
$\SL(2,\real)$\index{SL(2,R)*$\SL(2,\real)$}.
 \begin{enumerate}
 \item Show that $\{u^t\}$ is conjugate to~$\mathbb{U}_2$. 
 \item Suppose $a^s$ is a nontrivial,
\index{subgroup!one-parameter!hyperbolic}one-parameter, hyperbolic subgroup
of $\SL(2,\real)$\index{SL(2,R)*$\SL(2,\real)$} that
\index{normalizer}normalizes~$\{u^t\}$. Show there exists $h \in
\SL(2,\real)$, such that $h^{-1} \{a^s\} h = \mathbb{D}_2$ and $h^{-1}
\{u^t\} h = \mathbb{U}_2$.
 \end{enumerate}

\item Verify the calculations of Eg.~\ref{sl2brackets}.

\item Show that if $a$ is a hyperbolic element of $\SL(\ell,\real)$, and
$V$ is an $a$-invariant\index{subspace, invariant} subspace
of~$\real^\ell$, then $V$~has an $a$-invariant \term[complement (of a
subspace)]{complement}. That is, there is an $a$-invariant subspace~$W$
of~$\real^\ell$, such that $\real^\ell = V \oplus W$.
 \hint{A subspace of~$\real^\ell$ is $a$-invariant if and only if it is
a sum of subspaces of \index{eigenspace}eigenspaces.}

\item Show that if $L$ is a \index{subgroup!Levi}Levi subgroup of~$G$, then
$L \cap \Rad G$ is \index{subgroup!discrete}discrete.
 \hint{$L \cap \Rad G$ is a closed, solvable,
\index{subgroup!normal}normal subgroup of~$L$.}

\item \label{FromRellToR}
 Show that every continuous \index{homomorphism!of Lie groups}homomorphism
$\rho \colon \real^k \to \real$ is a linear map.
 \hint{Every \index{homomorphism!of Lie groups}homomorphism is
$\rational$-linear. Use continuity to show that $\rho$ is $\real$-linear.}

\item \label{RegularCharacter}
 Suppose $T$ is a \index{subgroup!connected|)}connected Lie subgroup of
$\mathbb{D}_\ell$, and $\rho \colon T \to \real^+$ is a $C^\infty$
\index{homomorphism!of Lie groups}homomorphism.
 \begin{enumerate}
 \item \label{RegularCharacter-smooth}
 Show there exist real numbers $\alpha_1,\cdots, \alpha_\ell$, such that
 $$ \rho(t) = t_{1,1}^{\alpha_1} \cdots t_{\ell,\ell}^{\alpha_\ell} $$
 for all $t \in T$.
 \item Show that if $\rho$ is \index{homomorphism!polynomial}polynomial,
then $\alpha_1,\cdots, \alpha_\ell$ are integers.
 \end{enumerate}
 \hint{\pref{RegularCharacter-smooth}~Use Exer.~\ref{FromRellToR}.}

\item \label{HighWtGen}
 Suppose $\module$ and~$\rho$ are as in Prop.~\ref{sl2Rmodules}.
 Show:
 \begin{enumerate}
 \item \label{HighWtGen-gen}
 No proper $\rho \bigl( \LieSL(2,\real)
\bigr)$-invariant\index{subspace, invariant} subspace of~$\module$ contains
$\ker \lie u^\rho$.
 \item \label{HighWtGen-VinW}
 If
 $V$ and~$W$ are $\rho \bigl( \LieSL(2,\real) \bigr)$-invariant
subspaces of~$\module$,
 such that
 $V \subsetneq W$,
 then $V \cap \ker \lie u^\rho \subsetneq W \cap \ker \lie u^\rho$.
 \end{enumerate}
\hint{\pref{HighWtGen-VinW}~Apply \pref{HighWtGen-gen} with $W$ in the
place of~$\module$.}

\item 
 Suppose 
 $$\bigset{ w_{i,j} 
 }{
 \begin{matrix}
 1 \le i \le n, \\
  0 \le j \le \lambda_i
 \end{matrix}
 }$$
 is a basis of a real vector space\/~$\module$, for some sequence
$\lambda_1,\ldots,\lambda_n$ of natural numbers. Show that the equations
\ref{sl2Rmodules}(\ref{sl2Rmodules-a},%
 \ref{sl2Rmodules-u},\ref{sl2Rmodules-v})
 determine linear transformations $\lie a^\rho$, $u^\rho$, and $\lie
v^\rho$ on~$\module$, such that the commutation relations
 (\ref{sl2brackets-(u,a)},\ref{sl2brackets-(v,a)},\ref{sl2brackets-(v,u)})
 of Eg.~\ref{sl2brackets} are satisfied. Thus, there is a Lie algebra
\index{homomorphism!of Lie algebras}homomorphism $\sigma \colon
\LieSL(2,\real) \to \LieSL(\module)$, such that
 $\sigma(\lie a) = \lie a^\rho$, 
  $\sigma(\lie u) = \lie u^\rho$, 
 and
  $\sigma(\lie v) = \lie v^\rho$.

 \end{exercises}
 \index{Lie algebra|)}
 \index{Lie group|)}

\begin{notes}

 The \index{algebraic group!real}algebraic groups that appear in these
lectures are defined over~$\real$, and our only interest is in their real
points.
 Furthermore, we are interested only in \index{Lie group!linear}linear
groups (that is, subgroups of $\SL(\ell,\real)$), not
\term[variety!abelian]{``abelian varieties."} Thus, our definitions and
terminology are tailored to this setting.

 There are many excellent textbooks on the theory of (linear)
\index{algebraic group!theory of}algebraic groups, including
\cite{Borel-AlgicGrps, Humphreys-AlgicGrps}, but they generally focus on
\index{algebraic group!over C*over~$\complex$}algebraic groups
over~$\complex$ (or some other algebraically closed field). The books of
V.~Platonov and A.~Rapinchuk \cite[Chap.~3]{PlatonovRapinchukBook} and
A.~L.~Onishchik and E.~B.~Vinberg \cite{OnishchikVinberg} are excellent
sources for information on \index{algebraic group!real}algebraic groups
over~$\real$.

\notesect{AlgicGrpsSect}

Standard textbooks discuss \term[variety]{varieties},
\term[Zariski!closed]{Zariski closed sets}, \term[algebraic
group!theory of]{algebraic groups}, \index{dimension!of a Zariski closed
set}dimension, and the \term{singular set}.

Whitney's Theorem \pref{Zar->AlmConn} appears in \cite[Cor.~1 of
Thm.~3.6, p.~121]{PlatonovRapinchukBook}.

%See \cite[\S AG.3]{Humphreys-AlgicGrps} for a discussion of the dimension
%of a Zariski closed set.

\notesect{ZarClosureSect}

The \term[Zariski!closure]{Zariski closure} is a standard topic.

The notion of being ``\term[Zariski!closed!almost]{almost Zariski closed}"
does not arise over an algebraically closed field, so it is not described
in most texts. Relevant material (though without using this terminology)
appears in \cite[\S3.2]{PlatonovRapinchukBook} and
\cite[\S3.1]{ZimmerBook}. References to numerous specific results on
almost-Zariski closed subgroups can be found in
\cite[\S3]{Witte-superlatt}.

Exercise~\ref{CentClosed} is a version of \cite[Prop.~8.2b,
p.~59]{Humphreys-AlgicGrps}.

\notesect{JordanDecSect}

\term[function!polynomial]{Polynomials}, \term[unipotent!element]{unipotent
elements}, and the \term[Jordan!decomposition]{Jordan decomposition} are
standard material. However, most texts consider the
\term[Jordan!decomposition]{Jordan decomposition} over~$\complex$,
not~$\real$. (\term[hyperbolic!element]{Hyperbolic elements} and
\term[elliptic element]{elliptic elements} are lumped together into a
single class of \index{semisimple element}{``semisimple" elements}.)

The \term[Jordan!decomposition]{real Jordan decomposition} appears in
\cite[Lem.~IX.7.1, p.~430]{Helgason-DGLGSS}, for example.

A solution of Exer.~\ref{norm-closed} appears in the proof of
\cite[Thm.~3.2.5, p.~42]{ZimmerBook}.

\notesect{StructureAlgGrpSect}

This material is standard, except for Thm.~\ref{Gcpct} and its corollary
(which do not occur over an algebraically closed field). 

The theory of \term{roots and weights} is described in many textbooks,
including \cite{Helgason-DGLGSS, Humphreys-LieAlgs, Serre-LieAlgs}. See
\cite[Table~V, p.~518]{Helgason-DGLGSS} for a list of the almost-simple
groups.

For the case of \index{Lie group!semisimple}semisimple groups, the
difficult direction ($\Leftarrow$) of Thm.~\ref{Gcpct} is immediate from
the ``\term{Iwasawa decomposition}" $G = K A N$, where $K$ is \index{Lie
group!compact}compact, $A$~is a hyperbolic \index{torus!hyperbolic}torus,
and $N$~is unipotent. This decomposition appears in \cite[Thm.~3.9,
p.~131]{PlatonovRapinchukBook}, or in many texts on \term[Lie group]{Lie
groups}.

The proof of \term[Theorem!Engel]{Engel's Theorem} in
Exer.~\ref{Unip->Lower} is taken from \cite[Thm.~17.5,
p.~112]{Humphreys-AlgicGrps}. The theorem of
\term[Theorem!Burnside]{Burnside} mentioned there (or the more general
\term[Theorem!Jacobson Density]{Jacobson Density Theorem}) appears in
graduate algebra texts, such as \cite[Cor.~3.4 of
Chap.~XVII]{Lang-Algebra}.

\notesect{ChevalleySect}

This is standard.

%\term[Theorem!Chevalley]{Chevalley's Theorem} \pref{ChevalleyFixLine}
%appears in \cite[Thm.~11.2]{Humphreys-AlgicGrps}, for example.

\notesect{AlmZarSect}

These results are well known, but do not appear in most texts on
\index{algebraic group!theory of}algebraic groups.

Proposition~\ref{Cpct->ZarClosed} is due to C.~Chevalley \cite[Prop.~2, \S
VI.5.2, p.~230]{Chevalley-Lie3}. A proof also appears in
\cite[\S8.6]{Bailey-AutoForms}.

See \cite[Thm.~8.1.1, p.~107]{Hochschild-AlgicGrps} for a proof of
Prop.~\ref{Unip->Closed}.

See \cite[Thm.~8.3.2, p.~112]{Hochschild-AlgicGrps} for a proof of
Thm.~\ref{SemiSimp->ZarClosed}.

The analogue of Prop.~\ref{AB-closed} over an algebraically closed field is
standard (e.g., \cite[Cor.~7.4, p.~54]{Humphreys-AlgicGrps}).  For a
derivation of Prop.~\ref{AB-closed} from this, see
\cite[Lem.~3.17]{Witte-superlatt}.

See \cite[Cor.~1 of Prop.~3.3, p.~113]{PlatonovRapinchukBook} for a proof
of Cor.~\ref{homo-img}.

Corollary~\ref{G'almZar} is proved in \cite[Thm.~15, \S II.14,
p.~177]{Chevalley-Lie2} and \cite[Thm8.3.3, p.~113]{Hochschild-AlgicGrps}.

Completing the proof of Cor.~\ref{ZarG'=G'} requires one to know that
$\Zar{G}/G$ is \index{group!abelian}abelian for every
\index{subgroup!connected}connected, \index{Lie
group!semisimple}semisimple subgroup~$G$ of $\SL(\ell,\real)$. In fact,
$\Zar{G}/G$ is trivial if $\Zar{G}$ is ``\term{simply connected}" as an
\index{algebraic group!simply connected}algebraic group \cite[Prop.~7.6,
p.~407]{PlatonovRapinchukBook}, and the general case follows from this by
using an exact sequence of \index{Galois!cohomology}Galois cohomology
groups: $\widetilde G_\real \to (\widetilde G/Z)_\real \to
\coho^1(\complex/\real, Z_\complex)$.

A proof of Lem.~\ref{poly+U=algic} appears in
\cite[\S2.2]{BorelPrasad-Sintegral}. (It is based on the analogous result
over an algebraically closed field, which is a standard result that appears
in \cite[Prop.~7.5, p.~55]{Humphreys-AlgicGrps}, for example.)

Exercise~\ref{unip-sc} is a version of \cite[Thm.~8.1.1,
p.~107]{Hochschild-Lie}.

\notesect{BDTSect}

This material is fairly standard in \term[ergodic!theory]{ergodic theory},
but not common in texts on \index{algebraic group!theory of}algebraic
groups.

The \index{Theorem!Borel Density}Borel Density Theorem \pref{BDT} was
proved for \index{Lie group!semisimple}semisimple groups in
\cite{Borel-BDT} \seeex{BDTorig}. (The theorem also appears in
\cite[Lem.~II.2.3 and Cor.~II.2.6, p.~84]{MargulisBook} and
\cite[Thm.~3.2.5, pp.~41--42]{ZimmerBook}, for example.)
 The generalization to all \term[Lie group]{Lie groups} is due to
S.~G.~Dani \cite[Cor.~2.6]{Dani-Quasi}. 

The \term[Theorem!Poincar\'e Recurrence]{Poincar\'e Recurrence Theorem}
\pref{Poincare} can be found in many textbooks on
\term[ergodic!theory]{ergodic theory}, including \cite[Cor.~I.1.8,
p.~8]{BekkaMayer}.

See \cite[Thm.~4.10, p.~205]{PlatonovRapinchukBook} for a solution of
Exer.~\ref{SLZarDensinSLREx}.

Exercise~\ref{CtblyManyS} is \cite[Cor.~A(2)]{RatnerOrbit}.

Exercise~\ref{NoOrbit->NoSubvar} is
\cite[Prop.~3.2]{MargulisTomanov-Ratner}.

\notesect{Def/QSect}

This material is standard in the theory of
``\term[group!arithmetic]{arithmetic groups}." (If $G$ is defined
over~$\rational$, then $G \cap
\SL(\ell,\integer)$\index{SL(l,Z)*$\SL(\ell,\integer)$} is said to be an
\defit[group!arithmetic]{arithmetic group}.) The book of Platonov and
Rapinchuk \cite{PlatonovRapinchukBook} is an excellent reference on the
subject. See \cite{Morris-ArithGrps} for an introduction. There are also
numerous other books and survey papers.\index{SL(l,R)*$\SL(\ell,\real)$|)}

Theorem~\ref{GZLatt} is due to A.~Borel and Harish-Chandra
\cite{BorelHarishChandra-GZLatt}. (Many special cases had previously been
treated by C.~L.~Siegel \cite{Siegel-Einheiten}.) Expositions can also be
found in \cite[Cor.~13.2]{Borel-ArithGrps} and
\cite[Thm.~4.13]{PlatonovRapinchukBook}. (A proof of only
Cor.~\ref{SLZLatt} appears in \cite[\S V.2]{BekkaMayer}.) These are based
on the reduction theory for \term[group!arithmetic]{arithmetic groups},
not unipotent flows.

The observation that Thm.~\ref{GZLatt} can be obtained from a variation of
Thm.~\ref{UnipNotNearInfty} is due to G.~A.~Margulis
\cite[Rem.~3.12(II)]{Margulis-LieGrpErgThy}.

\notesect{LieGrpSect}

 There are many textbooks on \term[Lie group]{Lie groups}, including
\cite{Helgason-DGLGSS, Hochschild-Lie, Varadarajan}. The expository
article of R.~Howe \cite{Howe-VeryBasic} provides an elementary
introduction.
 \index{Zariski!closed!almost|)}
  \index{Zariski!closed|)}
 \index{Zariski!closure|)}
 \index{hyperbolic!element|)}

 \index{function!polynomial|)}
 \end{notes}

%% file: RatnerProof.tex
\mychapter{Proof of the Measure-Classification Theorem}
\label{ProofChap}

% reset the equation counter manually (since no \section)
\setcounter{equation}{0}

In this chapter, we present the main ideas in a
proof of the following theorem. The reader is assumed to be familiar with
the concepts presented in Chap.~\ref{IntroChap}.

\begin{thm}[(Ratner)] \label{Thm}
 \index{Ratner's Theorems!Measure Classification}
 If
 \begin{itemize}
 \item $G$ is a closed, \index{subgroup!connected}connected subgroup of\/
$\SL(\ell,\real)$\index{SL(l,R)*$\SL(\ell,\real)$}, for some~$\ell$,
 \item $\Gamma$ is a \index{subgroup!discrete}discrete subgroup of~$G$,
 \item $u^t$ is a \term[unipotent!subgroup]{unipotent} one-parameter
subgroup of~$G$,
 and
 \item $\mu$ is an \term[ergodic!measure]{ergodic}
$u^t$-invariant\index{measure!invariant} probability measure on\/ $\Gamma
\backslash G$,
 \end{itemize}
 then $\mu$~is \defit[measure!homogeneous]{homogeneous}. 

More precisely, there exist
 \begin{itemize}
 \item a closed, \index{subgroup!connected}connected subgroup~$S$ of~$G$, 
 and
 \item a point~$x$ in\/ $\Gamma \backslash G$,
 \end{itemize}
 such that
 \begin{enumerate}
 \item $\mu$ is $S$-invariant\index{measure!invariant},
 and
 \item $\mu$ is \term[support!of a measure]{supported} on the orbit $xS$.
 \end{enumerate}

 \end{thm}
 \index{measure!invariant|(}
 \nindex{$\mu$ = ergodic $u^t$-invariant probability measure on $\Gamma
\backslash G$}

\begin{rem} \label{Haar}
 If we write $x = \Gamma g$, for some $g \in G$, and let $\Gamma_S =
(g^{-1} \Gamma g) \cap S$, then the conclusions imply that
 \begin{enumerate}
 \item under the natural identification of the orbit~$x S$ with
the \index{homogeneous!space}homogeneous space $\Gamma_S  \backslash S$,
the measure~$\mu$ is the \index{measure!Haar}Haar measure on $\Gamma_S 
\backslash S$, 
 \item $\Gamma_S$ is a \term{lattice} in~$S$,
 and
 \item $x S$ is closed
 \end{enumerate}
 \seeex{HaarEx}.
 \end{rem}

\begin{assump} \label{RatnerPfAssumeA}
 Later \seeassump{HyperAssump}, in order to simplify the details of
the proof while losing very few of the main ideas, we will make the
additional assumption that
 \begin{enumerate}
 \item \label{RatnerPfAssumeA-A}
 $\mu$ is invariant\index{measure!invariant} under a 
\term[subgroup!one-parameter!hyperbolic]{hyperbolic one-parameter
subgroup} $\{a^s\}$ that \index{normalizer}normalizes~$u^t$,
 and
 \item \label{RatnerPfAssumeA-SL2}
 $\langle a^s, u^t \rangle$ is contained in a subgroup
 $L = \langle u^t, a^s, v^r \rangle$ 
 that is \term{locally isomorphic} to
$\SL(2,\real)$\index{SL(2,R)*$\SL(2,\real)$}.
 \end{enumerate}
 See \S\ref{RatnerPfAssumeASect} for a discussion of the changes involved
in removing this hypothesis. The basic idea is that
Prop.~\ref{RatMeas-S=U} shows that we may assume \index{stabilizer!of a
measure}$\Stab_G(\mu)$ contains
a \term[subgroup!one-parameter]{one-parameter subgroup} that is not
unipotent. A more sophisticated version of this argument, using the
theory of \index{algebraic group!theory of}algebraic groups, shows that
slightly weakened forms of \pref{RatnerPfAssumeA-A} and
\pref{RatnerPfAssumeA-SL2} are true.  Making these assumptions from the
start simplifies a lot of the algebra, without losing \emph{any} of the
significant ideas from dynamics.
 \end{assump}

\begin{rem}
 Note that $G$ is \emph{not} assumed to be \index{Lie
group!semisimple}semisimple. Although the semisimple case is the most
interesting, we allow ourselves more freedom, principally because the
proof relies (at one point, in the proof of Thm.~\ref{SuppS}) on
induction on $\dim G$, and this induction is based on knowing the result
for all \index{subgroup!connected}connected subgroups, not only the
semisimple ones.
 \end{rem}

\begin{rem} \label{MayAssZar}
 There is no harm in assuming that $G$ is almost
\index{Zariski!closed!almost}Zariski closed \seeex{MayAssZarEx}. This
provides a slight simplification in a couple of places (see
Exer.~\ref{GminG0Uopen} and the proof of Thm.~\ref{SuppS}).
 \end{rem}

\begin{exercises}

\item \label{HaarEx}
 Prove the assertions of Rem.~\ref{Haar} from the conclusions of
Thm.~\ref{Thm}.

\item Show that Thm.~\ref{Thm} remains true without the assumption that
$G$ is \index{subgroup!connected}connected.
 \hint{$\mu$ must be \term[support!of a measure]{supported} on a single connected
\index{component!connected}component of $\Gamma \backslash G$. Apply
Thm.~\ref{Thm} with $G^\circ$ in the place of~$G$.}

\item \label{MayAssZarEx}
 Assume Thm.~\ref{Thm} is true under the additional hypothesis that $G$ is
\term[Zariski!closed!almost]{almost Zariski closed}. Prove that this
additional hypothesis can be eliminated.
 \index{SL(l,R)*$\SL(\ell,\real)$}
 \hint{$\Gamma \backslash G$ embeds in $\Gamma \backslash
\SL(\ell,\real)$.}

 \end{exercises}

\section{An outline of the proof} \label{OutlinePfSect}

Here are the main steps in the proof.
 \begin{enumerate}

 \item \emphit{Notation.}

 \begin{itemize}

 \item Let
 $S = \Stab_G(\mu)$\index{stabilizer!of a
measure}. We wish to show that $\mu$ is
\term[support!of a measure]{supported} on a single $S$-orbit.

 \item Let $\Lie G$ be the \term{Lie algebra} of~$G$ and $\Lie S$ be the
Lie algebra of~$S$.

 \item The expanding and contracting subspaces of~$a^s$ (for $s > 0$)
provide decompositions 
 $$ \text{$\Lie G = \Lie G_- + \Lie G_0 + \Lie G_+$ 
 \quad and \quad
 $\Lie S = \Lie S_- + \Lie S_0 + \Lie S_+$,} $$
 and we have corresponding subgroups $G_-$, $G_0$, $G_+$, $S_-$, $S_0$,
and~$S_+$ \seeNot{Notation}.

\item For convenience, let $U = S_+$. Note that $U$ is unipotent, and we
may assume $\{u^t\} \subset U$, so $\mu$ is
\term[ergodic!measure]{ergodic} for~$U$.

 \end{itemize}

 \item \label{Out-BasLem}
 We are interested in \defit[transverse divergence]{transverse}
divergence of nearby orbits. (We ignore relative motion \emphit{along}
the $U$-orbits, and project to $G \ominus U$.)
 The \term{shearing property} of unipotent flows implies, for a.e.\ 
 $x,y \in \Gamma \backslash G$, that if $x \approx y$, then the
transverse divergence of the $U$-orbits through $x$ and~$y$ is
fastest along some direction in~$S$ \seeprop{BasLemStrong}. Therefore, the
direction belongs to $G_- G_0$ \seecor{DivMinus}.

 \item We define a certain subgroup
 $$ \widetilde{S}_- = \{\, g \in G_- \mid \forall u \in U, \ u^{-1} g u
\in G_- G_0 U \,\} $$
 of~$G_-$ \cfdefn{StildeDefn}. Note that $S_- \subset \widetilde{S}_-$.

The motivation for this definition is that if $y \in x \widetilde{S}_-$,
then all of the \index{transverse divergence}transverse divergence belongs
to $G_- G_0$ --- there is no \index{component!G+*$G_+$-}$G_+$-component to
any of the transverse divergence. For clarity, we emphasize that this
restriction applies to \emph{all} transverse divergence, not only the
\emph{fastest} transverse divergence.

 \item \label{Out-inSmin}
 Combining \pref{Out-BasLem} with the dilation provided by the
translation~$a^{-s}$ shows, for a.e.\ $x,y \in \Gamma \backslash G$, that
if $y \in x G_-$, then $y \in x \widetilde{S}_-$ \seecor{GminInSmin}.

 \item \label{Out-NotS->G+Div}
 A Lie algebra calculation shows that if $y \approx x$, and $y = xg$,
with $g \in (G_- \ominus \widetilde{S}_-) G_0 G_+$, then the
\term{transverse divergence} of the $U$-orbits through $x$ and~$y$ is
fastest along some direction in~$G_+$ \seelem{NotS->G+Div}.

 \item \label{Out-NotSperp}
 Because the conclusions of \pref{Out-BasLem} and~\pref{Out-NotS->G+Div}
are contradictory, we see, for a.e.\ $x,y \in \Gamma \backslash G$, that 
 $$ \mbox{if $x \approx y$, then $y \notin x (G_- \ominus
\widetilde{S}_-) G_0 G_+$} $$
 \cfcor{MustS}. (Actually, a technical problem causes us obtain this
result only for $x$ and~$y$ in a set of measure $1 - \epsilon$.)

 \item \label{Out-Stilde=S}
 The relation between \index{entropy!vs.\ stretching}stretching and
entropy (Prop.~\ref{EntropyLemma}) provides bounds on the entropy
of~$a^s$, in terms of the the \index{Jacobian}Jacobian of~$a^s$ on~$U$
and (using \pref{Out-inSmin}) the \index{Jacobian}Jacobian of~$a^{-s}$
on~$\widetilde{S}_-$:
 $$ J(a^s, U) \le h_\mu(a^s) \le J \bigl( a^{-s}, \widetilde{S}_- \bigr)
.$$
 On the other hand, the structure of
$\LieSL(2,\real)$-modules\index{module} implies that
 $J(a^s, U) \ge J \bigl( a^{-s}, \widetilde{S}_- \bigr)$. 
 Thus, we conclude that $h_\mu(a^s) = J(a^{-s}, \widetilde{S}_-)$.  This
implies that $\widetilde{S}_- \subset \Stab_G(\mu)$\index{stabilizer!of a
measure}, so we must have
$\widetilde{S}_- = S_-$ \seeprop{Stilde=S}.

 \item \label{Out-SuppGU}
 By combining the conclusions of \pref{Out-NotSperp} and
\pref{Out-Stilde=S}, we show that $\mu(x S_- G_0 G_+)
> 0$, for some $x \in \Gamma \backslash G$ \seeprop{SuppGU}.

 \item By combining \pref{Out-SuppGU} with the (harmless) assumption
that $\mu$ is not \term[support!of a measure]{supported} on an orbit of any closed, proper
subgroup of~$G$, we show that $S_- = G_-$ (so $S_-$ is horospherical),
and then there are a number of ways to show that $S = G$ \seethm{SuppS}.
 \end{enumerate}

The following several sections expand this outline into a fairly complete
proof, modulo some details that are postponed to
\S\ref{PreciseSect}.

\section{Shearing and polynomial divergence} \label{ShearPolyPfSect}

 As we saw in Chap.~\ref{IntroChap}, \index{shearing property}shearing and
\index{polynomial!divergence}polynomial divergence are crucial ingredients
of the proof of Thm.~\ref{Thm}. Precise statements will be given in
\S\ref{PreciseSect}, but let us now describe them informally. Our goal
here is to prove that the direction of fastest divergence usually belongs
to the \index{stabilizer!of a measure}stabilizer of~$\mu$
(see Prop.~\ref{BasLem}, which follows Cor.~\ref{UnipCase}). This will
later be restated in a slightly more convenient (but weaker) form
\seecor{DivMinus}.

\begin{lem}[(Shearing)] \label{RatTransInN}
 If $U$ is any \index{subgroup!connected|(}connected, unipotent
subgroup\index{unipotent!subgroup} of~$G$, then the  \term{transverse
divergence} of any two nearby $U$-orbits is fastest along some direction
that is in the \index{normalizer|(}normalizer $N_G(U)$.
 \end{lem}

\begin{lem}[(Polynomial divergence)] \label{PolyDiv}
 If $U$ is a \index{subgroup!connected}connected,
\index{unipotent!subgroup}unipotent subgroup of~$G$, then any two nearby
$U$-orbits diverge at \index{polynomial!speed|indsee{speed,~polynomial}}
 \index{speed!polynomial}polynomial speed.

Hence, if it takes a certain amount of time for two nearby $U$-orbits to
diverge to a certain distance, then the amount {\upshape(}and
direction{\upshape)} of divergence will remain approximately the same for
a proportional length of time.
 \end{lem}

By combining these two results we will establish the following conclusion
\cfcor{FastestInStab}. It is the basis of the entire proof.

\begin{notation} \label{SDefn}
 Let $S = \Stab_G(\mu)^{\circ}$\index{stabilizer!of a
measure}. This is a closed subgroup of~$G$
\seeex{Stab(mu)Closed}.
 \nindex{$S$ = $\Stab_G(\mu)^{\circ}$}
 \end{notation}

\begin{prop} \label{BasLemStrong}
 If $U$ is any \index{subgroup!connected}connected,
\term[ergodic!action]{ergodic}, \index{unipotent!subgroup}unipotent
subgroup of~$S$, then there is a conull subset~$\Omega$ of $\Gamma
\backslash G$, such that, for all $x, y \in \Omega$, with $x \approx y$,
the $U$-orbits through $x$ and~$y$ \index{transverse divergence}diverge
fastest along some direction that belongs to~$S$.
 \end{prop}

This immediately implies the following interesting special case of
Ratner's Theorem \seeex{UnipCaseEx}, which was proved rather informally
in Chap.~\ref{IntroChap} \seeprop{RatMeas-S=U}.

\begin{cor} \label{UnipCase}
 If $U = \Stab_G(\mu)$\index{stabilizer!of a
measure} is \index{unipotent!subgroup}unipotent {\upshape(}and
\index{subgroup!connected}connected\/{\upshape)}, then $\mu$~is
\term[support!of a measure]{supported} on a single $U$-orbit.
 \end{cor}

Although Prop.~\ref{BasLemStrong} is true \seeex{BasLemTrueEx},
it seems to be very difficult to prove from scratch, so we will be
content with proving the following weaker version that does not yield a
conull subset, and imposes a restriction on the relation between $x$
and~$y$ \see{Cond*}.  (See Exer.~\ref{BasLemClassical} for
a \term{non-infinitesimal version} of the result.)

\begin{thmref}{BasLemStrong}
 \begin{prop} \label{BasLem}
 \label{BASLEM} %% to eliminate a Latex error
 For any
 \begin{itemize}
 \item \index{subgroup!connected}connected,
\term[ergodic!action]{ergodic}, \index{unipotent!subgroup}unipotent
subgroup~$U$ of~$S$, 
 and
 \item any $\epsilon > 0$,
 \end{itemize}
 there is a subset~$\Omega_\epsilon$ of\/ $\Gamma \backslash G$, such that
 \begin{enumerate}
 \item $\mu(\Omega_\epsilon) > 1 - \epsilon$,
 and 
 \item for all $x, y \in \Omega_\epsilon$, with $x \approx y$, and such
that a certain \index{technical assumption}technical assumption
\pref{TechAssump} is satisfied, the fastest \index{transverse
divergence}transverse divergence of the $U$-orbits through $x$ and~$y$ is
along some direction that belongs to~$S$.
 \end{enumerate}
 \end{prop}
 \end{thmref}

\begin{proof}[{Proof {\rm \cfcor{FastestInStab1D}}}]
 Let us assume that no \index{normalizer}$N_G(U)$-orbit has positive
measure, for otherwise it is easy to complete the proof
\cfex{NormalSubgrpEx}. Then, for a.e.~$x \in \Gamma \backslash G$, there
is a point $y \approx x$, such that 
 \begin{enumerate}
 \item $y \notin x \, N_G(U)$, 
 and
 \item $y$~is a \term[generic!point]{generic} point for~$\mu$
\seecor{GenericPtae}.
 \end{enumerate}

Because $y \notin x \, N_G(U)$, we know that the orbit $yU$ is not
parallel to~$x U$, so they diverge from each other. From
Lem.~\ref{RatTransInN}, we know that the direction of fastest
\index{transverse divergence}transverse divergence belongs to
$N_G(U)$\index{normalizer}, so there exist $u,u' \in U$, and $c \in N_G(U)
\ominus U$, such that 
 \begin{itemize}
 \item $yu' \approx (xu) c$,
 and
 \item $\|c\| \asymp 1$ (i.e., $\|c\|$ is finite, but not
\index{infinitesimal}infinitesimal).
 \end{itemize}
 Because $c \notin U = \Stab_G(\mu)$\index{stabilizer!of a
measure}, we know that $c_* \mu \neq \mu$.
Because $c \in N_G(U)$\index{normalizer}, this implies $c_* \mu \perp \mu$
\seeex{N(U)->Singular}, so there is a compact subset~$K$ with
$\mu(K) > 1 - \epsilon$ and $K \cap Kc = \emptyset$
\seeex{Sing->K}. 

We would like to complete the proof by saying that there are values
of~$u$ for which both of the two points $xu$ and~$yu'$ are arbitrarily
close to~$K$, which contradicts the fact that $d(K,Kc) > 0$. However,
there are two technical problems:
 \begin{enumerate}
 \item The set $K$~must be chosen before we know the value of~$c$. This
issue is handled by Lem.~\ref{KnotKc}. 
 \item The \index{Theorem!Pointwise Ergodic}Pointwise Ergodic Theorem
\pref{PtwiseErg} implies (for a.e.~$x$) that $xu$ is arbitrarily close
to~$K$ a huge proportion of the time. But this theorem does not apply
directly to~$y u'$, because $u'$ is a nontrivial function of~$u$. To
overcome this difficulty, we add an additional \index{technical
assumption}technical hypothesis on the element~$g$ with $y = xg$
\see{Cond*}. With this assumption, the result can be proved
\see{BasLemPf}, by showing that the \index{Jacobian}Jacobian of the
change of variables $u \mapsto u'$ is bounded above and below on some set
of reasonable size, and applying the uniform approximate version of the
\index{Theorem!Pointwise Ergodic!Uniform}Pointwise Ergodic Theorem
\seecor{UnifErg}. The uniform estimate is what requires us to restrict to
a set of measure $1 - \epsilon$, rather than a conull set.
 \qedhere
 \end{enumerate}
 \end{proof}

\begin{rem} \ 
 \begin{enumerate}
 \item The fact that $\Omega_\epsilon$ is not quite conull is not a
serious problem, although it does make one part of the proof more
complicated \cfprop{SuppGU}.
 \item We will apply Prop.~\ref{BasLem} only twice
(in the proofs of Cors.~\ref{GminInSmin} and~\ref{MustS}). In each case,
it is not difficult to verify that the \index{technical
assumption}technical assumption is satisfied
\seeexs{GminInSminEx}{NotS->G+DivTechEx}.
 \end{enumerate}
 \end{rem}

\begin{exercises}

\item \label{Stab(mu)Closed}
 Show that $\Stab_G(\mu)$\index{stabilizer!of a measure} is a closed
subgroup of~$G$.
 \hint{$g \in \Stab_G(\mu)$ if and only if $\int f(xg) \, d\mu(x) = \int
f \, d\mu$ for all continuous functions~$f$ with compact
\term[support!of a function]{support}.}

\item \label{NormUnimodEx}
 Suppose 
 \begin{itemize}
 \item $\nu$ is a (finite or infinite) Borel measure on~$G$,
 and
 \item $N$ is a \index{Lie group!unimodular}unimodular,
\index{subgroup!normal}normal subgroup of~$G$.
 \end{itemize}
 Show that if $\nu$~is right-invariant\index{measure!invariant} under~$N$
(that is, $\nu(A n) = \nu(A)$ for all $n \in N$), then $\nu$ is
left-invariant\index{measure!invariant} under~$N$.

\item \label{NormalSubgrpEx}
 Show that if 
 \begin{itemize}
 \item $N$ is a \index{Lie group!unimodular}unimodular,
\index{subgroup!normal}normal subgroup of~$G$,
 \item $N$ is contained in $\Stab_G(\mu)$\index{stabilizer!of a
measure},
 and
 \item $N$ is \term[ergodic!action]{ergodic} on $\Gamma \backslash G$,
 \end{itemize}
 then $\mu$ is \term[measure!homogeneous]{homogeneous}.
 \hint{Lift $\mu$ to an (infinite) measure~$\hat\mu$ on~$G$, such that
$\hat\mu$~is left invariant\index{measure!invariant} under~$\Gamma$, and
right invariant\index{measure!invariant} under~$N$.
Exercise~\ref{NormUnimodEx} implies that $\hat\mu$~is left invariant (and
\term[ergodic!measure]{ergodic}) under the closure~$H$ of $\Gamma N$.
Ergodicity implies that $\hat\mu$~is \term[support!of a
measure]{supported} on a single $H$-orbit.}

\item \label{UnipCaseEx}
 Prove Cor.~\ref{UnipCase} from Prop.~\ref{BasLemStrong} and
Exer.~\ref{NormalSubgrpEx}.
 \hint{If $\mu \bigl( x N_G(U) \bigr) > 0$\index{normalizer}, for some $x
\in \Gamma \backslash G$, then Exer.~\ref{NormalSubgrpEx} (with $N_G(U)$
in the place of~$G$) implies that $\mu$~is
\term[measure!homogeneous]{homogeneous}. Otherwise,
Prop.~\ref{BasLemStrong} implies that \index{stabilizer!of a
measure}$\Stab_G(\mu) \smallsetminus U \neq
\emptyset$.}

\item \label{BasLemTrueEx}
 Show that Thm.~\ref{Thm} implies Prop.~\ref{BasLemStrong}.

\item \label{N(U)->Singular}
 Show that if 
 \begin{itemize}
 \item $\mu$ is $U$-invariant\index{measure!invariant} and
\term[ergodic!measure]{ergodic},
 and 
 \item $c \in N_G(U)$\index{normalizer},
 \end{itemize}
 then 
 \begin{enumerate}
 \item $c_* \mu$ is $U$-invariant\index{measure!invariant} and
\term[ergodic!measure]{ergodic},
 and
 \item either $c_* \mu = \mu$ or $c_* \mu \perp \mu$.
 \end{enumerate}

\item \label{Sing->K}
 Suppose
 \begin{itemize}
 \item  $\epsilon > 0$, 
 \item $\mu$ is $U$-invariant\index{measure!invariant} and
\term[ergodic!measure]{ergodic},
 \item $c
\in N_G(U)$\index{normalizer}, 
 and 
 \item $c_* \mu \perp \mu$.
 \end{itemize}
 Show that there is a compact subset~$K$ of~$\Gamma \backslash G$, such
that 
 \begin{enumerate}
 \item $\mu(K) > 1 - \epsilon$,
 and
 \item $K \cap Kc = \emptyset$.
 \end{enumerate}

 \end{exercises}

\section{Assumptions and a restatement of~\ref{BasLem}}

\begin{assump} \label{HyperAssump}
  \label{HYPERASSUMP} %% to eliminate a Latex error
 As mentioned in Assump.~\ref{RatnerPfAssumeA}, we assume there exist
 \begin{itemize}
 \item a closed subgroup~$L$ of~$G$
 and
 \item a (nontrivial) \term[subgroup!one-parameter]{one-parameter
subgroup}~$\{a^s\}$ of~$L$,
 \end{itemize}
 such that
 \begin{enumerate}
 \item $\{u^t\} \subset L$,
 \item $\{a^s\}$ is hyperbolic, and
\index{normalizer}normalizes~$\{u^t\}$,
 \item $\mu$ is invariant under $\{a^s\}$, and
 \item $L$ is locally isomorphic to
$\SL(2,\real)$\index{SL(2,R)*$\SL(2,\real)$}.
 \nindex{$L$ = subgroup locally isomorphic to~$\SL(2,\real)$
containing~$u^t$}
 \nindex{$a^s$ = hyperbolic one-parameter subgroup in
$\Stab_{N_L(\{u^t\})}(\mu)$}
 \goodbreak
 \end{enumerate}
 \end{assump}

\begin{rem} \ 
 \begin{enumerate}
 \item Under an appropriate local isomorphism between $L$ and
$\SL(2,\real)$\index{SL(2,R)*$\SL(2,\real)$}, the subgroup $\langle
a^s, u^t \rangle$ maps to the group $\mathbb{D}_2 \mathbb{U}_2$ of
\index{lower-triangular matrices}lower triangular matrices in
$\SL(2,\real)$ \seeex{BorelIsDU}.
 \item Therefore, the parametrizations of $a^s$ and~$u^t$ can be chosen so
that $a^{-s} u^t a^s = u^{e^{2s} t}$ for all $s$ and~$t$.
 \item The \term{Mautner Phenomenon} implies that the measure~$\mu$ is
\term[ergodic!measure]{ergodic} for~$\{a^s\}$ \seecor{MautnerErgSL2}.
 \end{enumerate}
 \end{rem}

\begin{notation} \label{Notation} \ 
 \begin{itemize}
 \item For a (small) element~$g$ of~$G$, we use $\lie g$ to denote the
corresponding element $\log g$ of the \index{Lie algebra}Lie algebra
$\Lie G$.
 \nindex{$\Lie G$, $\Lie S$, $\Lie L$ = Lie algebra of~$G$, $S$, $L$}
 \nindex{$\lie g$ = $\log g$}
 \item Recall that $S = \Stab_G(\mu)^\circ$\index{stabilizer!of a
measure} \seeNot{SDefn}.
 \item By renormalizing, let us assume that
$[\lie u, \lie a] = 2 \lie u$ (where $a = a^1$ and $u = u^1$).
 \item Let $\{v^r\}$ be the (unique) one-parameter
\index{unipotent!subgroup!one-parameter}unipotent subgroup of~$L$, such
that $[\lie v, \lie a] = -2 \lie v$ and $[\lie v, \lie u] = \lie a$
\seeeg{sl2brackets}.
 \nindex{$v^r$ = unipotent subgroup of~$L$ opposite to
$\{u^t\}$}
 \item Let $\bigoplus_{\lambda \in \integer} \Lie G_\lambda$ be the
decomposition of~$\Lie G$ into \term[weight space]{weight spaces} of~$\lie
a$: that is,
 $$ \Lie G_\lambda = \bigset{ \lie g \in \Lie G }{ [\lie g, \lie a] =
\lambda \lie g } .$$
 \item Let $\Lie G_+ = \bigoplus_{\lambda > 0} \Lie G_\lambda$,
 $\Lie G_- = \bigoplus_{\lambda < 0} \Lie G_\lambda$,
 $\Lie S_+ = \Lie S \cap \Lie G_+$,
 $\Lie S_- = \Lie S \cap \Lie G_-$,
 and
 $\Lie S_0 = \Lie S \cap \Lie G_0$.
 Then
 $$ \mbox{$\Lie G = \Lie G_- + \Lie G_0 + \Lie G_+$
 \qquad and \qquad
 $\Lie S = \Lie S_- + \Lie S_0 + \Lie S_+$}
 .$$
 These are direct sums of vector spaces, although they are not direct
sums of \index{Lie algebra}Lie algebras.
 \nindex{$\Lie G_-, \Lie G_0, \Lie G_+$ = subspaces of~$\Lie G$
(determined by $a^s$)}
 \nindex{$\Lie S_-, \Lie S_0, \Lie S_+$ = subspaces of~$\Lie S$
(determined by $a^s$)}
 \item Let $G_+, G_-, G_0, S_+, S_-, S_0$ be the
\index{subgroup!connected|)}connected subgroups of~$G$ corresponding to
the Lie subalgebras $\Lie G_+, \Lie G_-, \Lie G_0, \Lie S_+, \Lie S_-,
\Lie S_0$, respectively \seeex{GplusSubalg}.
 \nindex{$G_-$, $G_0$, $G_+$, $S_-$, $S_0$, $S_+$
 = corresponding subgroups}
 \item Let $U = S_+$ (and let $\Lie U$ be the \index{Lie algebra}Lie
algebra of~$U$).
 \nindex{$U$ = $S_+$}
 \nindex{$\Lie U$ = Lie algebra of~$U$}
 \end{itemize}
 \end{notation}

Because $S_- S_0 U = S_- S_0 S_+$ contains a neighborhood of~$e$ in~$S$
\seeex{SminS0SplusOpen}, Prop.~\ref{BasLem} states that the direction of
fastest \index{transverse divergence}transverse divergence belongs to
$S_- S_0$. The following corollary is \emphit{a~priori} weaker (because
$G_-$ and~$G_0$ are presumably larger than~$S_-$ and~$S_0$), but it is
the only consequence of Lem.~\ref{KnotKc} or Lem.~\ref{RatTransInN} that
we will need in our later arguments.

\begin{cor} \label{DivMinus}
 For any $\epsilon > 0$, there is a subset $\Omega_\epsilon$ of\/ $\Gamma
\backslash G$, such that
 \begin{enumerate}
 \item $\mu(\Omega_\epsilon) > 1 - \epsilon$, 
 and
 \item \label{DivMinus-Div}
 for all $x, y \in \Omega_\epsilon$, with $x \approx y$, and such that a
certain \index{technical assumption}technical assumption
\pref{TechAssump} is satisfied, the fastest \index{transverse
divergence}transverse divergence of the $U$-orbits through $x$ and~$y$ is
along some direction that belongs to~$G_- G_0$.
 \end{enumerate}
 \end{cor}

\begin{exercises}

\item \label{GplusSubalg}
 Show $\Lie G_+$, $\Lie G_-$, and~$\Lie G_0$ are \index{Lie
algebra}subalgebras of~$\Lie G$.
 \hint{$[\Lie G_{\lambda_1}, \Lie G_{\lambda_2}] \subset \Lie
G_{\lambda_1 + \lambda_2}$.}

\item \label{SminS0SplusOpen}
 Show $S_- S_0 S_+$ contains a neighborhood of~$e$ in~$S$.
 \hint{Because $\Lie S_- + \Lie S_0 + \Lie S_- + = \Lie S$, this follows
from the \term[Theorem!Inverse Function]{Inverse Function Theorem}.}

 \end{exercises}

\section{Definition of the subgroup $\widetilde{S}$}
\label{DefineStildeSect}

 To exploit Cor.~\ref{DivMinus}, let us introduce some notation. The
corollary states that orbits diverge \emph{fastest} along some direction
in $G_- G_0$, but it will be important to understand when \emph{all} of
the \index{transverse divergence}transverse divergence, not just the
fastest part, is along $G_- G_0$. 
 More precisely, we wish to understand the elements~$g$ of~$G$, such that
if $y = xg$, then the orbits through $x$ and~$y$ diverge transversely only
along directions in $G_- G_0$: the
\index{component!G+*$G_+$-}$G_+$-component of the relative motion should
belong to~$U$, so the $G_+$-component of the divergence is trivial.
Because the divergence is measured by $u^{-1} g u$ (thought of as an
element of $G/U$), this suggests that we wish to understand
 $$ \bigset{ g \in G  }{
 \begin{matrix}
 u^{-1} g u \in G_- G_0 U, \\
 \forall u \in \mbox{some neighborhood of~$e$ in~$U$} 
 \end{matrix} 
 }.$$
 This is the right idea, but replacing $G_- G_0 U$ with its
\term[Zariski!closure]{Zariski closure} $\Zar{G_- G_0
U}$\index{Zariski!closure} yields a slightly better theory. (For example,
the resulting subset of~$G$ turns out to be a subgroup!) Fortunately,
when $g$~is close to~$e$ (which is the case we are usually interested
in), this alteration of the definition makes no difference at all
\seeex{StildeNoZar}. (This is because $G_- G_0 U$ contains a neighborhood
of~$e$ in $\Zar{G_- G_0 U}$\index{Zariski!closure} \seeex{GminG0Uopen}.)
Thus, the non-expert may wish to think of $\Zar{G_- G_0 U}$ as simply
being $G_- G_0 U$, although this is not strictly correct.

\begin{defn} \label{StildeDefn}
 Let
 $$\widetilde{S} = \bigset{
 g \in G 
 }{
 \mbox{$u^{-1} g u \in
\Zar{G_- G_0 U}$\index{Zariski!closure}, for all $u \in U$} 
 }
 $$
 \index{Zariski!closure}
 \nindex{$\widetilde{S}$ = $\bigset{g \in G}
 { \mbox{$u^{-1} g u \in \Zar{G_- G_0 U}$\index{Zariski!closure},
 for all $u \in U$}  }$}
 and \nindex{$\widetilde{S}_-$ = $\widetilde{S} \cap G_-$}
 $$ \widetilde{S}_- = \widetilde{S} \cap G_- .$$
 It is more or less obvious that $S \subset \widetilde{S}$
\seeex{SinStildeEx}. Although this is much less obvious, it should also
be noted that $\widetilde{S}$ is a closed subgroup of~$G$
\seeex{StildeSubgrp}. 
 \end{defn}

\begin{rem}  \label{StildeRem}
 Here is an alternate approach to the definition of $\widetilde{S}$, or,
at least, its identity \index{component!identity}component.
 \begin{enumerate}
 \item \label{StildeRem-Lie}
 Let 
 $$ \widetilde{\Lie S} = \bigset{ \lie g \in \Lie G }{
 \lie g (\ad \lie u)^k  \in \Lie G_- + \Lie G_0 + \Lie U,
 \ \forall k \ge 0, \ \forall \lie u \in \Lie U } 
 .$$
 Then $\widetilde{\Lie S}$ is a \index{Lie algebra}Lie subalgebra of~$\Lie
G$ \seeex{StildeLieAlgEx}, so we may let $\widetilde{S}^\circ$ be the
corresponding \index{subgroup!connected|(}connected Lie subgroup of~$G$.
(We will see in~\pref{StildeRem-same} below that this agrees with
Defn.~\ref{StildeDefn}.)
 \item \label{StildeRem-max}
 From the point of view in~\pref{StildeRem-Lie}, it is not difficult to
see that $\widetilde{S}^\circ$ is the unique maximal
\index{subgroup!connected}connected subgroup of~$G$, such that
 \begin{enumerate}
 \item \label{StildeRem-max-U} 
 %% want single letter name for this subitem
 %% should use \label mechanism ???
 \xdef\StildeRemmaxU{\alph{enumii}}
 $\widetilde{S}^\circ \cap G_+ = U$,
 and
 \item \label{StildeRem-max-norm} 
 %% want single letter name for this subitem
 %% should use \label mechanism ???
 \xdef\StildeRemmaxnorm{\alph{enumii}}
 $\widetilde{S}^\circ$ is \index{normalizer}normalized by~$a^t$
 \end{enumerate}
 \seeexs{Stilde+=UEx}{StildeMaxEx1}.
 This makes it obvious that $S \subset \widetilde{S}^\circ$. It is also
easy to verify directly that $\Lie S \subset \widetilde{\Lie S}$
\seeex{SinStilde}.
 \item \label{StildeRem-same}
 It is not difficult to see that the identity
\index{component!identity}component of the subgroup defined in
Defn.~\ref{StildeDefn} is also the subgroup characterized in
\pref{StildeRem-max} \seeex{StildeMaxEx2}, so this alternate approach
agrees with the original definition of~$\widetilde{S}$.
 \end{enumerate}
 \end{rem}

\begin{eg} \label{StildeEg}
 Remark~\ref{StildeRem} makes it easy to calculate $\widetilde{S}^\circ$.
 \begin{enumerate}
 \item \label{StildeEg-S=G}
 We have $\widetilde{S} = G$ if and only if
$U = G_+$ \seeex{StildeEg-S=GEx}.

 \item \label{StildeEg-topcorner}
 If \index{SL(3,R)*$\SL(3,\real)$}
 $$ G = \SL(3,\real),
 \quad
 \lie a = 
 \begin{bmatrix}
 1 & 0 & 0 \\
 0 & 0 & 0 \\
 0 & 0 & -1 
 \end{bmatrix}
 ,
 \quad \mbox{and} \quad
 \Lie U =
  \begin{bmatrix}
 0 & 0 & 0 \\
 0 & 0 & 0 \\
 * & 0 & 0 
 \end{bmatrix}
 ,$$
 then
 $$ \Lie G_+ =
 \begin{bmatrix}
 0 & 0 & 0 \\
 * & 0 & 0 \\
 * & * & 0  
 \end{bmatrix}
 \quad \mbox{and} \quad
 \widetilde{\Lie S} = 
 \begin{bmatrix}
 * & 0 & * \\
 0 & * & 0 \\
 * & 0 & * 
 \end{bmatrix}
 $$
 \seeex{StildeEg-topcornerEx}.

 \item \label{StildeEg-2D} 
 If \index{SL(3,R)*$\SL(3,\real)$}
 $$ G = \SL(3,\real),
 \quad
 \lie a = 
 \begin{bmatrix}
 1 & 0 & 0 \\
 0 & 0 & 0 \\
 0 & 0 & -1 
 \end{bmatrix}
 ,
 \quad \mbox{and} \quad
 \Lie U =
  \begin{bmatrix}
 0 & 0 & 0 \\
 * & 0 & 0 \\
 * & 0 & 0 
 \end{bmatrix}
 ,$$
 then
 $$
 \Lie G_+ =
 \begin{bmatrix}
 0 & 0 & 0 \\
 * & 0 & 0 \\
 * & * & 0  
 \end{bmatrix}
 \quad \mbox{and} \quad
 \widetilde{\Lie S} = 
 \begin{bmatrix}
 * & 0 & * \\
 * & * & * \\
 * & 0 & * 
 \end{bmatrix}
 $$
 \seeex{StildeEg-2DEx}.

 \item \label{StildeEg-SO12}
 If \index{SL(3,R)*$\SL(3,\real)$}
 $$ G = \SL(3,\real),
 \quad
 \lie a = 
 \begin{bmatrix}
 2 & 0 & 0 \\
 0 & 0 & 0 \\
 0 & 0 & -2 
 \end{bmatrix}
 ,
 \quad \mbox{and} \quad
 \Lie U =
  \real \begin{bmatrix}
 0 & 0 & 0 \\
 1 & 0 & 0 \\
 0 & 1 & 0 
 \end{bmatrix}
 ,$$
 then
 $$ 
 \Lie G_+ =
 \begin{bmatrix}
 0 & 0 & 0 \\
 * & 0 & 0 \\
 * & * & 0  
 \end{bmatrix}
 \quad \mbox{and} \quad
 \widetilde{\Lie S} = 
 \real \begin{bmatrix}
 0 & 1 & 0 \\
 0 & 0 & 1 \\
 0 & 0 & 0 
 \end{bmatrix}
 + 
 \begin{bmatrix}
 * & 0 & 0 \\
 0 & * & 0 \\
 0 & 0 & * 
 \end{bmatrix}
 + \Lie U
 $$
 \seeex{StildeEg-SO12Ex}.

 \item \label{StildeEg-SL2xSL2}
 If \index{SL(2,R)*$\SL(2,\real)$}
 $$G = \SL(2,\real) \times \SL(2,\real),
 \quad
 \lie a = 
 \left(
 \begin{bmatrix}
 1 & 0 \\
 0 & -1
 \end{bmatrix}
 ,
 \begin{bmatrix}
 1 & 0 \\
 0 & -1
 \end{bmatrix}
 \right)
 , $$
 and
 $$ \Lie U =
  \real  \left(
 \begin{bmatrix}
 0 & 0 \\
 1 & 0
 \end{bmatrix}
 ,
 \begin{bmatrix}
 0 & 0 \\
 1 & 0
 \end{bmatrix}
 \right)
 ,$$
 then
 $$ 
 \Lie G_+ =
  \begin{bmatrix}
 0 & 0 \\
 * & 0
 \end{bmatrix}
  \times
 \begin{bmatrix}
 0 & 0 \\
 * & 0
 \end{bmatrix}
  \mbox{ and } 
 \widetilde{\Lie S} = 
  \real  \left(
 \begin{bmatrix}
 0 & 1 \\
 0 & 0
 \end{bmatrix}
 ,
 \begin{bmatrix}
 0 & 1 \\
 0 & 0
 \end{bmatrix}
 \right)
 + 
 \real \lie a
 + \Lie U
 $$
 \seeex{StildeEg-SL2xSL2Ex}.

 \end{enumerate}
 \end{eg}

\begin{exercises}

\item \label{Norm->Conn}
 Show that if 
 \begin{itemize}
 \item $V$ is any subgroup of~$G_+$ (or of~$G_-$),
 and
 \item $V$ is \index{normalizer}normalized by~$\{a^t\}$,
 \end{itemize}
 then $V$ is \index{subgroup!connected}connected.
 \hint{If $v \in G_+$, then $a^{-t} v a^t \to e$ as $t \to -\infty$.}

\item \label{H-H0H+DenseEx}
 Show that if $H$ is a \index{subgroup!connected}connected subgroup
of~$G$, and $H$~is \index{normalizer}normalized by~$\{a^t\}$, then $H
\subset \Zar{H_- H_0 H_+}$\index{Zariski!closure}.
 \hint{$\dim \Zar{H_-} \, \Zar{H_0} \, \Zar{H_+} = \dim \Zar{H}$. Use
Exer.~\ref{SubvarSmallerDim}.}

\item Show, directly from Defn.~\ref{StildeDefn}, that $N_{G_-}(U)
\subset \widetilde{S}_-$\index{normalizer}.

\item \label{SinStildeEx}
 Show, directly from Defn.~\ref{StildeDefn}, that $S \subset
\widetilde{S}$.
 \hint{Use Exer.~\ref{H-H0H+DenseEx}.}

\item \label{StildeNeedClosure}
 Let \index{SL(2,R)*$\SL(2,\real)$}
 $G = \SL(2,\real)$,
 $a^t = \begin{bmatrix}
 e^t & 0 \\
 0 & e^{-t}
 \end{bmatrix}$
 and
 $U = G_+ = \begin{bmatrix}
 1 & 0 \\
 * & 1
 \end{bmatrix}$.
 \begin{enumerate}
 \item Show that $G_- G_0 G_+ \neq G$. 
 \item For $g \in G$, show that if $u^{-1} g u \in G_- G_0 U$, for all $u
\in U$, then $g \in G_0 U$.
 \item Show, for all $g \in G$, and all $u \in U$, that $u^{-1} g u \in
\Zar{G_- G_0 U}$\index{Zariski!closure}. Therefore $\widetilde{S} = G$.
 \end{enumerate}
 \hint{Letting $v = (0,1)$, and considering the usual
\term{representation} of~$G$ on~$\real^2$, we have $U =
\Stab_G(v)$\index{stabilizer!of a vector}. Thus, $G/U$ may be
identified with $\real^2 \smallsetminus \{0\}$. This identifies $G_- G_0
U/U$ with $\{(x,y) \in \real^2 \mid x > 0 \}$.}

\item \label{GminG0Uopen}
 Assume $G$ is almost \index{Zariski!closed!almost}Zariski closed
\seerem{MayAssZar}.
 Define the \index{function!polynomial}polynomial $\psi \colon G_- G_0
\times G_+ \to G_- G_0 G_+$ by $\psi(g,u) = gu$. (Note that $G_- G_0 G_+$
is an open subset of~$G$ \cfex{SminS0SplusOpen}.) Assume the inverse
of~$\psi$ is \index{function!rational}rational (although we do not prove
it, this is indeed always the case, cf.~Exer.~\ref{StildeNeedClosure}).

 Show that $G_- G_0 U$ is an open subset of $\Zar{G_- G_0
U}$\index{Zariski!closure}.
 \hint{$G_- G_0 U$ is the inverse image of~$U$ under a
\index{function!rational}rational map $\psi^{-1}_+ \colon G_- G_0 G_+ \to
G_+$.}

\item \label{VginVEx}
 \begin{enumerate}
 \item
 Show that if 
 \begin{itemize}
 \item $V$ is a \index{Zariski!closed}Zariski closed subset of
$\SL(\ell,\real)$\index{SL(l,R)*$\SL(\ell,\real)$}, 
 \item $g \in \SL(\ell,\real)$,
 and
 \item $V g \subset V$,
 \end{itemize}
 then $Vg = V$.
 \item \label{VginVEx-subgrp}
 Show that if $V$ is a \index{Zariski!closed}Zariski closed subset
of $\SL(\ell,\real)$\index{SL(l,R)*$\SL(\ell,\real)$}, then
 $$ \{\, g \in \SL(\ell,\real) \mid V g \subset V \,\} $$
 is a closed subgroup of $\SL(\ell,\real)$.
 \item Construct an example to show that the conclusion of
\pref{VginVEx-subgrp} can fail if $V$ is assumed only to be closed, not
\index{Zariski!closed}Zariski closed.
 \index{counterexample}
 \end{enumerate}
 \hint{Use Exer.~\ref{SubvarSmallerDim}.}

\item \label{StildeSubgrp}
 Show, directly from Defn.~\ref{StildeDefn}, that $\widetilde{S}$~is a
subgroup of~$G$.
 \hint{Show that 
 $$\widetilde{S} = \bigset{
 g \in G 
 }{
 \Zar{G_- G_0 U} \, g \subset \Zar{G_- G_0 U}
 }
 ,$$
 \index{Zariski!closure}and apply Exer.~\ref{VginVEx-subgrp}.}

\item \label{NotStilde->smallU}
 Show that if $g \notin \widetilde{S}$, then \index{Zariski!closure}
 $$ \{\, u \in U \mid u^{-1} g u \in \Zar{G_- G_0 U} \,\} $$
 is \term{nowhere dense} in~$U$. That is, its closure does not contain
any open subset of~$U$.
 \hint{It is a \index{Zariski!closed}Zariski closed, proper subset
of~$U$.}

\item \label{StildeNoZar}
 Show that there is a neighborhood~$W$ of~$e$ in~$G$, such that
 $$ \widetilde{S} \cap W
 = \bigset{ g \in W }{
 \begin{matrix}
  u^{-1} g u \in G_- G_0 U, \\
 \forall u \in
\mbox{some neighborhood of~$e$ in~$U$} 
 \end{matrix}
 } .$$
 \hint{Use Exers.~\ref{NotStilde->smallU} and \ref{GminG0Uopen}.}

\item \label{StildeLieAlgEx}
 Show, directly from the definition \see{StildeDefn}, that
 \begin{enumerate}
 \item $\widetilde{\Lie S}$~is invariant under $\ad \lie a$,
 and
 \item $\widetilde{\Lie S}$~is a Lie subalgebra of~$\Lie G$.
 \end{enumerate}
 \hint{If ${\lie g}_1 \in \widetilde{\Lie S}_{\lambda_1}$, ${\lie g}_2 \in
\widetilde{\Lie S}_{\lambda_2}$, $\lie u \in \Lie U_{\lambda_3}$, and
$\lambda_1 + \lambda_2 + (k_1 + k_2) \lambda_3 > 0$, then ${\lie g}_i (\ad
\lie u)^{k_i} \in \Lie U$, for some $i \in \{1,2\}$, so
 $$[{\lie g}_1 (\ad \lie u)^{k_1}, {\lie g}_2 (\ad \lie u)^{k_2}] \in \Lie
G_- + \Lie G_0 + \Lie U ,$$
 and it follows that $\widetilde{\Lie S}$~is a Lie subalgebra.}

\item \label{Stilde+=UEx}
 Show, directly from Defn.~\ref{StildeDefn}, that 
 \begin{enumerate}
 \item \label{Stilde+=UEx-=U}
 $\widetilde{S} \cap G_+ = U$,
 and
 \item \label{Stilde+=UEx-norm}
 $\widetilde{S}$ is \index{normalizer}normalized by~$\{a^t\}$.
 \end{enumerate}
 \hint{It suffices to show that $\widetilde{\Lie S}_+ = \Lie U$
\seeex{Norm->Conn}, and that $\widetilde{\Lie S}$ is $(\Ad_G
a^t)$-invariant.}

\item \label{StildeMaxEx1}
 Show, directly from Defn.~\ref{StildeDefn}, that if $H$ is any
\index{subgroup!connected}connected subgroup of~$G$, such that
 \begin{enumerate}
 \item $H \cap G_+ = U$,
 and
 \item $H$ is \index{normalizer}normalized by~$\{a^t\}$,
 \end{enumerate}
 then $H \subset \widetilde{S}$.
 \hint{It suffices to show that $\Lie H \subset \widetilde{\Lie S}$.}

\item \label{SinStilde}
 Show, directly from the definition of $\widetilde{\Lie S}$
in Rem.~\fullref{StildeRem}{Lie}, that $\Lie S \subset \widetilde{\Lie
S}$.

\item \label{StildeMaxEx2}
 Verify, directly from Defn.~\ref{StildeDefn} (and assuming that
$\widetilde{S}$ is a subgroup), 
 \begin{enumerate}
 \item that $\widetilde{S}$ satisfies conditions (\StildeRemmaxU)
and~(\StildeRemmaxnorm) of Rem.~\fullref{StildeRem}{max},
 and
 \item conversely, that if $H$ is a \index{subgroup!connected}connected
subgroup of~$G$, such that $H \cap G_+ = U$ and $H$~is
\index{normalizer}normalized by~$\{a^t\}$, then $H \subset \widetilde{S}$.
 \end{enumerate}

\item \label{StildeEg-S=GEx}
 Verify Eg.~\fullref{StildeEg}{S=G}.

\item \label{StildeEg-topcornerEx}
 Verify Eg.~\fullref{StildeEg}{topcorner}.

\item \label{StildeEg-2DEx}
 Verify Eg.~\fullref{StildeEg}{2D}.

\item \label{StildeEg-SO12Ex}
 Verify Eg.~\fullref{StildeEg}{SO12}.

\item \label{StildeEg-SL2xSL2Ex}
 Verify Eg.~\fullref{StildeEg}{SL2xSL2}.

 \end{exercises}

\section{Two important consequences of shearing} \label{ShearConseq}
 \index{shearing property}

Our ultimate goal is to find a conull subset~$\Omega$ of $\Gamma
\backslash G$, such that if $x,y \in \Omega$, then $y \in x S$. In this
section, we establish two consequences of Cor.~\ref{DivMinus} that
represent major progress toward this goal
\seecors{GminInSmin}{MustS}. These results deal with $\widetilde{S}$,
rather than~$S$, but that turns out not to be a very serious problem,
because $\widetilde{S} \cap G_+ = S \cap G_+$ \fullseerem{StildeRem}{max}
and $\widetilde{S} \cap G_- = S \cap G_-$ \seeprop{Stilde=S}.

\begin{notation}
 Let 
 \begin{itemize}
 \item $\Lie G_+ \ominus \Lie U$ be an $a^s$-invariant \term[complement
(of a subspace)]{complement} to $\Lie U$ in~$\Lie G_+$,
 \nindex{$\Lie G_+ \ominus \Lie U$
 = complement to $\Lie U$ in~$\Lie G_+$}

 \item $\Lie G_- \ominus \widetilde{\Lie S}_-$ be an $a^s$-invariant
\term[complement (of a subspace)]{complement} to $\widetilde{\Lie S}_-$
in~$\Lie G_-$,
 \nindex{$\Lie G_- \ominus \widetilde{\Lie S}_-$
 = complement to $\widetilde{\Lie S}_-$ in~$\Lie G_-$}

 \item $G_+ \ominus U = \exp (\Lie G_- \ominus \Lie U)$,
 \nindex{$G_+ \ominus U$ = $\exp(\Lie G_+ \ominus \Lie U)$}

 and
 \item $G_- \ominus \widetilde{S}_- = \exp (\Lie G_- \ominus
\widetilde{\Lie S}_-)$.
 \nindex{$G_- \ominus \widetilde{S}_-$ = 
 $\exp(\Lie G_- \ominus \widetilde{\Lie S}_-)$}

 \end{itemize}
 Note that the natural maps $(G_+ \ominus U) \times U \to G_+$ and $(G_-
\ominus \widetilde{S}_-) \times \widetilde{S}_- \to G_-$ (defined by
$(g,h) \mapsto g h$) are \index{diffeomorphism}diffeomorphisms
\seeex{ProdDiffeo}.
 \end{notation}

\begin{cor} \label{GminInSmin}
 There is a conull subset~$\Omega$ of\/ $\Gamma \backslash G$, such that
if $x,y \in \Omega$, and $y \in x G_-$, then $y \in x \widetilde{S}_-$.
 \end{cor}

\begin{proof}
 Choose $\Omega_0$ as in the conclusion of Cor.~\ref{DivMinus}. From
the \index{Theorem!Pointwise Ergodic}Pointwise Ergodic Theorem
\pref{PtwiseErgThmFlow}, we know that
 $$ \Omega = \bigset{
 \vphantom{\Bigl(}
 x \in \Gamma \backslash G 
 }{ 
 \mbox{$\{\, t \in \real^+\mid x a^t \in \Omega_0 \,\}$ is unbounded}
 }$$
 is conull \seeex{UnbddHitsAE}.

 We have $y = xg$, for some $g \in G_-$. Because $a^{-t} g a^t \to e$ as
$t \to \infty$, we may assume, by replacing $x$ and~$y$ with $x a^t$ and
$y a^t$ for some \index{infinitely!large}infinitely large~$t$, that $g$~is
\index{infinitesimal}infinitesimal (and that $x,y \in \Omega_0$). (See
Exer.~\ref{GminInSminClassical} for a \term{non-infinitesimal version} of
the proof.)

Suppose $g \notin \widetilde{S}_-$ (this will lead to a
contradiction). From the definition of~$\widetilde{S}_-$, this means there
is some $u \in U$, such that $u^{-1} g u \notin G_- G_0 U$: write $u^{-1}
g u = h c u'$ with $h \in G_- G_0$, $c \in G_+ \ominus U$, and $u' \in U$.
We may assume $h$~is \index{infinitesimal}infinitesimal (because we could
choose $u$~to be finite, or even infinitesimal, if desired
\seeex{NotStilde->smallU}). Translating again by an
(\index{infinitely!large}infinitely large) element of~$\{a^t\}$, with $t
\ge 0$, we may assume $c$ is infinitely large. Because $h$~is
\index{infinitesimal}infinitesimal, this clearly implies that the orbits
through $x$ and~$y$ diverge fastest along a direction in~$G_+$, not a
direction in $G_- G_0$. This contradicts Cor.~\ref{DivMinus}. (See
Exer.~\ref{GminInSminEx} for a verification of the \index{technical
assumption}technical assumption \pref{TechAssump} in that corollary.)
 \end{proof}

An easy calculation (involving only algebra, not dynamics) establishes the
following. (See Exer.~\ref{NotS->G+DivClassical} for a
\term{non-infinitesimal version}.)

\begin{lem} \label{NotS->G+Div}
 If 
 \begin{itemize}
 \item $y = xg$
 with
 \item $g \in (G_- \ominus \widetilde{S}_-) G_0 G_+$,
 and
 \item $g \approx e$,
 \end{itemize}
 then the \index{transverse divergence}transverse divergence of the
$U$-orbits through $x$ and~$y$ is fastest along some direction in~$G_+$.
 \end{lem}

\begin{proof}
 Choose $s > 0$ (\index{infinitely!large}infinitely large), such that
$\hat g = a^s g a^{-s}$ is finite, but not
\index{infinitesimal}infinitesimal, and write $\hat g = \hat g_- \hat g_0
\hat g_+$, with $\hat g_- \in G_-$, $\hat g_0 \in G_0$, and $\hat g_+ \in
G_+$. (Note that $\hat g_0$ and $\hat g_+$ are
\index{infinitesimal}infinitesimal, but $\hat g_-$ is not.) Because $\hat
g_- \in G_- \ominus \widetilde{S}_-$, we know that $\hat g$ is not
\index{infinitesimally close}infinitely close to~$\widetilde S_-$, so
there is some finite $u \in U$, such that $u^{-1} \hat g$ is \emph{not}
infinitesimally close to~$G_- G_0 U$. 
 \index{infinitely!close|indsee{infinitesimally~close}}

 Let 
 $\hat u = a^{-s} u a^s$,
 and consider
 $\hat u^{-1} g \hat u = a^{-s} ( u^{-1} \hat g u \bigr) a^s$. 
 \begin{itemize}
 \item Because $u^{-1} \hat g u$ is finite (since $u$ and $\hat g$ are
finite), we know that each of $(u^{-1} \hat g u)_-$ and $(u^{-1} \hat g
u)_0$ is finite. Therefore $(\hat u^{-1} g \hat u)_-$ and $(\hat u^{-1} g
\hat u)_0$ are finite, because conjugation by~$a^s$ does not expand $G_-$
or~$G_0$.
 \item On the other hand, we know that $(\hat u^{-1} g \hat u)_+$ is
\index{infinitely!far}infinitely far from~$U$, because the distance
between $u^{-1} \hat g u$ and~$U$ is not
\index{infinitesimal}infinitesimal, and conjugation by~$a^s$
expands~$G_+$ by an infinite factor.
 \end{itemize}
 Therefore, the fastest \index{transverse divergence}divergence is clearly
along a direction in~$G_+$.
 \end{proof}

The conclusion of the above lemma contradicts the conclusion of
Cor.~\fullref{DivMinus}{Div} (and the \index{technical
assumption}technical assumption \pref{TechAssump} is automatically
satisfied in this situation \seeex{NotS->G+DivTechEx}), so we have the
following conclusion:

\begin{cor} \label{MustS}
  For any $\epsilon > 0$,
 there is a subset\/~$\Omega_\epsilon$ of\/ $\Gamma \backslash G$, such
that
 \begin{enumerate}
 \item $\mu(\Omega_\epsilon) > 1 - \epsilon$,
 and 
 \item for all $x, y \in \Omega_\epsilon$, with $x \approx y$, we have $y
\notin x (G_- \ominus \widetilde{S}_-) G_0 G_+$.
 \end{enumerate}
 \end{cor}

This can be restated in the following \term[non-infinitesimal
version]{non-infinitesimal terms} \seeex{MustSClassicalPf}:

\begin{thmref}{MustS}
\begin{cor} \label{MustSClassical}
  For any $\epsilon > 0$,
 there is a subset\/~$\Omega_\epsilon$ of\/ $\Gamma \backslash G$, and
some $\delta > 0$, such that
 \begin{enumerate}
 \item $\mu(\Omega_\epsilon) > 1 - \epsilon$,
 and 
 \item for all $x, y \in \Omega_\epsilon$, with $d(x,y) < \delta$, we
have $y \notin x (G_- \ominus \widetilde{S}_-) G_0 G_+$.
 \end{enumerate}
 \end{cor}
 \end{thmref}

 \onlyoneexercise
 \begin{exercises}

\item \label{ProdDiffeo}
 Show that if $\Lie V$ and~$\Lie W$ are two \index{complement
(of a subspace)}complementary $a^t$-invariant subspaces of~$\Lie G_+$,
then the natural map $\exp \Lie V \times \exp \Lie W \to G_+$, defined by
$(v,w) \mapsto v w$, is a \index{diffeomorphism}diffeomorphism.
 \hint{The \term[Theorem!Inverse Function]{Inverse Function Theorem}
implies that the map is a \term[diffeomorphism!local]{local
diffeomorphism} near~$e$. Conjugate by~$a^s$ to expand the good
neighborhood.}

 \end{exercises}
 \endgroup

\section{Comparing $\widetilde{S}_-$ with~$S_-$}
 \label{CompareSSect}

 \index{entropy!of a dynamical system|(}
 We will now show that $\widetilde{S}_- = S_-$ \seeprop{Stilde=S}. To do
this, we use the following lemma on the \term{entropy} of translations on
\term[homogeneous!space]{homogeneous spaces}.
 Corollary~\ref{GminInSmin} is what makes this possible, by verifying the
hypotheses of Lem.~\fullref{EntLem}{GminInW}, with $W = \widetilde{S}_-$. 

\begin{thmref}{EntropyLemma}
 \begin{lem} \label{EntLem}
 Suppose $W$ is a closed, \index{subgroup!connected|)}connected subgroup
of~$G_-$ that is \index{normalizer|)}normalized by~$a$, and let 
 $$J(a^{-1}, W) = \det \bigl( (\Ad a^{-1})|_{\Lie W} \bigr)$$
 be the \index{Jacobian}Jacobian of~$a^{-1}$ on~$W$. 
 \begin{enumerate}

 \item \label{EntLem-Inv->}
 If $\mu$ is $W$-invariant, then $h_\mu(a) \ge \log J(a^{-1},
W)$\index{Jacobian}.

 \item \label{EntLem-GminInW}
 If there is a conull, Borel subset~$\Omega$ of $\Gamma \backslash
G$, such that $\Omega \cap x G_- \subset x W$, for every $x \in \Omega$,
then $h_\mu(a) \le \log J(a^{-1}, W)$\index{Jacobian}. 

 \item \label{EntLem-equal}
 If the hypotheses of \pref{EntLem-GminInW} are satisfied, and
equality holds in its conclusion, then $\mu$~is
$W$-invariant\index{measure!invariant}.

 \end{enumerate}
 \end{lem}
 \end{thmref}

\begin{prop} \label{Stilde=S}
 We have $\widetilde{S}_- = S_-$.
 \end{prop}

\begin{proof}[Proof {\rm (cf.\ proofs of Cors.~\ref{Entropy->SL2inv} and
\ref{FinFib->FinCov})}]
 We already know that $\widetilde{S}_- \supset S_-$
\fullseerem{StildeRem}{max}. Thus, because $\widetilde{S}_- \subset G_-$,
it suffices to show that $\widetilde{S}_- \subset S$. That is, it suffices
to show that $\mu$~is $\widetilde{S}_-$-invariant.

 From Lem.~\fullref{EntLem}{Inv->}, with $a^{-1}$ in the role of~$a$,
and $U$~in the role of~$W$, we have
 \index{entropy!of a dynamical system}
 \index{Jacobian}
 $$ h_\mu(a^{-1}) \ge \log J(a, U) .$$
 From Cor.~\ref{GminInSmin} and Lem.~\fullref{EntLem}{GminInW}, we
have 
 $$ h_\mu(a) \le \log J(a^{-1}, \widetilde{S}_-) .$$
 Combining these two inequalities with the fact that $h_\mu(a) =
h_\mu(a^{-1})$ \seeex{h(T)=h(Tinv)Ex}, we have
 $$ \log J(a, U) \le h_\mu(a^{-1}) = h_\mu(a) \le \log J(a^{-1},
\widetilde{S}_-) .$$
 Thus, if we show that 
 \begin{equation} \label{Stilde=SPf-J<J}
 \log J(a^{-1}, \widetilde{S}_-) \le \log J(a, U) ,
 \end{equation}
 then we must have equality throughout, and the desired conclusion will
follow from Lem.~\fullref{EntLem}{equal}.

Because $\lie u$ belongs to the \index{Lie algebra}Lie algebra~$\Lie L$
of~$L$ \seeNot{Notation}, the structure of
$\LieSL(2,\real)$-modules\index{module} implies, for each $\lambda \in
\integer^+$, that the restriction  $(\ad_{\Lie G} \lie u)^\lambda|_{\Lie
G_{-\lambda}}$ is a bijection from the \index{weight space}weight
space~$\Lie G_{-\lambda}$ onto the weight space~$\Lie G_\lambda$
\seeex{sl2WtBij}. If $\lie g \in \widetilde {\Lie S}_- \cap \Lie
G_{-\lambda}$, then Rem.~\fullref{StildeRem}{Lie} implies $\lie g
(\ad_{\Lie G} \lie u)^\lambda \in (\Lie G_- + \Lie G_0 + \Lie U) \cap
\Lie G_\lambda = \Lie U \cap \Lie G_\lambda$, so we conclude that
$(\ad_{\Lie G} \lie u)^\lambda|_{\widetilde{\Lie S}_- \cap \Lie
G_{-\lambda}}$ is an embedding of~$\widetilde{\Lie S}_- \cap \Lie
G_{-\lambda}$ into~$\Lie U \cap \Lie G_\lambda$. So 
 $$\dim (\widetilde{\Lie S}_- \cap \Lie G_{-\lambda}) \le \dim (\Lie U
\cap \Lie G_\lambda) .$$
 The \index{eigenvalue}eigenvalue of $\Ad_G a = \exp( \ad_{\Lie G} \lie
a)$ on~$\Lie G_\lambda$ is~$e^\lambda$, and the
\index{eigenvalue}eigenvalue of $\Ad_G a^{-1}$ on~$\Lie G_{-\lambda}$ is
also~$e^\lambda$ \seeex{Wt(aInv)}.
 Hence, 
 \index{Jacobian}
 \begin{align*}
 \log J(a^{-1}, \widetilde{S}_-) 
 &= \log \det \bigl( \Ad_G a^{-1} \bigr) |_{\widetilde{\Lie S}_-} \\
 &= \log \prod_{\lambda \in \integer^+} (e^\lambda)^{\dim (\widetilde {\Lie
S}_- \cap \Lie G_{-\lambda})} \\
 &= \sum_{\lambda \in \integer^+}
\bigl( \dim (\widetilde {\Lie S}_- \cap \Lie
G_{-\lambda}) \bigr) \cdot \log e^\lambda \\
 &\le \sum_{\lambda \in \integer^+} (\dim \Lie
U \cap \Lie G_\lambda) \cdot \log e^\lambda \\
 &= \log J(a, U)
 ,
 \end{align*}
 as desired.
 \end{proof}

 \index{entropy!of a dynamical system|)}

\begin{exercises}

\item \label{sl2WtBij}
  Suppose 
 \begin{itemize}
 \item $\module$ is a finite-dimensional\index{dimension!of a vector
space} real vector space,
 and
 \item $\rho \colon \LieSL(2,\real) \to \LieSL(\module)$ is a \index{Lie
algebra}Lie algebra \index{homomorphism!of Lie algebras}homomorphism.
 \end{itemize}
 Show, for every $m \in \integer^{\ge 0}$, that
 $(\lie u^\rho)^m$ is a bijection from $\module_{-m}$
to~$\module_m$.\index{weight space}
 \hint{Use Prop.~\ref{sl2Rmodules}. If $-\lambda_i \le 2j - \lambda_i \le
0$, then $w_{i,j} (\lie u^\rho)^{\lambda_i - 2j}$ is a nonzero multiple
of~$w_{i,\lambda_i-j}$,
 and $w_{i,\lambda_i-j} (\lie v^\rho)^{\lambda_i - 2j}$ is a nonzero
multiple of~$w_{i,j}$.}

\item \label{Wt(aInv)}
 In the notation of the proof of Prop.~\ref{Stilde=S}, show that the
\index{eigenvalue}eigenvalue of $\Ad_G a^{-1}$ on~$\Lie
G_{-\lambda}$\index{weight space} is the same as the the eigenvalue of
$\Ad_G a$ on~$\Lie G_\lambda$.

 \end{exercises}

\section{Completion of the proof} \label{EndOfPfSect}

 We wish to show, for some $x \in \Gamma \backslash G$, that $\mu(xS) >
0$. In other words, that $\mu(x S_- S_0 S_+) > 0$. The following
weaker result is a crucial step in this direction.

\begin{prop} \label{SuppGU}
 For some $x \in \Gamma \backslash G$, we have $\mu(x S_- G_0 G_+) > 0$.
 \end{prop}

\begin{proof}
 Assume that the desired conclusion fails. (This will lead to a
contradiction.) Let $\Omega_\epsilon$ be as in Cor.~\ref{MustS}, with
$\epsilon$~sufficiently small.

 Because the conclusion of the proposition is assumed to fail, there
exist $x,y \in \Omega_\epsilon$, with $x \approx y$ and $y = x g$, such
that $g \notin S_- G_0 G_+$. (See Exer.~\ref{SuppGUClassical} for a
\term[non-infinitesimal version]{non-infinitesimal} proof.) Thus, we may
write 
 $$ \mbox{$g = v w h$ with $v \in S_-$, $w \in (G_- \ominus S_-)
\smallsetminus \{e\}$, and $h \in G_0 G_+$} .$$
 For simplicity, let us pretend that $\Omega_\epsilon$~is
$S_-$-invariant. (This is not so far from the truth, because $\mu$~is
$S_-$-invariant\index{measure!invariant} and $\mu(\Omega_\epsilon)$ is
very close to~$1$, so the actual proof is only a little more complicated
\seeex{SuppGUEx}.)
 Then we may replace $x$ with~$x v$, so that $g = w h \in (G_- \ominus
S_-) G_0 G_+$. This contradicts the definition of~$\Omega_\epsilon$.
 \end{proof}

We can now complete the proof (using some of the theory of
\index{algebraic group!theory of}algebraic groups).

 \begin{thm} \label{SuppS}
 $\mu$ is \term[support!of a measure]{supported} on a single $S$-orbit.
 \end{thm}

\begin{proof}
 There is no harm in assuming that $G$ is almost
\index{Zariski!closed!almost}Zariski closed \seerem{MayAssZar}. By
\label{InductionSpot}\term[induction on $\dim G$]{induction on $\dim G$},
we may assume that there does \emph{not} exist a subgroup~$H$ of~$G$,
such that 
 \begin{itemize}
 \item $H$ is almost \index{Zariski!closed!almost}Zariski closed,
 \item $U \subset H$, 
	and 
 \item some $H$-orbit has full measure.
 \end{itemize}
 Then a short argument \seeex{NoOrbit->NoSubvar} implies, for all $x \in
\Gamma \backslash G$, that
 \begin{equation} \label{Variety}
 \begin{matrix}
 \mbox{if $V$ is any subset of~$G$,} \\
 \mbox{such that $\mu(x V) > 0$, then $G \subset \Zar{V}$.}
 \end{matrix}
 \end{equation}
 \index{Zariski!closure} \index{variety}
 This hypothesis will allow us to show that $S = G$.

\begin{claim}
 We have $S_- = G_-$.
 \end{claim}
 Prop.~\ref{SuppGU} states that $\mu \bigl(x S_- G_0 G_+) > 0$, so,
from~\pref{Variety}, we know that $G \subset \Zar{S_- G_0
G_+}$\index{Zariski!closure}. This implies that $S_- G_0 G_+$ must
contain an open subset of~$G$ \seeex{GminG0Uopen}. Therefore 
 $$\dim S_- \ge \dim G - \dim(G_0 G_+) = \dim G_- .$$
 Because $S_- \subset G_-$, and $G_-$ is
\index{subgroup!connected}connected, this implies that $S_- = G_-$, as
desired.

 The subgroup $G_-$ is a \term[subgroup!horospherical]{horospherical}
subgroup of~$G$ \seerem{horospherical}, so we have shown that $\mu$~is
invariant under a horospherical subgroup of~$G$.

There are now at least three ways to complete the argument. 
 \begin{enumerate} \renewcommand{\theenumi}{\alph{enumi}}

 \item \label{Pf-FinishSymm}
 We showed that $\mu$ is $G_-$-invariant. By going through the same
argument, but with $v^r$ in the place of~$u^t$, we could show that
$\mu$~is $G_+$-invariant. So $S$~contains $\langle L, G_+, G_- \rangle$,
which is easily seen to be a (\index{Lie group!unimodular}unimodular)
\index{subgroup!normal}normal subgroup of~$G$ \seeex{HoroGenNorm}. Then
Exer.~\ref{NormalSubgrpEx} applies.

 \item \label{Pf-FinishByEntropy}
 By using considerations of \index{entropy!of a dynamical system}entropy,
much as in the proof of Prop.~\ref{Stilde=S}, one can show that
$G_+ \subset S$ \seeex{EntGminToGplus}, and then
Exer.~\ref{NormalSubgrpEx} applies, once again.

 \item \label{Pf-FinishByHoro}
 If we assume that $\Gamma \backslash G$ is compact (and in some other
cases), then a completely separate proof of the theorem is known for
measures that are invariant under a
\index{subgroup!horospherical}horospherical subgroup. (An example of an
argument of this type appears in Exers.~\ref{ParabUniqErg}
and~\ref{HoroUniqErg}.) Such special cases were known several years
before the general theorem. \qedhere

 \end{enumerate}
 \end{proof}

\begin{exercises}

\item \label{SuppGUEx}
 Prove Prop.~\ref{SuppGU} (without assuming $\Omega_\epsilon$ is
$S_-$-invariant).
 \hint{Because $\Omega_\epsilon$ contains 99\% of the $S_-$-orbits of both
$x$ and~$y$, it is possible to find $x' \in x S_- \cap \Omega_\epsilon$
and $y' \in y S_- \cap \Omega_\epsilon$, such that $y' \in x' (G_-
\ominus S_-) G_0 G_+$.}

\item \label{SuppGUClassical}
 Prove Prop.~\ref{SuppGU} \term[non-infinitesimal version]{without using
infinitesimals}. 
 \hint{Use Cor.~\ref{MustSClassical}.}

\item \label{HoroGenNorm}
 Show that $\langle G_-, G_+ \rangle$ is a \index{subgroup!normal}normal
subgroup of~$G^\circ$.
 \hint{It suffices to show that it is \index{normalizer}normalized by
$G_-$, $G_0$, and~$G_+$.}

\item \label{EntGminToGplus}
 \begin{enumerate}
 \item Show that $J(a^{-1}, G_-) = J(a, G_+)$\index{Jacobian}.
 \item Use Lem.~\ref{EntLem} (at the beginning of
\S\ref{CompareSSect}) to show that if $\mu$ is $G_-$-invariant, then it is
$G_+$-invariant.
 \end{enumerate}

%\item \label{horominimal} 
% Show that if $\Gamma \backslash G$ is compact, and $G_+$ has a dense
%orbit on $\Gamma \backslash G$, then $G_+$ is minimal on $\Gamma
%\backslash G$; that is, every $G_+$ orbit is dense.

\item \label{HoroUniqErg}
 Let 
 \begin{itemize}
 \item $G$ be a \index{subgroup!connected}connected, \index{Lie
group!semisimple}semisimple subgroup of
$\SL(\ell,\real)$\index{SL(l,R)*$\SL(\ell,\real)$},
 \item $\Gamma$ be a \term{lattice} in~$G$, such that $\Gamma \backslash
G$ is compact,
 \item $\mu$ be a probability measure on $\Gamma \backslash G$,
 \item $a^s$ be a nontrivial 
\term[subgroup!one-parameter!hyperbolic]{hyperbolic one-parameter
subgroup} of~$G$,
 and
 \item $G_+$ be the corresponding expanding
\term[subgroup!horospherical]{horospherical subgroup} of~$G$.
 \end{itemize}
 Show that if $\mu$ is $G_+$-invariant, then $\mu$ is the
\index{measure!Haar}Haar measure on $\Gamma \backslash G$.
 \hint{Cf.~hint to Exer.~\ref{ParabUniqErg}. (Let $U_\epsilon \subset
G_+$, $A_\epsilon \subset G_0$, and $V_\epsilon \subset G_-$.) Because
$\mu$ is not assumed to be $a^s$-invariant, it may not be possible to
choose a \term[generic!point]{generic} point~$y$ for~$\mu$, such that $y
a^{s_k} \to y$. Instead, show that the \index{mixing}mixing property
\pref{MixingRem} can be strengthened to apply to the compact family of
subsets $\{\, y U_\epsilon A_\epsilon V_\epsilon \mid y \in \Gamma
\backslash G \,\}$.}

 \end{exercises}

\section{Some precise statements} \label{PreciseSect}

Let us now state these results more precisely, beginning with the
statement that \index{function!polynomial}polynomials stay near their
largest value for a proportional length of time.

\begin{lem}
 For any~$d$ and~$\epsilon$, and any \index{averaging!sequence|(}averaging
sequence $\{E_n\}$ of open sets in any \index{unipotent!subgroup}unipotent
subgroup~$U$ of~$G$, there is a ball~$B$ around~$e$ in~$U$, such that if 
 \begin{itemize}
 \item $f \colon U \to \real^m$ is any
\index{function!polynomial}polynomial of degree $\le d$,
 \item $E_n$ is an averaging set in the averaging sequence $\{E_n\}$, 
 and
 \item $\sup_{u \in E_n} \|f(u)\| \le 1$,
 \end{itemize}
 then 
 $ \| f(v_1 u v_2) - f(u) \| < \epsilon$, for all $u \in E_n$, and all
$v_1, v_2 \in B_n = a^{-n} B a^n$.
 \end{lem}

\begin{rem}
 Note that $\nu_U(B_n) / \nu_U(E_n) = \nu_U(B)/\nu_U(E)$
is independent of~$n$; thus, $B_n$~represents an amount of time
proportional to~$E_n$.
 \end{rem}

\begin{proof}
 The set $\poly^d$ of real \index{function!polynomial}polynomials of
degree $\le d$ is a finite-dimensional\index{dimension!of a vector space}
vector space, so 
 $$ \bigset{
 f \in \poly^d 
 }{
 \sup_{u \in E} \|f(u)\| \le 1 
 } $$
 is compact. Thus, there is a ball~$B$ around~$e$ in~$U$, such that the
conclusion of the lemma holds for $n = 0$.  Rescaling by~$a^n$
then implies that it must also hold for any~$n$.
 \end{proof}

As was noted in the previous chapter, if $y = x g$, then the relative
displacement between $xu$ and~$yu$ is $u^{-1} g u$. For each fixed~$g$,
this is a \index{function!polynomial}polynomial function on~$U$, and the
degree is bounded independent of~$g$. The following observation makes a
similar statement about the \index{transverse divergence}transverse
divergence of two $U$-orbits. It is a formalization of
Lem.~\ref{PolyDiv}.

\begin{rem} \label{PolyDivInG/U}
 Given $y = xg$, the relative displacement between $xu$ and~$yu$ is
$u^{-1} g u$. To measure the part of this displacement that is transverse
to the $U$-orbit, we wish to multiply by an element~$u'$ of~$U$, to make
$(u^{-1} g u) u'$ as small as possible: equivalently, we can simply
think of $u^{-1} g u$ in the quotient space $G/U$. That is, 
 \begin{quote} 
 \emphit{the transverse distance between the two $U$-orbits {\upshape(}at
the point $xu${\upshape)} is measured by the position of the point
$(u^{-1} g u) U = u^{-1} g U$ in the \term[homogeneous!space]{homogeneous
space} $G/U$.}
 \end{quote}
 \index{transverse divergence}
 Because $U$ is \index{Zariski!closed}Zariski closed
\seeprop{Unip->Closed}, we know, from \index{Theorem!Chevalley}Chevalley's
Theorem \pref{ChevalleyFixLine}, that, for some~$m$, there is 
 \begin{itemize}
 \item a polynomial \index{homomorphism!polynomial}homomorphism $\rho
\colon G \to \SL(m,\real)$\index{SL(l,R)*$\SL(\ell,\real)$},
 and
 \item a vector $w \in \real^m$,
 \end{itemize} such that (writing our linear transformations on the left)
we have
 $$ U = \{\, u \in G \mid \rho(u) w = w \,\} .$$
 Thus, we may identify $G/U$ with the orbit $w G$, and, because $\rho$
is a \index{function!polynomial}polynomial, we know that $u \mapsto
\rho(u^{-1} g u) w $ is a polynomial function on~$U$. Hence, the
transverse distance between the two $U$-orbits is completely described by
a polynomial function.
 \end{rem}

We now make precise the statement in Lem.~\ref{RatTransInN} that the
direction of \index{transverse divergence}fastest divergence is in the
direction of the \index{normalizer|(}normalizer. (See
Rem.~\ref{DivergeInNClassical} for a \term{non-infinitesimal version} of
the result.)

 \begin{prop} \label{RatTransInNformal}
 Suppose 
 \begin{itemize}
 \item $U$ is a \index{subgroup!connected}connected,
\index{unipotent!subgroup}unipotent subgroup of~$G$,
 \item $x,y \in \Gamma \backslash G$, 
 \item $y = x g$, for some $g \in G$, with $g \approx e$,
 and 
 \item ${E}$~is an\/ {\upshape(}infinitely large\/{\upshape)}
\index{averaging!set}averaging set,
 \end{itemize}
 such that 
 \begin{itemize}
 \item $g^{E} = \{\, u^{-1} g u \mid u \in {E} \,\}$ has finite
diameter in $G/U$.
 \end{itemize}
 Then each element of~$g^{E} U/U$ is \index{infinitesimally
close}infinitesimally close to some element of
$N_G(U)/U$\index{normalizer}.
 \end{prop}

\begin{proof}
 Let $g' \in g^{E}$.
 Note that
 $$ N_G(U)/U = \{\, x \in G/U \mid \mbox{$ux = x$, for all $u \in U$} \,\}
.$$
 \index{normalizer}
 Thus, it suffices to show that $u g' U \approx g' U$, for each finite $u
\in U$.

We may assume $g' U$ is a finite (not \index{infinitesimal}infinitesimal)
distance from the base point~$e U$, so its distance is comparable to the
farthest distance in~$g^{E} U/U$. It took
\index{infinitely!long}infinitely long to achieve this distance, so
polynomial divergence implies that it takes a proportional, hence
infinite, amount of time to move any additional finite distance. Thus, in
any finite time, the point $g' U$ moves only
\index{infinitesimal}infinitesimally. Therefore, $u g' U \approx g'
U$, as desired.
 \end{proof}

\begin{rem} \label{DivergeInNClassical}
 The above statement and proof are written in terms of
\index{infinitesimal}infinitesimals. To obtain a
\index{non-infinitesimal version}non-infinitesimal version, replace 
 \begin{itemize}
 \item $x$ and~$y$ with convergent sequences $\{x_k\}$ and $\{y_k\}$, such
that $d(x_k,y_k) \to e$,
\item $g$ with the sequence $\{g_k\}$, defined by $x_k g_k = y_k$,
 and
 \item $E$ with an \index{averaging!sequence}averaging sequence~$E_n$,
such that $g_k^{E_{n_k}}$ is bounded in $G/U$ (independent of~$k$).
 \end{itemize}
 The conclusion is that if $\{g'_k\}$ is any sequence, such that
 \begin{itemize}
 \item $g'_k \in g_k^{E_{n_k}}$ for each~$k$, 
 and
 \item $g'_k U/U$ converges,
 \end{itemize}
 then the limit is an element of \index{normalizer}$N_G(U)/U$
\seeex{DivergeInNClassicalEx}.
 \end{rem}

\begin{lem} \label{KnotKc}
 If 
 \begin{itemize}
 \item $C$ is any compact subset of \index{normalizer}$N_G(U)
\smallsetminus \Stab_G(\mu)$\index{stabilizer!of a
measure},
 and
 \item $\epsilon > 0$,
 \end{itemize}
 then there is a compact subset~$K$ of\/ $\Gamma \backslash G$, such that
 \begin{enumerate}
 \item $\mu(K) > 1 - \epsilon$
 and
 \item $K \cap K c = \emptyset$, for all $c \in C$.
 \end{enumerate}
 \end{lem}

\begin{proof}
 Let $\Omega$ be the set of all points in $\Gamma \backslash G$ that are
\term[generic!point]{generic} for~$\mu$
(see Defn.~\ref{GenericDefn} and Thm.~\ref{PtwiseErg}).
 It suffices to show that $\Omega \cap \Omega c = \emptyset$, for all $c
\in N_G(U) \smallsetminus \Stab_G(\mu)$\index{stabilizer!of a
measure}\index{normalizer}, for then we may choose $K$~to be
any compact subset of~$\Omega$ with $\mu(K) > 1 - \epsilon$.

Fix $c \in C$. We choose a compact subset $K_c$ of~$\Omega$ with $K_c
\cap K_c c = \emptyset$, and $\mu(K_c) > 1 - \delta$, where
$\delta$~depends only on~$c$, but will be specified later.

Now suppose $x, xc \in \Omega$. Except for a proportion~$\delta$ of the
time, we have $xu$ very near to~$K_c$ (because $x \in \Omega_c$). Thus,
it suffices to have $x u c$ very close to~$K_c$ more than a
proportion~$\delta$ of the time. That is, we wish to have $(xc) (c^{-1} u
c)$ very close to~$K_c$ a significant proportion of the time. 

We do have $(xc) u$ very close to~$K_c$ a huge proportion of the time. 
 Now $c$ acts on~$U$ by conjugation, and the \index{Jacobian}Jacobian of
this \index{diffeomorphism}diffeomorphism is constant (hence bounded), as
is the maximum \index{eigenvalue}eigenvalue of the derivative. Thus, we
obtain the desired conclusion by choosing $\delta$~sufficiently small
(and $E$ to be a nice set) \seeex{KnotKcEx}. 
 \end{proof}

\begin{numproof}[Completing the proof of Prop.~\ref{BasLem}]
\label{BasLemPf}
 Fix a set $\Omega_0$ as in the Uniform \index{Theorem!Pointwise
Ergodic}Pointwise Ergodic Theorem \pref{UnifErg}.
 Suppose $x,y \in \Omega_0$ with $x \approx y$, and write $y = x g$.
Given an (\index{infinitely!large}infinite)
\index{averaging!set}averaging set $E_n = a^{-n} E a_n$, such that
$g^{E_n}$ is bounded in $G/U$, and any $v \in E$, we wish to show that
$(a^{-n} v a^n) g U$ is \index{infinitesimal}infinitesimally close to
$\Stab_G(\mu)/U$\index{stabilizer!of a measure}. The proof of
Prop.~\ref{BasLem} will apply if we show that $y u'$ is close to~$K$ a
significant proportion of the time. 

To do this, we make the additional \index{technical assumption}technical
assumption that%
 \begin{equation} \label{Cond*}
 \mbox{$g^* = a^n g a^{-n}$ is finite (or
\index{infinitesimal}infinitesimal).}
 \end{equation}
 Let us assume that $E$ is a ball around~$e$. Choose a small
neighborhood~$B$ of~$v$ in~$E$, and define 
 $$ \mbox{$\sigma \colon B \to U$ by $u g^* \in (G \ominus U) \cdot
\sigma(u)$, \qquad for $u \in B$}
 ,$$
 so 
 $$u' = a^{-n} \, \sigma(u) a^n .$$

 The \index{Jacobian}Jacobian of~$\sigma$ is bounded (between~$1/J$
and~$J$, say), so we can choose $\epsilon$ so small that 
 $$ (1 - J^2 \epsilon) \cdot \nu_U(B) > \epsilon \cdot \nu_U(E).$$
 (The compact set~$K$ should be chosen with $\mu(K) > 1 - \epsilon$.)

By applying Cor.~\ref{UnifErg} to the \index{averaging!sequence}averaging
sequence $\sigma(B)_n$ (and noting that $n$~is
\index{infinitely!large}infinitely large), and observing that 
 $$y (a^{-n} u a^n) = x (a^{-n} g^* u a^n) ,$$
 we see that 
 $$ \nu_U \bigl( \{\, u \in \sigma(B) \mid x (a^{-n} g^* u a^n)
\not\approx K \,\} \bigr) \le \epsilon \, \nu_U(B_n) .$$
 Therefore, the choice of~$\epsilon$ implies
 $$ \nu_U \bigl( \{\, u \in B \mid x \bigl( a^{-n} g^* \, \sigma(u) \,
a^n \bigr) \approx K \,\} \bigr) > \epsilon \, \nu_U(E) .$$
 Because 
 \begin{align*}
 x \bigl( a^{-n} g^* \, \sigma(u) \, a^n \bigr)
 &= x \, ( a^{-n} g^* a^n) \, \bigl( a^{-n} \sigma(u) \, a^n \bigr)
 \\&= x g u' 
 \\&= y u' ,
 \end{align*}
 this completes the proof. 
 \end{numproof}
 \vskip-\lastskip
 \goodbreak
 \medskip

\begin{techassump} \label{TechAssump} \ 
 \begin{enumerate} 
 \item The \index{technical assumption}technical assumption \pref{Cond*}
in the proof of Prop.~\ref{BasLem} can be stated in the following
explicit form if $g$~is \index{infinitesimal}infinitesimal: there are 
 \begin{itemize}
 \item an (\index{infinitely!large}infinite) integer~$n$, 
 and
 \item a finite element~$u_0$ of~$U$,
 \end{itemize}
 such that
 \begin{enumerate}
 \item $a^{-n} u_0 a^n g \in G_- G_0 G_+$,
 \item $a^{-n} u_0 a^n g U$ is not
\index{infinitesimally close}infinitesimally close to $eU$ in~$G/U$,
 and
 \item $a^n g a^{-n}$ is finite (or \index{infinitesimal}infinitesimal).
 \end{enumerate}
 \item \label{TechAssump-Classical}
 In \term[non-infinitesimal version]{non-infinitesimal terms}, the
assumption on $\{g_k\}$ is: there are 
 \begin{itemize}
 \item a sequence $n_k \to \infty$,
 and
 \item a bounded sequence~$\{u_k\}$ in~$U$,
 \end{itemize}
 such that
 \begin{enumerate}
 \item $a^{-n_k} u_k a^{n_k} g_k \in G_- G_0 G_+$,
 \item no subsequence of $a^{-n_k} u_k a^{n_k} g_k U$ converges to $eU$
in~$G/U$,
 and
 \item $a^{n_k} g_k a^{-n_k}$ is bounded.
 \end{enumerate}
 \end{enumerate}
 \end{techassump}

\begin{exercises}

\item \label{GminInSminEx}
 Show that if $g = a^{-t} v a^t$, for some standard $v \in G_-$, and
$g$~is \index{infinitesimal}infinitesimal, then either
 \begin{enumerate}
 \item $g$~satisfies the \index{technical assumption}technical assumption
\pref{TechAssump},
 or 
 \item $g \in \widetilde{S}_-$.
 \end{enumerate}

\item \label{NotS->G+DivTechEx}
 Show that if $g$ is as in Lem.~\ref{NotS->G+Div}, then $g$~satisfies the
\index{technical assumption}technical assumption \pref{TechAssump}.
 \hint{Choose $n>0$ so that $a^n g a^{-n}$ is finite, but not
\index{infinitesimal}infinitesimal. Then $a^n g a^{-n}$ is not
\index{infinitesimally close}infinitesimally close to $\widetilde{S}$, so
there is some (small) $u \in U$, such that $u (a^n g a^{-n})$ is not
infinitesimally close to $G_- G_0 U$. Conjugate by~$a^n$.}

 \item \label{DivergeInNClassicalEx}
 Provide a (\term[non-infinitesimal version]{non-infinitesimal}) proof of
Rem.~\ref{DivergeInNClassical}.

\item \label{KnotKcEx}
 Complete the proof of Lem.~\ref{KnotKc}, by showing that if $E$ is
a convex neighborhood of~$e$ in~$U$, and $\delta$~is sufficiently small,
then, for all~$n$ and every subset~$X$ of~$E_n$\index{averaging!sequence}
with $\nu_U(X) \ge (1 - \delta) \nu_U(E_n)$, we have 
 $$ \nu_U \bigl( \{\, u \in E_n \mid c^{-1} u c \in X \,\} \bigr)
 > \delta \, \nu_U(X) .$$
 \hint{There is some $k > 0$, such that 
 $$ \mbox{$c^{-1} E_{n-k} c \subset E_n$, for all~$n$} .$$
 Choose $\delta$ small enough that 
 $$ \nu_U(E_{n-k}) > (J + 1) \delta \nu_U(E_n) ,$$
 where $J$ is the \index{Jacobian}Jacobian of the conjugation
\index{diffeomorphism}diffeomorphism.}

\item \label{BasLemClassical} 
 Prove the \term{non-infinitesimal version} of Prop.~\ref{BasLem}:
 For any $\epsilon > 0$, there is a compact subset $\Omega_\epsilon$ of\/
$\Gamma \backslash G$, with $\mu(\Omega_\epsilon) > 1 - \epsilon$, and
such that if
 \begin{itemize}
 \item $\{x_k\}$ and $\{y_k\}$ are convergent sequences in
$\Omega_\epsilon$, 
 \item $\{g_k\}$ is a sequence in~$G$ that satisfies
\fullref{TechAssump}{Classical},
 \item $x_k g_k = y_k$,
 \item $g_k \to e$,
 \item $\{E_n\}$ is an \index{averaging!sequence}averaging sequence,
and $\{n_k\}$ is a sequence of natural numbers, such that $g_k^{E_{n_k}}$
is bounded in $G/U$ (independent of~$k$),
 \item $g'_k \in g_k^{E_{n_k}}$, 
 and
 \item $g'_k U/U$ converges,
 \end{itemize}
 then the limit of $\{g'_k U\}$ is an element of $S/U$.

\item \label{GminInSminClassical}
 Prove Cor.~\ref{GminInSmin} \term[non-infinitesimal version]{without
using infinitesimals}.

\item \label{NotS->G+DivClassical}
 Prove the \term{non-infinitesimal version} of Lem.~\ref{NotS->G+Div}:
 If 
 \begin{itemize}
 \item $\{g_n\}$ is a sequence in~$(G_- \ominus \widetilde{S}_-) G_0 G_+$,
 \item $\{E_n\}$ is an \index{averaging!sequence|)}averaging sequence,
and $\{n_k\}$ is a sequence of natural numbers, such that $g_k^{E_{n_k}}$
is bounded in $G/U$ (independent of~$k$),
 \item $g'_k \in g_k^{E_{n_k}}$, 
 and
 \item $g'_k U/U$ converges,
 \end{itemize}
 then the limit of $\{g'_k U\}$ is an element of $G_+/U$.

\item \label{MustSClassicalPf}
 Prove Cor.~\ref{MustSClassical}.

\end{exercises}

\section{How to eliminate Assumption~\ref{HyperAssump}}
 \label{RatnerPfAssumeASect}

Let $\hat U$ be a maximal \index{subgroup!connected}connected,
\index{unipotent!subgroup!maximal}unipotent subgroup of~$S$, and assume
$\{u^t\} \subset \hat U$.

From Rem.~\ref{PolyDivInG/U}\index{polynomial!divergence}, we know, for
$x,y \in \Gamma \backslash G$ with $x \approx y$, that the transverse
\index{component!transverse}component of the relative position between
$xu$ and~$yu$ is a \index{function!polynomial}polynomial function of~$u$.
 (Actually, it is a \index{function!rational}rational function
\cfex{BigCellSL2Ratl-ratl}, but this technical issue does not cause any
serious problems, because the function is unbounded on~$U$, just like a
polynomial would be.)
 Furthermore, the \index{transverse divergence}transverse
\index{component!transverse}component belongs to~$S$ (usually) and
\index{normalizer|)}normalizes~$\hat U$
\seeprops{BasLem}{RatTransInNformal}. Let $\hat S$ be the closure of the
subgroup of \index{normalizer}$N_S(\hat U)$ that is generated by the
image of one of these polynomial maps, together with~$\hat U$. Then $\hat
S$ is (almost) \index{Zariski!closed!almost}Zariski closed
\seelem{poly+U=algic}, and the maximal
\index{unipotent!subgroup!maximal}unipotent subgroup~$\hat U$ is
\index{subgroup!normal}normal, so the structure theory of
\index{algebraic group!theory of}algebraic groups implies that there is a
hyperbolic \index{torus!hyperbolic}torus~$T$ of~$\hat S$, and a compact
subgroup~$C$ of~$\hat S$, such that 
 $ \hat S = T C \hat U $ (see Thm.~\ref{AlgGrpTLU}
and Cor.~\fullref{AlgGrpCpct}{semi}).
 Any nonconstant \index{function!polynomial}polynomial is unbounded, so
(by definition of~$\hat S$), we see that $\hat S/\hat U$ is not compact;
 thus, $T$ is not compact. Let 
 \begin{itemize}
 \item $\{a^s\}$\index{subgroup!one-parameter!hyperbolic} be a noncompact
one-dimensional\index{dimension!of a Lie group} subgroup of~$T$,
 and
 \item $U = S_+$.
 \end{itemize}

This does not establish \pref{HyperAssump}, but it comes close:
 \begin{itemize}
 \item $\mu$ is invariant under $\{a^s\}$,
 and
 \item $\{a^s\}$ is hyperbolic\index{subgroup!one-parameter!hyperbolic},
and \index{normalizer}normalizes~$U$.
 \end{itemize}
 We have not constructed a subgroup~$L$, isomorphic to
$\SL(2,\real)$\index{SL(2,R)*$\SL(2,\real)$}, that contains~$\{a^s\}$,
but the only real use of that assumption was to prove that
$J(a^{-s},\widetilde S_-) \le J(a^s, S_+)$\index{Jacobian}
\see{Stilde=SPf-J<J}. Instead of using the theory of
$\SL(2,\real)$\index{SL(2,R)*$\SL(2,\real)$}-modules, one shows, by using
the theory of \index{algebraic group!theory of}algebraic groups, and
choosing $\{a^s\}$ carefully, that $J(a^s, H) \ge 1$, for every
\index{Zariski!closed}Zariski closed subgroup~$H$ of~$G$ that is
\index{normalizer}normalized by $a^s \, U$
 \seeex{J(a,H)>1Ex}.

An additional complication comes from the fact that $a^s$ may not act
\term[ergodic!flow]{ergodically} (w.r.t.~$\mu$): although $u^t$ is
\term[ergodic!measure]{ergodic}, we cannot apply the \term{Mautner
Phenomenon}, because $U = S_+$ may not contain~$\{u^t\}$ (since $(u^t)_-$
or $(u^t)_0$ may be nontrivial). Thus, one works with
\term[ergodic!component]{ergodic components} of~$\mu$. The key point is
that the arguments establishing Prop.~\ref{Stilde=S} actually show that
each \term[ergodic!component]{ergodic component} of~$\mu$ is $\widetilde
S_-$-invariant. But then it immediately follows that $\mu$ itself is
$\widetilde S_-$-invariant, as desired, so nothing was lost.

\begin{exercises}

\item \label{BigCellSL2Ratl}
 Let $B = \left\{ \begin{bmatrix} * & * \\ 0 & * \end{bmatrix}
\right\} \subset \SL(2,\real)$\index{SL(2,R)*$\SL(2,\real)$},
 and define $\psi \colon \mathbb{U}_2 \times B \to
\SL(2,\real)$\index{SL(2,R)*$\SL(2,\real)$} by $\psi( u, b) = u b$.
 \begin{enumerate}
 \item Show $\psi$ is a polynomial\index{function!polynomial}.
 \item Show $\psi$ is injective.
 \item Show the image of~$\psi$ is a dense, open subset~$\neigh$ of
$\SL(2,\real)$\index{SL(2,R)*$\SL(2,\real)$}.
 \item \label{BigCellSL2Ratl-ratl}
 Show $\psi^{-1}$ is a 
\index{function!rational}rational function on~$\neigh$.
 \end{enumerate}
 \hint{Solve $\begin{bmatrix} \mathsf{x} & \mathsf{y} \\ \mathsf{z} & \mathsf{w}
\end{bmatrix}
 = \begin{bmatrix} 1 & 0 \\ \mathsf{u} & 1 \end{bmatrix}
 \begin{bmatrix} \mathsf{a} & \mathsf{b} \\ 0 & 1/\mathsf{a}
\end{bmatrix}$
 for $\mathsf{a}$, $\mathsf{b}$, and~$\mathsf{u}$.}

 \item \label{J(a,H)>1Ex}
 Show there is a (nontrivial)
\term[subgroup!one-parameter!hyperbolic]{ hyperbolic one-parameter
subgroup} $\{a^s\}$ of~$\hat S$, such that \index{Jacobian}$J(a^s, H) \ge
1$, for every \index{Zariski!closed!almost}almost-Zariski closed
subgroup~$H$ of~$G$ that is \index{normalizer}normalized by $\{a^s\} \,
U$, and every $s > 0$.
 \hint{Let $\phi \colon U \to \hat S$ be a
\index{function!polynomial}polynomial, such that $\langle \phi(U), U
\rangle = \hat S$. For each $H$, we have $J(u,H) = 1$ for all $u \in U$,
and the function $J \bigl( \phi(u), H \bigr)$ is a polynomial on~$U$.
Although there may be infinitely many different possibilities for~$H$,
they give rise to only finitely many different polynomials, up to a
bounded error. Choose $u \in U$, such that $|J(\phi(u), H)|$ is large for
all~$H$, and let $a^1 = \phi(u)_h$ be the hyperbolic part in the
\index{Jordan!decomposition}Jordan decomposition of~$\phi(u)$.}

 \end{exercises}

\begin{notes}

Our presentation in this chapter borrows heavily from the original proof
of M.~Ratner \cite{Ratner-Solvable, Ratner-SS, Ratner-Meas}, but its
structure is based on the approach of G.~A.~Margulis and G.~M.~Tomanov
\cite{MargulisTomanov-Ratner}.
% (A complete proof of
%\index{Ratner's Theorems}Ratner's Theorem can be found in either of these
%sources.) 
 The two approaches are similar at the start, but, instead of employing
the \index{entropy!of a dynamical system}entropy calculations of
\S\ref{CompareSSect} to finish the proof, Ratner \cite[Lem.~4.1,
Lem.~5.2, and proof of Lem.~6.2]{Ratner-SS} bounded the number of small
\index{counting rectangular boxes}rectangular boxes needed to cover
certain subsets of $\Gamma \backslash G$.
 This allowed her to show \cite[Thm.~6.1]{Ratner-SS} that the
measure~$\mu$ is supported on an orbit of a subgroup~$L$, such that
 \begin{itemize}
 \item $L$ contains both $u^t$ and $a^s$, 
 \item the Jacobian $J(a^s,\Lie L)$ of~$a^s$ on the Lie algebra of~$L$
is~$1$,
 and
 \item $L_0 L_+ \subset \Stab_G(\mu)$.
 \end{itemize}
 Then an elementary argument \cite[\S7]{Ratner-SS} shows that $L \subset
\Stab_G(\mu)$.

The proof of Margulis and Tomanov is shorter, but less elementary,
because it uses more of the theory of \index{algebraic group!theory
of}algebraic groups.

\notesect{ShearPolyPfSect}

Lemma~\ref{RatTransInN} is a version of \cite[Thm.~3.1
(``R-property")]{Ratner-Solvable}\index{R-property|indsee{property,~R-}}\index{property!R-}
and (the first part of) \cite[Prop.~6.1]{MargulisTomanov-Ratner}.

Lemma~\ref{PolyDiv} is implicit in \cite[Thm.~3.1]{Ratner-Solvable}
and is the topic of \cite[\S5.4]{MargulisTomanov-Ratner}.

Proposition~\ref{BasLem} is 
\cite[Lem.~7.5]{MargulisTomanov-Ratner}. It is also implicit in the work
of M.~Ratner (see, for example, \cite[Lem.~3.3]{Ratner-SS}).

\notesect{DefineStildeSect}

The definition \pref{StildeDefn} of~$\widetilde{S}$ is based on
\cite[\S8.1]{MargulisTomanov-Ratner} (where $\widetilde{S}$ is denoted
$\mathcal{F}(s)$ and $\widetilde{S}_-$ is denoted $U^-(s)$).

\notesect{ShearConseq}

Corollary~\ref{GminInSmin} is \cite[Cor.~8.4]{MargulisTomanov-Ratner}.

Lemma~\ref{NotS->G+Div} is a special case of the last sentence of
\cite[Prop.~6.7]{MargulisTomanov-Ratner}.

\notesect{CompareSSect}

Proposition~\ref{Stilde=S} is \cite[Step~1 of
10.5]{MargulisTomanov-Ratner}.

\notesect{EndOfPfSect}

The proof of Prop.~\ref{SuppGU} is based on
\cite[Lem.~3.3]{MargulisTomanov-Ratner}.

The \term[Claim]{\emphit{Claim}} in the proof of Thm.~\ref{SuppS} is
\cite[Step~2 of 10.5]{MargulisTomanov-Ratner}.

The use of \index{entropy!of a dynamical system}entropy to prove that if
$\mu$ is $G_-$-invariant\index{measure!invariant}, then it is
$G_+$-invariant (alternative~\pref{Pf-FinishByEntropy} on
p.~\pageref{Pf-FinishByEntropy}) is due to Margulis and Tomanov
\cite[Step~3 of 10.5]{MargulisTomanov-Ratner}.

References to results on invariant measures for
\index{subgroup!horospherical}horospherical subgroups
(alternative~\pref{Pf-FinishByHoro} on p.~\pageref{Pf-FinishByHoro}) can
be found in the historical notes at the end of Chap.~\ref{IntroChap}.

Exercise~\ref{NoOrbit->NoSubvar} is
\cite[Prop.~3.2]{MargulisTomanov-Ratner}.

\notesect{PreciseSect}

The \index{technical assumption}technical assumption \pref{TechAssump}
needed for the proof of \pref{BasLem} is based on the condition~($*$) of
\cite[Defn.~6.6]{MargulisTomanov-Ratner}. (In Ratner's approach, this role
is played by \cite[Lem.~3.1]{Ratner-SS} and related results.)

 Exers.~\ref{GminInSminEx} and~\ref{NotS->G+DivTechEx} are special cases
of \cite[Prop.~6.7]{MargulisTomanov-Ratner}.

\notesect{RatnerPfAssumeASect}

 That $\hat S/U$ is not compact is part of
\cite[Prop.~6.1]{MargulisTomanov-Ratner}. 

That $a^s$ may be chosen to satisfy the condition
\index{Jacobian}$J(a^s, H) \ge 1$ is
\cite[Prop.~6.3b]{MargulisTomanov-Ratner}.

The (possible) \term[ergodic!measure]{nonergodicity} of~$\mu$ is addressed
in \cite[Step~1 of 10.5]{MargulisTomanov-Ratner}.

\end{notes}
 \index{measure!invariant|)}